\documentclass[10pt,twoside,leqno]{amsart}
\usepackage{amssymb,amsbsy,amsmath,amsfonts,amssymb,amscd,epsfig,times,graphics,color}
\usepackage[latin1]{inputenc}
\newcommand{\B}                 {\mathbb{B}}
\newcommand{\C}                 {\mathbb{C}}
\newcommand{\I}                 {\mathbb{I}}

\newcommand{\LL}                {\mathbb{L}}
\newcommand{\R}                 {\mathbb{R}}
\newcommand{\N}                 {\mathbb{N}}

\newtheorem{theorem}{Theorem}
\newtheorem{lemma}{Lemma}
\newtheorem{proposition}{Proposition}
\newtheorem{corollary}{Corollary}

\begin{document}


\title[Characteristic foliations and envelopes of
holomorphy]{Characteristic foliations on maximally real submanifolds
of $\C^n$ and envelopes of holomorphy}

\author{Jo\"el Merker and Egmont Porten}

\address{
CNRS, Universit\'e de Provence, LATP, UMR 6632, CMI, 
39 rue Joliot-Curie, F-13453 Marseille Cedex 13, France}

\email{merker@cmi.univ-mrs.fr} 

\address{Humboldt-Universit\"at zu Berlin, 
Mathematisch-Naturwissenschaftenschaft\-liche Fa\-kult\"at II, 
Institut f\"ur Mathematik, Rudower Chaussee 25, 
D-12489 Berlin, Germany}

\email{egmont@mathematik.hu-berlin.de}


\date{\number\year-\number\month-\number\day}

\begin{abstract}
Let $S$ be an arbitrary real $2$-surface, with or without boundary,
contained in a hypersurface $M$ of $\C^2$, with $S$ and $M$ of class
$\mathcal{ C}^{ 2, \alpha}$, where $0< \alpha <1$. If $S$ is totally
real except at finitely many complex tangencies which are hyperbolic
in the sense of E.~Bishop and if the union of separatrices is a tree
of curves without cycles, we show that every compact $K$ of $S$ is
CR-, $\mathcal{ W}$- and $L^{\sf p}$-removable (Theorem~1.3). Our
purely local techniques enable us to formulate substantial
generalizations of this statement, for the removability of closed sets
in totally real $1$-codimensional submanifolds contained in generic
submanifolds of CR dimension $1$.
\end{abstract}

\maketitle

\begin{center}
\begin{minipage}[t]{11cm}
\baselineskip =0.35cm
{\scriptsize

\centerline{\bf Table of contents}

\smallskip

{\bf Part~I \dotfill 1.}

\ \
{\bf 1.~Introduction \dotfill 1.}

\ \
{\bf 2.~Description of the proof of Theorem~1.2 and organization 
of the paper \dotfill 8.}

\ \
{\bf 3.~Strategy per absurdum for the proofs of
Theorems~1.2' and~1.4 \dotfill 13.} 

\ \
{\bf 4.~Construction of
a semi-local half-wedge \dotfill 19.}

\ \
{\bf 5.~Choice of a special point of $C_{\rm nr}$ to be removed 
locally \dotfill 32.}

{\bf Part~II \dotfill 45.}

\ \
{\bf 6.~Three preparatory lemmas on H\"older spaces \dotfill 45.}

\ \
{\bf 7.~Families of analytic discs half-attached
to maximally real submanifolds \dotfill 47.}

\ \
{\bf 8.~Geometric properties of families of half-attached analytic
discs \dotfill 55.}

\ \
{\bf 9.~End of proof of Theorem~1.2': application of 
the continuity principle \dotfill 63.}

\ \
{\bf 10.~Three proofs of Theorem~1.4 \dotfill 69.}

\ \
{\bf 11.~$\mathcal{ W}$-removability implies
$L^{\sf p}$-removability \dotfill 73.}

\ \
{\bf 12.~Proofs of Theorem~1.1 and of Theorem~1.3
\dotfill 78.}

\ \
{\bf 13.~Applications to the edge of the wedge theorem \dotfill 87.}

\ \
{\bf 14.~An example of a nonremovable three-dimensional torus \dotfill 
91.}

\ \
{\bf 15.~References \dotfill 96.}

\smallskip

{\footnotesize\tt [With 24 figures]}

}\end{minipage}
\end{center}

\bigskip

\section*{\S1.~Introduction}

In the past fifteen years, remarkable progress has been made towards
the understanding of the holomorphic extendability properties of CR
functions. At the origin of this development, the most fundamental
achievement was the deep discovery, due to the effort of numerous
mathematicians, that the so-called {\sl CR orbits} are the adequate
underlying objects for the semi-local CR analysis on a general
embedded CR manifold. As an independent and now established theory in
several complex variables, one may find a precise correspondence
between such orbits and progressively attached analytic discs covering
a thick part of the envelope of holomorphy of CR manifolds, {\it cf.}
\cite{ b}, \cite{ trv}, \cite{ tr1}, \cite{ tu1}, \cite{ ber}, \cite{
tu2}, \cite{ m1}, \cite{ j2} and \cite{ p3} for a recent synthesis.

Within this framework, it became mathematically accessible to
endeavour the general study of removable singularities on embedded CR
manifolds $M\subset \C^n$ of arbitrary CR dimension and of arbitrary
codimension, not necessarily being the boundaries of (strictly)
pseudoconvex domains. With respect to their size or ``mass'', the
interesting singularities can be essentially ordered by
their codimension in $M$. For instance, provided it does not perturb
the fact that $M$ consists of a single CR orbit, an arbitrary closed
subset $C\subset M$ which is of vanishing 2-codimensional Hausdorff
content is {\it always}\, removable, as is shown in~\cite{ cs} in the
hypersurface case and in \cite{ mp3}, Theorem~1.1, in arbitrary
codimension. Hence one is left to study the removability of
singularities of codimension at most two. Since the general problem
of characterizing removability seems at the moment to be out of reach
(even for $M$ being a hypersurface), it is advisible to focus on
geometrically accessible singularities, namely singularities contained
in a CR submanifold of $M$. A complete study of the automatic
removability of two-codimensional singularities may be found in
Theorem~4 of~\cite{ mp1}. Having in mind the classical Painlev\'e
problem, we will mainly consider in this paper singularities which are
closed sets $C$ contained in a {\it codimension one}\, 
submanifold $M^1$ of $M$ which is {\it generic}\, in $\C^n$.

The known results on singularities of codimension one can be
subdivided into two strongly different groups according to the {\it CR
dimension} of $M$. If ${\rm CRdim}\, M \geq 2$, then a generic
hypersurface $M^1\subset M$ is itself a CR manifold of positive CR
dimension, and singularities $C\subset M^1$ can be understood on the
basis of the interplay between $C$ and the CR orbits of $M^1$. Deep
results in this direction were established when $M$ is a hypersurface
of $\C^n$ in~\cite{ j4}, \cite{ j5} and then generalized to CR
manifolds of arbitrary codimension in~\cite{ p1}: the 
geometric condition insuring automatic removability is simply
that $C$ does not contain any CR orbit of $M^1$.

On the other hand, if ${\rm CRdim}\, M=1$ the geometric situation
becomes highly different, as a generic hypersurface $M^1\subset M$ is
now (maximal) {\it totally real}. Fortunately, as a substitute for the
CR orbits of $M^1$, one can consider the so-called {\sl characteristic
foliation}\, of $M^1$, obtained by integrating the characteristic line
field $T^c M|_{M^1} \cap TM^1$. But removability theorems exploiting
this concept were only known for hypersurfaces in $\C^2$ and, until
very recently, only in the strictly pseudoconvex case. Furthermore,
a geometric condition insuring automatic removablitity
has not yet been clearly delineated.

Hence, with respect to the current state of the art, there was a
two-fold gap about codimension one removable singularities contained
in generic submanifolds $M$ of CR dimension one: firstly, to establish
a satisfying theory for non-pseudoconvex hypersurfaces in $\C^2$ and
secondly, to understand the situation in higher codimension. This
second main task was formulated as the first open problem in a list
p.~432 of \cite{j5} ({\it see}\, also the comments pp.~431--432 about
the relative geometric simplicity of the case ${\rm CRdim}\, M \geq
2$). {\it A priori}, it is not clear at all whether the two directions
of research are related somehow, but in the present work, we shall
fill in this two-fold gap by devising a new semi-local approach which
applies uniformly with respect to codimension.

For the detailed discussion of our result we have to introduce some
terminology which will be used throughout the article. Let $M$ be a
generic submanifold of $\C^n$ and let $C$ be a closed subset of $M$.
Recall from~\cite{mp3} that a {\sl wedgelike domain} attached to a
generic submanifold $M'\subset \C^n$ is a domain containing a local
wedge of edge $M'$ at every point of $M'$. Our wedgelike domains will
always be {\it nonempty}. Let us define three basic notions of
removability. Firstly, we say that $C$ is {\sl CR-removable} if there
exists a wedgelike domain $\mathcal{ W}$ attached to $M$ to which
every continuous CR function $f\in\mathcal{ C}_{CR}^0(M\backslash C)$
extends holomorphically. Secondly, as in~\cite{mp3}, p.~486, we say
that $C$ is {\sl $\mathcal{ W}$-removable} if for every wedgelike
domain $\mathcal{ W}_1$ attached to $M\backslash C$, there is a
wedgelike domain $\mathcal{ W}_2$ attached to $M$ and a wedgelike
domain $\mathcal{ W}_3 \subset \mathcal{ W}_1 \cap \mathcal{ W}_2$
attached to $M\backslash C$ such that for every holomorphic function
$f\in\mathcal{ O}(\mathcal{ W}_1)$, there exists a holomorphic
function $F\in\mathcal{ O}(\mathcal{ W}_2)$ which coincides with $f$
in $\mathcal{ W}_3$. Thirdly, with ${\sf p} \in\R\cup \{+\infty\}$
satisfying ${\sf p} \geq 1$, we say that $C$ is {\it $L^{\sf
p}$-removable} if every locally integrable function $f\in L^{\sf
p}_{loc}(M)$ which is CR in the distributional sense on $M\backslash
C$ is in fact CR on all of $M$.

The first notion of removability is a generalization of the kind of
removability considered in most of the pioneering papers \cite{cs},
\cite{d}, \cite{fs}, \cite{j1}, \cite{kr}, \cite{l}, \cite{lu},
\cite{st} about removable singularities in boundaries of domains $D
\subset \subset \C^n$. We observe that a wedgelike open set attached
to a hypersurface $M$ is just a (global) one-sided neighborhood of
$M$, namely a domain $\omega$ with $\overline{\omega}\supset M$ such
that for every point $p\in M$, the domain $\omega$ contains the
intersection of a neighborhood of $p$ in $\C^n$ with one side of
$M$. If now a closed set $C$ contained in a $\mathcal{ C}^1$-smooth
bounded boundary $\partial D$ is CR-removable, then an application of
the Hartogs-Bochner theorem shows that CR functions on $\partial
D\backslash C$ can be holomorphically extended to $D$. The second
notion of removability is a way to isolate the part of the question
related to envelopes of holomorphy. The third notion of removability
has the advantage of being completely intrinsic with respect to $M$
and may be relevant in the study of non-embeddable CR manifolds.

To avoid confusion, we state precisely our submanifold notion: $Y$ is
a {\sl submanifold of $X$} if $Y$ and $X$ are equipped with a manifold
structure, if there exists an immersion $i$ of $Y$ into $X$ and if the
manifold topology of $Y$ and the topology of $i(Y)$ inherited from the
topology of $X$ coincide, so that one may identify the submanifold $Y$
with the subset $i(Y)\subset X$. Furthermore, our submanifolds will
always be connected.

Let us now enter the discussion of the case $n=2$. Here we shall
denote the submanifold $M^1\subset M$, which is a surface in $\C^2$,
by $S$. In~\cite{b}, E.~Bishop showed that a two-dimensional surface
in $\C^2$ of class at least $\mathcal{ C}^2$ having an isolated
complex tangency at one of its points $p$ may be represented by a
complex equation of the form $w=z\bar z+\lambda(z^2+\bar z^2)+{\rm
o}(\vert z\vert^2)$, in terms of local holomorphic coordinates $(z,w)$
vanishing at $p$, where the real parameter $\lambda\in [0,\infty)$ is
a biholomorphic invariant of $S$. The point $p$ is said to be {\sl
elliptic} if $\lambda\in [0,\frac{1}{2})$, {\sl parabolic} if $\lambda
= \frac{1}{2}$ and {\sl hyperbolic} if $\lambda \in (\frac{1}{2},
\infty)$. Recall that $M$ is called {\sl globally minimal} if it
consists of a single CR orbit ({\it cf.} \cite{tr1}, \cite{tr2};
\cite{mp1}, pp.~814--815; and~\cite{j4}, pp.~266--269). Throughout
this paper, we shall work in the $\mathcal{ C}^{ 2, \alpha}$-smooth
category, where $0 < \alpha < 1$. Our first main new result is as
follows.

\def\thetheorem{{\bf 1.1}}\begin{theorem}
Let $M$ be a globally minimal $\mathcal{ C}^{ 2, \alpha}$-smooth
hypersurface in $\C^2$ and let $D \subset M$ be a $\mathcal{ C}^{ 2,
\alpha}$-smooth surface which is
\begin{itemize}
\item[{\bf (a)}]
$\mathcal{C}^{2,\alpha}$-diffeomorphic to the unit
$2$-disc of $\R^2$ and 
\item[{\bf (b)}]
totally real outside a discrete subset of isolated complex tangencies
which are hyperbolic in the sense of E.~Bishop.
\end{itemize}
Then every compact subset $K$ of $D$ is CR-, $\mathcal{ W}$- and
$L^{\sf p}$-removable.
\end{theorem}

As a corollary, one obtains a corresponding result about holomorphic
extension from $\partial\Omega\backslash K$ for the case that $M$ is
the boundary of a relatively compact domain $\Omega\subset\C^2$. Note
that $\partial \Omega$ is automatically globally minimal
(\cite{ j4}, Section~2). We will first recall the historical
background of Theorem 1.1 and explain afterwards on this basis the
main ideas and techniques necessary for the proof. 

In 1988, applying a global
version of the Kontinuit\"atssatz, B.~J\"oricke~\cite{j1} established
a remarkable theorem: every compact subset of a totally real
$\mathcal{ C}^2$-smooth $2$-disc lying on the boundary of the unit
ball in $S^3=\partial \B_2\subset \C^2$ is CR-removable. This
discovery motivated the work~\cite{fs} by
F.~Forstneri$\check{\text{\rm c}}$ and E.L.~Stout, where it is shown
that every $\mathcal{ C}^2$-smooth compact $2$-disc contained in a
strictly pseudoconvex $\mathcal{ C}^2$-smooth boundary $\partial
\Omega$ contained in a $2$-dimensional Stein manifold $\mathcal{ M}$
which is totally real except at a finite number of hyperbolic complex
tangencies is removable; the proof mainly relies on a previous work by
E.~Bedford and W.~Klingenberg about the hulls of $2$-spheres contained
in such strictly pseudoconvex boundaries $\Omega \subset \mathcal{
M}$, which may be filled by Levi-flat $3$-spheres after a generic
small perturbation (\cite{bk}, Theorem~1). Indirectly, it followed
from~\cite{j1} and~\cite{fs} that such compact totally real $2$-discs
$D\subset \partial \Omega$ (possibly having finitely many hyperbolic
complex tangencies) are $\mathcal{ O}(\overline{\Omega})$-convex and
in particular polynomially convex if $D=\B_2$ and $\mathcal{ M}=\C^2$,
thanks to a previous work~\cite{st} by E.L.~Sout, where it is shown
(Theorem~II.10) that a compact subset $K$ of a $\mathcal{ C}^2$-smooth
strictly pseudoconvex boundary $\partial \Omega$ in a Stein manifold
is removable {\it if and only if}\, $K$ is $\mathcal{
O}(\overline{\Omega})$-convex. It is also established in~\cite{fs}
that a neighborhood of an isolated hyperbolic complex tangency in
$\C^2$ is polynomially convex. These papers have been followed by the
work~\cite{d}, where the question of $\mathcal{
O}(\overline{\Omega})$-convexity of {\it arbitrary compact surfaces}\,
$S$ (with or without boundary, not necessarily diffeomorphic to a
$2$-disc) contained in a $\mathcal{ C}^2$-smooth strictly pseudoconvex
domain $\Omega\subset \C^2$ is dealt with directly. Using K.~Oka's
characterization of the envelope of a compact, J.~Duval shows that the
essential hull $\widehat{ K}_{\rm ess}:= \overline{ \widehat{ K}_{
\mathcal{ O }( \overline{ \Omega} )} \backslash K}$ must cross
every leaf of the characteristic foliation on the totally real part of
$S$ and he deduces that a compact $2$-disc having only hyperbolic
complex tangencies is $\mathcal{ O }(\overline{ D })$-convex.

All the above proofs heavily rely on strong pseudoconvexity, in
contrast to the experience, familiar at least in the case ${\rm
CRdim}\, M \geq 2$, that removability should depend rather on the
structure of CR orbits than on Levi curvature. The first theorem for
the non-pseudoconvex situation was established by the second author in
\cite{p2}. He proved that every compact subset of a totally real disc
embedded in a globally minimal $\mathcal{ C}^\infty$-smooth
hypersurface in $\C^2$ is always CR-removable. We would like to point
out that, seeking theorems without any assumption of pseudoconvexity
leads to substantial open problems, because one loses almost all of
the strong interweavings between function-theoretic tools and
geometric arguments which are valid in the pseudoconvex realm, for
instance: Hopf Lemma, plurisubharmonic exhaustions, envelopes of
function spaces, local maximum modulus principle, Stein neighborhood
basis, {\it etc.}

To discuss the main elements of our approach, let us briefly explain
the geometric setup of the proof of Theorem 1.1. The characteristic
foliation has isolated singularities at the hyperbolic points, where it
looks like the phase diagram of a saddle point. In particular there
are four local separatrices accumulating orthogonally at each
hyperbolic point. Hence we can decompose the $2$-disc $D$ as a union
$D=T_D \cup D_o$, where $T_D$ consists of the union of the hyperbolic
points of $D$ together with the separatrices issuing from them, and
where $D_o:= D \backslash T_D$ is the remaining open submanifold of
$D$, contained in the totally real part of $D$. By
H.~Poincar\'e and I.~Bendixson's theory, $T_D$ is a tree of $\mathcal{
C}^{2,\alpha}$-smooth curves which contains no subset homeomorphic to
the unit circle, {\it cf.}~\cite{ d}. Accordingly, we decompose $K:=
K_{ T_D} \cup C_o$, where $K_{T_D}:= K \cap T_D$ is a proper closed
subset of the tree $T_D$ and where $C_o := K \cap D_o$ is a relatively
closed subset of $D_o$.

The hard part of the proof, which was actually the starting point of
the whole paper, will consist in removing the closed subset $C_o$ of
the $2$-dimensional surface $S:=D_o$ lying in $M \backslash
T_D$. Thereafter the removal of the remaining part $K_{T_D}$ will
be done by means of an investigation of the behaviour of
the CR orbits near $T_D$, close in spirit to our previous
methods in ~\cite{mp1} ({\it see}\, Section~12 below for the details).

Let us formulate the first crucial part of the above argument as an
independent theorem about the removal of closed subsets contained in a
totally real surface $S$. We point out that now $S$ may have {\it
arbitrary topology}.

\def\thetheorem{{\bf 1.2}}\begin{theorem}
Let $M$ be a $\mathcal{ C}^{ 2, \alpha }$-smooth globally minimal
hypersurface in $\C^2$, let $S \subset M$ be a $\mathcal{ C}^{2,
\alpha}$-smooth surface, open or closed, with or without boundary,
which is totally real at every point. Let $C$ be a proper closed
subset of $S$ and assume that the following topological condition
holds{\rm :}

\smallskip
\begin{itemize}
\item[$\mathcal{ F }_S^c \{C\}:$] For every closed subset $C' \subset
C$, there exists a simple $\mathcal{ C}^{2,\alpha}$-smooth curve
$\gamma: [-1, 1] \to S$, whose range is contained in a single leaf of
the characteristic foliation $\mathcal{ F }_S^c$ {\rm (}obtained
by integrating the characteristic line field $T^cM \vert_S \cap 
TS${\rm )}, with $\gamma (-1)
\not \in C'$, $\gamma (0) \in C'$ and $\gamma (1) \not \in C'$, such
that $C'$ lies completely in one closed side of $\gamma [-1, 1]$
with respect to the topology of $S$ in a neighborhood of
$\gamma[-1,1]$.
\end{itemize}

\smallskip
\noindent
Then $C$ is CR-, $\mathcal{ W}$- and $L^{\sf p}$-removable.
\end{theorem}

The condition $\mathcal{ F}_S^c \{C \}$ is a {\it common condition}\,
on $C$ and on the characteristic foliation $\mathcal{ F}_S^c$, namely
on the relative disposition of $\mathcal{ F}_S^c$ with respect to $C$,
not only on $S$; an illustration may be found in {\sc Figure~2} below.
In the strictly pseudoconvex context, this condition appeared
implicitly during the course of the proofs given in \cite{d}. Note
that the relevance of the characteristic foliation had earlier been
discovered in contact geometry, {\it cf.} \cite{be}, \cite{e}. It is
interesting to notice that it re-appears in the situation of Theorem
1.2, where the underlying distribution $T^c M$ is allowed to be very
far from contact.

As is known, it follows from a subcase of H.~Poincar\'e and
I.~Bendixson's theory that if $S$ is diffeomorphic to a real $2$-disc
or if $S=D_o$ as above, then $\mathcal{ F}_S^c \{C\}$ is automatically
satisfied for an arbitry compact subset $C$ of $S$. On the contrary, it
may be not satisfied when for instance $S$ is an annulus equipped with
a radial foliation together with $C$ containing a continuous closed
curve around the hole of $S$. Crucially, it is elementary to construct
an example of such an annulus which is truly nonremovable. Indeed, the
small closed curve $C$ which consists of the transversal intersection
of a strictly convex boundary $\partial D$ with a complex line close
to a boundary point may be enlarged as a thin maximally real strip $S
\subset \partial D$ which is diffeomorphic to an annulus; in this
setting, $C$ is obviously nonremovable and {\it the characteristic
foliation is everywhere transversal to $C$}. Consequently, the
geometric condition $\mathcal{ F}_S^c\{C\}$ is the optimal one
insuring automatic removability for all choices of $M$, $S$ and
$C$. Further examples of closed subsets in surfaces with arbitrary
genus equipped with such foliations may be exhibited.

In the proof of Theorem~1.2, after some contraction $C'$ of $C$, we
may assume that no point of $C'$ is locally removable ({\it see}\,
Sections~2 and~3 below). Then the very assumption $\mathcal{
F}_S^c\{C\}$ yields the existence of a characteristic segment
$\gamma[-1,1]$, such that $C'$ lies on one side of
$\gamma[-1,1]$. Reasoning by contradiction, our aim is to show that
there exists at least one special point $p\in C'\cap \gamma(-1,1)$,
which is {\it locally}\, CR-, $\mathcal{ W}$- and $L^{\sf
p}$-removable. The choice of such a point $p$, achieved in Section~5
below, will be nontrivial.

The strategy for the local removal of $p$ is to construct an analytic
disc $A$ such that a segment of its boundary $\partial A$ is attached
to $S$ and touches $C'$ in only one point $p$. Several geometrical
assumptions have to be met to ensure that a sufficiently rich family
of deformations of $A$ have boundaries disjoint from $C'$, that
analytic extension along these discs is possible (i.e.~appropriate
moment conditions are satisfied), and that the union of these good
discs is large enough to give analytic extension to a one-sided
neighborhood of $M$: this is where the (semi)localization and the
choice of the special point $p\in C'$ will be key ingredients.
Let us explain why localization is crucial.

Working globallly, the second author produced in~\cite{ p2} a
convenient disc by applying the powerful E.~Bedford and W.~Klingenberg
theorem to an appropriate 2-sphere containing a neighborhood of the
{\it entire}\, singularity $C'$. This method requires global
properties of $S$ like $S$ being a totally real $2$-disc, which
ensures the existence of a nice Stein neighborhood basis of
$C'$. Already for real discs with isolated hyperbolic points, it is
not clear whether this argument can be generalized (however, we would
like to mention that recent results of M.~Slapar in~\cite{ sl}
indicate that this could be possible at least if the geometry near the
hyperbolic points satisfies some additional assumptions). In the case
where $M$ is an arbitrary globally minimal hypersurface, where $S$ has
arbitrary topology and has complex tangencies, the reduction to
E.~Bedford and W.~Klingenberg's theorem seems impossible, {\it cf.}
the example of an unknotted nonfillable $2$-sphere in $\C^2$
constructed by J.E.~Forn{\ae}ss and D.~Ma in~\cite{fm}. Also, to the
authors' knowledge, the possibility of filling by Levi-flat
$3$-spheres the (not necessarily generic) $2$-spheres lying on a {\it
nonpseudoconvex}\, hypersurface is a delicate open problem. In
addition, for the higher codimensional generalization of Theorem~1.2,
the idea of global filling seems to be irrelevant at present times,
because no analog of the E.~Bedford and W.~Klingenberg theorem is
known in dimension $n\geq 3$. As we aim to deal with surfaces $S$
having arbitrary topology and to generalize these results in arbitrary
codimension, we shall endeavour to firmly {\it localize the
removability arguments}, using only small analytic discs.

Thus, our way to overcome these obstacles is to consider {\it local}
discs $A$ which are only partially attached to $S$. The delicate point
is that we have at the same time (i) to control the geometry of
$\partial A$ near $p\in C'$ and (ii) to guarantee that the rest of the
boundary stays in the region where holomorphic extension is already
known. In fact, (ii) will be incorporated in our very special and
tricky choice of $p\in C'$. For (i), we have to sharpen known existence
theorems about partially attached analytic discs and to combine it
with a careful study of the complex/real geometry of the pair
$(M,S)$. Importantly, our construction of such analytic discs is
achieved elementarily in a self-contained way. A precise description
of the proof in the hypersurface case (only) may be found in Section~2
below. With some substantial extra work, we shall generalize this
purely local strategy of proof to higher codimension.

To conclude with the removal of surfaces, let us formulate a more
general version of Theorem~1.1, whitout the restricted topological
assumption that $S$ be diffeomorphic to a real disc. Applying
Theorem~1.2 for the removal of $K\cap (S\backslash T_S)$ and a slight
generalization of Theorem~4 (ii) in~\cite{mp1} for the removal of
$K\cap T_S$ (more precisions will be given in \S13 below), we shall
obtain the following statement, implying Theorem~1.1 as a 
direct corollary.

\def\thetheorem{1.3}\begin{theorem}
Let $M$ be a $\mathcal{ C}^{ 2, \alpha }$-smooth globally minimal
hypersurface in $\C^2$, let $S \subset M$ be a $\mathcal{ C}^{2,
\alpha}$-smooth totally real surface, open or closed, with or without
boundary, which is totally real outside a discrete subset of isolated
complex tangencies which are hyperbolic in the sense of E.~Bishop.
Let $T_S$ be the union of hyperbolic points of $S$ together with all
separatrices issued from hyperbolic points and assume that $T_S$ does
not contain any subset which is homeomorphich to the unit circle. Let
$K$ be a proper compact subset of $S$ and assume that $\mathcal{ F
}_{S\backslash T_S}^c \{K \cap (S \backslash T_S)\}$ holds.

Then $K$ is CR-, $\mathcal{ W}$- and $L^{\sf p}$-removable.
\end{theorem}

As was already emphasized, our main motivation for this work was to
{\it devise a local strategy of proof for Theorems~1.1, 1.2 and~1.3 in
order to generalize them to higher codimension}. In fact, we will
realize the program sketched above for generic submanifolds of CR
dimension $1$ and of arbitrary codimension. Thus, let $M$ be a
$\mathcal{C}^{2,\alpha}$-smooth globally minimal
generic submanifold of codimension $(n-1)$ in $\C^n$, hence of CR
dimension $1$, where $n\geq 2$. Let $M^1$ be a maximally real
$\mathcal{C}^{2,\alpha}$-smooth one-codimensional submanifold of $M$
which is generic in $\C^n$. As in the surface case, $M^1$ carries a
{\sl characteristic foliation} $\mathcal{ F}_{M^1}^c$, whose leaves
are the integral curves of the line distribution $TM^1\cap
T^cM\vert_{M^1}$. Of course the assumption that the singularity lies
on one side of some characteristic segment is no longer reasonable.
We will generalize it as a condition requiring (approximatively
speaking) that there be always a characteristic segment which is
accessible from the complement of $C$ in $M^1$ along one direction
normal to the characteristic segment.

The generalization of Theorem~1.2, which is our
principal result in this paper, is as follows.

\def\thetheorem{{\bf 1.2'}}\begin{theorem}
Let $M$, $M^1$, $\mathcal{F}_{M^1}^c$ be as above and let $C$ be a
proper closed subset of $M$. Assume that the following topological
condition, meaning that $C$ is not transversal to the characteristic
foliation, holds{\rm :}
\smallskip
\begin{itemize}
\item[$\mathcal{ F}_{ M^1 }^c \{C\}:$] For every closed subset $C'
\subset C$, there exists a simple $\mathcal{ C}^{2, \alpha}$-smooth
curve $\gamma: [-1,1] \to M^1$ whose range $\gamma[-1,1]$ is
contained in a single leaf of the characteristic foliation
$\mathcal{F}_{ M^1}^c$ with $\gamma(-1)\not\in C'$, $\gamma (0) \in
C'$ and $\gamma (1) \not \in C'$, there exists a local $(n-
1)$-dimensional transversal $R^1 \subset M^1$ to $\gamma$ passing
through $\gamma(0)$ and there exists a thin elongated
open neighborhood $V_1$ of
$\gamma [-1,1]$ in $M^1$ such that if $\pi_{ \mathcal{F}_{ M^1}^c}:
V_1 \to R^1$ denotes the semi-local projection parallel to the leaves
of the characteristic foliation $\mathcal{F}_{ M^1}^c$, then
$\gamma(0)$ lies on the boundary, relatively to the topology of $R^1$,
of $\pi_{\mathcal{F}_{M^1}^c}(C' \cap V_1)$.
\end{itemize}
\smallskip
\noindent 
Then $C$ is CR-, $\mathcal{ W}$- and $L^{\sf p}$-removable.
\end{theorem}

The condition $\mathcal{ F}_{ M^1}^c \{C \}$, which is of course
independent of the choice of the transversal $R^1$ and of the thin
neighborhood $V_1$, is illustrated in {\sc Figure~8} of \S5.1 below;
clearly, in the case $n =2$, it means that $C' \cap V_1$ lies
completely in one side of $\gamma [-1, 1]$, with respect to the
topology of $M^1$, as written in the statement of
Theorem~1.2. Applying some of our previous results in this direction
(\cite{mp1}, \cite{mp3}), we shall provide in the end of Section~13
below some formulations of applications of Theorem~1.2', close to
being analogs of Theorem~1.3 in higher codimension.

Importantly, in order to let the geometric condition $\mathcal{
F}_{M^1}^c\{C\}$ appear less mysterious and to argue that it provides
the adequate generalization of Theorem~1.2 to higher codimension, in
the last Section~14 below, we shall describe an example of $M$, $M^1$
and $C$ in $\C^3$ violating the condition $\mathcal{ F}_{M^1}^c\{C\}$,
such that $C$ is transversal to the characteristic foliation and is
truly nonremovable. This example will be analogous in some sense to
the example of a nonremovable annulus discussed after the statement of
Theorem~1.2. Since there is no H.~Poincar\'e and I.~Bendixson 
theorem for
foliations of $3$-dimensional balls by curves, it will be even
possible to insure that $M$ and $M^1$ are diffeomorphic to real balls
of dimension $4$ and $3$ respectively. We may therefore conclude that
Theorem~1.2' provides the desirable answer to the (already cited {\em
supra}) Problem~2.1 raised by B.~J\"oricke in~\cite{ j5}, p.~432.

To pursue the presentation of our results, let us comment the
assumption that $M$ be of codimension $(n-1)$. Geometrically speaking,
the study of closed singularities $C$ lying in a one-codimensional
generic submanifold $M^1$ of a generic submanifold $M\subset \C^n$
which is of CR dimension $m\geq 2$ is more simple. Indeed, thanks to
the fact that $M^1$ is of CR dimension $m-1\geq 1$, there exist local
Bishop discs {\it completely attached to $M^1$}, and this helps much in
describing the envelope of holomorphy of a wedge attached to
$M\backslash C$. On the contrary, in the case where $M$ is of CR
dimension $1$, small analytic discs attached to a maximally real $M^1$
are (trivially) inexistent. This is why, in the proof of our main
Theorems~1.2 and~1.2', we shall deal only with small analytic discs
whose boundary is in part (only) contained in $M^1$. Such discs are
known to exist; we would like to mention that historically speaking,
the first construction of discs partially attached to maximally real
submanifolds was exhibited by S.~Pinchuk in~\cite{p}, who developed
the ideas of E.~Bishop~\cite{b}.

Finally we will test our main Theorem~1.2' in applications. First of
all, we clarify its relation to known removability results in CR
dimension greater than one. Here the motivation is simply that most
questions of CR geometry should be reducible to CR dimension $1$ by
slicing. It turns out that the main known theorems about removable
singularities, due to E.~Chirka, E.~L.~Stout, and B.~J\"oricke for
hypersurfaces, and by the authors in higher codimension (\cite{p1},
Theorem~1 about $L^{\sf p}$-removability; \cite{m2}, Theorem~3 about
CR- and $\mathcal{ W}$- removability) are all a rather direct
consequence of Theorem 1.2'. Since these results have not yet been
published in complete form, we take the occasion of including them in
the present paper, {\it as a corollary of Theorem~1.2'}, yet devising
a new geometric approach.

\def\thetheorem{{\bf 1.4}}\begin{theorem} 
Let $M$ be a $\mathcal{ C}^{2, \alpha}$-smooth globally minimal
generic submanifold of $\C^n$ of CR dimension $m\geq 2$ and of
codimension $d=n- m\geq 1$, let $M^1 \subset M$ be a $\mathcal{ C}^{2,
\alpha}$-smooth one-codimensional submanifold which is generic in
$\C^n$ and let $C \subset M^1$ be a proper closed subset of
$M$. Assume that the following
condition holds{\rm :}
\begin{itemize}
\item[$\mathcal{ O}_{M^1}^{CR}\{C\}:$] The closed subset
$C$ does not contain any CR orbit of $M^1$.
\end{itemize}
Then $C$ is CR-, $\mathcal{ W}$- and $L^{\sf p}$-removable.
\end{theorem}

Notice the difference with the case $m=1$, where the analog of CR
orbits would consist of characteristic curves: the condition
$\mathcal{ F}_{M^1}^c\{C\}$ does {\it not}\, say that $C$ should not
contain any maximal characteristic curve. In fact, we observe that
there cannot exist a uniform removability statement covering both the
case $m=1$ and the case $m\geq 2$, whence Theorem~1.2' is stronger
than Theorem~1.4. Indeed, the elementary example of a nonremovable
circle in an annulus contained in the boundary of a strictly convex
domain of $\C^2$ shows that $C$ may be truly nonremovable whereas it
does not contain any characteristic curve. In the strictly
pseudoconvex hypersurface setting, it is well known that Hopf's Lemma
implies that boundaries of Riemann surfaces contained in $C$ (and also
the track on $C$ of its essential hull, {\it cf.}~\cite{d}) should be
everywhere transversal to the characteristic foliation. Of course, this
implies conversely that $C$ cannot contain such boundaries (unless
they are empty) if $\mathcal{ F}_{M^1}^c\{C\}$ is satisfied. The
reason why $\mathcal{ F}_{M^1}^c\{C\}$ implies that $C$ is removable
also in the nonpseudoconvex setting and in arbitrary codimension will
be appearant later. Finally, we mention that the $L^{\sf
p}$-removability of $C$ in Theorem~1.4 holds more generally with no
assumption of global minimality on $M$, as already noticed in~\cite{
j5}, \cite{ p1}, \cite{mp1}. However, since the case where $M$ is not
globally minimal essentially reduces to the consideration of its CR
orbits, which are globally minimal by definition, we shall only
deal with globally minimal generic submanifolds $M$ throughout this
paper. 

As a final comment, we point out that, {\it because the previously known
proofs of Theorem~1.4 were of local type, it is satisfactory to bring
in this paper a purely local framework for the treatment of
Theorems~1.1, 1.2, 1.3 and~1.2'}.

Our second group of applications concerns the classical
edge of the wedge theorem. Typically one considers a maximally real
edge $E$ to which an open double wedge $(\mathcal{ W}_1, \mathcal{
W}_2)$ is attached from opposite directions. One may interprete this
configuration as a partial thickening of a generic CR manifold
$M\subset E\cup \mathcal{ W}_1\cup \mathcal{ W}_2$ containing $E$ as a
generic hypersurface. The classical edge of the wedge theorem states
that functions which are continuous on $ \mathcal{ W}_1 \cup \mathcal{
W}_2 \cup E$ and holomorphic in $\mathcal{ W}_1 \cup \mathcal{ W}_2$
extend holomorphically to a neighborhood of $E$. Theorem 1.2' implies
that it suffices to assume continuity outside a removable
singularities of $E$. This allows us to derive an edge of the wedge
theorem for meromorphic extension (Section~13 below).

This paper is divided in two parts: Part I contains the strategy per
absurdum for the proof of Theorem~1.2', the construction of what we
call a semi-local half-wedge and the choice of a special point to be
removed locally. Part~II contains the explicit construction of
families of half-attached analytic discs, the end of proof of
Theorem~1.2' and the proofs of the various applications. The reader
will find a more detailed description of the content of the paper in
\S2.16 below.

\subsection*{1.5.~Acknowledgements}
The authors would like to thank B.~J\"oricke for several valuable
scientific exchanges. They acknowledge generous support from the
European TMR research network ERBFMRXCT 98063 and they also thank the
universities of Berlin (Humboldt), of G\"oteborg (Chalmers), of
Marseille (Provence) and of Uppsala for invitations which provided
opportunities for fruitful mathematical research.

\section*{\S2.~Description of the proof of Theorem~1.2
and organization of the paper}

The main part of this paper is devoted to the proof of Theorem~1.2',
which will occupy Sections~3, 4, 5, 6, 7, 8 and~9 below. In this
preliminary section, we shall summarize the hypersurface version of
Theorem~1.2', namely Theorem~1.2. Our goal is to provide a
conceptional description of the basic geometric constructions, which
should be helpful to read the whole paper. Because precise, complete
and rigorous formulations will be developed in the next sections, we
allow here the use of a slightly informal language.

\subsection*{2.1.~Strategy per absurdum} Let $M$, $S$, and $C$ 
be as in Theorem~1.2. It is essentially known that both the CR- and
the $L^{\sf p}$-removability of $C$ are a (relatively mild)
consequence of the $\mathcal{ W}$-removability of $C$ ({\it see}\,
\S3.14 and Section~11 below). Thus, we shall describe in this section
only the $\mathcal{ W}$-removability of $C$. 

First of all, as $M$ is globally minimal, it may be proved that for
every closed subset $C' \subset C$, the complement $M\backslash C'$ is
also globally minimal ({\it see}\, Lemma~3.5 below). As $M$ is of
codimension one in $\C^2$, a wedge attached to $M\backslash C$ is
simply a connected one-sided neighborhood of $M\backslash C$ in
$\C^2$. Let us denote such a one-sided neighborhood by $\omega_1$. The
goal is to prove that there exists a one-sided neighborhood $\omega$
attached to $M$ to which holomorphic functions in $\omega_1$ extend
holomorphically. By the definition of $\mathcal{ W}$-removability,
this will show that $C$ is $\mathcal{ W}$-removable.

Reasoning by contradiction, we shall denote by $C_{\rm nr}$ the
smallest {\sl nonremovable} subpart of $C$. By this we mean that
holomorphic functions in $\omega_1$ extend holomorphically to a
one-sided neighborhood $\omega_2$ of $M\backslash C_{\rm nr}$ in
$\C^2$ and that $C_{\rm nr}$ is the smallest subset of $C$ such that
this extension property holds. If $C_{\rm nr}$ is empty, the
conclusion of Theorem~1.2 holds, gratuitously: nothing has to be
proved. If $C_{ \rm nr}$ is nonempty, to come to an absurd, it
suffices to show that at least one point of $C_{\rm nr}$ is {\sl
locally removable}. By this, we mean that there exists a local
one-sided neighborhood $\omega_3$ of at least one point of $C_{\rm
nr}$ such that holomorphic functions in $\omega_2$ extend
holomorphically to $\omega_3$. In fact, the choice of such a point
will be the most delicate and the most tricky part of the proof.

In order to be in position to apply the continuity principle, we now
deform slightly $M$ inside the one-sided neighborhood $\omega_2$,
keeping $C_{\rm nr}$ fixed, getting a hypersurface $M^d$ (with $d$
like ``{\sl d}eformed'') satisfying $M^d\backslash C_{\rm nr} \subset
\omega_2$. We notice that a local one-sided neighborhood of $M^d$ at
one point $p$ of $C_{\rm nr}$ always contains a local one-sided
neighborhood of $M$ at $p$ (the reader may draw a figure), so we may
well work on $M^d$ instead of working on $M$ (however, the analogous
property about wedges over deformed generic submanifolds is untrue in
codimension $\geq 2$, {\it see}\, \S3.16 below, where supplementary
arguments are needed).

Replacing the notation $C_{\rm nr}$ by the notation $C$, the notation
$M^d$ by the notation $M$ and the notation $\omega_2$ by the notation
$\Omega$, we see that Theorem~1.2 is reduced to the following main
proposition, whose formulation is essentially analogous to that
of Theorem~1.2, except that it
suffices to remove at least one special point.

\def\theproposition{2.2}\begin{proposition}
Let $M$ be a $\mathcal{ C}^{ 2, \alpha}$-smooth globally minimal
hypersurface in $\C^2$, let $S \subset M$ be a $\mathcal{ C}^{ 2,
\alpha}$-smooth surface which is totally real at every point. Let $C$
be a {\rm nonempty} proper closed subset of $S$ and assume that the
nontransversality condition $\mathcal{ F }_S^c \{C\}$ formulated in
Theorem~1.2 holds. Let $\Omega$ be an arbitrary neighborhood of $M
\backslash C$ in $\C^n$. Then there exists a special point $p_{\rm sp}
\in C$ and there exists a local one-sided neighborhood $\omega_{ p_{
\rm sp}}$ of $M$ in $\C^2$ 
at $p_{\rm sp}$ such that holomorphic functions in
$\Omega$ extend holomorphically to $\omega_{ p_{ \rm sp}}$.
\end{proposition}

\subsection*{2.3.~Holomorphic extension to a half-one-sided 
neighborhood of $M$} The choice of the special point $p_{\rm sp}$ will
be achieved in two main steps. According to the nontransversality
assumption $\mathcal{ F}_S^c\{C\}$, there exists a characteristic
segment $\gamma: [-1, 1] \to S$ with $\gamma(-1)\not\in C$, with
$\gamma(0)\in C$ and with $\gamma(1)\not \in C$ such that $C$ lies in
one (closed, semi-local) side of $\gamma$ in $S$. As $\gamma$ is a
Jordan arc, we may orient $S$ in $M$ along $\gamma$, hence we may
choose a semi-local open side $(S_\gamma)^+$ of $S$ in $M$ along
$\gamma$. In the first main step (to be conducted in Section~4 in the
context of the general codimensional case Theorem~1.2'), we shall
construct what we call a {\sl semi-local half-wedge $\mathcal{
HW}_\gamma^+$ attached to $(S_\gamma)^+$ along $\gamma$}. By this, we
mean the ``half part'' of a wedge attached to a neighborhood of the
characteristic segment $\gamma$ in $M$, which yields a wedge attached
to the semi-local one-sided neighborhood $(S_\gamma)^+$. For an
illustration, {\it see}\, {\sc Figure~6} below, in which one should
replace the notation $M^1$ by the notation $S$. Such a half-wedge may
be interpreted as a wedge attached to a neighborhood of $\gamma$ in
$S$ which is not arbitrary, but should satisfy a further property:
locally in a neighborhood of every point of $\gamma$, either the
half-wedge contains $(S_\gamma)^+$ or one of its two ribs contains
$(S_\gamma)^+$, as illustrated in {\sc Figure~6} below. Importantly
also, the cones of this attached half-wedge should vary continuously
as we move along $\gamma$, {\it cf.} again {\sc Figure~6}.

The way how we will construct this half-wedge $\mathcal{ HW
}_\gamma^+$ is as follows. As illustrated in {\sc Figure~1} just
below, we shall first construct a string of analytic discs $Z_{ r :s }
( \zeta)$, where $r$ is the approximate radius of $Z_{r: s}( \partial
\Delta)$, whose boundaries are contained in $(S_\gamma )^+ \subset M$
and which touch the curve $\gamma$ only at the point $\gamma (s)$, for
every $s \in [-1, 1]$, namely $Z_{r: s} (1) = \gamma (s)$ and $Z_{r
:s}\left(\partial \Delta \backslash \{ 1\} \right) \subset (S_\gamma
)^+$.

\bigskip
\begin{center}
\input disc-string.pstex_t
\end{center}

From now on, we fix a small radius $r_0$. By deforming the discs
$Z_{r_0:s}(\zeta)$ in $\Omega$ near their opposite points
$Z_{r_0:s}(-1)$, which lie at a positive distance from the singularity
$C$, we shall construct in Section~4 a family of analytic disc
$Z_{r_0,t:s}(\zeta)$, where $t\in \R$ is a small parameter, so that
the disc boundaries $Z_{r_0,t:s}(\partial \Delta)$ are pivoting
tangentially to $S$ at the point $\gamma(s)\equiv Z_{r_0,t:s}(1)$,
which remains fixed as $t$ varies. Precisely, we mean that
$\frac{\partial Z_{r_0,t:s}}{ \partial \theta}(1)\in T_{\gamma(s)}S$
and that the mapping $t\longmapsto \frac{\partial Z_{r_0,t:s}}{
\partial \theta}(1)$ is of rank $1$ at $t=0$. This construction and
the next ones will be achieved thanks to the solvability Bishop's
equation. Furthermore, we may add a small translation parameter $\chi
\in \R$, getting a family $Z_{r_0, t,\chi :s} (\zeta)$ with the
property that the mapping $(\chi, s) \longmapsto Z_{ r_0,t, \chi
:s}(1)\in S$ is a diffeomorphism onto a neighborhood of $\gamma( [-1,
1])$ in $S$, still with the property that the point $Z_{r_0,t,\chi
:s}(1)$ is fixed equal to the point $Z_{r_0, 0,\chi:s}(1)$ as $t$
varies. Finally, we may add a small translation parameter $\nu\in \R$
with $\nu >0$, getting a family $Z_{r_0,t,\chi,\nu:s}(\zeta)$ with
$Z_{r_0,t,\chi,0:s}(\zeta)\equiv Z_{r_0,t,\chi:s}(\zeta)$, such that
the mapping $(\chi,\nu,s)\longmapsto Z_{r_0,t,\chi,\nu:s}(1)$ is a
diffeomorphism onto the semi-local one-sided neighborhood
$(S_\gamma)^+$ of $S$ along $\gamma$ in $M$, provided $\nu >0$. Then
the semi-local attached half-wedge may be defined as
\def\theequation{2.4}\begin{equation}
\mathcal{ HW}_\gamma^+:=\left\{
Z_{r_0,t,\chi,\nu:s}(\rho): \
\vert t \vert < \varepsilon, \
\vert \chi \vert < \varepsilon, \ 
0< \nu < \varepsilon, \
1-\varepsilon < \rho < 1, \
-1\leq s \leq 1
\right\},
\end{equation}
for some small $\varepsilon >0$. In the first main technical step (to
be conducted in Section~4 below in the context of Theorem~1.2'), we
shall show that every holomorphic function $f \in \mathcal{ O} (
\Omega)$ extends holomorphically to $\mathcal{ HW }_\gamma^+$. To
prove Proposition~2.2, we shall find a special point $p_{\rm sp} \in
C$ such that there exists a local one-sided neighborhood $\omega_{
p_{\rm sp}}$ at $p_{\rm sp}$ such that holomorphic functions in
$\Omega \cup \mathcal{ HW }_\gamma^+$ extend holomorphically to
$\omega_{ p_{\rm sp}}$.

\subsection*{2.5.~Field of cones on $S$}
We continue the description of the proof of Theorem~1.2 with the full
family of analytic discs $Z_{ r_0, t, \chi, \nu:s} (\zeta)$. Thanks to
a technical application of the implicit function theorem, we can
arrange from the beginning that the vectors $\frac{ \partial Z_{r_0,
t, \chi,0: s} }{ \partial \theta } (1)$ are tangent to $S$ at the
point $Z_{ r_0, 0,\chi,0:s} (1) \in S$ when $t$ varies, for all fixed
$s$. Then by construction, the disc boundaries
$Z_{r_0, t, \chi, 0:s}( \partial \Delta)$ are pivoting tangentially to
$S$ at the point $Z_{ r_0, t, \chi, 0 :s} (1) \equiv Z_{r_0, 0, \chi,
0:s} (1)$. It follows that when $t$ varies, the oriented half-lines
$\R^+ \cdot \frac{ \partial Z_{ r_0,t,\chi,0 : s} }{ \partial \theta}
(1)$ describe an open infinite oriented cone in the tangent space to
$S$ at the point $Z_{r_0, 0, \chi, 0:s}(1)$. Consequently, we may
define a {\sl field of cones} $p \mapsto {\sf C}_p$ as
\def\theequation{2.6}\begin{equation}
{\sf C}_p:= \left\{ \R^+ \cdot 
\frac{ \partial Z_{r_0,t,\chi,0: s} }{
\partial \theta} (1) : \ \vert
t \vert < \varepsilon \right\},
\end{equation}
at every point $p=Z_{ r_0, 0, \chi, 0: s} (1) \in S$ of a neighborhood
of $\gamma$ in $S$. The following figure provides an intuitive
illustration. One should think that the small cones are generated when
the small discs boundaries of {\sc Figure~1} pivote tangentially to
$S$.

\bigskip
\begin{center}
\input field-cones-hypersurface.pstex_t
\end{center}

After having defined this field of cones, we shall {\sl fill} the
cones as follows. Remind that a neighborhood of $\gamma$ in $S$ is
foliated by characteristic segments, which are approximatively
parallel to $\gamma$. In {\sc Figure~2} above, one should think that
the characteristic foliation is horizontal. So there exists a nowhere
vanishing vector field $p \mapsto X_p$ defined in a neighborhood of
$\gamma$ whose integral curves are characteristic segments. We define
the {\sl filled cone} ${\sf FC}_p$ by
\def\theequation{2.7}\begin{equation}
{\sf FC}_p:=\left\{
\lambda \cdot X_p + (1-\lambda) \cdot v_p : \ 
0 \leq \lambda < 1, \
v_p\in {\sf C}_p
\right\}.
\end{equation}
Geometrically, we rotate every half-line $\R^+ \cdot v_p$ towards the
characteristic half-line $\R^+ \cdot X_p$ and we call the result the
{\sl filling} of ${\sf C}_p$. In {\sc Figure~2} above, all the cones
${\sf C}_p$ coincide with their fillings. Thus we have constructed a
field of filled cones $p \longmapsto {\sf FC}_p$ over a neighborhood
of $\gamma$ in $S$.

\subsection*{2.8.~Small analytic discs half-attached to $S$}
The next main observation is that small analytic discs which are
half-attached to $S$ are essentially contained in the half-wedge
$\mathcal{ HW }_\gamma^+$, provided that they are approximatively
directed by the filled cone ${\sf FC}_p$. Let us be more precise. Let
$\partial^+ \Delta:= \{\zeta \in \partial \Delta: \ {\rm Re} \, \zeta
\geq 0 \}$ denote the {\sl positive half part} of the unit circle. We
say that an analytic disc $A: \overline{ \Delta } \to \C^2$ is {\sl
half-attached} to $S$ if $A( \partial^+ \Delta)$ is contained in
$S$. Here, $A$ is at least of class $\mathcal{ C}^1$ over
$\overline{\Delta}$ and holomorphic in $\Delta$. In addition, we shall
always assume that our discs $A$ are embeddings of $\overline{\Delta}$
into $\C^2$. We shall say that $A$ is {\sl approximatively straight}
(in an informal sense) if $A(\Delta)$ is close in $\mathcal{
C}^1$-norm to an open subset of the complex line generated by the
complex vector $\frac{\partial A}{ \partial \zeta}(1)$. Finally, we
say that $A$ is {\sl approximatively directed by the filled cone ${\sf
FC}_p$ at $p=A(1)$}, if the vector $\frac{\partial A}{\partial
\theta}(1)\in T_pS$ belongs to ${\sf FC}_p$. Although this
terminology will not be re-employed in the next sections, we may
formulate a crucial geometric observation as follows.

\def\thelemma{2.9}\begin{lemma}
A sufficiently small approximatively straight analytic disc $A:
\overline{ \Delta} \to \C^2$ of class at least $\mathcal{ C}^1$ which
is half-attached to $S$ and which is approximatively directed by the
filled cone ${\sf FC}_p$ at $p=A(1)\in S$, necessarily satisfies
\def\theequation{2.10}\begin{equation}
A\left(\overline{ \Delta } \backslash \partial^+ \Delta\right) 
\subset \mathcal{ HW}_\gamma^+.
\end{equation}
\end{lemma}

In the context of the general Theorem~1.2',
this property (with more precisions) will be checked in 
Section~8 below. 

\subsection*{2.11.~Choice of a special point}
In the second main step of the proof (to be conducted in Section~5 in
the context of Theorem~1.2'), we shall choose the desired special
point $p_{ \rm sp}$ of Proposition~2.2 to be removed locally as
follows. Since we shall remove $p_{ \rm sp}$ by means of half-attached
analytic discs (applying the continuity principle), we want to find a
special point $p_{ \rm sp}\in C$ so that the following two conditions
hold true:
\begin{itemize}
\item[{\bf (i)}]
There exists a small approximatively straight analytic disc $A:
\overline{\Delta} \to \C^2$ with $A(1)= p_{\rm sp}$ which is
half-attached to $S$ such that $A$ is approximatively directed by the
filled cone ${\sf FC}_{ p_{ \rm sp}}$ (so that the conclusion of
Lemma~2.9 above holds true).
\item[{\bf (ii)}]
The same disc satisfies $A \left( \partial^+ \Delta \backslash \{1\}
\right) \subset S \backslash C$.
\end{itemize}

In particular, since $M \backslash C$ is contained in $\Omega$, it
follows from these two conditions that the blunt disc boundary
$A\left( \partial \Delta \backslash \{ 1 \} \right)$ is contained in
the open subset $\Omega \cup \mathcal{ HW }_\gamma^+$, a property that
will be appropriate for the application of the continuity principle,
as we shall explain in Section~9 below.

To fulfill conditions {\bf (i)} and {\bf (ii)} above, we first
construct a supporting real segment at a special point of the nonempty
closed subset $C\subset S$.

\def\thelemma{2.12}\begin{lemma}
There exists at least one special point $p_{ \rm sp } \in C$
arbitrarily close to $\gamma$ in a neighborhood of which the following
two properties hold true{\rm :}
\begin{itemize}
\item[{\bf (i')}]
There exists a small $\mathcal{ C}^{2, \alpha}$-smooth open segment
$H_{p_{\rm sp}}\subset S$ passing through $p_{ \rm sp}$ such that an
oriented tangent half-line to $H_{p_{\rm sp}}$ at $p_{\rm sp}$ is
contained in the filled cone ${\sf FC}_{p_{\rm sp}}$, as illustrated
in {\sc Figure~3} below.
\item[{\bf (ii')}]
The same segment is a supporting segment in the following
sense{\rm :} locally in a
neighborhood of $p_{\rm sp}$, the set $C \backslash \{p_{\rm sp}\}$ is
contained in one open side $(H_{ p_{\rm sp}})^-$ if $H_{p_{\rm sp}}$
in $S$, as illustrated in {\sc Figure~3} just below.
\end{itemize}
\end{lemma}

\bigskip
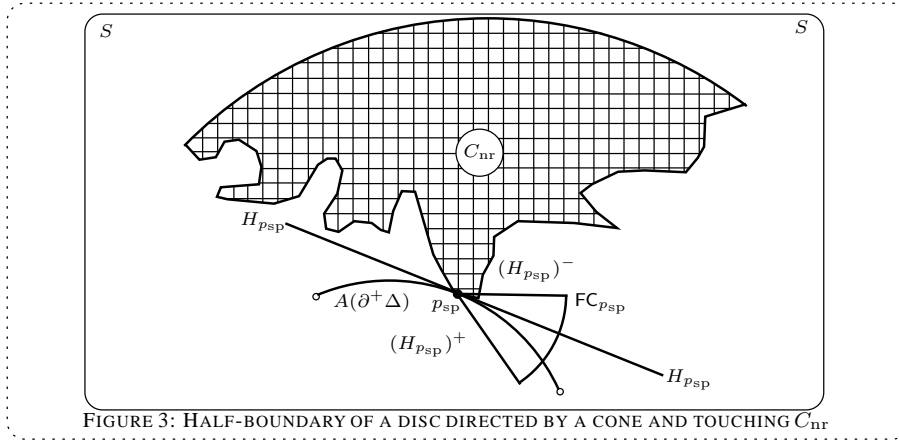
\begin{center}
\input one-touch.pstex_t
\end{center}

The way how we prove Lemma~2.12 is illustrated intuitively in {\sc
Figure~2} above. For $\lambda\in \R$ with $0 \leq \lambda < 1$ very
close to $1$, the vector field $p \longmapsto v_p^\lambda := \lambda
\cdot X_p + (1-\lambda) \cdot v_p$ is very close to the characteristic
vector field $p\mapsto X_p$. By construction, this vector field runs
into the filled field of cones $p\mapsto {\sf FC}_p$. In {\sc
Figure~2}, the integral curves of $p\mapsto v_p^\lambda$ are drawn as
dotted lines, which are almost horizontal if $\lambda$ is very close
to $1$. If we choose the first dotted integral curve from the lower
part of {\sc Figure~2} which touches $C$ at one special point $p_{\rm
sp}\in C$ and if we choose for $H_{p_{\rm sp}}$ a small segment of
this first dotted integral curve, we may check that properties {\bf
(i)} and {\bf (ii)} are satisfied, modulo some mild technicalities. A
rigorous complete proof of Lemma~2.12 will be provided in Section~5
below.

\subsection*{2.13.~Construction of analytic discs
half-attached to $S$} Small analytic discs which are half-attached to
a $\mathcal{ C}^{2,\alpha}$-smooth maximally real submanifold $M^1$ of
$\C^n$ and which are approximatively straight will be constructed in
Section~7 below. For this, we shall use the solution of Bishop's
equation with parameters in H\"older spaces, obtained by A.~Tumanov
in~\cite{ tu3} with an optimal loss of regularity. In Section~8, we
shall check that it follows from the general constructions of
Section~7 that there exists a small analytic disc $A$ half-attached to
$S$ with $A(1)= p_{\rm sp}$ whose half boundary is tangent to
$H_{p_{\rm sp}}$ at $p_{\rm sp}$ and which satisfies property {\bf
(ii)} above, as drawn in {\sc Figure~3} above. Thus, the two geometric
properties {\bf (i')} and {\bf (ii')} satisfied by the real segment
$H_{p_{\rm sp}}$ may be realized by the half-boundary of a
half-attached analytic disc.

\subsection*{2.14.~Translation of half-attached and continuity 
principle} By means of the results of Section~7, we shall see that we
may include the disc $A (\zeta)$ is a parametrized family $A_{ x,v}
(\zeta)$ of analytic discs half-attached to $S$, where $x\in \R^2$ and
$v \in \R$ are small, so that the mapping $x \mapsto A_{ x,0} (1) \in
S$ is a local diffeomorphism onto a neighborhood of $p_{ \rm sp}$ in
$S$ and so that the mapping $v \mapsto \frac{ \partial A_{ 0,v }}{
\partial \theta} (1)$ is of rank $1$ at $v=0$. Furthermore, we
introduce a new parameter $u \in \R$ in order to ``translate'' the
totally real surface $S$ in $M$ by means of a family of $S_u \subset
M$ with $S_0=S$ and $S_u \subset ( S_\gamma )^+$ for $u>0$. Thanks to
the tools developed in Section~7, we deduce that there exists a
deformed family of analytic discs $A_{ x, v, u} (\zeta)$ which are
half-attached to $S_u$ and which satisfies $A_{x, v,0} (\zeta) \equiv
A_{ x,v} (\zeta)$. In particular, this family covers a local
one-sided neighborhood $\omega_{p_{\rm sp}}$ of $M$ at $p_{\rm sp}$
defined by
\def\theequation{2.15}\begin{equation}
\omega_{p_{\rm sp}}:= \left\{
A_{x,v,u}(\rho): \ 
\vert x \vert < \varepsilon, \ 
\vert v \vert < \varepsilon, \
\vert u \vert < \varepsilon, \ 
1-\varepsilon < \rho < 1
\right\},
\end{equation}
for some $\varepsilon>0$.

In the third and last main step of the proof (to be conducted in
Section~9 below), we shall prove that every disc $A_{x, v,u} (\zeta)$
with $u \neq 0$ is analytically isotopic to a point with the boundary
of every disc of the isotopy being contained in $\Omega \cup \mathcal{
HW}_\gamma^+$. Thanks to the continuity principle, we shall deduce
that every holomorphic function $f \in \mathcal{ O}\left( \Omega \cup
\mathcal{ HW}_\gamma^+ \right)$ extends holomorphically to $\omega_{
p_{\rm sp}}$ minus a certain thin closed subset $\mathcal{ C}_{p_{\rm
sp}}$ of $\omega_{ p_{\rm sp}}$. Finally, we shall conclude both the
proof of Proposition~2.2 and the proof of Theorem~1.2 by checking that
the thin closed set $\mathcal{ C}_{p_{\rm sp}}$ is in fact removable
for holomorphic functions defined in $\omega_{ p_{\rm sp }} \backslash
\mathcal{ C}_{ p_{ \rm sp}}$.

\subsection*{2.16.~Organization of the paper} 
As was already announced, Sections~3, 4, 5, 6, 7, 8 and~9 below will be
entirely devoted to the proof of Theorem~1.2', which will be
endeavoured directly in arbitrary codimension, without any further
reference to the hypersurface version. Only in Section~3 shall we also
consider the beginning of the proof of Theorem~1.4. During the
development of the proof of Theorem~1.2', in comparison to the quick
description of the proof of Theorem~1.2 achieved just above, we shall
unavoidably encounter some supplementary technical complications caused
by the codimension being 
$\geq 2$, namely technicalities which are absent in
codimension $1$. We would like to mention that the crucial geometric
argument which enables us to choose the desired special point will be
conducted in the central Section~5 below.

Then Section~10 is devoted to summarize three geometrically distinct
proofs of Theorem~1.4. In Section~11, we check that both the CR- and
the $L^{\sf p}$-removability of $C$ are a consequence of the
$\mathcal{ W}$-removability of $C$. In Section~12, we
provide the proof of Theorem~1.1, of Theorem~1.3 and of further
applications. This Section~12 may be read before entering the proof of
Theorem~1.2'. Finally, in Section~13, we provide some applications of
our removability results to the edge of the wedge theorem for
meromorphic functions.

\section*{\S3.~Strategy per absurdum
for the proofs of Theorems~1.2' and~1.4} 

\subsection*{3.1.~Preliminary}
For the proof of Theorems~1.2', as in \cite{cs}, \cite{m2},
\cite{mp1}, \cite{mp3}, \cite{p1}, we shall proceed by contradiction.
This strategy possesses a considerable advantage: it will not be
necessary to control the size of the local subsets of $C$ that are
progressively removed, which simplifies substantially the presentation
and the understandability of the reasonings. We shall explain how to
reduce CR- and $L^{\sf p}$-removability of $C$ to its $\mathcal{
W}$-removability. Also, it may be argued that the $\mathcal{
W}$-removability of $C$ can be reduced to the simpler case where the
functions which we have to extend are even holomorphic in a
neighborhood of $M\backslash C$ in $\C^n$. Whereas such a strategy is
essentially carried out in detail in previous references (with some
variations), we shall for completeness recall the complete reasonings
briefly here, in \S3.2 and in \S3.16 below.

\subsection*{3.2.~Global minimality of some complements} 
For background notions about CR orbits in CR manifolds, we refer the
reader to~\cite{su}, \cite{j2}, \cite{mp1}, \cite{j4}. We just recall
a few standard facts: if $p$ belongs to a generic submanifold of
$\C^n$ of class at least $\mathcal{ C}^2$, a point $q\in M$ belongs to
the CR orbit $\mathcal{ O}_{CR}(M,p)$ if and only if there exists a
piecewise smooth curve $\lambda: [0,1] \to M$ with $\lambda(0)=p$,
$\lambda(1)=q$ such that $d\lambda(s)/ds\in T_{\gamma(s)}^cM
\backslash \{0\}$ at every $s\in [0,1]$ at which $\lambda$ is
differentiable; CR orbits make a partition of $M$; CR orbits are
immersed $\mathcal{ C}^{1,\alpha}$-smooth submanifolds of $M$,
according to H.J.~Sussmann's Theorem~4.1 in~\cite{su} specialized in
the CR category; Every maximal $T^cM$-integral immersed submanifold of
$M$ must contain the CR orbit of each of its points; and finally, a
trivial, but often useful fact: if $N$ is a $T^cM$-integral
submanifold of $M$, namely $T_pN \subset T_p^cM$ for every point $p\in
N$, then the local flow of every $T^cM$-tangent vector field on $M$
stabilizes locally $N$.

In the two geometric situations of Theorems~1.2' and~1.4, we shall
apply the following two ~Lemmas~3.3 and~3.5 about the CR structure of
the complement $M\backslash C'$, where $C'\subset C\subset M^1$ is an
arbitrary proper closed subset of $C$.

\def\thelemma{3.3}\begin{lemma} Let $M$ be a $\mathcal{ C}^{2,
\alpha}$-smooth generic submanifold of $\C^n$ $(n\geq 2)$ of
codimension $d\geq 1$ and of CR dimension $m:=n-d\geq 1$, let $M^1$ be
a $\mathcal{ C}^{2,\alpha}$-smooth 
one-codimensional submanifold of $M$ which is generic in $\C^n$ and
let $C'$ be an arbitrary proper closed subset of $M^1$. If either
\begin{itemize}
\item[{\bf (1)}]
$M$ is of CR dimension $m=1$ and the condition 
$\mathcal{ F}_{M^1}^c\{C'\}$ of Theorem~1.2' holds; or if
\item[{\bf (2)}]
$M$ is of CR dimension $m\geq 2$ and the condition 
$\mathcal{ O}_{M^1 }^{CR }\{C'\}$ of Theorem~1.4 holds,
\end{itemize}
then for each point $p'\in C'$, there exists a piecewise $\mathcal{
C}^{2, \alpha}$-smooth curve $\gamma: [0,1]\to M^1$ satisfying $d
\gamma(s)/ds\in T_{\gamma(s) }M^1 \cap T_{\gamma(s)}^c M\backslash
\{0\}$ at every $s\in [0,1]$ at which $\lambda$ is differentiable,
such that $\gamma(0)=p'$ and $\gamma(1)$ does not belong to $C'$.
\end{lemma}

\proof
In the case $m=1$, we proceed by contradiction and we suppose that
there exists a point $p'\in C'$ such that all piecewise $\mathcal{
C}^{2,\alpha}$-smooth curves $\gamma: [0,1]\to M^1$ with $d\gamma(s)
/ds\in T_{\gamma(s)}M^1\cap T_{\gamma(s)}^cM\backslash \{0\}$ at every
$s\in [0,1]$ at which $\gamma$ is differentiable, which have origin
$p'$ are entirely contained in $C'$. Notice that since the bundle
$TM^1\cap T^cM\vert_{M^1}$ is of real dimension one and of class
$\mathcal{ C}^{1,\alpha}$, such curves $\gamma$ are in fact $\mathcal{
C}^{2,\alpha}$-smooth at every point. It follows immediately there
cannot exist a curve $\gamma: [-1,1] \to M^1$ contained in a single
leaf of the characteristic foliation $\mathcal{ F}_{M^1}^c$ with
$\gamma(0)\in C'$ and $\gamma(-1), \gamma(1)\not \in C'$, which
contradicts the condition $\mathcal{ F}_{M^1}^c\{C'\}$.

In the case $m\geq 2$, by genericity of $M^1$, the complex tangent
bundle $T^cM^1$, which is of real dimension $(2m-2)$, is a
one-codimensional subbundle of the $(2m-1)$-dimensional bundle
$TM^1\cap T^cM\vert_{M^1}$, namely $T^cM^1\subset TM^1 \cap
T^cM\vert_{M^1}$. Let $p\in C'$ be an arbitrary point. By the
assumption $\mathcal{ O}_{M^1}^{CR}\{C'\}$, the CR orbit of $p'$ is
not contained in $C'$. Equivalently, there exists a piecewise
$\mathcal{ C}^{2,\alpha}$-smooth curve $\gamma: [0,1] \to M^1$
satisfying
\def\theequation{3.4}\begin{equation}
d\gamma(s)/ds\in T_{\gamma(s)}^cM^1\backslash \{0\}\subset
T_{\gamma(s)}M^1 \cap T_{\gamma(s)}^cM\backslash \{0\} 
\end{equation} 
at each $s\in [0,1]$ at which $\gamma$ is
differentiable,
such that $\gamma(0)=p'$ and $\gamma(1)$ does not belong to $C'$.
Hence the conclusion of Lemma~3.3 is immediately satisfied.
This completes the proof.
\endproof

As an application, we deduce that under the respective assumptions
$\mathcal{ F}_{M^1}^c\{C\}$ and $\mathcal{ O}_{M^1}^{CR}\{C\}$ of
Theorems~1.2' and of Theorem~1.4, the complement $M \backslash C'$ is
globally minimal, for every closed subset $C'\subset C$.

\def\thelemma{3.5}\begin{lemma} 
Let $M$ be a $\mathcal{ C}^{2, \alpha}$-smooth generic submanifold of
$\C^n$ $(n \geq 2)$ of codimension $d\geq 1$ and of CR dimension
$m:=n-d \geq 1$, let $M^1$ be a $\mathcal{ C}^{2,\alpha}$-smooth
one-codimensional submanifold of $M$ which is generic in $\C^n$ and
let $C'$ be an arbitrary {\rm nonempty} proper closed subset of
$M^1$. Assume that for each point $q'\in C'$, there exists a piecewise
$\mathcal{ C}^{2,\alpha}$-smooth curve $\gamma: [0,1] \to M^1$ with
$d\gamma(s) / ds\in T_{ \gamma(s)}M^1\cap T_{ \gamma(s)}^cM\backslash
\{0\}$ at every $s\in [0,1]$ at which 
$\gamma$ is differentiable, such that $\gamma(0)=q'$ and $\gamma(1)$
does not belong to $C'$. Then the CR orbit in $M\backslash C'$ of
every point $p\in M\backslash C'$ coincides with its CR orbit in $M$
minus $C'$, namely
\def\theequation{3.6}\begin{equation}
\mathcal{ O}_{CR}(M\backslash C', p)= 
\mathcal{ O}_{CR}(M,p)\backslash C'.
\end{equation}
In particular, as a corollary, if $M$ is globally minimal, then
$M\backslash C'$ is also globally minimal.
\end{lemma}

\proof
Let us first explain the last sentence, which applies
to the situations considered in both Theorems~1.2' and~1.4: 
if $\mathcal{ O}_{CR}(M,p)=M$, then by~\thetag{3.6}, 
$\mathcal{ O}_{CR}(M\backslash C',p)\equiv
M\backslash C'$, which proves that
$M\backslash C'$ is globally minimal.

To establish~\thetag{3.6}, we shall need the following crucial lemma,
deserving an illustration: {\sc Figure~4} below.

\def\thelemma{3.7}\begin{lemma}
Under the assumptions of Lemma~3.5, for every point $q'\in C' \subset
M^1$, there exists a $\mathcal{ C}^{1,\alpha}$-smooth locally embedded
submanifold $\Omega_{q'}$ of $M$ passing through $q'$ which is
transverse to $M^1$ at $q'$ in $M$, namely which satisfies
$T_{q'}\Omega_{q'}+T_{q'}M^1= T_{q'}M$, such that
\begin{itemize}
\item[{\bf (1)}]
$\Omega_{q'}$ is a $T^cM$-integral submanifold, namely 
$T_p^cM\subset T_p\Omega_{q'}$, for every point 
$p\in \Omega_{q'}$.
\item[{\bf (2)}]
$\Omega_{q'}\backslash C'$ is contained in 
a single CR orbit of $M$.
\item[{\bf (3)}] $\Omega_{q'}$ is
also contained in a single CR orbit of $M\backslash C'$.
\end{itemize}
\end{lemma}

\proof
So, let $q'\in C'\subset M^1$. Since $M^1$ is generic in $\C^n$, there
exists a $\mathcal{ C}^{1,\alpha}$-smooth vector field $Y$ defined in
a neighborhood of $q'$ which is complex tangent to $M$ but locally
transversal to $M^1$, {\it cf.} {\sc Figure~4} just below (for easier
readability, we have erased the hatching of $C'$ in a neighborhood of
$q'$).

\bigskip
\begin{center}
\input orbit-complement.pstex_t
\end{center}

\noindent
Following the integral curve of $Y$ issued from $q'$, we can define a
point $q_\epsilon'$ in an $\epsilon$-neighborhood of $q'$ which does
not belong to $M^1$. By assumption, there exists a piecewise smooth
curve $\gamma: [0,1]\to M^1$ with $d\gamma(s) /ds \in T_{\gamma(s)}M^1
\cap T_{\gamma(s)}^cM\backslash \{0\}$ at every $s\in [0,1]$ at which
$\gamma$ is differentiable, such that $\gamma(0)=q'$ and $\gamma(1)$
does not belong to $C'$. For simplicity, we shall assume that
$\gamma$ consists of a single smooth piece, the case where $\gamma$
consists of finitely many smooth pieces being treated in a completely
similar way. With this assumption (which will simplify slightly the
technicalities), it follows that there exists a vector field $X$
defined in a neighborhood of the curve $\gamma([0,1])$ in $M$ which is
complex tangent to $M$ and whose restriction to $M^1$ is a semi-local
section of $TM^1\cap T^cM\vert_{M^1}$, such that $\gamma$ is an
integral curve of $X$ and such that $\gamma(1)=\exp(X)(q')\in
M^1\backslash C'$. In addition, we can assume that the vector field
$Y$ is defined in the same neighborhood of $\gamma([0,1])$ in $M$ and
everywhere transversal to $M^1$. If $\epsilon$ is sufficiently small,
{\it i.e.} if $q_\epsilon'$ is sufficiently close to $q'$, the point
$r_\epsilon':= \exp(X)(q_\epsilon')$ is still very close to
$M^1$. Thus, we can define a new point $r'\in M^1$ to be the unique
intersection of the integral of $Y$ issued from $r_\epsilon'$ with
$M^1$. By choosing $\epsilon$ small enough, the point $r_\epsilon'$
will be arbitrarily close to $\gamma(1)\not \in C'$, and consequently,
we can assume that $r'$ also does not belong to $C'$, as drawn in {\sc
Figure~4} above. Notice that the integral curve of $X$ from
$q_\epsilon'$ to $r_\epsilon'$ is contained in $M\backslash M^1$,
since the flow of $X$ stabilizes $M^1$, whence the two points
$r_\epsilon'$ and $r'$ belong to the CR orbit $\mathcal{
O}_{CR}(M\backslash C',q_\epsilon')$.

Let $\Omega_{r'}$ denote a small piece of the immersed submanifold
$\mathcal{ O}_{CR}(M\backslash C',r')$ passing through $r'$. By the
standard properties of CR orbits, we can assume that $\Omega_{r'}$ is
an embedded $\mathcal{ C}^{1,\alpha}$-smooth submanifold of
$M\backslash C'$ of the same CR dimension as $M\backslash C'$ and we
can in addition assume that $\Omega_{r'}$ contains $r_\epsilon'$,
provided $\epsilon$ is small enough. Since $Y$ is a vector field
complex tangent to $M$, the submanifold $\Omega_{r'}$ is necessarily
stretched along the flow lines of $Y$, hence it is transversal to
$M^1$.

We then define the submanifold $\exp(-X)(\Omega_{r'})$, close to the
point $q'$ (we shall argue in a while that it passes in fact through
$q'$). Since the flow of $X$ stabilizes $M^1$, it follows that
$\exp(-X)(\Omega_{r'})$ is transversal to $M^1$ and that
$\exp(-X)(\Omega_{r'})$ is divided in two sides by its
one-codimensional $\mathcal{ C}^{1,\alpha}$-smooth submanifold
$M^1\cap \exp(-X)(\Omega_{r'})$. Furthermore, the flow of $X$
stabilizes the two sides of $M^1$ in $M$, semi-locally in a
neighborhood of $\gamma([0,1])$, {\it see} again {\sc Figure~4} above.
Consequently, every integral curve of $X$ issued from every point in
$\Omega_{r'}\backslash M^1$ stays in $M\backslash M^1$, hence in
$M\backslash C'$ and it follows that the submanifold
\def\theequation{3.8}\begin{equation}
\exp(-X)(\Omega_{r'})\backslash M^1,
\end{equation}
consisting of two connected pieces, is contained in the single CR
orbit $\mathcal{ O}_{CR}(M\backslash C', p')$. By the characteristic
property of a CR orbit, this means that the two connected pieces of
$\exp(-X)(\Omega_{r'})\backslash M^1$ are CR submanifolds of
$M\backslash C'$ of the same CR dimension as $M\backslash C'$.
Furthermore, since the intersection $M^1\cap \exp(-X)(\Omega_{r'})$ is
one-codimensional, it follows by continuity that {\it the $\mathcal{
C}^{1,\alpha}$-smooth submanifold $\exp(-X) (\Omega_{r'})$ is in fact a
CR submanifold of $M$ of the same CR dimension as $M$}. This 
proves property {\bf (1)}. 

Since $q_\epsilon'$ belongs to $\exp (-X) (\Omega_{r'})$ and since the
flow of the complex tangent vector field $Y$ necessarily stabilizes
the $T^cM$-integral submanifold $\exp( -X) (\Omega_{ r'})$, the point
$q'$ which belongs to an integral curve of $Y$ issued from
$q_\epsilon'$, must belong to the submanifold $\exp(-X)( \Omega_{
r'})$, which we can now denote by $\Omega_{ q'}$, as in {\sc
Figure~4} above.

Observe that locally in a neighborhood of $q'$, the integral curves of
$Y$ are transversal to $M^1$ and meet $M^1$ only at one point.
Shrinking if necessary $\Omega_{q'}$ a little bit and using positively
or negatively oriented integral curves of $Y$ with origin all points
in $\Omega_{q'} \backslash M^1$ lying in both sides, we deduce that
$\Omega_{q'}\backslash C'$ is contained in the single CR orbit
$\mathcal{ O}_{CR}(M\backslash C', r')$, which proves property {\bf
(3)}. Using again $Y$ to join points of $C'\cap \Omega_{q'}$, we
deduce also that $\Omega_{q'}$ is contained in the single CR orbit
$\mathcal{ O}_{CR}(M,r')$, which proves property {\bf (2)}.

The proof of Lemma~3.7 is complete.
\endproof

We can now prove Lemma~3.5. It suffices to establish that for every two
points $p \in M \backslash C'$ and $q \in \mathcal{ O}_{ CR} (M, p)$
with $q\not\in C'$, the point $q$ belongs in fact to the CR orbit of
$p$ in $M\backslash C'$, namely $q$ belongs to 
$\mathcal{ O}_{CR}(M\backslash C', p)$.

Since $q$ belongs to the CR orbit of $p$ in $M$, there exists a
piecewise $\mathcal{ C}^{2,\alpha}$-smooth curve $\lambda: [0,1] \to
M$ with $\lambda(0)=p$, $\lambda(1)=q$ and $d\lambda(s)/ds\in
T_{\lambda(s)}^cM\backslash \{0\}$ at every $s\in [0,1]$ at which
$\lambda$ is differentiable. For every $s$ with $0\leq s \leq 1$, we
define a local $\mathcal{ C}^{1,\alpha}$-smooth submanifold
$\Omega_{\lambda(s)}$ of $M$ passing through $\lambda(s)$ as follows:
\begin{itemize}
\item[{\bf (1)}]
If $\lambda(s)$ does not belong to $C'$, choose for $\Omega_{
\lambda(s)}$
a piece of the CR orbit of $\lambda(s)$ in $M\backslash C'$.
\item[{\bf (2)}]
If $\lambda(s)$ belongs to $C'$, choose for $\Omega_{\lambda(s)}$
the submanifold constructed in Lemma~3.7 above.
\end{itemize}
By Lemma~3.7, for each $s$, the complement
$\Omega_{\lambda(s)}\backslash C'$ is contained in a single CR orbit
of $M\backslash C'$. Since each $\Omega_{\lambda(s)}$ is a
$T^cM$-integral submanifold, for each $s\in [0,1]$, a neighborhood of
$\lambda(s)$ in the arc $\lambda([0,1])$ is necessarily contained
$\Omega_{\gamma(s)}$. By the Borel-Lebesgue covering lemma, we can
therefore find an integer $k\geq 1$ and real numbers
\def\theequation{3.9}\begin{equation}
0=s_1<r_1<t_1<s_2<r_2<t_2<\cdots\cdots <
s_{k-1} < r_{k-1} < t_{k-1} < s_k=1,
\end{equation}
such that
$\lambda([0,1])$ is covered by $\Omega_{\lambda(0)} \cup
\Omega_{\lambda(s_2)} \cup \cdots \cup \Omega_{\lambda(s_{k-1})}\cup
\Omega_{\lambda(1)}$ and such that in addition, 
$\lambda([r_j,t_j])\subset \Omega_{\lambda(s_j)}\cap \,
\Omega_{\lambda(s_{j+1})}$ for $j=1,\dots,k-1$.

\def\thelemma{3.10}\begin{lemma}
The following union minus $C'$
\def\theequation{3.11}\begin{equation}
\left(
\Omega_{\lambda(0)} \cup \Omega_{\lambda(s_2)}\cup \cdots\cdots
\cup \Omega_{\lambda(s_{k-1})}\cup
\Omega_{\lambda(1)}\right)\backslash C'
\end{equation}
is contained in a single CR orbit of $M\backslash C'$.
\end{lemma}

\proof
It suffices to prove that for every $j=1,\dots,k-1$, the union $\left(
\Omega_{\lambda(s_j)}\cup \Omega_{\lambda(s_{j+1})}\right) \backslash
C'$ minus $C'$ is contained in a single CR orbit of $M\backslash C'$.

Two cases are to be considered. Firstly, assume that
$\lambda([r_j,t_j])$ is not contained in $C'$, namely there exists
$u_j$ with $r_j\leq u_j\leq t_j$ such that
\def\theequation{3.12}\begin{equation}
\gamma(u_j)\in \left(\Omega_{\lambda(s_j)}\cap 
\Omega_{\lambda(s_{j+1})}\right)\backslash C'.
\end{equation}
Because $\Omega_{\lambda(s_j)}\backslash C'$ and
$\Omega_{\lambda(s_{j+1})} \backslash C'$ are both contained in a
single CR orbit of $M\backslash C'$, it follows from~\thetag{3.12}
that they are contained in the same CR orbit of $M\backslash C'$, as
desired.

Secondly, assume that $\lambda([r_j,t_j])$ is contained in
$C'$. Choose $u_j$ arbitrary with $r_j\leq u_j \leq t_j$. By
construction, $\lambda(u_j)$ belongs to $\Omega_{\lambda(s_j)}\cap
\Omega_{\lambda(s_{j+1})}$ and both $\Omega_{\lambda(s_j)}$ and
$\Omega_{\lambda(s_{j+1})}$ are $T^cM$-integral submanifolds of $M$
passing through the point $\lambda(u_j)$. Let $Y$ be a local section
of $T^cM$ defined in a neighborhood of $\lambda(u_j)$ which is not
tangent to $M^1$ at $\lambda(u_j)$. On the integral curve of $Y$
issued from $\lambda(u_j)$, we can choose a point
$\lambda(u_j)_\epsilon$ arbitrarily close to $\lambda(u_j)$ which does
not belong to $C'$. Since $Y$ is a section of $T^cM$, it is tangent to
both $\Omega_{\lambda(s_j)}$ and $\Omega_{\lambda(s_{j+1})}$, hence
we deduce that
\def\theequation{3.13}\begin{equation}
\gamma(u_j)_\epsilon\in \left(\Omega_{\lambda(s_j)}\cap 
\Omega_{\lambda(s_{j+1})}\right)\backslash C'.
\end{equation}
Consequently, as in the first case, it follows that
$\Omega_{\lambda(s_j)}\backslash C'$ and $\Omega_{\lambda(s_{j+1})}
\backslash C'$ are both contained in the same CR orbit of $M\backslash
C'$, as desired. This completes the proof of Lemma~3.10.
\endproof

Since $p$ and $q$ belong to the set~\thetag{3.11}, we
deduce that the points $p=\lambda(0)\in M\backslash C'$ and
the point $q=\lambda(1)\in \mathcal{ O}_{CR}(M,p)\backslash C'$
belong to the same CR orbit of $M\backslash C'$, as desired. 
This completes the proof of Lemma~3.5.
\endproof

\subsection*{3.14.~Reduction of CR- and of $L^{\sf p}$-removability 
to $\mathcal{ W}$-removability} First of all, we remind that it
follows from successive efforts of numerous mathematicians ({\it
cf.}~\cite{a}, \cite{j2}, \cite{m1}, \cite{tr2}, \cite{tu1},
\cite{tu2}, \cite{tu3}) that for every $\mathcal{
C}^{2,\alpha}$-smooth globally minimal submanifold $M'$ of $\C^n$,
there exists a wedge $\mathcal{ W}'$ attached to $M$ to which all
continuous CR functions on $M$ extend holomorphically. It follows that
the CR-removability of the closed subset $C\subset M^1$ claimed in
Theorems~1.2' and~1.4 is an immediate consequence of its $\mathcal{
W}$-removability. Based on the construction of analytic discs
half-attached to $M^1$ which will be achieved in Section~7 below, we
shall also be able to settle the reduction of $L^{\sf p}$-removability
in the end of the paper, and we formulate a convenient lemma, whose
proof is postponed to \S11 below.

\def\thelemma{3.15}\begin{lemma}
Under the assumptions of Theorem~1.2' and of Theorem~1.4, if the closed subset
$C\subset M^1$ is $\mathcal{ W}$-removable, then it is
$L^{\sf p}$-removable, for all ${\sf p}$ with 
$1\leq {\sf p}\leq \infty$.
\end{lemma}

\subsection*{3.16.~Strategy per absurdum: removal of a single point 
of the residual non-removable subset} Thus, it suffices to demonstrate
that the closed subsets $C$ of Theorems~1.2' and~1.4 are $\mathcal{
W}$-removable, {\it cf.} the definition given in Section~1. Let us fix
a wedgelike domain $\mathcal{ W}_1$ attached to $M\backslash C$ and
remind that all our wedgelike domains are assumed to be nonempty. Our
precise goal is to establish that there exists a wedgelike domain
$\mathcal{ W}_2$ attached to $M$ (including $C$) and a wedgelike
domain $\mathcal{ W}_3\subset \mathcal{ W}_1 \cap \mathcal{ W}_2$
attached to $M\backslash C$ such that for every holomorphic function
$f\in \mathcal{ O}(\mathcal{ W}_1)$, there exists a holomorphic
function $F\in \mathcal{ O}(\mathcal{ W}_2)$ which coincides with $f$
in $\mathcal{ W}_3$. At first, we need some more definition.

Let $C'$ be an arbitrary closed subset of $C$. We shall say that
$M\backslash C'$ {\sl enjoys the wedge extension property} if there
exist a wedgelike domain $\mathcal{ W}_2'$ attached to $M\backslash
C'$ and a wedgelike subdomain $\mathcal{ W}_3 \subset \mathcal{
W}_1\cap \mathcal{ W}_2'$ attached to $M\backslash C$ such that, for
every function $f\in\mathcal{ O}(\mathcal{ W}_1)$, there exists a
function $F'\in \mathcal{ O}(\mathcal{ W}_2')$ which coincides with
$f$ in $\mathcal{ W}_3$.

The notion of wedge removability can be localized as follows. Let
again $C'\subset C$ be arbitrary. We shall say that a point $p' \in
C'$ is {\sl locally $\mathcal{ W}$-removable with respect to $C'$} if
for every wedgelike domain $\mathcal{ W}_1'$ attached to $M\backslash
C'$, there exists a neighborhood $U'$ of $p'$ in $M$, there exists a
wedgelike domain $\mathcal{ W}_2'$ attached to $(M \backslash C') \cup
U'$ and there exists a wedgelike subdomain $\mathcal{ W}_3' \subset
\mathcal{ W}_1' \cap \mathcal{ W}_2'$ attached to $M\backslash C'$
such that for every holomorphic function $f \in \mathcal{ O}
(\mathcal{ W}_1')$, there exists a holomorphic function $F'\in
\mathcal{ O}(\mathcal{ W}_2')$ which coincides with $f$ in $\mathcal{
W}_3'$.

Supppose now that $M\backslash C_1'$ and $M\backslash C_2'$ enjoy the
wedge extension property, for some two closed subsets $C_1', C_2'
\subset C$. Using Ayrapetian's version of the edge of the wedge
theorem (\cite{a}, \cite{tu1}, \cite{tu2}), the two wedgelike domains
attached to $M\backslash C_1'$ and to $M\backslash C_2'$ can be glued
(after appropriate shrinking) to a wedgelike domain $\mathcal{ W}_1$
attached to $M\backslash (C_1'\cap C_2')$ in such a way that
$M\backslash (C_1'\cap C_2')$ enjoys the $\mathcal{ W}$-extension
property. Also, if $M\backslash C'$ enjoys the wedge extension
property and if $p'\in C'$ is locally $\mathcal{ W}$-removable with
respect to $C'$, then again by means of the edge of the wedge theorem,
it follows that there exists a neighborhood $U'$ of $p'$ in $M$ such
that $(M\backslash C') \cup U'$ enjoys the wedge extension property.

Based on these preliminary remarks, we may define the following set of
closed subsets of $C$:
\def\theequation{3.17}\begin{equation}
\mathcal{ C}:=\{C'\subset C \ \hbox{closed} \ ;
\ M\backslash C' \ \text{\rm enjoys the 
$\mathcal{ W}$-extension property} \}.
\end{equation}
Then the residual set 
\def\theequation{3.18}\begin{equation}
C_{\text{\rm nr}}:= \bigcap_{C' \in\mathcal{ C}} C'
\end{equation}
is a closed subset of $M^1$ contained in $C$. It follows from the
above (abstract nonsense) considerations that $M\backslash C_{\rm nr}$
enjoys the wedge extension property and that no point of $C_{\rm nr}$
is locally $\mathcal{ W}$-removable with respect to $C_{\rm nr}$.
Here, we may think that the letters ``nr'' abbreviate ``{\sl n}on-{\sl
r}emova\-ble'', because by the very definition of $C_{\rm nr}$, none
of its points should be locally $\mathcal{ W}$-removable. Notice also
that $M\backslash C_{\rm nr}$ is globally minimal, thanks to
Lemma~3.5.

Clearly, to establish Theorem~1.1, it is enough to show that
$C_{\text{\rm nr}}=\emptyset$.

We shall argue indirectly (by contradiction) and assume that
$C_{\text{\rm nr}}\neq \emptyset$. In order to derive a
contradiction, it clearly suffices to show that there exists at least
one point $p\in C_{\rm nr}$ which is in fact locally $\mathcal{
W}$-removable with respect to $C_{\rm nr}$.

At this point, we notice that the main assumptions $\mathcal{
F}_{M^1}^c\{C\}$ and $\mathcal{ O}_{M^1}^{CR}\{C\}$ of Theorem~1.2'
and of Theorem~1.4 imply trivially that for every closed subset $C'$
of $C$, then the condition $\mathcal{ F}_{M^1}^c\{C'\}$ and the
condition $\mathcal{ O}_{M^1}^{CR}\{C'\}$ also hold true: these two
assumptions are obviously stable by passing to closed subsets. In
particular, $\mathcal{ F}_{M^1}^c\{C_{\rm nr}\}$ and $\mathcal{
O}_{M^1}^{CR}\{C_{\rm nr}\}$ hold true. Consequently, by following a
{\it per absurdum} strategy, we are led to prove two statements wich
are totally similar to Theorem~1.2' and to Theorem~1.4, except that we
now have only to establish that {\it a single point of $C_{\rm nr}$}
is locally $\mathcal{ W}$-removable with respect to $C_{\rm nr}$. This
preliminary logical consideration will simplify substantially the
whole architecture of the two proofs. Another important advantage of
this strategy which will not be appearant until the very end of the
two proofs in Sections~9 and~10 below is that we are even allowed to
select a special point $p_{\rm sp}$ of $C_{\rm nr}$ by requiring
some nice geometric disposition of $C_{\rm nr}$ in a neighborhood of
$p_{\rm sp}$ before removing it. Sections~4 and~5 below are devoted to
such a selection.

So we are led to show that for every wedgelike domain $\mathcal{ W}_1$
attached to $M\backslash C_{\rm nr}$, there exists a special point
$p_{\rm sp}\in C_{\rm nr}$, there exists a neighborhood $U_{p_{\rm
sp}}$ of $p_{\rm sp}$ in $M$, there exists a wedgelike domain
$\mathcal{ W}_2$ attached to $\left(M\backslash C_{\rm nr}\right)\cup
U_{p_{\rm sp}}$ and there exists a wedgelike domain $\mathcal{
W}_3\subset \mathcal{ W}_1 \cap \mathcal{ W}_2$ attached to
$M\backslash C_{\rm nr}$ such that for every holomorphic function
$f\in \mathcal{ O}(\mathcal{ W}_1)$, there exists a function
$F\in\mathcal{ O}(\mathcal{ W}_2)$ which coincides with $f$ in
$\mathcal{ W}_3$.

A further convenient simplification of the task may be achieved by
deforming slightly $M$ inside the wedge $\mathcal{ W}_1$ attached to
$M\backslash C_{\rm nr}$. Indeed, by means of a partition of unity, we
may perform arbitrarily small $\mathcal{ C}^{2,\alpha}$-smooth
deformations $M^d$ of $M$ leaving $C_{\text{\rm nr}}$ fixed and moving
$M\backslash C_{\text{\rm nr}}$ inside the wedgelike domain $\mathcal{
W}_1$. Furthermore, we can make $M^d$ to depend on a single small
real parameter $d\geq 0$ with $M^0=M$ and $M^d\backslash C_{\rm
nr}\subset \mathcal{ W}_1$ for all $d>0$. Now, {\it the wedgelike
domain $\mathcal{ W}_1$ becomes a neighborhood of $M^d$ in
$\C^n$}. Let us denote by $\Omega$ this neighborhood. After some
substantial technical work has been achieved, at the end of the proofs
of Theorem~1.2' and~1.4 to be conducted in Sections~9 and~10 below, we
shall construct a local wedge $\mathcal{ W}_{p_{\rm sp}}^d$ of edge
$M^d$ at $p_{\rm sp}$ by means of small Bishop analytic discs glued to
$M^d$, to $\Omega$ and to another subset (which we will call a {\sl
half-wedge}, {\it see} Section~4 below) such that every holomorphic
function $f\in \mathcal{ O}(\Omega)$ extends holomorphically to
$\mathcal{ W}_{p_{\rm sp}}^d$. Using the stability of Bishop's
equation under perturbation, we shall argue in \S9.27
below that our constructions are stable under such small deformations,
whence in the limit $d\rightarrow 0$, the wedges $\mathcal{ W}_{p_{\rm
sp}}^d$ tend smoothly to a local wedge $\mathcal{ W}_{p_{\rm sp}}:=
\mathcal{ W}_{p_{\rm sp}}^0$ of edge a neighborhood $U_{p_{\rm sp}}$
of $p_{\rm sp}$ in $M^0 \equiv M$ (notice that in codimension $\geq
2$, a wedge of edge a deformation $M^d$ of $M$ does not contain in
general a wedge of edge $M$, hence such an argument is needed). In
addition, we shall derive univalent holomorphic extension to
$\mathcal{ W}_{p_{\rm sp}}$. Finally, using again the
edge of the wedge theorem to fill the space between $\mathcal{ W}_1$
and $\mathcal{W}_{p_{\rm sp}}$, possibly after appropriate
contractions of these two wedgelike domains, we may construct a
wedgelike domain $\mathcal{ W}_2$ attached to $\left(M\backslash
C\right) \cup U_{p_{\rm sp}}$ and a wedgelike domain $\mathcal{
W}_3\subset \mathcal{ W}_1 \cap \mathcal{ W}_{p_{\rm sp}}$ attached to
$M\backslash C$ such that for every holomorphic function $f\in
\mathcal{ O}(\mathcal{ W}_1)$, there exists a function $F\in\mathcal{
O}(\mathcal{ W}_2)$ which coincides with $f$ in $\mathcal{ W}_3$. In
conclusion, we thus reach the desired contradiction to the definition
of $C_{\text{\rm nr}}$.

As a summary of the above discussion, we have essentially shown that
it suffices to prove Theorems~1.2' and~1.4 with the following two
extra simplifying assumptions:
\begin{itemize}
\item[{\bf 1)}] 
Instead of functions which are holomorphic in a
wedgelike domain attached to $M\backslash C_{\rm nr}$, we consider
functions which are holomorphic in a neighborhood
$\Omega$ of $M\backslash
C_{\rm nr}$ in $\C^{m+n}$. 
\item[{\bf 2)}]
Proceeding by contradiction, we
have argued that it suffices to remove at least one point of $C_{\rm
nr}$.
\end{itemize}
 
Consequently, we may formulate the local statement that remains to
prove: after replacing $C_{\rm nr}$ by $C$ and $M^d$ by $M$, we are
led to establish the following main assertion, to which Theorems~1.2'
and~1.4 are essentially reduced.

\def\thetheorem{3.19}\begin{theorem}
Let $M$ be a $\mathcal{ C}^{ 2,\alpha}$-smooth globally minimal
generic submanifold of $\C^n$ of CR dimension $m\geq 1$ and of
codimension $d=n- m\geq 1$, let $M^1 \subset M$ be a $\mathcal{ C}^{2,
\alpha}$-smooth one-codimensional submanifold which is generic in
$\C^n$, and let $C$ be a {\rm nonempty} proper closed subset of $M^1$.
\begin{itemize}
\item[{\bf (i)}]
If $m=1$, assume that the condition $\mathcal{ F}_{M^1}^c\{C\}$
holds.
\item[{\bf (ii)}]
If $m\geq 2$, assume that the condition 
$\mathcal{ O}_{M^1}^{CR}\{C\}$ holds.
\end{itemize}
Let $\Omega$ be an arbitrary neighborhood of $M \backslash C$ in
$\C^n$. Then there exist a special point $p_{ \rm sp} \in C$, there
exists a local wedge $\mathcal{ W}_{ p_{ \rm sp}}$ of edge $M$ at
$p_{\rm sp}$ and there exists a subneighborhood $\Omega' \subset
\Omega$ of $M \backslash C$ in $\C^n$ with $\mathcal{ W}_{p_{ \rm sp
}} \cap \Omega'$ connected such that for every holomorphic function $f
\in \mathcal{ O}( \Omega)$, there exists a holomorphic function $F\in
\mathcal{ O} \left( \mathcal{ W}_{p_{ \rm sp}}\cup \Omega' \right)$
which coincides with $f$ in $\mathcal{ W}_{ p_{\rm sp }}\cap \Omega'$.
\end{theorem}

However, we remind the necessity of some supplementary arguments about
the stability of our constructions under deformation. The proof of our
main Theorem~3.19 will occupy Sections~4, 5, 6, 7, 8, 9 and~10 below
and the deformation arguments will appear lastly in \S9.27.
From now on, the main question is: {\it How to choose
the special point $p_{\rm sp}$ to be removed locally}~?

\subsection*{3.20.~Choice of a special point $p\in C$ in the CR 
dimension $m\geq 2$ case} In the case of CR dimension $m\geq 2$,
essentially all points of $C$ can play the role of the special point
$p_{\rm sp}$. However, since we want to devise a new proof of
Theorem~1.4 which differs from the two proofs given in~\cite{ j4},
~\cite{p1} and in \cite{m2}, it will be convenient to choose a special
point $p_1 \in M^1$ which has the property that locally in a
neighborhood of $p_1$, the singularity $C\subset M^1$ lies behind a
generic ``wall'' $H^1$ contained in $M^1$ and of codimension one in
$M^1$. Notice that as $m\geq 2$, the dimension of a two-codimensional
submanifold $H^1$ of $M$ is equal to $2m+d-2\geq n$, whence $M^1$ may
perfectly be generic in $\C^n$.

\def\thelemma{3.21}\begin{lemma}
Let $M$ be a $\mathcal{ C}^{2, \alpha}$-smooth globally minimal
generic submanifold of $\C^n$ of CR dimension $m \geq 2$ and of
codimension $d = n-m \geq 1$, let $M^1\subset M$ be a $\mathcal{
C}^{2, \alpha}$-smooth one-codimensional submanifold which is generic
in $\C^n$ and let $C\subset M^1$ be a {\rm nonempty} proper closed
subset which does not contain any CR orbit of $M^1$. Then there
exists a point $p_1\in C$ and a $\mathcal{ C}^{2,\alpha}$-smooth
one-codimensional submanifold $H^1\subset M^1$ passing through $p_1$
which is {\rm generic} in $\C^n$ such that $C\backslash \{p_1\}$ lies,
in a neighborhood of $p_1$, in one open side $(H^1)^-$ of $H^1$ in
$M^1$.
\end{lemma}

\proof
The proof is completely similar to the proof of
Lemma~2.1 in~\cite{mp3}, {\it see} especially 
{\sc Figure~1}, p.~490. 
\endproof

With Lemma~3.21 at hand, we can now state a more precise version 
of case {\bf (ii)} of Theorem~3.19, which will be the
main removability proposition. 

\def\theproposition{3.22}\begin{proposition}
Let $M$ be a $\mathcal{ C}^{2, \alpha}$-smooth globally minimal
generic submanifold of $\C^n$ of CR dimension $m \geq 2$ and of
codimension $d = n-m \geq 1$, let $M^1\subset M$ be a $\mathcal{
C}^{2, \alpha}$-smooth one-codimensional submanifold which is generic
in $\C^n$, let $p_1 \in M^1$, let $H^1\subset M^1$ be a $\mathcal{
C}^{2, \alpha}$-smooth one-codimensional submanifold of $M^1$ passing
through $p_1$ which is also generic in $\C^n$ {\rm (}this is possible,
thanks to the assumption $m \geq 2${\rm )} and let $(H^1 )^-$ denote
an open local one-sided neighborhood of $H^1$ in $M^1$. Suppose that
$C \subset M^1$ is a {\rm nonempty} proper closed subset of $M^1$ with
$p_1 \in C$ which is situated, locally in a neighborhood of $p_1$,
only in one side of $H^1$, namely $C \subset (H^1)^- \cup \{
p_1\}$. Let $\Omega$ be a neighborhood of $M \backslash C$ in
$\C^n$. Then there exists a local wedge $\mathcal{ W}_{p_1}$ of edge
$M$ at $p_1$ with $\Omega \cap \mathcal{ W}_{p_1}$ connected
(shrinking $\Omega$ if necessary) such that for every holomorphic
function $f\in\mathcal{ O} (\Omega)$, there exists a holomorphic
function $F\in \mathcal{ O} \left( \mathcal{ W}_{p_1} \cup \Omega
\right)$ with $F \vert_\Omega =f$.
\end{proposition}

In the CR dimension $m=1$ case, the choice of a special point $p\in C$
is much more delicate and will be done in the next two Sections~4
and~5 below, where the analog of Proposition~3.22 appears as the main
removability proposition~5.12. In the case $m=1$, the submanifold
$M^1$ is of real dimension equal to $n$ and it is not difficult to
generalize Lemma~3.21, obtaining a submanifold $H^1$ which is of real
dimension $(n-1)$ and totally real (but not generic) in
$\C^n$. However, in general, such a point $p_1\in C\cap H^1$ is {\it not
locally $\mathcal{ W}$-removable}\, in general. For instance, in the
hypersurface case $n=2$, locally in a neighborhood of $p_1$, the
closed set $C\subset (H^1)^-\cup \{p_1\}$ may coincide with the
intersection of $M$ with a local complex curve transverse to $M$ at
$p_1$, hence $C$ is not locally removable. In this (trivial) example,
the condition $\mathcal{ F}^c\{C\}$ is not satisfied and this
justifies a more refined geometrical analysis to chase a suitable
special point $p_{\rm sp}\in C$ to be removed locally.
This is the main purpose of Sections~4 and~5 below.

\section*{\S4.~Construction of a semi-local half wedge}

\subsection*{4.1.~Preliminary}
Let $M$ be a $\mathcal{ C}^{2, \alpha}$ globally minimal generic
submanifold of $\C^n$ of CR dimension $m=1$ and of codimension $d=n-m
\geq 1$ and let $M^1$ be a $\mathcal{ C}^{2,\alpha}$-smooth
one-codimensional submanifold which is generic in $\C^n$. Let
$\gamma:[-1,1] \to M^1$ be a $\mathcal{ C}^{2, \alpha}$-smooth curve,
embedding the segment $[-1,1]$ into $M$. Later, in Section~5 below, we
shall exploit the geometric condition formulated in Theorem~2.1' that
such a characteristic curve should satisfy, but in the present
Section~4, we shall not at all take account of this geometric
condition. Our goal is to construct a semi-local half-wedge attached
to a one-sided neighborhood of $M^1$ along $\gamma$ with the property
that holomorphic functions in the neighborhood $\Omega$ of
$M\backslash C$ in $\C^n$ do extend holomorphically to this
half-wedge. First of all, we need to define what we understand by the
term ``half-wedge''. Although all the geometric considerations of
this section may be generalized to the CR dimension $m\geq 2$ case
with slight modifications, we choose to endeavour the exposition in
the case $m=1$, because our constructions are essentially needed only
for this case (however in Section~10 below, since we aim to present a
new proof of Theorem~1.4, we shall also use the notion of a local
half-wedge, without any characteristic curve $\gamma$).

\subsection*{4.2.~Three equivalent definitions of attached 
half-wedges} First of all, we define the notion of a local
half-wedge. We shall denote by $\Delta_n (p,\delta)$ the open polydisc
centered at $p \in \C^n$ of radius $\delta >0$. Let $p_1\in M^1$, and
let ${\sf C}_1$ be an open infinite cone in the normal space $T_{p_1}
\C^n/ T_{ p_1}M$. Classically, a {\sl local wedge of edge $M$ at
$p_1$} is a set of the form: $\mathcal{ W}_{p_1}:= \{p+ {\sf c}_1: \,
p \in M, \, {\sf c}_1 \in { \sf C}_1\}\cap \Delta_n (p_1, \delta_1 )$,
for some $\delta_1 >0$. Sometimes, we shall use the following
terminology: if $v_1$ is a nonzero vector in $T_{p_1} \C^n /T_{
p_1}M$, we shall say that $\mathcal{ W}_{p_1}$ is a {\sl local wedge
at $(p_1, v_1)$}. This definition seems to be misleading in the sense
that different vectors $v_1$ seem to yield local wedges with different
directions, however, there is a concrete geometric meaning in this
definition that should be reminded: the positive half-line $\R^+ \cdot
v_1$ generated by the vector $v_1$ is locally contained in the wedge
$\mathcal{ W}_{p_1}$.

For us, a {\sl local half-wedge of edge $M$ at $p_1$} will be a set of
the form
\def\theequation{4.3}\begin{equation}
\mathcal{ HW}_{p_1}^+:= \left\{p+{\sf c}_1: \,
p\in U_1\cap (M^1)^+, \, {\sf c}_1 \in {\sf C}_1\right\}
\cap \Delta_n(p_1,\delta_1). 
\end{equation}
This yields a first definition and we shall delineate two further
definitions. Let $\Delta$ denote the unit disc in $\C$, let $\partial
\Delta$ denote its boundary, the unit circle and let $\overline{
\Delta}= \Delta \cup \partial \Delta$ denote its closure.
Throughout this paper, we shall denote by $\zeta=\rho e^{i\theta}$
the variable of $\overline{\Delta}$ with $0 \leq \rho \leq 1$ and
with $\vert \theta \vert \leq \pi$.

In fact, our local half-wedges (to be constructed in this section)
will be defined by means of a $\mathcal{ C}^{2, \alpha-0}$-smooth
$\C^n$-valued (remind $\mathcal{ C}^{2,\alpha-0}\equiv \bigcap_{\beta
< \alpha} \mathcal{ C}^{2,\beta}$) mapping $(t, \chi, \nu, \rho)
\longmapsto \mathcal{ Z}_{t,\chi, \nu} (\rho)$, which comes from
parametrized family of analytic discs of the form $\zeta \mapsto
\mathcal{ Z}_{ t,\chi,\nu} (\zeta)$, where the parameters $t \in
\R^{n-1}$, $\chi \in \R^n$, $\nu \in \R$ satisfy $\vert t \vert<
\varepsilon$, $\vert \chi \vert < \varepsilon$, $\vert \nu \vert <
\varepsilon$ for some small $\varepsilon >0$, and where $\mathcal{
Z}_{ t,\chi,\nu} (\zeta)$ is holomorphic with respect to $\zeta$ in
$\Delta$. This mapping will satisfy the following three properties:

\smallskip
\begin{itemize}
\item[{\bf (i)}]
$(\chi,\nu) \mapsto \mathcal{ Z}_{0,\chi,\nu}(1)$ is a 
diffeomorphism onto a neighborhood of $p_1$ in $M$,
the mapping $\chi \mapsto \mathcal{ Z}_{0,\chi,0}(1)$ is a 
diffeomorphism onto a neighborhood of $p_1$ in $M^1$
and $(M^1)^+$ corresponds to $\nu>0$ in the diffeomorphism
$(\chi,\nu) \mapsto \mathcal{ Z}_{0,\chi,\nu}(1)$.
\item[{\bf (ii)}]
$\mathcal{ Z}_{t,0,0}(1)=p_1$ and the half-boundary
$\mathcal{ Z}_{t, \chi, \nu}\left(\left\{ e^{i\theta} : \ 
\vert \theta \vert \leq \frac{\pi}{2} \right\}\right)$ is contained in
$M$ for all $t$, all $\chi$ and
all $\nu$.
\item[{\bf (iii)}]
The vector $v_1:= \frac{ \partial \mathcal{ Z}_{0, 0,0}}{ \partial
\theta} (1)\in T_{p_1 }\C^n$ is nonzero and belongs to $T_{p_1}M^1$.
Furthermore, the rank of the $\R^{ n-1}$-valued $\mathcal{
C}^{1,\alpha-0}$-smooth mapping 
\def\theequation{4.4}\begin{equation}
\R^{n-1}\ni t \longmapsto 
\frac{\partial \mathcal{ Z}_{t,0,0}}{\partial \theta}(1)
\in T_{p_1}M^1
\ {\rm mod} \ \left( T_{p_1}M^1\cap T_{p_1}^cM
\right) \cong \R^{n-1}
\end{equation}
is maximal equal to $(n-1)$ at $t=0$.
\end{itemize}

\smallskip
\noindent
By holomorphicity of the map $\zeta\mapsto \mathcal{
Z}_{t,\chi,\nu}(\zeta)$, we have $\frac{\partial \mathcal{
Z}_{t,\chi,\nu}}{ \partial \theta}(1) =J\cdot \frac{\partial \mathcal{
Z}_{t,\chi,\nu}}{ \partial \rho}(1)$, where $J$ denotes the complex
structure of $T\C^n$. Consequently, because $J$ induces an isomorphism
from $T_{p_1}M/T_{p_1}^cM\to T_{ p_1}\C^n/T_{p_1}M$, it follows from
property {\bf (iii)} above that the vectors $\frac{\partial \mathcal{
Z}_{t,0,0}}{\partial \rho}(1)$ cover an open cone containing $Jv_1$ in
the quotient space $T_{p_1}M/ T_{p_1}^cM$, as $v$ varies. Then a {\sl
local half-wedge of edge $(M^1)^+$ at $p_1$} will be a set of the form
\def\theequation{4.5}\begin{equation}
\mathcal{ HW}_{p_1}^+:=\left\{
\mathcal{ Z}_{t,\chi,\nu}(\rho)\in\C^n: \ 
\vert t \vert < \varepsilon, \ 
\vert \chi \vert < \varepsilon, \
0 < \nu <\varepsilon, \ 
1-\varepsilon < \rho < 1
\right\}.
\end{equation}
We notice that a complete local wedge of edge $M$ at $p_1$ is also
associated to such a family $\mathcal{ Z}_{t,\chi, \nu}(\zeta)$ and may
be defined as $\mathcal{ W}_{p_1}:= \left\{ \mathcal{
Z}_{t,\chi,\nu}(\rho): \ \vert t \vert < \varepsilon, \ \vert \chi
\vert < \varepsilon, \ \vert \nu \vert < \varepsilon, \ 1-\varepsilon
< \rho < 1 \right\}$.

As may be checked, this second definition of a half-wege is {\sl
essentially equivalent} to the first one, in the sense that a half
wedge in the sense of the first definition always contains a
half-wedge in the sense of the second definition, and vice versa,
after appropriate shrinkings of open neighborhoods and cones.

Furthermore, we may distinguish two cases: either the vector $v_1= 
\frac{ \partial \mathcal{ Z}_{0, 0,0}}{ \partial
\theta} (1)$ is
not complex tangent to $M$ at $p_1$ or it is complex tangent to $M$ at
$p_1$. In the first case, after possibly shrinking $\varepsilon>0$, it
may be checked that a local half-wedge of edge $(M^1 )^+$ coincides
with the intersection of a (full) local wedge $\mathcal{ W}_{ p_1}$ of
edge $M$ at $p_1$ with a one-sided neighborhood $(N^1)^+$ of a local
hypersurface $N^1$ which intersects $M$ locally transversally along
$M^1$ at $p_1$, as drawn in the left hand side of the following
figure, where $M$ is of codimension two.

\bigskip
\begin{center}
\input half-wedge.pstex_t
\end{center}

In the second case, the vector $v_1= \frac{ \partial \mathcal{ Z}_{0,
0,0}}{ \partial \theta} (1)$ is complex tangent to $M$ at $p_1$, hence
belongs to the characteristic direction $T_{p_1}M^1\cap T_{p_1}^cM$,
so the vector $-Jv_1$ which is interiorly tangent to the disc
$\mathcal{ Z}_{0,0,0}(\Delta)$, is tangent to $M$ at $p_1$, is not
tangent to $M^1$ at $p_1$, but points towards $(M^1)^+$ at $p_1$. It
may then be checked that a local half-wedge of edge $(M^1)^+$
coincides with the intersection of a local wedge $\mathcal{
W}_{p_1}^1$ of edge $M^1$ at $(p_1,-Jv_1)$ which contains the side
$(M^1)^+$, as drawn in the right hand side of {\sc Figure~5} above, in
which $M$ is of codimension one. This provides the third and the most
intuitive definition of the notion of local half-wedge.

Finally, we may define the desired notion of a semi-local attached
half-wedge. Let $\gamma: [-1,1]\to M^1$ be an embedded $\mathcal{
C}^{2, \alpha}$-smooth segment in $M^1$. Since the normal bundle to
$M^1$ in $M$ is trivial, we can choose a coherent family of one-sided
neighborhoods $(M_\gamma^1)^+$ of $M^1$ in $M$ along $\gamma$. A {\sl
half-wedge attached to a one-sided neighborhood $(M_\gamma^1)^+$ of
$M^1$ along $\gamma$} is a domain $\mathcal{ HW}_\gamma^+$ which
contains a local half-wedge of edge $(M^1)^+$ at $\gamma(s)$ for every
$s\in [-1,1]$. Another essentially equivalent definition is to
require that we have a family $\mathcal{ Z}_{t,\chi, \nu:s}(\rho)$ of
mappings smoothly varying with the parameter $s$ such that at each
point $\gamma(s)=\mathcal{ Z}_{t,\chi, \nu:s}(1)$, the three
conditions {\bf (i)}, {\bf (ii)} and {\bf (iii)} introduced above to
define a local half-wedge are satisfied. Intuitively speaking, the
direction of the cones defining the local half wedge at the point
$\gamma(s)$ are smoothly varying with respect to $s$.

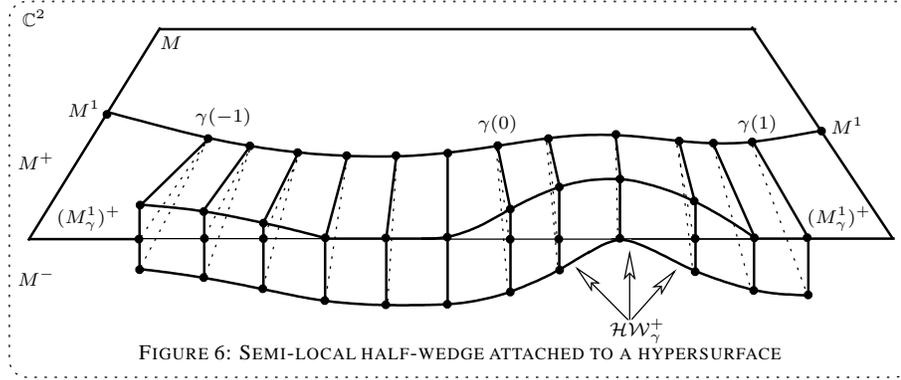
\bigskip
\begin{center}
\input attached-half-wedge.pstex_t
\end{center}

We can now state the main proposition of this section, which will be
of crucial use for the proof of Theorem~3.19 {\bf (i)}. Forgetting
for a while the complete content of the geometric condition $\mathcal{
F}_{M^1}^c\{C\}$ formulated in the assumptions of Theorem~1.2', which
we will analyze thoroughly in Section~5 below, we shall only assume
that we are given a characteristic segment $\gamma : [-1,1] \to M^1$
in the following proposition, whose proof is the main goal of this
Section~4.

\def\theproposition{4.6}\begin{proposition} 
Let $M$ be a $\mathcal{ C}^{2,\alpha}$ globally minimal generic
submanifold of $\C^n$ of CR dimension $m=1$ and of codimension
$d=n-m\geq 1$ and let $M^1$ be a $\mathcal{ C}^{2,\alpha}$-smooth
one-codimensional submanifold which is generic in $\C^n$. Let $\gamma:
[-1,1]\to M^1$ be an arbitrary $\mathcal{ C}^{2,\alpha}$-smooth
curve. Then there exist a neighborhood $V_\gamma$ of $\gamma[-1,1]$ in
$M$ and a semi-local one-sided neighborhood $(M_\gamma^1)^+$ of $M^1$
in $M$ along $\gamma$ which is the intersection of $V_\gamma$ with a
side $(M_\gamma^1)^+$ of $M^1$ along $\gamma$ and there exists a
semi-local half-wedge $\mathcal{ HW }_\gamma^+$ attached to
$(M_\gamma^1 )^+ \cap V_\gamma$ with $\Omega \cap \mathcal{
HW}_\gamma^+$ connected {\rm (}shrinking $\Omega$ if necessary{\rm )}
such that for every holomorphic function $f\in \mathcal{ O} (\Omega)$,
there exists a holomorphic function $F\in \mathcal{ O} \left(
\mathcal{ HW }_\gamma^+ \cup \Omega \right)$ with $F\vert_\Omega=f$.
\end{proposition}

To build $\mathcal{ HW}_\gamma^+$, we shall construct families of
analytic discs with boundaries in $(M_\gamma^1)^+$. First of all, we
need to formulate a special, adapted version of the so-called {\sl
approximation theorem} of~\cite{ bt}.

\subsection*{4.7.~Local approximation theorem}
As noted in~\cite{m2}, ~\cite{mp1} and~\cite{mp3}, when dealing with
some natural geometric assumptions on the singularity to be
removed\,--\,for instance, a two-codimensional singularity $N\subset
M$ with $T_pN \supset T_p^cM$ for some points $p\in N$ or metrically
thin singularities $E\subset M$ with $H^{2m+d-2}(E)=0$\,--\,it is
impossible to show {\it a priori}\, that continuous CR functions on
$M$ minus the singularity are approximable by polynomials, which
justifies the introduction of deformations and the use of the
continuity principle in~\cite{m2}, \cite{mp1}, \cite{mp3}. On the
contrary the genericity of the submanifold $M^1$ containing the
singularity $C$ enables us to get an approximation Lemma~4.8 just
below, in the spirit of~\cite{bt}. Together with the existence of
Bishop discs attached to $M^1$, the validity of this approximation
lemma on $M\backslash M^1$ is the second main reason for the relative
simplicity of the geometric proofs of Theorem~1.4 provided
in~\cite{j4}, \cite{m2}, \cite{p1}, in comparison the 
proof of Theorem~1.2' to be conducted in this paper.

\def\thelemma{4.8}\begin{lemma}
Let $M$ be a $\mathcal{ C}^{2, \alpha}$-smooth globally minimal
generic submanifold of $\C^n$ of CR dimension $m \geq 1$ and of
codimension $d = n-m \geq 1$, let $M^1 \subset M$ be a $\mathcal{
C}^{2, \alpha}$-smooth one-codimensional submanifold which is generic
in $\C^n$ which divides locally $M$ in two open sides $(M^1)^-$ and
$(M^1)^+$ and let $p_1\in M^1$. Then there exist two neighborhoods
$U_1$ and $V_1$ of $p_1$ in $M$ with $V_1 \subset \subset U_1$ such
that for every continuous CR function $f \in \mathcal{ C}_{CR}^0
\left(( M^1 )^+\cap U_1 \right)$, there exists a sequence of
holomorphic polynomials $(P_\nu)_{\nu \in\N}$ wich converges uniformly
towards $f$ on $(M^1)^+\cap V_1$. Of course, the same property holds
in the side $(M^1)^-$ instead of $(M^1)^+$.
\end{lemma}

\proof
The proof is a slight modification of~Proposition~5B in~\cite{m2} and
we summarize it, taking for granted that the reader is acquainted with
the approximation theorem proved in~\cite{bt} ({\it see}
also~\cite{j5}). Let $L_0^1$ be a maximally real submanifold passing
through $p_1$ and contained in $M^1\cap U_1$, for a sufficiently small
neighborhood $U_1$ of $p_1$ in $M^1$, possibly to be shrunk later. In
coordinates $z=(z_1,\dots,z_n)=x+iy\in\C^n$ vanishing at $p_1$, we can
assume that the tangent plane to $L_0^1$ at $p_1$ identifies with
$\R^n=\{y=0\}$. As the codimensions of $L_0^1$ in $M^1$ and in $M$ are
equal to $(d-1)$ and to $d$, we can include $L^1$ in a $d$-parameter
family of submanifolds $L_t^1$, where $t\in \R^d$ is small, so that
$L_t^1\cap V_1$ makes a foliation of $M\cap V^1$ by maximally real
$\mathcal{ C}^{2,\alpha}$-smooth submanifolds, for some neighborhood
$V_1\subset\subset U^1$ of $p^1$ in $M$, such that $L_t^1$ is
contained in $M^1$ for $t=(t_1,\dots, t_{d-1},0)$, {\it i.e.} for
$t_d=0$, such that $L_t^1\cap V_1$ is contained in $(M^1)^+$ for
$t_d>0$ and such that $L_t^1 \cap V_1$ is contained in $(M^1)^-$ for
$t< 0$. In addition, we can assume that all the $L_t^1$ coincide with
$L_0^1$ in a neighborhood of $\partial U_1$.

We shall first treat the case where $f$ is of class $\mathcal{
C}^1$. Thus, let $f$ be a $\mathcal{ C}^1$-smooth CR function on
$\left( M \backslash M^1\right) \cap U_1$, let $\tau \in \R$ with
$\tau > 0$, fix $t_0\in \R^d$ small with $t_{d;0} >0$, whence
$L_{t_0}^1\cap V^1$ is contained in $(M^1)^+$, let $\widehat{z} \in (
M^1)^+ \cap V_1$ be an arbitrary point and consider the following
integral which consists of the convolution of $f$ with the Gauss
kernel:
\def\theequation{4.9}\begin{equation}
G_\tau f(\widehat{z}):=
\left(\frac{\tau}{\pi}\right)^{n/2}
\int_{U_1\cap L_{t_0}} \, 
e^{-\tau(z-\widehat{z})^2}\, 
f(z) \, dz,
\end{equation}
where $(z- \widehat{z})^2:= (z_1- \widehat{ z}_1)^2 + \cdots+
(z_n-\widehat{ z}_n)^2$ and $dz := dz_1 \wedge \cdots \wedge d z_n$.
The point $\widehat{ z}$ belongs to some maximally real submanifold
$L_{\widehat{ t}}$ with $\widehat{ t}_d>0$. We now claim that the
value of $G_\tau f( \widehat{ z})$ is the same if we replace the
integration on the fixed submanifold $U_1 \cap L_{ t_0}$ in the
integral~\thetag{4.9} by an integration over $U_1 \cap L_{ \widehat{
t}}$. This key argument will follows from Stoke's theorem, from the
fact that $f$ is CR on $M\cap U_1$ and from the important fact that
between $L_{ t_0}$ and $L_{\widehat{ t}}$, we can construct a $(n+
1)$-dimensional submanifold $\Sigma$ with boundary $\partial \Sigma =
L_{t_0} - L_{ \widehat{ t}}$ {\it which is entirely contained in
$(M^1)^+\cap U_1$}, thanks to the fact that $M^1$ is generic and of
codimension $1$; indeed, we may compute
\def\theequation{4.10}\begin{equation}
\left\{
\aligned
G_\tau f(\widehat{ z})=
& \
\left(\frac{\tau}{\pi}\right)^{n/2}
\int_{U_1\cap L_{\widehat{ t}}}
e^{-\tau(z-\widehat{ z})^2}\, 
f(z) \, dz+\left(\frac{\tau}{\pi}\right)^{n/2}
\int_{\Sigma}\, 
d\left(
e^{-\tau (z-\widehat{ z})^2}\, f(z) \, dz
\right) \\
=
& \
\left(\frac{\tau}{\pi}\right)^{n/2}
\int_{U_1\cap L_{\widehat{ t}}}
e^{-\tau(z-\widehat{ z})^2}\, 
f(z) \, dz,
\endaligned\right.
\end{equation}
noticing that the second integral in the right hand side of the first
line vanishing, because $f$ and $e^{-(z-\widehat{ z})^2}$ being CR and
of class at least $\mathcal{ C}^1$, one has $d\left( e^{-\tau (z-
\widehat{ z})^2}\, f(z) \, dz \right)=0$. This proves the claim.

By analyzing the real and the imaginary part of the phase function
$-\tau(z-\widehat{ z})^2$ on $L_{\widehat{ t}}$, one can show by means
of a standard argument (convolution with Gauss' kernel is an
approximation of the Dirac mass) that the integral over $U_1 \cap
L_{\widehat{ t}}$ tends towards $f(\widehat{ z})$ as $\tau$ tends to
$\infty$, provided that the submanifold $U_1\cap L_{\widehat{ t}}$ is
sufficiently close to the real plane $\R^n$ in $\mathcal{ C}^1$
norm. Finally, by developing in power series and truncating the
exponential in the first expression~\thetag{4.9} which defines $G_\tau
f(\widehat{ z})$ and by integrating termwise, one constructs the
desired sequence of polynomials $(P_\nu (z))_{\nu\in\N}$.

The case where $f$ is only continuous follows from standard arguments
from the theory of distributions. This completes the proof of
Lemma~4.8.
\endproof

\subsection*{4.11.~A family of straightenings} 
Our main goal in the remainder of this section is to construct a
semi-local half-wedge attached to a one-sided neighborhood
$(M_\gamma^1)^+$ of $M^1$ in $M$ along $\gamma$ which consists of
analytic discs attached to $(M_\gamma^1)^+$. First of all, we need a
convenient family of normalizations of the local geometries of $M$ and
of $M^1$ along the points $\gamma(s)$ of our characteristic curve
$\gamma$, for all $s$ with $-1\leq s \leq 1$.

Let $\Omega$ be a thin neighborhood of $\gamma( [-1, 1])$ in $\C^n$,
say a union of polydiscs of fixed radius centered at the points
$\gamma (s)$. Then there exists $n$ real valued $\mathcal{ C}^{2,
\alpha}$-smooth functions $r_1 (z, \bar z), \dots, r_n (z,\bar z)$
defined in $\Omega$ such that $M \cap \Omega$ is given by the $(n-1)$
Cartesian equations $r_2 (z,\bar z)= \cdots= r_n (z,\bar z)=0$ and
such that moreover, $M^1\cap \Omega$ is given by the $n$ Cartesian
equations $r_1(z,\bar z)=r_2(z,\bar z)=\cdots=r_n(z,\bar z)=0$. We
first center the coordinates at $\gamma(s)$ by setting $z':=
z-\gamma(s)$. Then the defining functions centered at $z'=0$ become
\def\theequation{4.12}\begin{equation}
r_j\left(z'+
\gamma(s), \bar z'+\overline{\gamma(s)}\right)-r_j
\left(\gamma(s),\overline{\gamma(s)}\right)=:
r_j'(z',\bar z':s),
\end{equation}
for $j=1,\dots,n$, and they are parametrized by $s\in [-1,1]$. Now,
we drop the primes on coordinates and we denote by $r_j(z,\bar z:s)$,
$j=1,\dots,n$, the defining equations for the new 
$M_s$ and $M_s^1$, which
correspond to the old $M$ and $M^1$ locally in a neighborhood of
$\gamma(s)$. Next, we straighten the tangent planes by using the
linear change of coordinates $z'=A_s\cdot z$, where the $n\times n$
matrix $A_s$ is defined by $A_s:= 2i \, \left(\frac{\partial
r_j}{\partial z_k}(0,0:s) \right)_{1\leq j,k\leq n}$. Then the
defining equations for the two transformed 
$M_s'$ and for $M_s^{1'}$ are given by
\def\theequation{4.13}\begin{equation}
r_j'(z',\bar z':s):= r_j\left(A_s^{-1} \cdot z', \overline{A}_s^{-1}
\cdot \overline{ z}': s\right),
\end{equation}
and we check immediately that the matrix $\left( \frac{ \partial
r_j'}{ \partial z_k} (0,0:s) \right)_{1 \leq j,k \leq n}$ is equal to
$2i$ times the $n\times n$ identity matrix, whence $T_0 M_s'= \{y_2'=
\cdots = y_n' =0\}$ and $T_0{ M_s' }^1=\{y_1'=y_2'= \cdots =
y_n'=0\}$. It is important to notice that the matrix $A_s$ only
depends $\mathcal{ C}^{1,\alpha}$-smoothly with respect to $s$.
Consequently, if we now drop the primes on coordinates, the defining
equations for $M_s$ and for $M_s^1$ are of class $\mathcal{ C}^{2,
\alpha}$ with respect to $(z,\bar z)$ and only of class $\mathcal{
C}^{1, \alpha}$ with respect to $s$.

Applying then the $\mathcal{ C}^{2,\alpha}$-smooth
implicit, we deduce that there exist $(n-1)$ functions
$\varphi_j(x,y_1:s)$, $j=2,\dots,n$, which are all of class $\mathcal{
C}^{2,\alpha}$ with respect to $(x,y_1)$ in a real cube
$\I_{n+1}(2\rho_1):=\{(x,y_1)\in \R^n \times \R: \ \vert x \vert <
2\rho_1, \ \vert y_1 \vert < 2\rho_1\}$, for some $\rho_1>0$, which
are uniformly bounded in $\mathcal{ C}^{2,\alpha}$-norm as the
parameter $s$ varies in $[-1,1]$, which are of class $\mathcal{
C}^{1,\alpha}$ with respect to $s$, such that $M_s$ may be represented
in the polydisc $\Delta_n (\rho_1)$ by the $(n-1)$ graphed equations
\def\theequation{4.14}\begin{equation}
y_2=\varphi_2(x,y_1: s), \dots\dots, \, 
y_n=\varphi_n(x,y_1:s),
\end{equation}
or more concisely $y'=\varphi'(x,y_1:s)$, if we denote the coordinates
$(z_2,\dots,z_n)$ simply by $z'=x'+iy'$. Here, by construction, we
have the normalization conditions $\varphi_j(0:s)=\partial_{x_k}
\varphi_{j}(0:s)=\partial_{y_1} \varphi_{j}(0:s)=0$, for $j=2,\dots,n$
and $k=1,\dots,n$. Sometimes in the sequel, we shall use the notation
$\varphi_j(z_1,x':s)$ instead of $\varphi_j(x,y_1:s)$. Similarly,
again by means of the implicit function theorem, we obtain $n$
functions $h_k(x:s)$, for $k=1,\dots,n$ , which 
are of class $\mathcal{ C}^{2,\alpha}$
in the cube $\I_n(2\rho_1)$ (after possibly shrinking $\rho_1$)
enjoying the same regularity property with respect 
to $s$, such that $M_s^1$ is
represented in the polydisc $\Delta_n(\rho_1)$ by the $n$ graphed
equations
\def\theequation{4.15}\begin{equation}
y_1=h_1(x:s), \
y_2=h_2(x:s), \
\dots\dots, \ y_n=h_n(x:s).
\end{equation} 
In addition, we can assume that
\def\theequation{4.16}\begin{equation}
h_j(x:s)\equiv
\varphi_j\left(x,h_1(x:s):s\right), \ \ \ \ \ \ \
j=2,\dots,n.
\end{equation}
Here, by construction, we have the normalization conditions
$h_k(0:s)=\partial_{x_l}
h_k(0:s)=0$ for $k,l=1,\dots,n$.

In the sequel, we shall denote by $\widehat{ z}=\Phi_s(z)$ the
final change of coordinates which is centered at $\gamma(s)$ and
which straightens simultaneously the tangent planes to 
$M$ at $\gamma(s)$ and to $M^1$ at $\gamma(s)$ and we shall
denote by $M_s$ and by $M_s^1$ the transformations of
$M$ and of $M^1$.

Also, we shall remind that the following regularity
properties hold for the functions $\varphi_j (x,y_1:s)$ and
$h_k(x:s)$:
\begin{itemize}
\item[{\bf (a)}]
For fixed $s$, they are of class $\mathcal{ C}^{2,\alpha}$ with
respect to their principal variables, namely excluding the parameter
$s$.
\item[{\bf (b)}]
They are of class $\mathcal{ C}^{1,\alpha}$ with respect to all their
variables, including the parameter $s$.
\item[{\bf (c)}]
Each of their first order partial derivative with respect to one of
their principal variables is of class $\mathcal{ C}^{1,\alpha}$ with
respect to all their variables, including the parameter $s$.
\end{itemize}
Indeed, these properties are clearly satisfied for the
functions~\thetag{4.13} and they are inherited after the two
applications of the implicit function theorem which yielded the
functions $\varphi_j(x,y_1:s)$ and $h_k(x:s)$.

\subsection*{4.17.~Contact of a small ``round'' analytic
disc with $M^1$} Let $r\in \R$ with $0\leq r \leq r_1$, where $r_1$ is
small in comparison with $\rho_1$. Then the ``round'' analytic disc
$\overline{\Delta} \ni \zeta \to \widehat{ Z}_{1;r}(\zeta):=
ir(1-\zeta)\in \C$ with values in the complex plane equipped with the
coordinate $z_1=x_1+iy_1$ is centered at the point $ir$ of the
$y_1$-axis, is of radius $r$ and is contained in the open upper half
plane $\{z_1\in \C: y_1> 0\}$, except its boundary point $\widehat{
Z}_{1;r}(1)=0$. In addition, the tangent direction
$\frac{\partial}{\partial \theta} \widehat{ Z}_{1;r}(1)=r$ is directed
along the positive $x_1$-axis, {\it see} in advance {\sc Figure~7}
below.

As in~\cite{tu2}, \cite{tu3}, \cite{mp1}, \cite{mp3}, we denote by
$T_1$ the Hilbert transform (harmonic conjugate operator) on $\partial
\Delta$ vanishing at $1$, namely $(T_1 X) (1)=0$, whence
$T_1(T_1(X))=-X+X(1)$.

By lifting this disc contained in the complex tangent space
$T_{p_1}^cM\equiv \C_{z_1}\times \{0\}$, we may define an analytic
disc parametrized by $r$ and $s$ which is attached to $M$ 
of the form
\def\theequation{4.18}\begin{equation}
\widehat{ Z}_{r:s}(\zeta)=\left(
ir(1-\zeta), \widehat{ Z}_{r:s}'(\zeta)
\right)\in \C\times \C^{n-1}
\end{equation}
where the real part $\widehat{ X}_{r:s}'(\zeta)$ of $\widehat{
Z}_{r:s}'(\zeta)$ satisfies the following Bishop type equation on
$\partial \Delta$
\def\theequation{4.19}\begin{equation}
\widehat{ X}_{r:s}'(\zeta)= -\left[T_1
\varphi'\left(\widehat{ Z}_{1;r}(\cdot), 
\widehat{ X}_{r:s}'(\cdot):s\right)\right](\zeta), \ \ \ \ \ \
\zeta\in\partial\Delta.
\end{equation}
By~\cite{tu1}, \cite{tu3}, if $r_1$ is sufficiently small, there
exists a solution which is $\mathcal{ C}^{2,\alpha-0}$-smooth with
respect to $(r,\zeta)$, but only $\mathcal{ C}^{1, \alpha-0}$-smooth
with respect to $(r,\zeta,s)$. Notice that for $r=0$, the disc
$\widehat{ Z}_{1;0}(e^{i\theta })$ is constant equal to $0$ and by
uniqueness of the solution of~\thetag{4.19}, it follows that
$\widehat{ Z}_{0:s}'\left(e^{i\theta}\right) \equiv 0$. It follows
trivially that $\partial_\theta \widehat{ X}_{0:s}\left(e^{i\theta}
\right)\equiv 0$ and that $\partial_\theta\partial_\theta \widehat{
X}_{0:s} \left(e^{i\theta}\right)\equiv 0$, which 
will be used in a while. Notice also that
$\widehat{ X}_{r:s}(1)=0$ for all $r$ and all $s$.

On the other hand, since by assumption, we have $h_1(0:s)=0$ and
$\partial_{x_k}h_1(0:s)=0$ for $k=1,\dots,n$, it follows from the
chain rule that if we set
\def\theequation{4.20}\begin{equation}
F(r,\theta:s):= h_1\left(
\widehat{ X}_{r:s}\left(e^{i\theta}\right):s\right)
\end{equation}
where $\theta$ satisfies $0\leq \vert \theta \vert \leq \pi$,
then
the following four equations hold
\def\theequation{4.21}\begin{equation}
F(0,\theta:s)\equiv 0, \ \ \ \ \
F(r,0:s)\equiv 0, \ \ \ \ \
\partial_\theta F(0,\theta:s)\equiv 0, \ \ \ \ \ 
\partial_\theta F(r,0:s)\equiv 0.
\end{equation}
We claim that there exists a constant $C>0$ such that
the following five inequalities hold
for $0\leq \vert \theta \vert \leq \pi$, for
$0\leq r \leq r_1$, for $s\in [-1,1]$ and
for $\vert x \vert \leq \rho_1$:
\def\theequation{4.22}\begin{equation}
\left\{
\aligned
{}
& \
\left\vert
\widehat{ X}_{r:s}\left(e^{i\theta}\right) 
\right\vert\leq C\cdot r, \\
& \
\left\vert 
\partial_\theta \widehat{ X}_{r:s}\left(e^{i\theta}\right)
\right\vert \leq C\cdot r, \\
& \
\left\vert
\partial_\theta \partial_\theta \widehat{ X}_{r:s}\left(e^{i\theta}
\right) 
\right\vert \leq C\cdot r^{\frac{\alpha}{2}}, \\
& \
\sum_{k=1}^n 
\left\vert
\partial_{x_k} h_1(x) 
\right\vert\leq C\cdot \vert x \vert, \\
& \
\sum_{k_1,k_2=1}^n\, 
\left\vert
\partial_{x_{k_1}}\partial_{x_{k_2}}
h_1(x)
\right\vert \leq C.
\endaligned\right.
\end{equation}
The best constants for each inequality are {\it a priori} distinct,
but we simply take for $C$ the largest one. Indeed, the first, the
second and the third inequalities are elementary consequences of the
(uniform with respect to $s$) $\mathcal{ C}^{2,
\frac{\alpha}{2}}$-smoothness of $\widehat{
X}_{r:s}\left(e^{i\theta}\right)$ with respect to $(r,\theta)$, and of
the normalization conditions~\thetag{4.21}. The fourth and the fifth
inequalities are consequences of the $\mathcal{
C}^2$-smoothness of $h_1$ and of its first order normalizations
(complete argument may be easily be provided by imitating the
reasonings of the elementary Section~6 below). 

Computing now the second derivative
of $F(r,\theta:s)$ with respect to $\theta$, we obtain
\def\theequation{4.23}\begin{equation}
\small
\left\{
\aligned
\partial_\theta\partial_\theta
F(r,\theta:s)
& 
=
\sum_{k=1}^n\, 
\partial_{x_k} h_1 \left(
\widehat{ X}_{r:s}\left(e^{i\theta}\right):s
\right)\cdot
\partial_\theta \partial_\theta
\widehat{ X}_{k;r:s}\left(e^{i\theta}\right) + \\
&
+
\sum_{k_1,k_2=1}^n\, 
\partial_{x_{k_1}}\partial_{x_{k_2}} 
h_1\left(
\widehat{ X}_{r:s}\left(e^{i\theta}\right)\right)\cdot 
\partial_\theta \widehat{ X}_{k_1;r,s}\left(e^{i\theta}\right)\cdot 
\partial_\theta \widehat{ X}_{k_2;r,s}\left(e^{i\theta}\right),
\endaligned\right.
\end{equation}
and we may apply the majorations~\thetag{4.22} to get
\def\theequation{4.24}\begin{equation}
\left\{
\aligned
\left\vert
\partial_\theta \partial_\theta
F(r,\theta:s)
\right\vert\leq
& \
C\cdot \left\vert
\widehat{ X}_{r:s}(e^{i\theta}) \right\vert
\cdot C \cdot r^{\frac{\alpha}{2}}
+C\cdot (C\cdot r)^2 \\
\leq
& \
r\cdot C^3\left[
r^{\frac{\alpha}{2}}+r^2
\right].
\endaligned\right.
\end{equation}

We can now state and prove a lemma which shows that the disc
boundaries $\widehat{ Z}_{r:s}(\partial\Delta)$ touches $M^1$ only at
$p_1$ and lies in $(M^1)^+\cup \{p_1\}$.

\def\thelemma{4.25}\begin{lemma}
If $r_1 \leq \min\left( 1, \ \left( \frac{1}{4C^3 \pi^2}
\right)^{\frac{2}{\alpha}}\right)$, then $\widehat{ Z}_{r:s}(\partial
\Delta \backslash \{1\})$ is contained in $(M_s^1)^+$ for all $r$ with
$0 < r \leq r_1$ and all $s$ with $-1\leq s \leq 1$.
\end{lemma}

\proof
In the polydisc $\Delta_n(\rho_1)$, the positive half-side $(M_s^1)^+$
in $M$ is represented by the single equation $y_1>h_1(x:s)$, hence we
have to check that $\widehat{ Y}_{1;r}\left(e^{i\theta}\right) >
\left\vert h_1\left(\widehat{ X}_{r:s}\left(e^{i\theta}\right):s\right)
\right\vert$, for all $\theta$ with $0< \vert \theta \vert \leq \pi$.
According to~\thetag{4.18}, the $y_1$-component $\widehat{
Y}_{1;r}\left(e^{i\theta}\right)$ of $\widehat{
Z}_{r:s}\left(e^{i\theta}\right)$ is given by
$r\left(1-\cos \theta\right)$.

On the first hand, we observe the elementary minoration $r(1-\cos
\theta) \geq r\cdot \theta^2\cdot \frac {1}{\pi^2}$, valuable for
$0\leq \vert \theta \vert \leq \pi$.

On the second hand, taking account of the second and fourth
relations~\thetag{4.21}, Taylor's integral formula now yields
\def\theequation{4.26}\begin{equation} 
F(r,\theta:s)=\int_0^\theta \,
(\theta-\theta') \cdot \partial_\theta 
\partial_\theta F\left(r,\theta':s\right) \cdot
d\theta'.
\end{equation}
Observing that $r^2\leq r^{\frac{ \alpha}{2}}$, since $0<r \leq r_1\leq
1$, and using the majoration~\thetag{4.24}, we may estimate, taking
account of the assumption on $r_1$ written in the statement of the
lemma:
\def\theequation{4.27}\begin{equation}
\left\vert
F(r,\theta:s) 
\right\vert\leq r\cdot \frac{\theta^2}{2}
\cdot C^3 [2r^{\frac{\alpha}{2}}]
\leq r\cdot \theta^2 
\cdot \frac{1}{4\pi^2}.
\end{equation}
The desired inequality $r(1-\cos\theta) > \left\vert
F(r,\theta:s)
\right\vert$ for all $0< \vert \theta \vert \leq \pi$ is
proved.
\endproof

We now fix once for all a radius $r_0$ with $0< r_0 \leq r_1$. In the
remainder of the present Section~4, we shall deform the disc
$\widehat{ Z}_{r_0:s}(\zeta)$ by adding many more parameters. We
notice that for all $\theta$ with $0\leq \vert \theta \vert \leq
\frac{\pi}{4}$, we have the trivial minoration $\partial_\theta
\partial_\theta \widehat{ Y}_{1;r_0}\left(e^{i\theta}\right) = r_0
\cos \theta\geq \frac{ r_0}{\sqrt{2}}$. Also, by~\thetag{ 4.24} and
by the inequality on $r_1$ written in the statement of Lemma~4.25, we
deduce $\left\vert \partial_\theta \partial_\theta h_1\left( \widehat{
X}_{r_0:s}\left(e^{i\theta}\right) \right) \right\vert\leq \frac{
r_0}{2\pi^2}$ for all $\theta$ with $0\leq \vert \theta \vert \leq
\pi$. Since we shall need a generalization of Lemma~4.25 in
Lemma~4.51 below, let us remember these two interesting inequalities,
valid for $0\leq \vert \theta \vert \leq \frac{\pi}{4}$:
\def\theequation{4.28}\begin{equation}
\left\{
\aligned
\partial_\theta\partial_\theta
\widehat{ Y}_{1;r_0}\left(e^{i\theta}\right) \geq 
& \
\frac{r_0}{\sqrt{2}}, \\
\left\vert
\partial_\theta \partial_\theta h_1\left(
\widehat{ X}_{r_0:s}\left(e^{i\theta}\right)
\right)
\right\vert\leq 
& \
\frac{r_0}{2\pi^2},
\endaligned\right.
\end{equation}
noticing of course that $\frac{r_0}{2\pi^2}< \frac{r_0}{\sqrt{2}}$.

\subsection*{4.29.~Normal deformations of the disc 
$\widehat{ Z}_{r:s}(\zeta)$} So, we fix $r_0$ small with $0 < r_0 \leq
r_1$ and we consider the disc $\widehat{ Z}_{r_0:s}(\zeta)$ for
$\zeta\in\overline{\Delta}$. Then the point $\widehat{
Z}_{r_0:s}(-1)$ belongs to $(M_s^1)^+$ for each $s$ and stays at a
positive distance from $M_s^1$ as $s$ varies in $[-1,1]$. It follows
that we can choose a subneighborhood $\omega_s$ of $\widehat{
Z}_{r_0:s}(-1)$ in $\C^n$ which is contained in $\Omega$ and whose
diameter is uniformly positive with respect to $s$.

\bigskip
\begin{center}
\input normal-deformations.pstex_t
\end{center}

Following~\cite{tu2} and~\cite{mp1}, we shall introduce {\sl normal
deformations} of the analytic discs $\widehat{ Z}_{r_0 :s}(\zeta)$
parametrized by $s$ as follows. Let $\kappa: \R^{ n-1} \to \R^{n-1}$
be a $\mathcal{ C}^{2, \alpha}$-smooth mapping fixing the origin and
satisfying $\partial_{ x_k} \kappa_{j}(0)=\delta_j^k$ for $j,k=1,
\dots, n-1$, where $\delta_j^k$ denotes Kronecker's symbol. For
$j=2,\dots,n$, let $\eta_j = \eta_j( z_1, x':s)$ be a real-valued
$\mathcal{ C}^{2, \alpha}$-smooth function compactly supported in a
neighborhood of the point of $\R^{n+1}$ with coordinates $\left(
\widehat{ Z}_{ 1;r_0:s }(-1), \widehat{ X}_{ r_0:s}'(-1) \right)$ and
equal to $1$ at this point. We then define the $\mathcal{ C}^{2,
\alpha}$-smooth deformed generic submanifold $M_{ s,t}$ of equations
\def\theequation{4.30}\begin{equation}
\left\{
\aligned
y'
& \
=\varphi'(z_1,x': s) +\kappa(t) \cdot 
\eta'(z_1,x': s) \\
& \
=:\Phi'(z_1,x',t:s).
\endaligned\right.
\end{equation}
Notice that $M_{ s,0} \equiv M_s$ and that $M_{ s,t}$ coincides with
$M_s$ in a small neighborhood of the origin, for all $t$. If $\mu=
\mu \left( e^{i\theta}: s\right)$ is a real-valued nonnegative
$\mathcal{ C}^{2, \alpha}$-smooth function defined for $e^{i \theta}
\in \partial \Delta$ and for $s \in [-1,1]$ whose support is
concentrated near the segment $\{-1\} \times [-1,1]$, then applying
the existence Theorem~1.2 of~\cite{tu3}, for each fixed $s\in [-1,1]$,
we deduce the existence of a $\mathcal{ C}^{2,\alpha-0}$-smooth
solution of the Bishop type equation
\def\theequation{4.31}\begin{equation}
\widehat{ X}_{r_0,r:s}'\left(e^{i\theta}\right)=
-\left[T_1
\Phi'\left(
\widehat{ Z}_{1;r_0:s}(\cdot), 
\widehat{ X}_{r_0,t:s}'(\cdot), 
t\mu(\cdot : s):s\right)
\right]\left(e^{i\theta}\right),
\end{equation}
which enable us to construct a deformed family of analytic disc
\def\theequation{4.32}\begin{equation}
\widehat{ Z}_{r_0,t:s}\left(e^{i\theta}\right):= 
\left(
\widehat{ Z}_{1;r_0:s}(e^{i\theta}), 
\widehat{ X}_{r_0,t:s}'(e^{i\theta})+i
T_1\left[\widehat{ X}_{r_0,t:s}'(\cdot)\right]
\left(e^{i\theta}\right)
\right)
\end{equation}
whose boundaries are contained in $M\cup \omega_s$, by construction.
By an inspection of Theorem~1.2 in~\cite{tu3}, taking account of the
regularity properties {\bf (a)}, {\bf (b)} and {\bf (c)} stated
after~\thetag{4.16}, one
can show that
the general solution $\widehat{ Z}_{r_0,t:s}(\zeta)$ enjoys regularity
properties which are completely similar:
\begin{itemize}
\item[{\bf (a)}]
For fixed $s$, it is of class $\mathcal{ C}^{2,\alpha-0}$ with respect
to $(t,\zeta)$.
\item[{\bf (b)}]
It is of class $\mathcal{ C}^{1,\alpha-0}$ with respect to 
all the variables $(t,\zeta,s)$.
\item[{\bf (c)}]
Each of its first order partial derivative
with respect to the principal variables
$(t,\zeta)$ is of class $\mathcal{ C}^{1,\alpha-0}$
with respect to all the variables $(t,\zeta,s)$.
\end{itemize}

Since the solution is $\mathcal{ C}^{1,\alpha-0}$-smooth
with respect to $s$, it crucially follows that 
the vector
\def\theequation{4.33}\begin{equation}
v_{1:s}:= - \frac{\partial \widehat{ Z}_{r_0,t:s}}{\partial \rho}(1),
\end{equation}
which points inside the analytic disc, varies continuously with
respect to $s$. The following key proposition may be established (up
to a change of notations) just by reproducing the proof of Lemma~2.7
in~\cite{mp1}, taking account of the uniformity of all
differentiations with respect to the curve parameter $s$.

\def\thelemma{4.34}\begin{lemma}
There exists a real-valued nonnegative $\mathcal{
C}^{2,\alpha}$-smooth function $\mu=\mu(e^{i\theta}:s)$ defined for
$e^{i\theta}\in\partial \Delta$ and $s\in [-1,1]$ of support
concentrated near $\{-1\}\times [-1,1]$ such that the
mapping
\def\theequation{4.35}\begin{equation}
\R^{n-1}\ni t \longmapsto 
\left.\frac{\partial \widehat{ X}_{r_0,t:s}'}{\partial 
\theta}\left(e^{i\theta}\right)\right\vert_{
\theta=0}\in \R^{n-1}
\end{equation}
is maximal equal to $(n-1)$ at $t=0$.
\end{lemma}

Geometrically speaking, since the vector $\left. \frac{ \partial
\widehat{ X}_{1;r_0:s }}{ \partial \theta} \left(e^{i\theta} \right)
\right\vert_{ \theta=0}$ is nonzero, it follows that
when the parameter $t$ varies, then the set of lines
generated by the vectors $\left.\frac{ \partial \widehat{ X}_{ r_0,
t:s}}{\partial \theta}\left(e^{i \theta}\right) \right \vert_{
\theta=0}$ covers an open cone in the space $T_{ p_1}M^1\equiv \R^n$
equipped with coordinates $(x_1,x' )$, {\it see} again {\sc Figure~7}
above for an illustration.

\subsection*{4.36.~Adding pivoting and translation parameters}
Let $\chi=(\chi_1,\chi')\in \R\times \R^{n-1}$ 
and $\nu \in \R$ satisfying 
$\vert \chi \vert < \varepsilon$ and $\vert \nu \vert < \varepsilon$
for some small $\varepsilon>0$. Then the mapping
\def\theequation{4.37}\begin{equation}
\left\{
\aligned
\R^{n+1}\ni
(\chi_1,\chi',\nu)\longmapsto 
& \
\left( \,
\chi_1+i[h_1(\chi:s)+\nu], \chi'+i\varphi'
(\chi, h_1(\chi:s)+\nu:s) \,
\right) \\
& \ 
=: \widehat{p}(\chi,\nu:s)\in M_s
\endaligned\right.
\end{equation}
is a $\mathcal{ C}^{2,\alpha}$-smooth diffeomorphism
onto a neighborhood of the origin in $M_s$ with the
property that 
\begin{itemize}
\item[{\bf (a)}]
$\nu >0$ if and only if $\widehat{ p}(\chi,\nu:s)\in (M_s^1)^+$.
\item[{\bf (b)}]
$\nu=0$ if and only if $\widehat{ p}(\chi,\nu:s)\in M_s^1$.
\item[{\bf (c)}]
$\nu<0$ if and only if $\widehat{ p}(\chi,\nu:s)\in (M_s^1)^-$.
\end{itemize}
If $\tau\in \R$ with $\vert \tau \vert < \varepsilon$ is a
supplementary parameter, we may now define a crucial 
deformation of the first
component $\widehat{ Z}_{1; r_0:s}\left(e^{i \theta} \right)$ by setting
\def\theequation{4.38}\begin{equation}
\widehat{ Z}_{1; r_0,\tau,\chi,\nu:s}\left(e^{i\theta}\right):=
ir_0\left(1-e^{i\theta}\right)[1+i\tau]+\chi_1+
i[h_1(\chi:s)+\nu].
\end{equation}
Of course, we have $\widehat{ Z}_{ 1; r_0, 0,0,0:s} \left( e^{i
\theta} \right) \equiv \widehat{ Z}_{ 1;r_0 :s} \left( e^{i \theta }
\right)$. Geometrically speaking, this perturbation corresponds to
add firstly a small ``rotation parameter'' $\tau$ which rotates (and
slightly dilates) the disc $ir_0\left(1-e^{i\theta}\right)$ passing
through the origin in $\C_{ z_1}$, to add secondly a small
``translation parameter $(\chi_1, \chi')$ which will enable to cover a
neighborhood of the origin in $M_s^1$ and to add thirdly a small
translation parameter $\nu$ along the $y_1$-axis. Consequently, with
this first component $\widehat{ Z}_{1; r_0, \tau, \chi, \nu :s} \left(
e^{i \theta}\right)$, we can construct a $\C^n$-valued analytic disc
$\widehat{ Z}_{r_0, t, \tau, \chi, \nu: s} (\zeta)$ satisfying the
important property
\def\theequation{4.39}\begin{equation}
\widehat{ Z}_{r_0,t,\tau,\chi,\nu:s}(1)=\widehat{ p}(\chi,\nu:s),
\end{equation}
simply by solving the perturbed Bishop type equation
which extends~\thetag{ 4.31}
\def\theequation{4.40}\begin{equation}
\widehat{ X}_{r_0,t,\tau,\chi,\nu:s}'
\left(e^{i\theta}\right)=-
\left[
T_1\left(
\Phi'\left(\widehat{ Z}_{1;r_0,\tau,\chi,\nu:s}(\cdot), 
\widehat{ X}_{r_0,t,
\tau,\chi,\nu:s}'(\cdot), 
t\mu(\cdot:s): s\right)
\right)
\right]\left(e^{i\theta}\right).
\end{equation}
Of course, thanks to the sympathetic stability of Bishop's equation
under perturbation, the solution exists and satisfies smoothness
properties entirely similar to the ones stated after~\thetag{4.32}.
We can summarize the description of our final family of analytic discs
\def\theequation{4.41}\begin{equation}
\widehat{ Z}_{r_0,t,\tau,\chi,\nu:s}(\zeta):
\left\{
\aligned
r_0 = 
& \ 
\text{\rm approximate radius}. \\
t =
& \
\text{\rm normal deformation parameter}. \\
\tau =
& \
\text{\rm pivoting parameter}. \\
\chi =
& \
\text{\rm parameter of translation along} \ M^1. \\
\nu =
& \
\text{\rm parameter of translation in} \ M \ 
\text{\rm transversally to} \ M^1. \\
s =
& \
\text{\rm parameter of the characteristic curve} \ \gamma. \\
\zeta = 
& \ 
\text{\rm unit disc variable}.
\endaligned\right.
\end{equation}

For every $t$ and every $\chi$, we now want to adjust the pivoting
parameter $\tau$ in order that the disc boundary $\widehat{ Z}_{r_0,
t, \tau, \chi, 0:s} \left( e^{i\theta}\right)$ for $\nu=0$ is tangent
to $M_s^1$. This tangency condition will be useful in order to derive
the crucial Lemma~4.51 below.

\def\thelemma{4.42}\begin{lemma}
Shrinking $\varepsilon$ if necessary, there exists a unique $\mathcal{
C}^{1,\alpha-0}$-smooth map $(t,\chi,s) \mapsto \tau(t, \chi:s)$
defined for $\vert t \vert <\varepsilon$, for $\vert \chi \vert <
\varepsilon$ and for $s \in [-1,1]$ satisfying $\tau(0,
0:s)=\partial_{ t_j} \tau( 0,0:s)= \partial_{ \chi_k}\tau (0,0:s)=0$
for $j=1, \dots, n-1$ and $k=1, \dots,n$, such that the vector
\def\theequation{4.43}\begin{equation}
\left.
\frac{\partial}{\partial \theta}\right\vert_{\theta=0}
\widehat{ Z}_{r_0,t,\tau(t,\chi:s),\chi,0:s}\left(e^{i\theta}\right)
\end{equation}
is tangent to $M_s^1$ at the point $\widehat{ Z}_{ r_0, t, \tau(t,
\chi:s), \chi, 0:s}(1)= \widehat{ p} ( \chi, 0:s) \in M_s^1$.
\end{lemma}

\proof
We remind that $M_s$ is represented by the $(n-1)$ scalar equations
$y'= \varphi'( x, y_1: s)$ and that $M_s^1$ is represented by the $n$
equations $y_1=h_1(x:s)$ and $y'=\varphi '(x, h_1 (x: s):s)\equiv
h'(x':s)$. We can therefore compute the Cartesian equations of the
tangent plane to $M_s^1$ at the point $\widehat{p}(\chi,0:s)=
\chi+ih(\chi:s))$:
\def\theequation{4.44}\begin{equation}
\left\{
\aligned
{\sf Y}_1 -h_1(\chi:s)=
& \
\sum_{k=1}^n\, 
\partial_{x_k}h_{1}(\chi:s) \, \left[
{\sf X}_x-\chi_k\right], \\
{\sf Y}'-\varphi'(\chi,h_1(\chi:s):s) =
& \
\sum_{k=1}^n\, 
\left(
\partial_{x_k}\varphi'+
\partial_{y_1}\varphi'\cdot
\partial_{x_k}h_{1}
\right)
\left[
{\sf X}_k-\chi_k
\right].
\endaligned\right.
\end{equation}
On the other hand, we observe that the tangent vector
\def\theequation{4.45}\begin{equation}
\left. \frac{\partial }{\partial \theta}\right\vert_{\theta=0}
\widehat{ Z}_{r_0,t,\tau,\chi,0:s}\left(e^{i\theta}\right)= \left(
r_0[1+i\tau], \ \left. \frac{\partial }{\partial
\theta}\right\vert_{\theta=0} \widehat{ Z}_{r_0,t,\tau,\chi,0:s}'
\left(e^{i\theta}\right) \right)
\end{equation}
is already tangent to $M_s$ at the point $\widehat{ p}(\chi,0:s)$, 
because $M_{s,t}\equiv M_s$ in a neighborhood of the origin.
More precisely, since $\Phi'\equiv \varphi'$
in a neighborhood of the origin, 
we may differentiate with respect to 
$\theta$ at $\theta=0$ the relation
\def\theequation{4.46}\begin{equation}
\widehat{ Y}_{r_0,t,\tau,\chi,0:s}'\left(e^{i\theta}\right)\equiv 
\varphi'\left(\widehat{ X}_{r_0,\tau,\chi,0:s}\left(e^{i\theta}\right), 
\widehat{ Y}_{1;r_0,\tau,\chi,0:s}\left(e^{i\theta}\right),
\right)
\end{equation}
which is valid for $e^{i\theta}$ close to
$1$ in $\partial \Delta$, noticing in advance that it
follows immediately from~\thetag{4.38} that
\def\theequation{4.47}\begin{equation}
\left.
\frac{\partial }{\partial 
\theta} \right\vert_{\theta=0}
\widehat{ X}_{1;r_0,\tau,\chi,0:s}
\left(e^{i\theta}\right)=r_0 \ \ \ \ \ 
{\rm and} \ \ \ \ \ 
\left.
\frac{\partial }{\partial 
\theta} \right\vert_{\theta=0}
\widehat{ Y}_{1;r_0,\tau,\chi,0:s}
\left(e^{i\theta}\right)=r_0\tau,
\end{equation}
hence
we obtain by a direct application of the chain rule
\def\theequation{4.48}\begin{equation}
\left\{
\aligned
\left.
\frac{\partial }{\partial \theta} 
\right\vert_{\theta=0}\widehat{ Y}_{r_0,t,
\tau,\chi,0:s}'\left(e^{i\theta}\right)=
& \
\partial_{y_1}\varphi'\cdot r_0\tau+ \\
& \
+
\sum_{k=1}^n\,
\partial_{x_k}\varphi'\cdot \left(
\left.
\frac{\partial }{\partial \theta} 
\right\vert_{\theta=0}
\widehat{ X}_{k;r_0,t,\tau,\chi,0:s}
\left(e^{i\theta}\right)
\right).
\endaligned\right.
\end{equation}
On the other hand, the vector~\thetag{4.45} belongs to the tangent
plane to $M_s^1$ at $\widehat{ p}(\chi,0:s)$ whose equations are
computed in~\thetag{4.44} if and only if the following two conditions
are satisfied
\def\theequation{4.49}\begin{equation}
\left\{
\aligned
{}
&
r_0\tau=
\sum_{k=1}^n\, 
\partial_{x_k}h_1(\chi:s)\left[
\left.
\frac{\partial }{\partial \theta} 
\right\vert_{\theta=0}
\widehat{ X}_{k;r_0,t,\tau,
\chi,0:s}\left(e^{i\theta}\right)
\right], \\
&
\left.
\frac{\partial }{\partial \theta} 
\right\vert_{\theta=0}
\widehat{ Y}_{r_0,t,\tau,\chi,0:s}'
\left(e^{i\theta}\right)= \\
&
\ \ \ \ \ \ \ \ \ \ \ \ 
=
\sum_{k=1}^n\, 
\left(
\partial_{x_k}\varphi'+\partial_{y_1}
\varphi'\cdot \partial_{x_k}
h_1
\right)\cdot \left[
\left.
\frac{\partial }{\partial \theta} 
\right\vert_{\theta=0}
\widehat{ X}_{k;r_0,t,\tau,\chi,0:s}
\left(e^{i\theta}\right) \right].
\endaligned\right.
\end{equation}
We observe that the first line of~\thetag{4.49} together with the
relation~\thetag{4.48} already obtained implies the second line
of~\thetag{4.49} by an obvious linear combination. Consequently, the
vector~\thetag{4.45} belongs to the tangent plane to $M_s^1$ at
$\widehat{ p}( \chi, 0:s)$ if and only if the first line
of~\thetag{4.49} is satisfied. As $r_0$ is nonzero, as the first
order derivative $\partial_{x_k}h_1 (\chi :s)$ are of class $\mathcal{
C}^{1, \alpha}$ and vanish at $x=0$ and as $\left. \frac{\partial
}{\partial \theta} \right \vert_{ \theta=0} \widehat{ X}_{k; r_0, t,
\tau, \chi,0:s} \left(e^{ i \theta } \right)$ is of class $\mathcal{
C}^{1, \alpha-0}$ with respect to all variables $(t, \tau, \chi,s)$,
it follows from the implicit function theorem that there exists a
unique solution $\tau= \tau(t, \chi: s)$ of the first line
of~\thetag{4.49} which satisfies in addition the normalization
conditions $\tau( 0,0: s)= \partial_{ t_j}\tau ( 0, 0:s)= \partial_{
\chi_k}\tau (0, 0:s)=0$ for $j=1, \dots, n-1$ and $k=1, \dots, n$.
This completes the proof of Lemma~4.42.
\endproof

We now define the analytic disc 
\def\theequation{4.50}\begin{equation}
\widehat{\mathcal{ Z}}_{t,\chi,\nu:s}(\zeta):=
\widehat{ Z}_{r_0,t,\tau(t,\chi:s),\chi,\nu:s}(\zeta).
\end{equation}

\def\thelemma{4.51}\begin{lemma}
Shrinking $\varepsilon$ if necessary, 
the following two properties are satisfied{\rm :}
\begin{itemize}
\item[{\bf (1)}]
$\widehat{\mathcal{ Z}}_{
t,\chi,0:s}(\partial \Delta \backslash \{1\}) \subset (M_s^1)^+$
for all $t$, $\chi$, $\nu$ and $s$ with $\vert t \vert< \varepsilon$,
with
$\vert \chi \vert < \varepsilon$, with
$\vert \nu \vert < \varepsilon$ and with
$-1\leq s \leq 1$.
\item[{\bf (2)}]
If $\nu$
satisfies $0< \nu < \varepsilon$, then 
$\widehat{\mathcal{ Z}}_{t,\chi,\nu:s}(\partial \Delta)\subset
(M_s^1)^+$ for all $t$, $\chi$ and $s$ with 
$\vert t \vert <\varepsilon$, with $\vert \chi \vert
< \varepsilon$ and with $-1\leq s \leq 1$.
\end{itemize}
\end{lemma}

\proof
To establish property {\bf (1)}, we first observe that the disc
$\widehat{\mathcal{ Z}}_{0, 0,0 :s} \left(e^{i \theta}\right)$
identifies with the disc $\widehat{ Z}_{ r_0 :s}(e^{i \theta})$
defined in \S4.29. According to Lemma~4.25, we know that
$\widehat{\mathcal{ Z}}_{0,0,0:s}(\partial\Delta \backslash\{1\})$ is
contained in $(M_s^1)^+$. By continuity, if $\varepsilon$ is
sufficiently small, we can assume that for all $t$ with $\vert t\vert
< \varepsilon$, for all $\chi$ with $\vert \chi \vert < \varepsilon$
and for all $\theta$ with $\frac{\pi}{4}\leq \vert \theta \vert \leq
\pi$, the point $\widehat{\mathcal{ Z}}_{ t,\chi,0:s}\left(
e^{i\theta}\right)$ is contained in $(M_s^1)^+$. It remains to
control the part of $\partial \Delta$ which corresponds to $\vert
\theta \vert \leq \frac{\pi}{4}$.

Since the disc $\widehat{\mathcal{ Z}}_{
t, \chi, \nu:s}\left(e^{i \theta}\right)$ is of
class $\mathcal{ C}^2$ with respect to all its principal variables
$\left(t, \chi, \nu,e^{i\theta}\right)$, 
if $\vert t \vert < \varepsilon$, if
$\vert \chi \vert < \varepsilon$ and if $0 \leq \vert \theta \vert
\leq \frac{\pi }{4}$, for sufficiently small $\varepsilon$, then the
inequalities~\thetag{4.28} are just perturbed a little bit, so we can
assume that
\def\theequation{4.52}\begin{equation}
\left\{
\aligned
\partial_\theta\partial_\theta
\widehat{\mathcal{ Y}}_{1; t,\chi,0:s}\left(e^{i\theta}\right) \geq
& \
r_0, \\
\left\vert
\partial_\theta \partial_\theta
h_1\left(
\widehat{\mathcal{ X}}_{t,\chi,0:s}\left(e^{i\theta}\right)
\right)\right\vert \leq 
& \
\frac{r_0}{2}.
\endaligned\right.
\end{equation}
We claim that the inequality 
\def\theequation{4.53}\begin{equation}
\widehat{\mathcal{ Y}}_{1;t,\chi,0:s}\left(e^{i\theta}\right) >
\left\vert
h_1\left(\widehat{\mathcal{ X}}_{t,\chi,0:s}\left(e^{i\theta}\right)
\right)
\right\vert
\end{equation}
holds for all $0<\vert \theta \vert \leq \frac{\pi}{4}$, which 
will complete the proof of property {\bf (1)}. 

Indeed, we first remind that the tangency to $M_s^1$ of the vector
$\left. \frac{ \partial }{\partial \theta} \right \vert_{ \theta=0}
\widehat{\mathcal{ Z}}_{ t,\chi,0:s}\left(e^{ i \theta}\right)$ at the
point $\widehat{ p} ( \chi,0:s)$ is equivalent to the first
relation~\thetag{4.49}, which may be rewritten in terms of the
components of the disc $\widehat{\mathcal{ Z}}_{ t,\chi,0:s}\left(
e^{i\theta}\right)$ as follows
\def\theequation{4.54}\begin{equation}
\partial_\theta \widehat{\mathcal{ Y}}_{1;t,\chi,0:s}(1)=
\sum_{k=1}^n\, [\partial_{x_k} \, 
h_1]\left(\widehat{\mathcal{ X}}_{t,\chi,0:s}(1)\right)\cdot
\partial_\theta \widehat{\mathcal{ X}}_{k;t,\chi,0:s}(1).
\end{equation}
Substracting this relation from~\thetag{4.53} and
substracting also the relation $\widehat{ Y}_{1;t,\chi,0:s}(1)=
h_1\left(\mathcal{ X}_{t,\chi,0:s}(1)
\right)$, we see that it suffices to establish that
for all $\theta$ with $0< \vert \theta \vert \leq \frac{\pi}{4}$, 
we have the strict inequality
\def\theequation{4.55}\begin{equation}
\left\{
\aligned
{}
&
\widehat{\mathcal{ Y}}_{1;t,\chi,0:s}\left(e^{i\theta}\right)-
\widehat{\mathcal{ Y}}_{1;t,\chi,0:s}(1)-
\theta \cdot \partial_\theta 
\widehat{\mathcal{ Y}}_{1;t,\chi,0:s}(1)
>\\
& \
>
\left\vert
h_1\left(
\widehat{\mathcal{ X}}_{t,\chi,0:s}\left(e^{i\theta}\right)
\right)
-h_1\left(\widehat{\mathcal{ X}}_{t,\chi,0:s}(1)
\right) - \right. \\
& 
\ \ \ \ \ \ \ \ \ \ \ \ \ \ \ \ \ \ \ \
\left. 
-\sum_{k=1}^n\, [\partial_{x_k} \, 
h_1]\left(\widehat{\mathcal{ X}}_{t,\chi,0:s}(1)\right)\cdot
\partial_\theta \widehat{\mathcal{ X}}_{k;t,\chi,0:s}(1)
\right\vert
\endaligned\right.
\end{equation}
However, by means of Taylor's integral
formula, this last inequality may be rewritten 
as
\def\theequation{4.56}\begin{equation}
\small
\int_0^\theta \, 
(\theta -\theta')\cdot
\partial_\theta\partial_\theta
\widehat{\mathcal{ Y}}_{1;t,\chi,0:s}
\left(e^{i\theta'}\right) \cdot d\theta' >
\left\vert
\int_0^\theta \, 
(\theta-\theta') \cdot
\partial_\theta \partial_\theta
\left[
h_1\left(
\widehat{\mathcal{ X}}_{t,\chi,0:s}\left(e^{i\theta'}\right)
\right)
\right]\cdot d\theta'
\right\vert
\end{equation}
and it follows immediately by means of~\thetag{4.52}.

Secondly, to check property {\bf (2)}, we observe that by
the definition~\thetag{4.37}, the parameter $\nu$ corresponds to a
translation of the $z_1$-component of the disc boundary $\widehat{
\mathcal{ Z}}_{ t, \chi, 0:s} (\partial \Delta)$ along the $y_1$ axis.
More precisely, we have
\def\theequation{4.57}\begin{equation}
\frac{\partial }{\partial \nu}\, 
\widehat{ \mathcal{ Y}}_{1; t,\chi,\nu:s}(\zeta) \equiv 1, 
\ \ \ \ \ \ \ \ \ \ \ \ \ \ \ \ 
\frac{\partial }{\partial \nu}
\widehat{ \mathcal{ X}}_{1; t,\chi,\nu:s}(\zeta)\equiv 0.
\end{equation}
On the other hand, differentiating Bishop's equation~\thetag{4.40},
and using the smallness of the function $\Phi'$, it may be checked
that
\def\theequation{4.58}\begin{equation}
\left\vert 
\frac{\partial }{\partial \nu}
\widehat{ Z}_{r_0,t,\tau,\chi,\nu:s}'(e^{i\theta})
\right\vert < < 1, 
\end{equation}
if $r_0$ and $\varepsilon$ are sufficiently small. It follows that
the disc boundary $\widehat{ \mathcal{ Z}}_{t, \chi, \nu:s}( \partial
\Delta)$ is globally moved in the direction of the $y_1$-axis as
$\nu>0$ increases, hence is contained in $(M_s^1 )^+$.

The proof of Lemma~4.51 is complete.
\endproof

\subsection*{4.59.~Local half-wedges}
As a consequence of Lemma~4.34, of~\thetag{4.39} and of property {\bf
(2)} of Lemma~4.51, we conclude that for every $s\in [-1,1]$, our
discs $\widehat{\mathcal{ Z}}_{ t,\chi,\nu:s}(\zeta)$ satisfy all the
requirements {\bf (i)}, {\bf (ii)} and {\bf (iii)} of \S4.2 insuring
that the set defined by
\def\theequation{4.60}\begin{equation}
\mathcal{HW}_s^+:=\left\{
\widehat{\mathcal{ Z}}_{t,\chi,\nu:s}(\rho):\, 
\vert t\vert < \varepsilon, \ 
\vert \chi \vert < \varepsilon, \
0 < \nu < \varepsilon, \
1-\varepsilon < \rho < 1 
\right\}
\end{equation}
is a local half-wedge of edge $(M_s^1)^+$ at the origin
in the $\widehat{ z}$-coordinates, which corresponds 
to the point $\gamma(s)$ in the $z$-coordinates. Coming back to 
the $z$ coordinates, we may define the family of analytic discs
\def\theequation{4.61}\begin{equation}
\mathcal{ Z}_{t,\chi,\nu:s}(\zeta):= 
\Phi_s^{-1}\left(\widehat{ \mathcal{ Z}}_{t,\chi,\nu:s}(\zeta)\right)
\end{equation}
and we shall construct the desired semi-local attached
half-wedge $\mathcal{ HW}_\gamma^+$ of Proposition~4.6.

\subsection*{4.62.~Holomorphic extension to an 
attached half-wedge} Indeed, we can now complete the proof of
Proposition~4.6. Let us denote by $\widehat{z}=\Phi_s(z)$ the
parametrized change of coordinates defined in \S4.11, where the point
$\gamma(s)$ in $z$-coordinates corresponds to the origin in
$\widehat{z}$-coordinates. Given an arbitrary holomorphic function
$f\in \mathcal{ O}(\Omega)$ as in Proposition~4.6, by the change of
coordinates $\widehat{ z}=\Phi_s(z)$ and by restriction to
$(M_s^1)^+$, we get a CR function $\widehat{f}_s\in \mathcal{
C}_{CR}^0\left((M_s^1)^+\cap U_1\right)$, for some small neighborhood
$U_1$ of the origin in $\C^n$, whose size is uniform with respect to
$s$. Thanks to an obvious generalization of the approximation
Lemma~4.8 with a supplementary parameter $s\in [-1,1]$, we know that
there exists a second uniform neighborhood $V_1\subset \subset U_1$ of
the origin in $\C^n$ such that every continuous CR function in
$\mathcal{ C}_{CR}^0\left( (M_s^1)^+ \cap U_1 \right)$ is uniformly
approximable by polynomials on $(M_s^1)^+ \cap V_1$. In particular,
this property holds for the CR function $\widehat{f}_s$. Furthermore,
choosing $r_0$ and $\varepsilon$ sufficiently small, we can insure
that all the discs $\widehat{\mathcal{ Z}}_{t, 
\chi,\nu:s}(\zeta)$ are attached
to $(M_s^1)^+ \cap V_1$. It then follows from the maximum principle
applied to the approximating sequence of polynomials that for each
$s\in [-1,1]$, the function $\widehat{f}_s$ extends holomorphically to
the half-wedge defined by~\thetag{4.60}. Finally, we deduce that the
holomorphic function $f\in \mathcal{ O}(\Omega)$ extends
holomorphically to the half-wedge attached to the one-sided
neighborhood $(M_\gamma^1)^+$ defined by
\def\theequation{4.63}\begin{equation}
\left\{
\aligned
\mathcal{HW}_\gamma^+:=\left\{
\mathcal{ Z}_{t,\chi,\nu:s}(\rho):\, 
\vert t\vert < \varepsilon, \ 
\vert \chi \vert < \varepsilon, \ 
0< \nu < \varepsilon, \right. \\
\left. 
1-\varepsilon < \rho < 1, \ 
-1\leq s \leq s
\right\}.
\endaligned\right.
\end{equation}
Without shrinking $\Omega$ near the points $\mathcal{ Z}_{t, \chi,
\nu:s}(-1)$ (otherwise, the crucial rank property of Lemma~4.34 would
degenerate), we can shrink the open set $\Omega$ in a very thin
neighborhood of the characteristic segment $\gamma$ in $M$ and we can
shrink $\varepsilon >0$ if necessary in order that the intersection
$\Omega \cap \mathcal{ HW}_\gamma^+$ is connected. By the principle of
analytic continuation, this implies that there exists a well-defined
holomorphic function $F\in \mathcal{ O}\left( \Omega \cup \mathcal{ HW
}_\gamma^+ \right )$ with $F \vert_\Omega = f$.

The proof of Proposition~4.6 is complete.
\qed

\subsection*{4.64.~Local half-wedge in CR dimension $m\geq 2$}
Repeating the above constructions in the simpler case where the
curve $\gamma$ degenerates to the point $p_1$ of Theorem~3.22 and
adding some supplementary parameters along the $(m-1)$ remaining
complex tangent directions of $M$ at $p_1$ for the constructions of
analytic discs, we can show that under the assumptions of
Theorem~3.22, there exists a local half-wedge $\mathcal{ HW}_{p_1}^+$ 
at $p_1$ to which (shrinking $\Omega$ if
necessary) every holomorphic function $f\in \mathcal{ O}(\Omega)$
extends holomorphically. We shall not write down the details. 

\subsection*{4.65.~Transition} In the
case $m=2$, using such a local half-wedge and applying the continuity
principle along analytic discs of whose one boundary part is contained
in $M$ and whose second boundary part is contained in the local
half-wedge $\mathcal{ HW}_{p_1}^+$, we shall establish that $p_1$ is
$\mathcal{ W}$-removable in Section~10 below. At present, in the more
delicate case $m=1$, we shall pursue in Section~5 below our geometric
constructions for the choice of a special point $p_{\rm sp}\in C$
which satisfies the conclusion of Theorem~3.19 {\bf (i)}.

\section*{\S5.~Choice of a special point of $C_{\rm nr}$ 
to be removed locally}

\subsection*{5.1.~Choice of a first supporting hypersurface} 
Continuing with the proof of Theorem~3.19 {\bf (i)}, we shall now
analyze and use the important geometric 
condition $\mathcal{ F}_{M^1}^c\{C\}$
defined in Theorem~1.2'. We first delineate a convenient geometric
situation.

\def\thelemma{5.2}\begin{lemma} 
Under the assumptions of Theorem~3.19 {\bf (i)}, there exists a
$\mathcal{ C}^{2, \alpha}$-smooth segment $\gamma: [-1,1]\to M^1$
running in characteristic directions, namely satisfying $d \gamma(s)/
ds\in T_{ \gamma(s)}^c M\cap T_{\gamma (s)}M^1 \backslash \{0\}$ such
that $\gamma (-1) \not \in C$, $\gamma(0) \in C$, $\gamma (1) \not \in
C$, and there exists a $\mathcal{ C}^{ 1,\alpha}$-smooth hypersurface
$H^1$ of $M^1$ with $\gamma\subset H^1$ which is foliated by
characteristic segments close to $\gamma$, such that locally in a
neighborhood of $H^1$, the closed subset $C$ is contained in
$\gamma\cup (H^1)^-$, where $(H^1)^-$ denotes an open 
one-sided neighborhood of $H^1$ in $M^1$.
\end{lemma}

\proof
By the assumption $\mathcal{ F}_{M^1}^c\{C\}$, there exists a
characteristic curve $\widetilde{\gamma}: [-1,1]\to M^1$ with
$\widetilde{\gamma}(-1)\not\in C$, $\widetilde{\gamma}(0)\in C$ and
$\widetilde{\gamma}(1)\not\in C$, there exists a neighborhood
$V_{\widetilde{\gamma}}^1$ of $\widetilde{\gamma}$ in $M^1$ and there
exists a local $(n-1)$-dimensional submanifold $R^1$ passing through
$\widetilde{\gamma}(0)$ which is transversal to $\widetilde{\gamma}$
such that the semi-local projection $\pi_{\mathcal{ F}_{M^1}^c}:
V_{\widetilde{\gamma}}^1\to R^1$ parallel to the characteristic curves
maps $C$ onto the closed subset $\pi_{\mathcal{ F}_{M^1}^c}(C)$ with
the property that $\pi_{\mathcal{ F}_{M^1}^c}(\widetilde{\gamma})$
lies on the {\it boundary of $\pi_{\mathcal{ F}_{M^1}^c}(C)$ with
respect to the topology of $R^1$}. This property is illustrated in the
right hand side of the following figure.

\bigskip
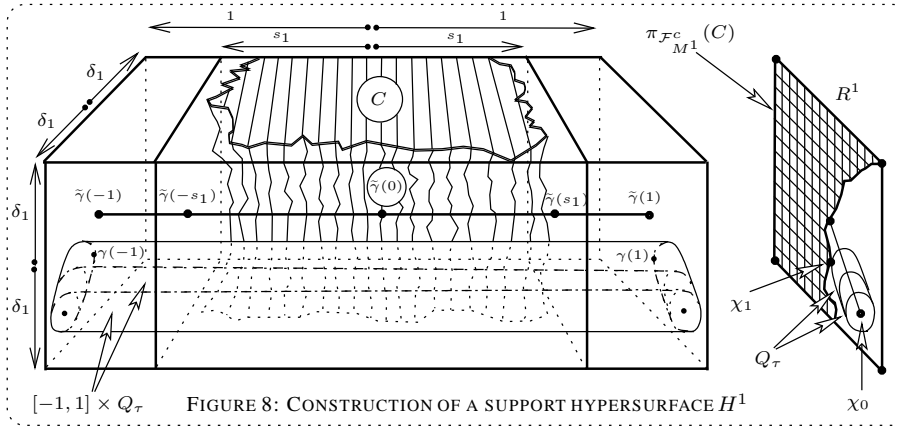
\begin{center}
\input tube-ellipse.pstex_t
\end{center}

However, we want in addition a foliated supporting hypersurface $H^1$,
which does not necessarily exist in a neighborhood of
$\widetilde{\gamma}$. To construct $H^1$, let us first straighten the
characteristic lines in a neighborhood of $\widetilde{ \gamma}$,
getting a product $[-1,1]\times [-\delta_1,\delta_1]^{n-1}$, for some
$\delta_1>0$, equipped with coordinates
$(s,\chi)=(s,\chi_2,\dots,\chi_n)\in \R\times \R^{n-1}$, so that the
lines $\{\chi={\rm ct.}\}$ correspond to characteristic lines. Such a
straightening is only of class $\mathcal{ C}^{1,\alpha}$, because the
line distribution $T^cM\vert_{M^1} \cap TM^1$ is only of class
$\mathcal{ C}^{1,\alpha}$. Clearly, we may assume that $\delta_1$ is
so small that there exists $s_1$ with $0< s_1 <1$ such that the two
cubes $[-1,-s_1]\times [-\delta_1,\delta_1]^{n-1}$ and $[s_1,1]\times
[-\delta_1,\delta_1]^{n-1}$ do not meet the singularity $C$, as drawn
in {\sc Figure~8} above.

We may identify the transversal $R^1$ with $[- \delta_1,
\delta_1]^{n-1}$; then the projection of $\widetilde{\gamma}$ 
is the origin of
$R^1$. By assumption, $\pi_{ \mathcal{ F}_{M^1}^c}(C)$ is a proper
closed subset of $R^1$ with the origin lying on its boundary. We can
therefore choose a point $\chi_0$ in the interior of $R^1$ lying
outside $\pi_{ \mathcal{ F}_{M^1 }^c }(C)$. Also, we can choose a
small open $(n- 1)$-dimensional ball $Q_0$ centered at this point
which is contained in the complement $R^1\backslash \pi_{ \mathcal{
F}_{M^1 }^c}(C)$. Furtermore, we can include this ball in a one
parameter family of $\mathcal{ C}^{1, \alpha}$-smooth domains $Q_\tau
\subset R^1$, for $\tau \geq 0$, which are parts of ellipsoids
stretched along the segment which joins the point $\chi_0$ with the
origin of $R^1$.

We then consider the tube domains $[-1,1] \times Q_\tau$ in $[-1, 1]
\times [- \delta_1, \delta_1]^{ n-1}$. Clearly, there exists the
smallest $\tau_1>0$ such that the tube $[-1,1]\times Q_{\tau_1}$ meets
the singularity $C$ on its boundary $[-1,1]\times \partial
Q_{\tau_1}$. In particular, there exists a point $\chi_1\in \partial
Q_{\tau_1}$ such that the characteristic segment $[-1,1]\times
\{\chi_1\}$ intersects $C$. Increasing a little bit the curvature of
$\partial Q_{\tau_1}$ in a neighborhood of $\chi_1$ if necessary, we
can assume that $\pi_{\mathcal{ F}_{M^1}^c}(C)\cap
\overline{Q_{\tau_1}}=\{\chi_1\}$ in a neighborhood of $\chi_1$.
Moreover, since by construction the two segments $[-1,-s_1]\times
\{\chi_1\}\cup [s_1,1]\times \{\chi_1\}$ do not meet $C$, we can
reparametrize the characteristic segment $[-1,1]\times \{\chi_1\}$ as
$\gamma: [-1,1]\to M^1$ with $\gamma(-1)\not\in C$, $\gamma(0)\in C$
and $\gamma(1)\not\in C$. Since all characteristic lines are
$\mathcal{ C}^{2,\alpha}$-smooth, we can choose the parametrization to
be of class $\mathcal{ C}^{2,\alpha}$. For the supporting hypersurface
$H^1$, it suffices to choose a piece of $[-1,1]\times \partial
Q_{\tau_1}$ near $[-1,1]\times \{\chi_1\}$. By construction, this
supporting hypersurface is only of class $\mathcal{ C}^{1,\alpha}$ and
we have that $C$ is contained in $\gamma \cup (H^1)^-$ semi-locally in
a neighborhood of $\gamma$, as desired. This completes the proof of
Lemma~5.2.
\endproof

\subsection*{5.3.~Field of cones on $M^1$}
With the characteristic segment $\gamma$ constructed in Lemma~5.2, by
an application of Proposition~4.6, we deduce that there exists a
semi-local half-wedge $\mathcal{ HW}_\gamma^+$ attached to
$(M_\gamma^1)^+\cap V_\gamma$, for some neighborhood $V_\gamma$ of
$\gamma$ in $M$, to which $\mathcal{ O}(\Omega)$ extends
holomorphically.

Then, we remind that by~\thetag{4.37}, \thetag{4.39}
and~\thetag{4.50}, for all $t$ with $\vert t \vert < \varepsilon$, the
point $\widehat{ \mathcal{ Z}}_{ t,\chi,0:s}(1)$ identifies with the
point $\widehat{ p}( \chi,0:s)\in M_s^1$ defined in~\thetag{4.37}
(which is independent of $t$) and the mapping $\chi \mapsto \widehat{
\mathcal{ Z}}_{t,\chi,0:s}(1)\in M_s^1$ is a local diffeomorphism.

Sometimes in the sequel, we shall denote the disc $\mathcal{ Z}_{t,
\chi, \nu:s} (\zeta) \equiv \Phi_s^{-1} \left( \widehat{ \mathcal{
Z}}_{t,\chi, \nu:s}(\zeta) \right)$ defined in~\thetag{4.61} by
$\mathcal{ Z}_{t, \chi_1, \chi',\nu:s}(\zeta)$, where $\chi'= (
\chi_2, \dots, \chi_n) \in \R^{n-1}$. Since the characteristic curve
is directed along the $x_1$-axis, which is transversal in $T_0M_s^1$
to the space $\{(0,\chi')\}$, it follows that the mapping $(s, \chi')
\longmapsto \mathcal{ Z}_{
t,0,\chi',0:s}(1)=\Phi_s^{-1}\left(\widehat{ p}(0,\chi',0:s)\right)$
is, independently of $t$, a diffeomorphism onto its image for $s\in
[-1,1]$ and for $\chi'$ close to the origin in $\R^{n-1}$. To fix
ideas, we shall let $\chi'$ vary in the {\it closed}\, cube
$[-\varepsilon,\varepsilon]^{n-1}$ (analogously to the fact that $s$
runs in the {\it closed}\, interval $[-1,1]$) and we shall denote by
$V_\gamma^1$ the closed subset of $M^1$ which is the image of this
diffeomorphism.

At every point $p:=\mathcal{ Z}_{t,0,\chi',0:s}(1)= \mathcal{
Z}_{0,0,\chi',0:s}(1)$ of this neighborhood $V_\gamma^1$, we define an
open infinite oriented cone contained in the $n$-dimensional linear
space $T_pM^1$ by
\def\theequation{5.4}\begin{equation}
{\sf C}_p:=\R^+\cdot 
\left\{
\frac{\partial \mathcal{ Z}_{
t,0,\chi',0:s}}{\partial \theta} 
(1): \ \vert t \vert < \varepsilon
\right\}.
\end{equation}
The fact that ${\sf C}_p$ is indeed an open cone follows
from Lemma~4.34, from~\thetag{4.61} and
from the fact that $\Phi_s^{-1}$ is a biholomorphism. 
This cone contains in its interior the nonzero vector 
\def\theequation{5.5}\begin{equation}
v_p^0:= \frac{\partial \mathcal{ Z}_{0,0,\chi',0:s}}{\partial \theta}
(1)\in {\sf C}_p\subset T_pM^1 \backslash \{0\}.
\end{equation}
We shall say that the cone ${\sf C}_p$ is the {\sl cone created at $p$
by the semi-local attached half-wedge $\mathcal{ HW}_\gamma^+$}
(more precisely, by the family of analytic discs which covers
this semi-local half-wedge).

As $p$ varies, the mapping $p\mapsto {\sf C}_p$ constitutes what we
shall call a {\sl field of cones over $V_\gamma^1$}, {\it see} {\sc
Figure~2} in Section~2 above and {\sc Figure~9} just below for
illustrations.

\bigskip
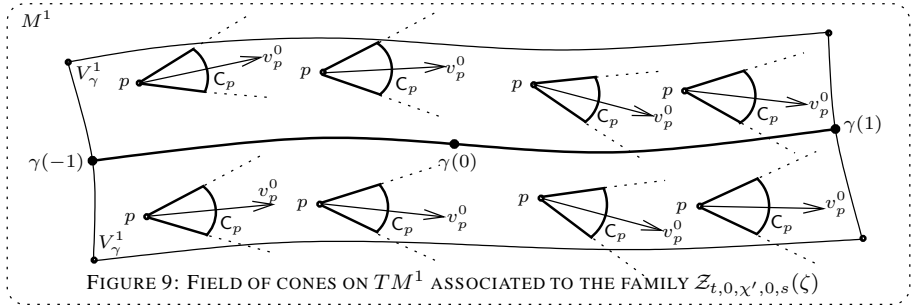
\begin{center} 
\input field-cones.pstex_t
\end{center}

The mapping $p\mapsto v_p^0$ defines a $\mathcal{
C}^{1,\alpha-0}$-smooth vector field tangent to $M^1$. This vector
field is contained in the field of cones $p\mapsto {\sf C}_p$. Over
$V_\gamma^1$, we can also consider a nowhere zero characteristic
vector field $X$ whose direction agrees with the orientation of
$\gamma$ and which satisfyies $\exp(sX)(\gamma(0))=\gamma(s)$ for all
$s\in [-1,1]$. Furthermore, for every $p\in V_\gamma^1$, we define the
{\sl filled cone}
\def\theequation{5.6}\begin{equation}
{\sf FC}_p:=\R^+\cdot \left\{
\lambda \cdot X_p+
(1-\lambda)\cdot v_p: \ 
0\leq \lambda < 1, \ 
v_p\in {\sf C}_p
\right\}.
\end{equation}
In the right hand side of {\sc Figure~10} just below, in the tangent
space $T_pM^1$ equipped with linear coordinates $(x_1,\dots,x_n)$ such
that the characteristic direction $T_p^cM\cap T_pM^1$ is directed
along the $x_1$-axis, we draw ${\sf C}_p$, its filling ${\sf FC}_p$
and its projection $\pi'({\sf C}_p)$ onto the $(x_2,,\dots,x_n)$-space
parallel to the $x_1$-axis.

\bigskip
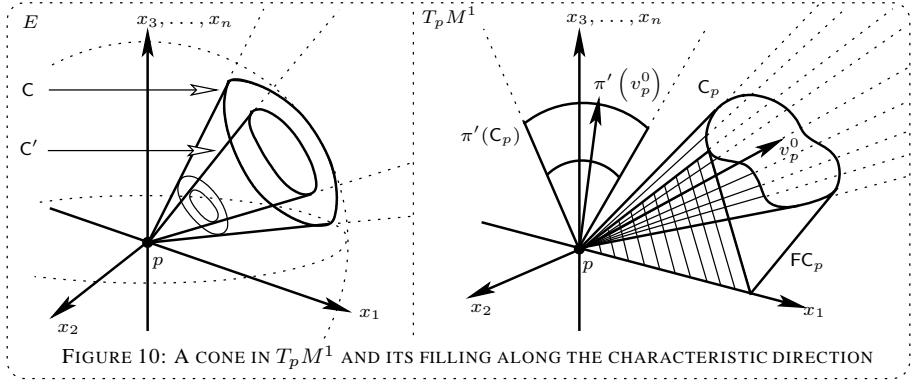
\begin{center}
\input filled-cone.pstex_t
\end{center}

Given an arbitrary nonzero vector $v_p \in {\sf C}_p$, where $p\in
V_\gamma^1$, it may be checked that a small neighborhood of the origin
in the positive half-line $\R^+\cdot Jv_p$ generated by $Jv_p$, where
$J$ denotes the complex structure of $T\C^n$, is contained in the
attached half-wedge $\mathcal{ HW}_\gamma^+$. More generally, the same
property holds for every nonzero vector $v_p$ which belongs to the
filled cone ${\sf FC}_p$. In fact, we shall establish that analytic
discs which are half-attached to $M^1$ at $p$, namely send
$\partial^+\Delta:=\{\zeta\in \partial \Delta: \ {\rm Re}\, \zeta \geq
0\}$ to $M^1$, which are closed to the complex line $v_p+Jv_p$ in
$\mathcal{ C}^1$ norm and which are sufficiently small are contained
in the attached half-wedge $\mathcal{ HW}_\gamma^+$. The interest of
this property and the reason why we have defined fields of cones and
their filling will be more appearant in Sections~8 and~9 below, where
we apply the continuity principle to achieve the proof of Theorem~3.19
{\bf (i)}. A precise statement, involving families of analytic discs
$A_c(\zeta)$ which will be constructed in Section~7 below, is as
follows. For the proof of a more precise statement involving families
of analytic discs $A_{x,v:c}^1(\zeta)$, we refer to Section~7 and
especially to Lemma~8.3 {\bf ($\mathbf{9_1}$)} below.

\def\thelemma{5.7}\begin{lemma}
Fix a point $p\in V_\gamma^1$ and a vector $v_p$ in the cone ${\sf
C}_p$ created by the semi-local attached half-wedge $\mathcal{
HW}_\gamma^+$ at $p$. Suppose that there exist two constants $c_1>0$
and $\Lambda_1 >1$ such that for every $c$ with $0< c \leq c_1$, there
exists a $\mathcal{ C}^{2,\alpha-0}$-smooth analytic disc $A_c(\zeta)$
with $A_c (\partial^+\Delta) \subset M^1$,
such that
\begin{itemize}
\item[{\bf (i)}]
The positive half-line generated by the boundary of $A_c$ at $\zeta=1$
coincides with the positive half-line generated by $v_p$, namely
$\R^+\cdot \frac{\partial A_c}{\partial \theta}(1) \equiv \R^+ \cdot
v_p$.
\item[{\bf (ii)}]
$\left\vert A_c (\zeta) \right\vert \leq c^ 2\cdot \Lambda_1$ for all
$\zeta \in \overline{ \Delta}$ and $c \cdot \frac{ 1}{ \Lambda_1} \leq
\left\vert \frac{ \partial A_c}{\partial \theta}(1) \right\vert \leq
c\cdot \Lambda_1$.
\item[{\bf (iii)}]
$\left\vert \frac{\partial A_c}{\partial \theta}\left( \rho
e^{i\theta}\right) - \frac{\partial A_c}{\partial \theta}(1)
\right\vert\leq c^2 \cdot \Lambda_1$ for all $\zeta=\rho e^{i\theta} \in
\overline{\Delta}$.
\end{itemize}
If $c_1$ is sufficiently small, then for every $c$ with $0< c \leq
c_1$, the set $A_c\left(\overline{\Delta} \backslash \partial^+\Delta
\right)$ is contained in the semi-local half-wedge $\mathcal{
HW}_\gamma^+$. Furthermore, the same conclusion holds if the nonzero
vector $v_p$ belongs to the filled cone ${\sf FC}_p$.
\end{lemma}

\subsection*{5.8.~Choice of the special point $p_{\rm sp}$ in 
the CR dimension $m=1$ case} We can now answer the question raised
after the statement of Theorem~3.19 {\bf (i)}, which was the main
purpose of the present Section~5: {\it How to choose the special point
$p_{\rm sp}$ to be removed locally}? In the following statement,
property {\bf (ii)} will be really crucial for the removal of $p_{\rm
sp}$, {\it see}\, in advance Proposition~5.12 below.

\def\thelemma{5.9}\begin{lemma}
Let $\gamma$ be the characteristic segment constructed in Lemma~5.2
above. Let $\mathcal{ HW }_\gamma^+$ be the semi-local attached
half-wedge of edge $(M_\gamma^1 )^+\cap V_\gamma$ constructed in
Proposition~4.6 above, and let $p\mapsto {\sf FC }_p$ be the filled
field of cones created in a closed neighborhood $V_\gamma^1$ of
$\gamma$ in $M^1$ by this semi-local
attached half-wedge $\mathcal{ HW }_\gamma^+$. Then
there exists a special point $p_{\rm sp}\in V_\gamma^1$ such that the
following two geometric properties are fulfilled{\rm :}
\begin{itemize}
\item[{\bf (i)}]
There exists a $\mathcal{ C}^{2, \alpha}$-smooth local supporting
hypersurface $H_{ \rm sp}$ of $M^1$ passing through $p_{\rm sp}$ such
that, locally in a neighborhood of $p_{ \rm sp}$, the closed subset
$C$ is contained in $(H_{ \rm sp})^-\cup \{p_{ \rm sp}\}$, where
$(H_{\rm sp})^-$ denotes an open one-sided neighborhood of
$H_{\rm sp}$ in $M^1$.
\item[{\bf (ii)}]
There exists a nonzero vector $v_{\rm sp}\in T_{p_{\rm sp}}H_{\rm sp}$
which belongs to the filled cone ${\sf FC}_{p_{\rm sp}}$.
\end{itemize}
\end{lemma}

\proof
According to Lemma~5.2, the singularity $C$ is contained
in $\gamma\cup (H^1)^-$, where $H^1$ is a $\mathcal{ C}^{
1,\alpha}$-smooth hypersurface containing $\gamma$ which is foliated
by characteristic segments.
If $\lambda\in [0,1)$ is very close to $1$, the vector field
over $V_\gamma^1$ defined by
\def\theequation{5.10}\begin{equation}
p\longmapsto 
v_p^\lambda:= \lambda\cdot X_p +(1-\lambda)\cdot v_p
\in T_pM^1
\end{equation}
is very close to the characteristic vector field $X_p$, so the
integral curves of $p\mapsto v_p^\lambda$ are very close to the
integral curves of $p\mapsto X_p$, which are the characteristic
segments. If $\lambda$ is sufficiently close to $1$, we can choose a
subneighborhood $V_\gamma^\lambda\subset V_\gamma^1$ of $\gamma$ which
is foliated by integral curves of $p\mapsto v_p^\lambda$. As in
Lemma~5.2, let us fix an $(n-1)$-dimensional submanifold $R^1$
transversal to $\gamma$ and passing through $\gamma(0)$. Since the
vector field $p\mapsto v_p^\lambda$ is very close to the
characteristic vector field, it follows that after projection onto
$R^1$ parallelly to the integral curves of $p\mapsto v_p^\lambda$, the
closed set $C\cap V_\gamma^\lambda$ is again a proper closed subset of
$R^1$. We notice that, by its very definition, the vector
$v_p^\lambda$ belongs to the filled cone ${\sf FC}_p$ for all $p
\in V_\gamma^\lambda$.

We can proceed exactly as in the
proof of Lemma~5.2 with the foliation of $V_\gamma^\lambda$ induced by
the integral curves of the vector field $p\mapsto v_p^\lambda$,
instead of the characteristic foliation, except that we want a
supporting hypersurface $H_{\rm sp}$ which is of class $\mathcal{
C}^{2,\alpha}$. Consequently, we first approximate the vector field
$p\mapsto v_p^\lambda$ by a new vector field $p\mapsto \widetilde{
v}_p^\lambda$ whose coefficients are of class $\mathcal{ C}^{2,
\alpha}$ (with respect to every local graphing function of $M^1$) and
which is very close to the vector field $p\mapsto v_p^\lambda$ in
$\mathcal{ C}^1$-norm. Again, we get a subneighborhood $\widetilde{
V}_\gamma^\lambda \subset V_\gamma^\lambda$ of $\gamma$ which is
foliated by integral curves of $p\mapsto \widetilde{ v}_p^\lambda$ and
a projection of $C\cap \widetilde{ V}_\gamma^\lambda$ which is a
proper closed subset of $R^1$. Moreover, if the approximation is
sufficiently sharp, we still have $\widetilde{ v}_p^\lambda \in {\sf
FC}_p$ for all $p \in \widetilde{ V}_\gamma^\lambda$. Then by
repeating the reasoning which yielded Lemma~5.2, using the foliation
of $\widetilde{ V}_\gamma^\lambda$ induced by $p \mapsto \widetilde{ v
}_p^\lambda$, we deduce that there exists an integral curve
$\widetilde{ \gamma}$ of the vector field $p\mapsto \widetilde{
v}_p^\lambda$ satisfying (after reparametrization) $\widetilde{\gamma}
(-1) \not \in C$, $\widetilde{\gamma}(0)\in C$ and
$\widetilde{\gamma}(1) \not\in C$, together with a $\mathcal{ C}^{
2,\alpha}$-smooth supporting hypersurface $\widetilde{ H}$ of
$\widetilde{ V}_\gamma^\lambda$ which contains $\widetilde{ \gamma}$
such that $C$ is contained in $\widetilde{ \gamma}\cup (\widetilde{
H})^-$. The fact that $\widetilde{ H}$ is of class $\mathcal{
C}^{2,\alpha}$ is due to the $\mathcal{ C}^{2,\alpha}$-smoothness of
the foliation on $\widetilde{ V}_\gamma^\lambda$ induced by the vector
field $p\mapsto \widetilde{ v }_\gamma^\lambda$. In {\sc Figure~11}
just below, suited for the case where $M^1$ is two-dimensional, we
have drawn as a dotted line the limiting integral curve $\widetilde{
\gamma}$ of $p\mapsto \widetilde{ v}_p^\lambda$ having the property
that $C$ lies in one closed side of $\widetilde{ \gamma}$ in
$\widetilde{ V}_\gamma^\lambda$.

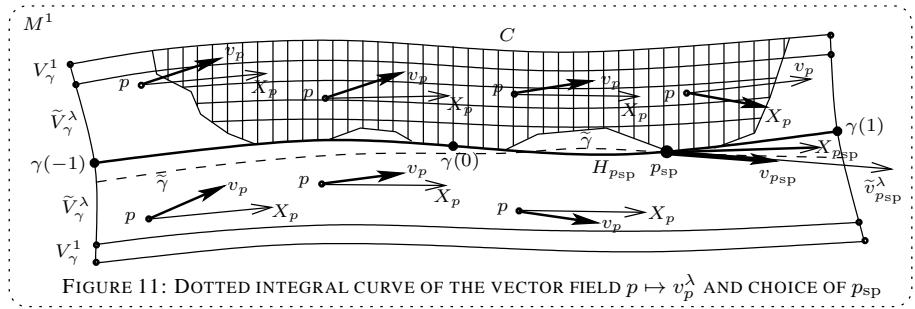
\bigskip
\begin{center}
\input perturbed-vector-field.pstex_t
\end{center}

To conclude the proof of Lemma~5.9, for the desired special point $p_{
\rm sp}$, it suffices to choose $\widetilde{ \gamma}(0)$. For the
desired local supporting hypersurface $H_{p_{\rm sp}}$, we cannot
choose directly a piece of $\widetilde{ H}$ passing through $p_{\rm
sp}$, because an open interval contained in $C\cap \widetilde{
\gamma}$ may well be contained in $\widetilde{ H}$ by the construction
in the proof of Lemma~5.2 that we have just reapplied. Fortunately,
since we know that locally in a neighborhood of $p_{\rm sp}$, the
closed subset $C$ is contained in $(\widetilde{ H})^-\cup \widetilde{
\gamma}$, it suffices to choose for the desired supporting
hypersurface $H_{p_{\rm sp}}\subset M^1$ a piece of a $\mathcal{
C}^{2,\alpha}$-smooth hypersurface passing through $p_1$, tangent to
$\widetilde{ H}$ at $p_1$ and satisfying $H_{p_{\rm sp}}\backslash
\{p_{\rm sp}\}\subset (\widetilde{ H})^+$ in a neighborhood of $p_{\rm
sp}$. Finally, for the nonzero vector $v_{\rm sp}$, it suffices to
choose any positive multiple of the vector $\widetilde{ v}_{ p_{\rm
sp}}^\lambda$. This completes the proof of Lemma~5.9.
\endproof

In Section~8 below, property {\bf (ii)} of Lemma~5.9 together with the
observation made in Lemma~5.7 will be crucial for the local $\mathcal{
W}$-removability of the special point $p_{\rm sp}$.

\subsection*{5.11.~Main removability proposition in the 
CR dimension $m=1$ case} We can now formulate the main removability
proposition to which Theorem~3.19 {\bf (i)} is now reduced. From now
on, we localize the situation at $p_{\rm sp}$, we denote this point
simply by $p_1$, we denote its supporting hypersurface simply by $H^1$
and we denote its associated vector simply by $v_1$. Furthermore, we
localize at $p_1$ the family of analytic discs considered in Section~4
for the construction of the semi-local attached half-wedge $\mathcal{
HW}_\gamma^+$, hence we get a family of analytic discs $\mathcal{
Z}_{t, \chi, \nu }( \zeta)$ which satisfy properties {\bf (i)}, {\bf
(ii)} and {\bf (iii)} of \S4.2 and which generate a local half-wedge
$\mathcal{ HW}_{p_1}^+ \subset \mathcal{ HW}_\gamma^+$ as defined
in~\thetag{4.5}. At present, we deduce from our constructions achieved
in Section~4 and in the beginning of Section~5 that Theorem~3.19 {\bf
(i)} is now a consequence of the following main Proposition~5.12 just
below, to the proof of which Sections~6, 7, 8 and~9 below are
devoted. From now on, all our geometric considerations will be
localized at the special point $p_1$.

\def\theproposition{5.12}\begin{proposition} 
Let $M$ be a $\mathcal{ C}^{2, \alpha}$-smooth generic submanifold of
$\C^n$ of CR dimension $m =1$ and of codimension $d = n-1 \geq 1$, let
$M^1\subset M$ be a $\mathcal{ C}^{2, \alpha }$-smooth
one-codimensional submanifold which is generic in $\C^n$, let $p_1 \in
M^1$, let $H^1\subset M^1$ be a $\mathcal{ C}^{2, \alpha}$-smooth
one-codimensional submanifold of $M^1$ passing through $p_1$ and let
$(H^1 )^-$ denote an open local one-sided neighborhood of $H^1$ in
$M^1$. Let $C \subset M^1$ be a {\rm nonempty} proper closed subset of
$M^1$ with $p_1 \in C$ which is situated, locally in a neighborhood of
$p_1$, only in one side of $H^1$, namely $C \subset (H^1)^- \cup \{
p_1\}$. Let $\Omega$ be a neighborhood of $M \backslash C$ in $\C^n$,
let $\mathcal{ HW}_{ p_1}^+$ be a local half-wedge of edge $(M^1)^+$
at $p_1$ generated by a family of analytic discs $\mathcal{
Z}_{t,\chi,\nu}(\zeta)$ satisfying the properties {\bf (i)}, {\bf
(ii)} and {\bf (iii)} of \S4.2, let ${\sf C}_{p_1}\subset T_{p_1}M^1$
be the cone created by $\mathcal{ HW}_{p_1}^+$ at $p_1$ and let ${\sf
FC}_{p_1}$ be its filling. {\rm As a main assumption}, suppose that
there exists a nonzero vector $v_1\in T_{p_1} H^1$ which belongs to
the filled cone ${\sf FC}_{p_1}$.
\begin{itemize}
\item[{\bf (I)}]
If $v_1$ does not belong to $T_{p_1}^cM$, then there exists a local
wedge $\mathcal{ W}_{ p_1}$ of edge $M$ at $(p_1,Jv_1)$ with
$\mathcal{ W}_{ p_1} \cap \left[\Omega \cup \mathcal{
HW}_{p_1}^+\right]$ connected {\rm (}shrinking $\Omega\cup \mathcal{
HW}_{p_1}^+$ if necessary{\rm )} such that for every holomorphic
function $f\in \mathcal{ O} \left( \Omega\cup \mathcal{ HW}_{p_1}^+
\right)$, there exists a holomorphic function $F \in \mathcal{ O}
\left( \mathcal{ W}_1 \cup \left[\Omega \cup \mathcal{ HW}_1^+
\right]\right)$ with $F\vert_{\Omega\cup \mathcal{ HW}_{p_1}^+} =f$.
\item[{\bf (II)}]
If $v_1$ belongs to $T_{p_1}^cM$, then there exists a neighborhood
$\omega_{p_1}$ of $p_1$ in $\C^n$ with $\omega_{p_1} \cap \left[\Omega
\cup \mathcal{ HW}_{p_1}^+\right]$ connected {\rm (}shrinking
$\Omega\cup \mathcal{ HW}_{p_1}^+$ if necessary{\rm )} such that for
every holomorphic function $f\in \mathcal{ O} \left( \Omega\cup
\mathcal{ HW}_{p_1}^+ \right)$, there exists a holomorphic function $F
\in \mathcal{ O} \left( \omega_{p_1} \cup \left[\Omega \cup \mathcal{
HW}_1^+ \right]\right)$ with $F\vert_{\Omega\cup \mathcal{
HW}_{p_1}^+} =f$.
\end{itemize}
\end{proposition}

In the CR dimension $m\geq 2$ case, we observe that an
analogous main removability proposition may be formulated simply by
adding to the assumptions of Proposition~3.22 a local half-wedge
$\mathcal{ HW}_{p_1}^+$, whose existence was established in \S4.64
above. The remainder of Section~5, and then Section~6, Section~7,
Section~8 and Section~9 below will be entirely devoted to the proof of
Proposition~5.12.

\subsection*{5.13.~A dichotomy}
Under the assumptions of Proposition~5.12, we shall indeed 
distinguish two cases: 
\begin{itemize}
\item[{\bf (I)}]
The nonzero vector $v_1$ does not belong to the characteristic
direction $T_{p_1}M^1\cap T_{p_1}^cM$.
\item[{\bf (II)}]
The nonzero vector $v_1$ belongs to the characteristic direction
$T_{p_1}M^1\cap T_{p_1}^cM$.
\end{itemize}

We must clarify the main assumption that $v_1$ belongs to the filling
${\sf FC}_{p_1}$
of the cone ${\sf C }_{ p_1}\subset T_{p_1} M^1$ created by the local
half-wedge $\mathcal{ HW}_{ p_1}^+$. As we have observed in \S4.2, in
the (generic) situation of Case {\bf (I)}, a local half-wedge may be
represented geometrically as the intersection of a (complete) local
wedge of edge $M$ at $p_1$, with a local one-sided neighborhood $(N^1
)^+$ of a hypersurface $N^1$ passing through $p_1$, which is
transversal to $M$ and which satisfies $N^1 \cap M \equiv M^1$ in a
neighborhood of $p_1$. The slope of the tangent space $T_{p_1} N^1$ to
$N^1$ at $p_1$ with respect to the tangent space $T_{p_1} M$ to $M$ at
$p_1$ may be understood in terms of the cone ${\sf C}_{p_1}$, as we
will now explain. Afterwards, we shall consider Case {\bf (II)}
separately.

\subsection*{5.14.~Cones, filled cones, subcones and 
local description of half-wedges in Case (I)} As in some of the
assumptions of Proposition~5.12, let $M$ be a $\mathcal{ C}^{2,
\alpha}$-smooth generic submanifold of $\C^n$ of CR dimension $m=1$
and of codimension $d=n-1 \geq 1$, let $p_1 \in M$, let $M^1$ be a
$\mathcal{ C}^{2, \alpha}$-smooth one-codimensional submanifold of $M$
passing through $p_1$. For the sake of concreteness, it will be
convenient to work in a holomorphic coordinate system $z=(z_1, \dots,
z_n)= (x_1+i y_1, \dots, x_n+ iy_n)$ centered at $p_1$ in which $T_{
p_1}M=\{y_2= \cdots =y_n =0\}$ and $T_{ p_1}M^1=\{y_1= y_2= \cdots=
y_n=0\}$ (the existence of such a coordinate system which straightens
both $T_{p_1}M$ and $T_{p_1}M^1$ is a direct consequence of the
considerations of \S4.11). Let $\pi': T_{p_1}M^1\to T_{p_1}M^1/
(T_{p_1}M^1\cap T_{p_1}^cM)$ denote the canonical projection, namey
$\pi'(x_1,x_2,\dots,x_n)=(x_2,\dots,x_n)$. Sometimes, we shall denote
the coordinates by $(z_1,z')= (x_1+iy_1,x'+iy')\in \C\times
\C^{n-1}$. In these coordinates, the characteristic direction is given
by the $x_1$-axis and we may assume that the tangent plane at $p_1$ of
the one-sided neighborhood $(M^1)^+$ is given by
$T_{p_1}(M^1)^+=\{y'=0, \ y_1>0\}$.

Let ${\sf C}_{ p_1} \subset T_{p_1} M^1$ be the infinite open cone
created by $\mathcal{ HW}_{ p_1}^+$ at $p_1$ and let ${\sf FC}_{ p_1}
\subset T_{ p_1} M^1$ be its filling. Let ${\sf C}_{ p_1}':= \pi'
\left( {\sf C}_{ p_1} \right)$ be its projection onto the $x'$-space,
which yields an $(n-1)$-dimensional infinite cone in the $x'$-space,
open with respect to its topology. Notice that, by the
definition~\thetag{ 5.6} of the filling (along the characteristic
direction) of a cone, the two projections $\pi'({\sf C}_{p_1})$ and
$\pi'({\sf FC}_{p_1})$ are identical. We now need to explain how
these three cones ${\sf C}_{p_1}$, ${\sf FC}_{p_1}$, ${\sf C}_{p_1}'$
and the nonzero vector $v_1\in {\sf FC}_{p_1}$ are disposed,
geometrically, {\it see} {\sc Figure~12} just below.

\bigskip
\begin{center}
\input cones-half-wedge.pstex_t
\end{center}

Because the disc $\mathcal{ Z}_{ t, \chi, \nu}$ of Proposition~5.12
(which is a localization in a neighborhood of the special point of the
discs constructed in Section~4) is small, the tangent vector
$\frac{\partial \mathcal{ Z}_{0, 0,0 }}{ \partial \theta} (1)$ is
necessarily close to the complex tangent plane $T_{p_1 }^cM$: this may
be checked directly by differentiating Bishop's equation~\thetag{4.40}
with respect to $\theta$, using the fact that the $\mathcal{
C}^1$-norm of $\Phi'$ is small. Moreover, since this vector $\frac{
\partial \mathcal{ Z}_{ 0,0,0 }}{ \partial \theta }(1)$ also belongs
to $T_{p_1} M^1$, it is in fact close to the positive
$x_1$-axis. Furthermore, since the vector $v_1$ belongs to the filling
of the open cone ${\sf C}_{ p_1}$ which contains the vector $\frac{
\partial \mathcal{ Z}_{ 0,0,0 }}{ \partial \theta }(1)$, and since in
the proof of Lemma~5.9 above we have chosen the special point, the
supporting hypersurface and the vector $v_1$ with a parameter
$\lambda$ very close to $1$, so that the vector field $p \mapsto
\widetilde{ v}_p^\lambda$ was very close to the characteristic vector
field $p \mapsto X_p$, it follows that the vector $v_1\equiv
\widetilde{ v}_{ p_{ \rm sp}}^\lambda$ is even closer to the positive
$x_1$-axis. However, we suppose in Case {\bf (I)} that $v_1$ is {\it
not}\, directed along the $x_1$-axis, so $v_1$ has coordinates
$(v_{1;1},v_{ 2;1}, \dots,v_{ n;1}) \in \R^n$ with $v_{1;1}>0$, with
$\vert v_{j;1}\vert < < v_{ 1;1}$ for $j=2, \dots,n$ and with at least
one $v_{j;1}$ being nonzero.

We need some general terminology. Let ${\sf C}$ be an open infinite
cone in a real linear subspace $E$ of dimension $q\geq 1$. We say
that ${\sf C}'$ is a {\sl proper subcone} and we write ${\sf
C}'\subset \subset {\sf C}$ ({\it see} the left hand side of {\sc
Figure~10} above for an illustration) if the intersection of ${\sf
C}'$ with the unit sphere of $E$ is a relatively compact subset of the
intersection of ${\sf C}$ with the unit sphere of $E$, this property
being independent of the choice of a norm on $E$. We say that ${\sf
C}$ is a {\sl linear cone} if it may be defined by ${\sf C}=\{x\in E:
\ \ell_1(x)>0,\dots,\ell_q(x)>0\}$ for some $q$ linearly independent
real linear forms $\ell_1,\dots,\ell_q$ on $E$.

In the $(x_2,\dots,x_n)$-space, we now choose an open infinite
strictly convex linear proper subcone ${\sf C}_1'\subset \subset {\sf
C }_{ p_1}'$ with the property that $v_1$ belongs to its filling ${\sf
FC}_1'$, {\it cf.} {\sc Figure~12} above. Here, we may assume that
${\sf C}_1'$ is described by $(n-1)$ strict inequalities $\ell_1'(x')
>0, \dots, \ell_{n-1}'(x')>0$, where the $\ell_j'( x')$ are linearly
independent linear forms. It then follows that there exists a linear
form $\sigma(x_1,x')$ of the form $\sigma(x_1, x')= x_1+a_2x_2+ \cdots
+a_nx_n$ such that the original filled cone ${\sf FC}_{p_1}$ is
contained in the linear cone
\def\theequation{5.15}\begin{equation}
{\sf C}_1:= \left\{
(x_1,x')\in \R^n : \
\ell_1'(x')>0,\dots,
\ell_{n-1}'(x')>0, \ 
\sigma(x_1,x')>0
\right\},
\end{equation}
which contains the vector $v_1$. This cone is automatically filled,
namely ${\sf C}_1\equiv {\sf FC}_1$.

We remind that by genericity of $M$, the complex structure $J$ of
$T\C^n$ induces an isomorphim $T_{p_1}M/T_{p_1}^cM \to
T_{p_1}\C^n/T_{p_1}M$. Hence $J{\sf C}_{p_1}'$ and $J{\sf C}_1'$ are
open infinite strictly convex linear proper cones in
$T_{p_1}\C^n/T_{p_1}M\cong \{(0,y')\in \C^n\}$. Since $J{\sf C}_1'$
is a proper subcone of $J{\sf C}_{p_1}'$ and since in the classical
definition of a wedge, only the projection of the cone on the quotient
space $T_{p_1}M/T_{p_1}^cM$ has a contribution to the wedge, it then
follows that the complete wedge $\mathcal{ W}_{p_1}$ associated to the
family $\mathcal{ Z}_{t,\chi,\nu}(\zeta)$ ({\it cf.} the paragraph
after~\thetag{4.5}) contains a wedge of the form
\def\theequation{5.16}\begin{equation}
\mathcal{ W}_1:=\left\{
p+{\sf c}_1': \ p\in M, \ 
{\sf c}_1' \in J{\sf C}_1'
\right\}\cap \Delta_n(p_1,\delta_1),
\end{equation}
for some $\delta_1$ with $0 < \delta_1 < \varepsilon$, where
$\varepsilon$ is as in \S4.2. Furthermore, as observed in \S4.2,
there exists a $\mathcal{ C}^{2, \alpha}$-smooth hypersurface $N^1$ of
$\C^n$ passing through $p_1$ with the property that $N^1 \cap M \equiv
M^1$ locally in a neighborhood of $p_1$ such that, shrinking $\delta_1
>0$ if necessary, the local half-wedge $\mathcal{ HW}_{ p_1}^+$
contains a local half-wedge $\mathcal{ HW}_1^+$ of edge $(M^1)^+$ at
$p_1$ which is described as the geometric intersection of the complete
wedge $\mathcal{ W}_{p_1}$ with a one-sided neighborhood $(N^1)^+$,
namely
\def\theequation{5.17}\begin{equation}
\mathcal{ HW}_1^+:=
\mathcal{ W}_1\cap (N^1)^+.
\end{equation}
An illustration for the case $n =2$ where $M \subset \C^2$ is a
hypesurface is provided in the left hand side of {\sc Figure~12}. In
addition, it follows from the definition of $\mathcal{ HW}_{p_1}^+$ by
means of the segments $\mathcal{ Z}_{t, \chi, \nu} \left ((1-
\varepsilon, 1) \right)$ that we can assume that
\def\theequation{5.18}\begin{equation}
T_{p_1}(N^1)^+=T_{p_1}M \oplus J(\Sigma^1)^+,
\end{equation}
where $(\Sigma_1 )^+$ is the hyperplane one-sided neighborhood $\{(
x_1, x'): \ \sigma (x_1, x')>0\} \subset T_{p_1} M^1$. Equivalently,
$T_{p_1}(N^1)^+$ is represented by the inequality
$y_1+a_2y_2+\cdots+a_ny_n>0$. Consequently, there exists a $\mathcal{
C}^{2,\alpha}$-smooth function $\psi(x,y')$ with
$\psi(0)=\partial_{x_k}\psi(0)= \partial_{y_j}\psi(0)=0$ for
$k=1,\dots,n$ and $j=2,\dots,n$ such that $N^1$ is represented by the
equation $y_1+a_2y_2+\cdots+a_ny_n=\psi(x,y')$ and $(N^1)^+$ by the
inequation $y_1+a_2y_2+\cdots+a_ny_n>\psi(x,y')$.

\subsection*{5.19.~Cones, filled cones, subcones and 
local description of half-wedges in Case (II)} Secondly, we assume
that the nonzero vector $v_1$ of Proposition~5.12 belongs to the
characteristic direction $T_{p_1 }M^1 \cap T_{p_1 }^cM$. In this
case, as observed in \S4.2, the half-wedge $\mathcal{ HW}_{p_1}^+$
coincides with a local wedge of edge $M^1$ at $(p_1,Jv_1)$. After a
real dilatation of the $z_1$-axis, we can assume that
$v_1=(1,0,\dots,0)$. Choosing an open infinite strictly convex linear
proper subcone ${\sf C}_2\subset \subset {\sf C}_{p_1}\subset
T_{p_1}M^1=\R_x^n$ defined by $n$ strict inequalities
$\ell_1(x)>0,\dots, \ell_n(x)>0$, where the $\ell_j(x)$ are linearly
independent real linear forms\,--\,of course with ${\sf C}_2$
containing the vector $v_1$\,--\,it follows that there exists
$\delta_1>0$ such that the half-wedge $\mathcal{ HW}_{p_1}^+$ contains
the following local wedge of edge $M^1$ at $p_1$:
\def\theequation{5.20}\begin{equation}
\mathcal{ W}_2:=\left\{
p+{\sf c}_2: \ 
p\in M^1, \
{\sf c}_2 \in J{\sf C}_2
\right\}\cap \Delta_n
(p_1,\delta_1).
\end{equation}
We remind that it was observed in \S4.2 ({\it cf.} especially the
right hand side of {\sc Figure}~5) that $\mathcal{ W}_2$ contains
$(M^1)^+$ locally in a neighborhood of $p_1$. In \S5.22 below, we
shall provide a more concrete representation of $\mathcal{ W}_2$
in an appropriate system of coordinates.

\subsection*{5.21.~A trichotomy}
Let us pursue this discussion more concretely by introducing further
normalizations. Our goal will now be to construct appropriate
normalized coordinate systems. Analyzing further the dichotomy 
introduced in \S5.13 by taking account of the 
presence of the one-codimensional 
submanifold $H^1\subset M^1$,
we shall distinguish three cases by dividing Case {\bf (I)} in two
subcases {\bf ($\mathbf{ I_1}$)} and {\bf ($\mathbf{ I_2}$)} as
follows:
\begin{itemize}
\item[{\bf ($\mathbf{I_1}$)}]
The nonzero vector $v_1$ does not belong to the characteristic
direction $T_{p_1}M^1\cap T_{p_1}^cM$ and
$\dim_\R \left(T_{p_1}H^1 \cap T_{p_1}^cM \right)=0$.
\item[{\bf ($\mathbf{I_2}$)}]
The nonzero vector 
$v_1$ does not belong to the characteristic direction 
$T_{p_1}M^1\cap T_{p_1}^cM$ and
$\dim_\R \left(T_{p_1}H^1 \cap T_{p_1}^cM \right)=1$ 
(this possibility can only occur when $n\geq 3$).
\item[{\bf (II)}]
The nonzero vector $v_1$ belongs to the
characteristic direction $T_{p_1}M^1\cap 
T_{p_1}^cM$.
\end{itemize}

In case {\bf ($\mathbf{I_1}$)}, we notice that the
assumption $T_{p_1}H^1\cap T_{p_1}^cM=\{0\}$ implies that
$v_1$ does not belong to the characteristic direction, 
because $v_1\in T_{p_1}H^1$. Also, 
in case {\bf (II)}, we notice that $\dim_\R \left(T_{ p_1} H^1 \cap
T_{ p_1}^c M \right)=1$ because $v_1 \in T_{ p_1} H^1$, because
$T_{p_1}H^1 \subset T_{ p_1}M^1$ and because the characteristic
direction $T_{p_1}M^1\cap T_{p_1}^cM$ is one-dimensional.

In each of the above three cases, it will be convenient in
Section~8 below to work with simultaneously normalized defining
(in)equations for $M$, for $M^1$, for $(M^1)^+$, for $H^1$, for $(H^1
)^+$, for ${\sf C }_1'$, for $v_1$, for ${\sf C }_1$, for $(N^1)^+$
and for $\mathcal{ HW}_{ p_1}^+$, in a single coordinate system
centered at $p_1$. In the next paragraphs, we shall set up further
elementary normalization lemmas in a {\it common system of
coordinates}, firstly for Case {\bf ($\mathbf{I_1}$)}, secondly for
Case {\bf ($\mathbf{I_2}$)} and thirdly for Case {\bf (II)}.

First of all, in the above coordinate system $(z_1,z')$ with $T_{p_1
}M=\{ y_2= \cdots= y_n =0\}$ and with $T_{p_1 }M^1=\{y_1 = y_2=
\cdots= y_n=0\}$, by means of the implicit function theorem, we can
represent locally $M$ by $(n-1)$ grahed equations of the form $y_2=
\varphi_2 (x,y_1),\dots, y_n= \varphi_n (x,y_1)$, where
the $\varphi_j$ are $\mathcal{ C}^{2, \alpha}$-smooth functions
satisfying $\varphi_j(0)= \partial_{ x_k} \varphi_j(0)= \partial_{y_1}
\varphi_j (0)= 0 $ for $j=2,\dots,n$, $k=1,\dots,n$ and
we can represent $M^1$ by $n$ graphed equations $y_1= h_1(x), y_2=
h_2(x), \dots,y_n= h_n(x)$, where the $h_j$ are $\mathcal{
C}^{2,\alpha}$-smooth functions satisfying
$h_j(0)=\partial_{x_k}h_j(0)=0$ for $j,k=1, \dots,n$.

\subsection*{5.22.~First order normalizations in Case 
{\bf ($\mathbf{I_1}$)}} Thus, let us deal first with Case {\bf
($\mathbf{I_1}$)}. After a possible permutation of coordinates, we can
assume that $T_{p_1}H^1$, which is a one-codimensional subspace of
$T_{p_1}M^1$, is given by the equations
\def\theequation{5.23}\begin{equation}
x_1=b_2x_2+\cdots+b_nx_n, \ 
y_1=0,y'=0,
\end{equation}
for some real numbers $b_2,\dots,b_n$. If we define the linear
invertible transformation $\widehat{z}_1:=z_1-b_2z_2-\cdots-b_nz_n$,
$\widehat{z}':=z'$, then the plane $T_{ p_1} H^1$ written
in~\thetag{5.23} clearly transforms to the plane $\widehat{
x}_1=\widehat{ y}_1= \widehat{ y}'=0$, and (fortunately) $T_{ p_1}M$
and $T_{ p_1} M^1$ are left unchanged, namely $T_{p_1} \widehat{M}=\{
\widehat{y}'=0\}$ and $T_{p_1}\widehat{ M^1}= \{\widehat{ y}_1=
\widehat{ y}'=0\}$.

Dropping the hats on coordinates, we have $T_{p_1}M=\{y'=0\}$,
$T_{p_1}M^1=\{y_1=y'=0\}$, $T_{p_1}H^1=\{x_1=y_1=y'=0\}$. Let ${\sf
C}_1'\subset \subset {\sf C}_{p_1}'$ be the open infinite strictly
convex linear cone introduced in \S5.14, which is contained in the
real $(n-1)$-dimensional space $\{(0,x')\}$ and which is defined by
$(n-1)$ strict inequalities $\ell_1'(x')>0,\dots,
\ell_{n-1}'(x')>0$. By means of a real linear invertible
transformation of the form $\widehat{ z}_1:= z_1$, $\widehat{ z}':=
A'\cdot z'$, where $A'$ is an $(n-1)\times (n-1)$ real matrix, we can
transform ${\sf C}_1'$ to a cone $\widehat{\sf C}_1'$ defined by the
simpler inequalities $\widehat{ x}_2>0,\dots, \widehat{ x}_n>0$.
Fortunately, this transformation stabilizes $T_{p_1}M$, $T_{p_1}M^1$
and $T_{p_1}H^1$.

Dropping the hats on coordinates, we now have $T_{p_1}M=\{y'=0\}$,
$T_{p_1}M^1=\{y_1=y'=0\}$, $T_{p_1}H^1=\{x_1=y_1=y'=0\}$ and ${\sf
C}_1'=\{(0,x'): \ x_2>0,\dots,x_n>0\}$. Then the nonzero vector
$v_1\in T_{p_1}H^1$ which belongs to ${\sf C}_1'$ has coordinates
$v_1=(0,v_{2;1},\dots, v_{n;1})\in\R^n$, where
$v_{2;1}>0,\dots,v_{n;1}>0$. By means of real dilatations or real
contractions of the real axes $\R_{x_2},\dots, \R_{x_n}$ (a
transformation which does not perturb the previously achieved
normalizations), we can assume that $v_1=(0,1,\dots,1)$ and that
$T_{p_1}(M^1)^+=\{y'=0, y_1>0\}$, $T_{p_1}(H^1)^+=\{y=0, \ x_1>0\}$.

Finally, the linear one-codimensional subspace $\sigma_1\subset
T_{p_1}M^1$ introduced in \S5.14 which does not contain the
characteristic direction $T_{p_1}M^1\cap T_{p_1}^cM \equiv \R_{x_1}$
may be represented by an equation of the form
$\sigma(x_1,x'):=x_1+a_2x_2+\cdots+a_nx_n=0$, for some real numbers
$a_2,\dots,a_n$. By~\thetag{5.15}, the vector $v_1$ belongs to the
cone ${\sf C}_1$, hence $a_2+\dots+a_n>0$. After a dilatation of the
$x_1$-axis, we can even assume that $a_2+\dots+a_n=1$. We remind that
by~\thetag{5.18}, the half-space $T_{p_1}(N^1)^+$ is given by
$y_1+a_2y_2+\cdots+a_ny_n>0$, hence there exists a $\mathcal{
C}^{2,\alpha}$-smooth function $\psi(x,y')$ with
$\psi(0)=\partial_{x_k}\psi(0)= \partial_{y_j}\psi(0)=0$ for
$k=1,\dots,n$ and $j=2,\dots,n$ such that $(N^1)^+$ is represented by
the inequation $y_1+a_2y_2+\cdots+a_ny_n>
\psi(x,y')$. Consequently, in this coordinate system, we may
represent concretely the local half-wedge $\mathcal{ HW}_1^+\subset
\mathcal{ HW}_{p_1}^+$ constructed in \S5.14 as
\def\theequation{5.24}\begin{equation}
\left\{
\aligned
{}
&
\mathcal{ HW}_1^+=
\left\{
(z_1,z')\in \C^n:
\vert z_1 \vert < \delta_1, \
\vert z' \vert < \delta_1, \right. \\
& \ \ \ \ \ \ \ \ \ \ \ \ \ \ \ \ \ \ \
y_1+a_2y_2+\dots+a_ny_n-\psi(x,y')>0,\\
& \ \ \ \ \ \ \ \ \ \ \ \ \ \ \ \ \ \
\left.
y_2-\varphi_2(x,y_1) > 0, \dots, 
y_n-\varphi_n(x,y_1) >0 \ 
\right\}.
\endaligned\right.
\end{equation}
For the continuation of the proof of Proposition~5.12, it will also be
convenient to proceed to further second order normalizations of the
totally real submanifolds $M^1$ and $H^1$. These normalizations will
all be tangent to the identity tranformation, hence they will leave the
previously achieved normalizations unchanged.

\subsection*{5.25.~Second order normalizations in Case 
($\mathbf{ I_1}$)}
Let us then perform a second order Taylor development of the defining
equations of $M^1$
\def\theequation{5.26}\begin{equation}
y=h(x)=\sum_{k_1,k_2=1}^n\, 
a_{k_1,k_2}\, x_{k_1} x_{k_2}+o(\vert x \vert^2),
\end{equation}
where the $a_{k_1,k_2}=\frac{1}{2}
\partial_{x_{k_1}}
\partial_{x_{k_2}}h (0)$ are vectors of $\R^n$.
If we define the quadratic invertible transformation
\def\theequation{5.27}\begin{equation}
\widehat{z}:= z-i\sum_{k_1,k_2=1}^n\, 
a_{k_1,k_2}\, z_{k_1} z_{k_2}=\Phi(z), 
\end{equation}
which is tangent to the identity mapping at the origin, 
then for $x+iy=x+ih(x)\in M^1$, 
we have by replacing~\thetag{5.26} in the imaginary 
part of $\widehat{z}$ given by~\thetag{5.27}
\def\theequation{5.28}\begin{equation}
\left\{
\aligned
\widehat{y}= 
& \
y - \sum_{k_1,k_2=1}^n\, 
a_{k_1,k_2}\, x_{k_1}x_{k_2}+
\sum_{k_1,k_2=1}^n\, a_{k_1,k_2}\, y_{k_1}y_{k_2} \\
= 
& \ 
{\rm o}(\vert x \vert^2) \\
= 
& \
{\rm o}\left(\left\vert
{\rm Re}\, \Phi^{-1}(\widehat{z})
\right\vert^2\right)= {\rm o}\left( \left\vert
(\widehat{x}, \widehat{y})
\right\vert^2 \right),
\endaligned\right.
\end{equation}
whence by applying the $\mathcal{ C}^{2,\alpha}$ implicit function
theorem to solve~\thetag{5.28} in terms of $\widehat{ y}$, we find
that $\widehat{ M^1}:=\Phi( M^1)$ may be represented by an equation of
the form $\widehat{ y}= \widehat{ h}(\widehat{ x})$, for some
$\R^n$-valued local $\mathcal{ C}^{2,\alpha}$-smooth mapping
$\widehat{ h}$ which satisfies $\widehat{ h}( \widehat{ x})= {\rm o}(
\vert \widehat{ x} \vert^2)$.

Finally, dropping the hats on coordinates, we can assume that the
functions $h_1, \dots, h_n$ vanish at the origin to second order.
Since $T_{p_1}H^1=\{y=0, \ x_1=0\}$, there exists a $\mathcal{
C}^{2,\alpha}$-smooth function $g(x')$ with $g(0)=\partial_{x_k}
g(0)=0$ for $k=2,\dots,n$ such that $(H^1)^+$ is given by the equation
$x_1 >g(x')$. We want to normalize also the defining equation
$x_1=g(x')$ of $H^1$. Instead of requiring, similarly as for
$h_1,\dots,h_n$, that $g$ vanishes to second order at the origin
(which would be possible), we shall normalize $g$ in order that
$g(x')=-x_1^2-\cdots - x_n^2+ {\rm o}\left(\vert x' \vert^{2+\alpha}
\right)$ (which will also be possible, thanks to the total reality of
$H^1$). The reason why we want $g$ to be strictly concave is a trick
that will be useful in Section~8 below.

Thus, we now perform a second
order Taylor development of the defining equations of $H^1$
\def\theequation{5.29}\begin{equation}
\left\{
\aligned
{}
&
x_1=g\left(x'\right)=\sum_{k_1,k_2=2}^n\, 
b_{k_1,k_2} \, x_{k_1}x_{k_2}+ o (\left\vert
x'
\right\vert^2),\\
&
y=h\left(g(x'),x'\right)=:k\left(x'\right)=
{\rm o}(\left\vert x' \right\vert^2),
\endaligned\right.
\end{equation}
where the $b_{k_1,k_2}=\frac{1}{2}
\partial_{x_{k_1}}
\partial_{x_{k_2}} g(0)$
are real numbers. If we define the quadratic invertible transformation
\def\theequation{5.30}\begin{equation}
\left\{
\aligned
{}
&
\widehat{z}_1:= z_1-\sum_{k_1,k_2=2}^n\, 
b_{k_1,k_2}\, z_{k_1}z_{k_2} -z_2^2-\cdots-z_n^2, \\
&
\widehat{z}':= z', 
\endaligned\right.
\end{equation}
which is tangent to the identity mapping, then 
for $\left(g(x')+ik_1(x'),x'+ik'(x')\right)\in H^1$,
we have by replacing~\thetag{5.29} in the 
real part of $\widehat{z}_1$, given by~\thetag{5.30}:
\def\theequation{5.31}\begin{equation}
\left\{
\aligned
\widehat{x}_1=
& \
x_1-\sum_{k_1,k_2=2}^n\, b_{k_1,
k_2}\, x_{k_1}x_{k_2}+\sum_{k_1,k_2=2}^n\, 
b_{k_1,k_2}\, y_{k_1}y_{k_2}-
\sum_{k=2}^n\, x_k^2+\sum_{k=2}^n\, 
y_k^2, \\
=
& \
-x_2^2-\cdots-x_n^2+{\rm o}\left(
\left\vert x' \right\vert^2\right) \\
=
& \
-\widehat{x}_2^2-\cdots-\widehat{x}_n^2-{\rm o}\left(\left\vert
(\widehat{x},\widehat{y})\right\vert^2\right).
\endaligned\right.
\end{equation}
Similarly (dropping the elementary computations), we may obtain for
the imaginary part of $\widehat{z}_1$ and for the imaginary part of
$\widehat{z}'$
\def\theequation{5.32}\begin{equation}
\widehat{y}_1={\rm o}\left(\left\vert
(\widehat{x},\widehat{y})\right\vert^2\right) \ \ \ \ \ 
{\rm and} \ \ \ \ \ 
\widehat{y}'={\rm o}\left(\left\vert
(\widehat{x},\widehat{y})\right\vert^2\right),
\end{equation}
whence by applying the $\mathcal{ C}^{2,\alpha}$ implicit function
theorem to solve the system~\thetag{5.31}, \thetag{5.32} in terms
of $\widehat{x}_1$, $\widehat{y}_1$ and $\widehat{y}'$, we find
that $\widehat{ H}^1:=\Phi(H^1)$ may be represented 
by equations
of the form
\def\theequation{5.33}\begin{equation}
\left\{
\aligned
\widehat{x}_1=
& \
\widehat{g}\left(\widehat{x}'\right)=
-\widehat{x}_2^2-\cdots-\widehat{x}_n^2+
{\rm o}\left(
\left\vert
\widehat{x}'
\right\vert^2
\right), \\
\widehat{y}=
& \
\widehat{k}\left(\widehat{x}'\right)=o\left(
\left\vert
\widehat{x}'
\right\vert^2
\right).
\endaligned\right.
\end{equation}
It remains to check that the above transformation has not perturbed
the previous second order normalizations of $h_1, \dots, h_n$ (this is
important), which is easy: replacing $y$ by $h(x) = {\rm o} (\vert x
\vert^2)$ in the imaginary parts of $\widehat{z}_1$ and of
$\widehat{z}'$ defined by the transformation~\thetag{5.30}, we get
firstly
\def\theequation{5.34}\begin{equation}
\left\{
\aligned
\widehat{y}_1 =
& \
y_1 -\sum_{k_1,k_2}^n\, 
b_{k_1,k_2}\, \left(
x_{k_1}y_{k_2}+y_{k_1}x_{k_2}
\right)-2\sum_{k=2}^n\, 
x_ky_k \\
=
& \
{\rm o}(\vert x\vert^2) \\
=
& \
{\rm o}\left(\left\vert
{\rm Re}\, \Phi^{-1}(\widehat{z})
\right\vert^2\right)= {\rm o}\left( \left\vert
(\widehat{x}, \widehat{y})
\right\vert^2 \right),
\endaligned\right.
\end{equation}
and similarly 
\def\theequation{5.35}\begin{equation}
\widehat{y}'=o\left( \left\vert
(\widehat{x}, \widehat{y})
\right\vert^2 \right),
\end{equation}
whence by applying the $\mathcal{ C}^{2,\alpha}$ implicit function
theorem to solve the system~\thetag{5.34}, \thetag{5.35} in terms of
$\widehat{y}$, we find that $\widehat{ M}^1 :=\Phi(M^1)$ may be
represented by equations of the form $\widehat{y}=\widehat{ h} \left(
\widehat{x} \right)= {\rm o}\left(\left \vert \widehat{ x}\right \vert^2
\right)$. Thus, after dropping the hats on coordinates, all the
desired normalizations are satisfied. We shall now summarize these
normalizations and we shall formulate just afterwards the analogous
normalizations for Cases {\bf ($\mathbf{I_2}$)} and {\bf (II)}.

\subsection*{5.36.~Simultaneous normalization lemma}
In the following lemma, the final choice of sufficiently small radii
$\rho_1>0$ and $\delta_1>0$ is made after that all the biholomorphic
changes of coordinates and all the applications of the implicit function
theorem are achieved.

\def\thelemma{5.37}\begin{lemma}
Let $M$, $M^1$, $p_1$, $H^1$, $(H^1)^+$, $\mathcal{ HW}_{p_1}^+$,
${\sf C}_{p_1}$ and ${\sf FC}_{p_1}$ be as in Proposition~5.12. Then
there exists a sub-half-wedge $\mathcal{ HW}_1^+$ contained in
$\mathcal{ HW}_{p_1}^+$ such that the following normalizations hold in
each of the three cases {\bf ($\mathbf{I_1}$)}, {\bf
($\mathbf{I_2}$)} and
{\bf (II)}{\rm :}

\smallskip
\begin{itemize}
\item[{\bf ($\mathbf{I_1}$)}]
If $\dim_\R \left( T_{ p_1} H^1 \cap T_{p_1 }^cM \right) =0$, then
there exists a system of holomorphic coordinates $z= (z_1, \dots,
z_n)=(x_1+i y_1,\dots, x_n+iy_n)$ vanishing at $p_1$ with the vector
$v_1$ equal to $(0,1,\dots,1)$, there exists positive numbers $\rho_1$
and $\delta_1$ with $0< \delta_1 < \rho_1$, there exist $\mathcal{
C}^{2,\alpha}$-smooth functions $\varphi_2,\dots, \varphi_n$, $h_1,
\dots, h_n$, $g$, $k_1, \dots,k_n,\psi$, all defined in real cubes of
edge $2 \rho_1$ and of the appropriate dimension, and there exist real
numbers $a_1,\dots,a_n$ with $a_2+\cdots+a_n=1$, such that, if we
denote $z':=(z_2, \dots, z_n)=x' +iy'$, then $M$, $M^1$, $(M^1)^+$,
$H^1$, $(H^1)^+$ and $N^1$ are represented in the polydisc of radius
$\rho_1$ centered at $p_1$ by the following graphed (in)equations and
the sub-half-wedge $\mathcal{ HW}_1^+ \subset \mathcal{ HW}_{p_1}^+$ is
represented in the polydisc of radius $\delta_1$ centered at $p_1$ by
the following inequations
\def\theequation{5.38}\begin{equation}
\left\{
\aligned
M: \ \ \ \ \ 
& \
y_2=\varphi_2(x,y_1), \dots\dots, \
y_n=\varphi_n(x,y_1), \\
M^1: \ \ \ \ \ 
& \ 
y_1=h_1(x),y_2=h_2(x),\dots\dots, \ 
y_n=h_n(x), \\
(M^1)^+: \ \ \ \ \ 
& \ 
y_1>h_1(x),y_2=\varphi_2(x,y_1),\dots\dots, \ 
y_n=\varphi_n(x,y_1), \\
H^1: \ \ \ \ \ 
& \
x_1=g(x'), \
y_1=k_1(x'), \dots\dots, \ 
y_n=k_n(x'), \\
(H^1)^+: \ \ \ \ \ 
& \
x_1>g(x'), \
y_1=h_1(x),y_2=h_2(x), \dots\dots, \ 
y_n=h_n(x), \\
N^1 : \ \ \ \ \ 
& \
y_1+a_2y_2+\cdots+a_ny_n=\psi(x,y'), \\
\mathcal{ HW}_1^+: \ \ \ \ \
& \
y_1+a_2y_2+\cdots+a_ny_n>\psi(x,y'), \\
& \ \ \ \ \ \ \ \ \ \ \ \ \ \ \ \ \ 
y_2>\varphi_2(x,y_1),\dots, y_n>\varphi_n(x,y_1), 
\endaligned\right.
\end{equation}
where we can assume that $M^1$ coincides with the intersection $M\cap
\{y_1=h_1(x)\}$, that $H^1$ coincides with the intersection
$M^1\cap \{x_1=g(x')\}$ and that $N^1$ contains $M^1$, which yields at
the level of defining equations the following three collections of
identities
\def\theequation{5.39}\begin{equation}
\left\{
\aligned
h_2(x)\equiv 
& \
\varphi_2(x,h_1(x)),\dots\dots, h_n(x)\equiv 
\varphi_n(x,h_1(x)), \\
k_1(x')\equiv
& \
h_1(g(x'),x'), \dots\dots,
k_n(x')\equiv h_n(g(x'),x'), \\
\psi(x,h'(x)) \equiv
& \
h_1(x)+a_2h_2(x)+\cdots+
a_nh_n(x),
\endaligned\right.
\end{equation}
and where the following normalizations hold (where $\delta_a^b$, equal
to $0$ if $a\neq b$ and to $1$ if $a=b$, denotes Kronecker's
symbol){\rm :}
\def\theequation{5.40}\begin{equation}
\left\{
\aligned
{}
&
\varphi_{j}(0)=\partial_{x_k}\varphi_j(0)
=\partial_{y_1}\varphi_j(0)=0, \ \ \ \ \ j=2,\dots,n, \
k=1,\dots,n,\\
&
h_j(0)=\partial_{x_k}
h_j(0)=\partial_{x_{k_1}}
\partial_{x_{k_2}}h_j(0)=0, \ \ \ \ \ \ \ \ \ \
j,k,k_1,k_2=1,\dots,n, \\
&
g(0)=\partial_{x_k}
g(0)=k_j(0)=\partial_{x_k}k_j(0)=0, \ \ \ \ \ 
j=1,\dots,n, \ k=2,\dots,n, \\
& \ \ \ \ \ \ \ \ \ \ \ \ \ \ \ \ \ \ \ \ \ \ \ \ \ \ \ \
\ \ \ \ \ \ \ \ \ \ \ \ \ \ \ \
\partial_{x_{k_1}}\partial_{x_{k_2}}
g(0)=-\delta_{k_1}^{k_2}, \ \ \ \ \ \ \ 
k_1,k_2=2,\dots,n, \\
&
\psi(0)=\partial_{x_k}\psi(0)=\partial_{y_j}\psi(0)=0, \ \ \ \ \
k=1,\dots,n, \ j=2,\dots,n.
\endaligned\right.
\end{equation}
In other words, $T_0M=\{y'=0\}$ (hence $T_0^cM$ coincides with the
complex $z_1$-axis), $T_0N^1=\{ y_1+a_2y_2+\cdots+a_ny_n=0\}$ and the
second order Taylor approximations of the defining equations of $M^1$,
of $H^1$ and of $(H^1)^+$ are the quadrics
\def\theequation{5.41}\begin{equation}
\left\{
\aligned
T_{p_1}^{(2)}M^1: \ \ \ \ \ 
& \
y_1=0, \dots\dots, \ y_n=0,
\\
T_{p_1}^{(2)}H^1: \ \ \ \ \ 
& \
x_1=-x_2^2-\cdots-x_n^2, \ y_1=0,\dots\dots, \ y_n=0, \\
T_{p_1}^{(2)}(H^1)^+: \ \ \ \ \ 
& \
x_1>-x_2^2-\cdots-x_n^2, \ y_1=0,\dots\dots, \ y_n=0.
\endaligned\right.
\end{equation}
\item[{\bf ($\mathbf{I_2}$)}]
Similarly, if $\dim_\R \left( T_{p_1}H^1\cap T_{p_1}^c M\right) = 1$
and if $v_1$ is not complex tangent to $M$ {\rm (}this possibility can
only occur in the case $n \geq 3${\rm )}, then there exists a system
of holomorphic coordinates $z= (z_1,\dots,z_n)=
(x_1+iy_1,\dots,x_n+iy_n)$ vanishing at $p_1$ with $v_1$ equal to
$(1,\dots,1,0)$, there exists positive numbers $\rho_1$ and $\delta_1$
with $0< \delta_1 < \rho_1$, there exist $\mathcal{ C}^{
2,\alpha}$-smoooth functions $\varphi_2,\dots,\varphi_n, h_1,\dots,
h_n, g, k_1,\dots,k_n,\psi$ all defined in real cubes of edge
$2\rho_1$ and of the appropriate dimension, such that if we denote
$z'':=(z_1,\dots,z_{n-1})=x''+iy''$ and $z'=(z_2,\dots,z_n)=
x'+iy'$, then $M$, $M^1$, $(M^1)^+$,
$H^1$, $(H^1)^+$ and $N^1$ are represented in the polydisc of radius
$\rho_1$ centered at $p_1$ by the following graphed (in)equations and
the sub-half-wedge $\mathcal{ HW}_1^+\subset \mathcal{ HW}_{p_1}^+$ is
represented in the polydisc of radius $\delta_1$ centered at $p_1$ by
the following inequations
\def\theequation{5.42}\begin{equation}
\left\{
\aligned
M: \ \ \ \ \ 
& \
y_2=\varphi_2(x,y_1), 
\dots\dots, \
y_n=\varphi_n(x,y_1), \\
M^1: \ \ \ \ \ 
& \ 
y_1=h_1(x),y_2=h_2(x),\dots\dots, \ 
y_n=h_n(x), \\
(M^1)^+: \ \ \ \ \ 
& \ 
y_1>h_1(x),y_2= \varphi_2
(x,y_1), \dots\dots, \ 
y_n=\varphi_n(x,y_1), \\
H^1: \ \ \ \ \ 
& \
x_n=g(x''), \
y_1=k_1(x''), \dots\dots, \ 
y_n=k_n(x''), \\
(H^1)^+: \ \ \ \ \ 
& \
x_n>g(x''), \
y_1=h_1(x),y_2=h_2(x), \dots\dots, \ 
y_n=h_n(x), \\
N^1 : \ \ \ \ \ 
& \
y_2+\cdots+y_{n-1}-y_n=
\psi(x,y'), \\
\mathcal{ HW}_1^+: \ \ \ \ \
& \
y_2+\cdots+y_{n-1}-y_n>
\psi(x,y'), \\
& \ \ \ \ \ \ \ \ \ \ \ \ \ \ \ \ \ 
y_1>\varphi_1(x,y_1),\dots, 
y_{n-1}>\varphi_{n-1}(x,y_1), 
\endaligned\right.
\end{equation}
where we can assume that $M^1$ coincides with the intersection $M\cap
\{y_1=h_1(x)\}$, that $H^1$ coincides with the intersection
$M^1\cap \{x_1=g(x')\}$ and that $N^1$ contains $M^1$, which yields at
the level of defining equations the following three collections of
identities
\def\theequation{5.43}\begin{equation}
\left\{
\aligned
h_2(x)\equiv 
& \
\varphi_2(x,h_1(x)),\dots\dots, h_n(x)\equiv 
\varphi_n(x,h_1(x)), \\
k_1(x'')\equiv
& \
h_1(x'',g(x'')), \dots\dots,
k_n(x'')\equiv h_n(x'',g(x'')), \\
\psi(x,h'(x)) \equiv
& \
h_1(x)+h_2(x)+\cdots+h_{n-1}(x)-
h_n(x),
\endaligned\right.
\end{equation}
and where the following normalizations hold{\rm :}
\def\theequation{5.44}\begin{equation}
\left\{
\aligned
{}
&
\varphi_{j}(0)=\partial_{x_k}\varphi_j(0)
=\partial_{y_1}\varphi_j(0)=0, \ \ \ \ \ j=2,\dots,n, \
k=2,\dots,n, \\
&
h_j(0)=\partial_{x_k}
h_j(0)=\partial_{x_{k_1}}
\partial_{x_{k_2}}h_j(0)=0, \ \ \ \ \ \ \ \ \ \
j,k,k_1,k_2=1,\dots,n, \\
&
g(0)=\partial_{x_k}
g(0)=k_j(0)=\partial_{x_k}k_j(0)=0, \ \ \ \ \ 
j=1,\dots,n, \ k=1,\dots,n-1, \\
& \ \ \ \ \ \ \ \ \ \ \ \ \ \ \ \ \ \ \ \ \ \ \ \ \ \ \ \
\ \ \ \ \ \ \ \ \ \ \ \ \ \ \ \
\partial_{x_{k_1}}\partial_{x_{k_2}}
g(0)=-\delta_{k_1}^{k_2}, \ \ \ \ \ \ \ 
k_1,k_2=1,\dots,n-1, \\
&
\psi(0)=\partial_{x_k}\psi(0)=\partial_{y_j}\psi(0)=0, \ \ \ \ \
k=1,\dots,n, \ j=2,\dots,n.
\endaligned\right.
\end{equation}
In other words, $T_0M=\{y'=0\}$ {\rm (}hence $T_0^cM$ coincides with
the complex $z_1$-axis{\rm )}, $T_0N^1=\{
y_1+y_2+\cdots+y_{n-1}-y_n=0\}$ and the second order Taylor
approximations of the defining equations of $M^1$, of $H^1$ and of
$(H^1)^+$ are the quadrics
\def\theequation{5.45}\begin{equation}
\left\{
\aligned
T_{p_1}^{(2)}M^1: \ \ \ \ \ 
& \
y_1=0, \dots\dots, \ y_n=0,
\\
T_{p_1}^{(2)}H^1: \ \ \ \ \ 
& \
x_n=-x_1^2-\cdots-x_{n-1}^2, \ y_1=0,\dots\dots, \ y_n=0, \\
T_{p_1}^{(2)}(H^1)^+: \ \ \ \ \ 
& \
x_n>-x_1^2-\cdots-x_{n-1}^2, \ y_1=0,\dots\dots, \ y_n=0.
\endaligned\right.
\end{equation}
\item[{\bf (II)}]
Finally, if $\dim_\R\left( T_{p_1}H^1\cap T_{p_1}^c M\right) = 1$ and
if $v_1$ is complex tangent to $M$ {\rm (}this possibility can occur
in all cases $n\geq 2${\rm )}, then there exists a system of holomorphic
coordinates $z=(z_1,\dots,z_n)= (x_1+iy_1,\dots,x_n+iy_n)$ vanishing
at $p_1$ with $v_1$ equal to $(1,0,\dots,0)$, there exist positive
numbers $\rho_1$ and $\delta_1$ with $0< \delta_1 < \rho_1$, there
exist $\mathcal{ C}^{ 2,\alpha}$-smoooth functions
$\varphi_2,\dots,\varphi_n, h_1,\dots, h_n, g, k_1,\dots,k_n$ all
defined in real cubes of edge $2\rho_1$ and of the appropriate
dimension, such that if we denote $z'':=(z_1,\dots,z_{n-1})=x''+iy''$
and $z'=(z_2,\dots,z_n)= x'+iy'$, then $M$, $M^1$, $(M^1)^+$, $H^1$
and $(H^1)^+$ are represented in the polydisc of radius $\rho_1$
centered at $p_1$ by the first five (in)equations of~\thetag{5.42}
together with the normalizations~\thetag{5.45} and such that the local
wedge $\mathcal{ W}_2 \subset \mathcal{ HW}_{p_1}^+$ of edge $M^1$
at $p_1$ is represented in the polydisc of radius $\delta_1$ centered
at $p_1$ by the following inequations
\def\theequation{5.46}\begin{equation}
\left\{
\aligned
\mathcal{ W}_2: \ \ \ \ \
& \
y_1-h_1(x)>-\left[y_2-h_2(x)\right], \dots\dots, \
y_1-h_1(x)>-\left[y_n-h_n(x)\right], \\
& \ \ \ 
y_1-h_1(x)>y_2-h_2(x)+\cdots+
y_n-h_n(x).
\endaligned\right.
\end{equation}
\end{itemize}
\end{lemma}

\subsection*{5.47.~Summarizing figure and proof of Lemma~5.37}
As an illustration for this technical lemma, by specifying the value
$n=3$, we have drawn in the following figure the cones ${\sf C}_1$ and
${\sf C}_2$ together with the vector $v_1$, the tangent plane
$T_{p_1}H^1$ and the hyperplane $\Sigma^1$ in the three cases {\bf
($\mathbf{I_1}$)}, {\bf ($\mathbf{I_2}$)} and {\bf (II)}. In the left
part of this figure, the cone ${\sf C}_1$ is given by $x_2>0$,
$x_3>0$, $x_1>-\frac{1}{2}x_2-\frac{1}{2}x_3$, namely we have chosen
the values $a_2=a_3=\frac{1}{2}$ for the drawing; in the central part,
the cone ${\sf C}_1$ is given by $x_1>0$, $x_2>0$, $x_2>x_3$; in the
right part, the cone ${\sf C}_2$ is given by $x_1>-x_2$, $x_1>-x_3$,
$x_1>x_2+x_3$.

\bigskip
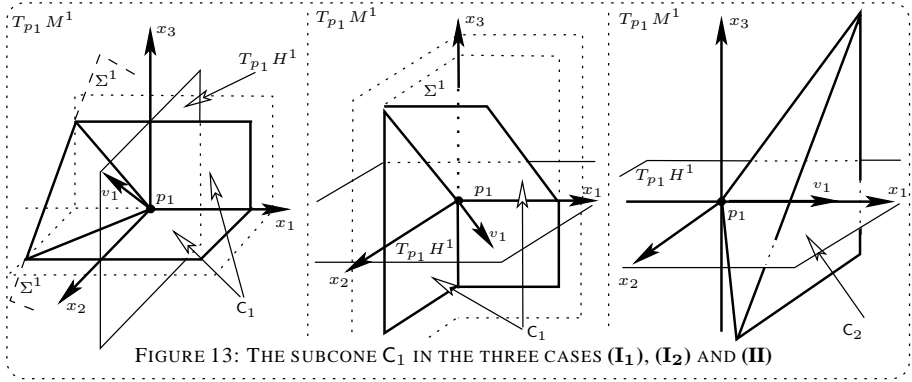
\begin{center}
\input three-cones.pstex_t
\end{center}

\proof
Case {\bf ($\mathbf{I_1}$)} has been completed before
the statement of Lemma~5.37. 

For Case {\bf ($\mathbf{I_2}$)}, we reason similarly, as follows. We
start with the normalizations $T_{p_1} M=\{y'=0\}$ and $T_{p_1} M^1=
\{y=0\}$ as in the end of \S5.21. By assumption, $T_{p_1}H^1$ contains
the characteristic direction, which coincides with the $x_1$-axis. By
means of an elementary real linear transformation of the form
$\widehat{ z}_1:= z_1$, $\widehat{ z}'=A'\cdot z'$, we may first
normalize $T_{p_1} H^1$ to be the hyperplane (after dropping the hats
on coordinates) $\{x_n=0, \ y=0\}$. 
Similarly, we may normalize $v_1$ to be
the vector $(1,1,\dots,1,0)$. Let again $\pi': (x_1,x') \mapsto x'$
denote the canonical projection on the $x'$-space. Then
$\pi'(v_1)=(1,\dots,1,0)$. Using again a real linear transformation of
the form $\widehat{ z}_1:= z_1$, $\widehat{ z}'=A'\cdot z'$, we can
assume that the proper subcone ${\sf C}_1'\subset \subset {\sf
C}_{p_1}'\equiv \pi'({\sf C}_{p_1})$ which contains the vector $v_1$ is
given (after dropping the hats on coordinates) by
\def\theequation{5.48}\begin{equation}
{\sf C}_1': \ \ \ \ \ 
x_2>0,\dots,x_{n-1}>0, \ 
x_2+\dots+x_{n-1}>x_n.
\end{equation}
Following \S5.14 ({\it cf.} {\sc Figure}~12), we choose a linear cone
${\sf C}_1\subset\subset {\sf FC}_{p_1}$ defined by the $(n-1)$
inequations of ${\sf C}_1'$ plus one inequation of the form
$x_1>a_2x_2+\cdots+a_nx_n$ with $1>a_2+\cdots+a_{n-1}$, since $v_1$
belongs to ${\sf C}_1$. Then by means of a real linear transformation
of the form $\widehat{ z}_1:=z_1 + a_2z_2+\cdots +a_nz_n$, $\widehat{
z}':= z'$, which stabilizes $\pi'(v_1)$ and the
inequations~\thetag{5.48} of ${\sf C}_1'$, we can assume that the
supplementary inequation for ${\sf C}_1$, namely the inequation for
$(\Sigma^1)^+$, is simply (after dropping the hats on coordinates)
$x_1>0$. Then the vector $v_1$ is mapped to the vector of coordinates
$(1-a_2-\cdots-a_n,1,\dots,1,0)$, which we map to the vector of
coordinates $(1,1,\dots,1,0)$ by an obvious positive scaling of the
$x_1$-axis. In conclusion, in the final system of coordinates,
the cone ${\sf C}_1$ is given by
\def\theequation{5.49}\begin{equation}
{\sf C}_1: \ \ \ \ \ 
x_1>0,x_2>0,\dots,x_{n-1}>0, \ 
x_2+\cdots+x_{n-1}-x_n>0.
\end{equation}
This implies that the half-wedge $\mathcal{ HW}_1^+\subset \mathcal{
HW}_{p_1}^+$ may be represented by the inequations of the last two line
of~\thetag{5.42}. To conclude the proof of Case {\bf
($\mathbf{I_2}$)} of Lemma~5.37, it suffices to observe that, as in
Case {\bf ($\mathbf{I_1}$)}, the further second order normalizations
do not perturb the previously achieved first order normalizations,
because the transformations are tangent to the identity mapping at the
origin.

Finally, we treat Case {\bf (II)} of Lemma~5.37, starting with the
system of coordinates $(z_1,\dots,z_n)$ of the end of \S5.21. After
an elementary real linear transformation stabilizing the characteristic
$x_1$-axis, we can assume that $v_1=(1,0,\dots,0)$ and that the
convex infinite linear cone ${\sf C}_2$ introduced in \S5.19 which
contains $v_1$ is given by the inequations
\def\theequation{5.50}\begin{equation}
x_1>-x_2,\dots\dots,\ x_1>-x_n,\ 
x_1>x_2+\cdots+x_n.
\end{equation}
This implies that the local wedge $\mathcal{ W}_2\subset \mathcal{
HW}_{p_1}^+$ of edge $M^1$ at $p_1$ introduced in \S5.19 may be
represented by the inequations~\thetag{5.46}. Finally, the second
order normalizations, which are tangent to the identity mapping, are
achieved as in the two previous cases {\bf ($\mathbf{I_1}$)} and {\bf
($\mathbf{I_2}$)}.

The proof of Lemma~5.37 is complete.
\endproof

\section*{\S6.~Three preparatory lemmas on H\"older spaces}

In this section, we first collect a few very elementary lemmas that
will be useful in our geometric construction of half-attached analytic
discs which will be achieved in Section~7 below. From now on, we
shall admit the convenient index notation $g_{x_k}$ for the partial
derivative which was denoted up to now by $\partial_{x_k}g$.

\subsection*{6.1.~Local growth of $\mathcal{ C}^{2,\alpha}$-smooth 
mappings} Let $n\in \N$, $n \geq 1$ and let $x= (x_1, \dots, x_n) \in
\R^n$. We shall use the norm $\vert x \vert := \max_{1 \leq k\leq n}\,
\vert x_k \vert$. If $g= g(x)$ is an $\R^n$-valued $\mathcal{
C}^1$-smooth mapping on the real cube $\{x\in \R^n: \ \vert x \vert <
2\rho_1\}$, for some $\rho_1>0$, and if $\vert x'\vert, \vert x''
\vert \leq \rho$, for some $\rho<2\rho_1$, we have the trivial
estimate
\def\theequation{6.2}\begin{equation}
\left\vert g(x')- g(x'') \right\vert \leq \left\vert 
x' -x'' \right\vert \cdot \left(
\sum_{k=1}^n \, \sup_{\vert x 
\vert \leq \rho} \, \vert g_{x_k}(x)
\vert \right),
\end{equation} 
where we denote by $g_{j,x_k}$ the partial derivative $\partial
g_j/\partial x_k$. Notice that by the definition of the norm $\vert
\cdot \vert$, we have in~\thetag{6.2} that $\vert g(x) \vert \equiv
\max_{1\leq k\leq n} \, \vert g_j (x) \vert$ and that $\vert g_{x_k}
(x)\vert \equiv \max_{1 \leq k \leq n}\, \vert g_{j,x_k}(x) \vert$.

Let $\alpha$ with $0 < \alpha < 1$ and let $h= h(x) = (h_1(x), \dots,
h_n(x ))$ be an $\R^n$-valued mapping which is of class $\mathcal{
C}^{2, \alpha}$ on the real cube $\{x\in \R^n: \ \vert x \vert <
2\rho_1\}$, for some $\rho_1 >0$. For every $\rho< 2 \rho_1$, we
consider the $\mathcal{ C}^{2, \alpha}$ norm of $h$ over $\{ \vert x
\vert \leq \rho\}$ which is defined precisely as:
\def\theequation{6.3}\begin{equation}
\left\{
\aligned
\vert \vert h \vert \vert_{\mathcal{ C}^{2,\alpha}(
\{\vert x \vert \leq \rho\})}:= 
& \
\sup_{\vert x \vert \leq \rho}\, 
\vert h(x)\vert+
\sum_{k=1}^n\, 
\sup_{\vert x \vert \leq \rho}\, 
\left\vert
h_{x_k}(x)\right\vert+
\sum_{k_1,k_2=1}^n\, 
\left\vert h_{x_{k_1}x_{k_2}}(x)\right\vert+\\
& \
+\sum_{k_1,k_2=1}^n\, 
\sup_{\left\vert x'\right\vert, \,
\left\vert x''\right\vert \leq \rho, \ 
x'\neq x''}\ 
\frac{
\left\vert h_{x_{k_1}x_{k_2}}(x')-
h_{x_{k_1}x_{k_2}}(x'')
\right\vert
}
{\left\vert x' - x'' \right\vert^\alpha}<\infty,
\endaligned\right.
\end{equation}
and which is finite. With these definitions at hand, the
following lemma can easily be established by means of~\thetag{6.2}.

\def\thelemma{{\bf 6.4}}\begin{lemma}
Under the above assumptions, let 
\def\theequation{6.5}\begin{equation}
K_1:= \vert \vert h \vert \vert_{
\mathcal{ C}^{2, \alpha}(\{ \vert x \vert \leq \rho_1\})}<\infty
\end{equation}
be the $\mathcal{ C}^{2, \alpha}$ norm of $h$ over the cube $\{\vert x
\vert \leq \rho_1\}$ and assume that $h_j (0)=0$, $h_{j, x_k} (0)=0$
and $h_{j, x_{k_1}x_{k_2}}(0)=0$, for all $j,k, k_1,
k_2=1,\dots,n$. Then the following three inequalities hold for $\vert
x \vert \leq \rho_1${\rm :}
\def\theequation{6.6}\begin{equation}
\left\{
\aligned
{\bf [1]:}
& \ \ \ \ \ 
\vert h(x)\vert 
\leq 
\vert x\vert^{2+\alpha} \cdot K_1,\\
{\bf [2]:} 
& \ \ \ \ \ 
\sum_{k=1}^n\, \left\vert
h_{x_k}(x) \right\vert 
\leq 
\vert x \vert^{1+\alpha} \cdot K_1, \\
{\bf [3]:} 
& \ \ \ \ \ \
\sum_{k_1,k_2=1}^n \, 
\left\vert
h_{x_{k_1}x_{k_2}}(x)
\right\vert 
\leq 
\vert x \vert^\alpha \cdot K_1.
\endaligned\right.
\end{equation} 
\end{lemma}

\subsection*{6.7.~A $\mathcal{ C}^{1,\alpha}$ estimate for
composition of mappings} Recall that $\Delta$ is the open unit disc in
$\C$ and that $\partial \Delta$ is its boundary, namely the unit
circle. We shall constantly denote the complex variable in $\overline{
\Delta}:= \Delta \cup \partial \Delta$ by $\zeta=\rho \, e^{i
\theta}$, where $0 \leq \rho \leq 1$ and where $\vert \theta \vert
\leq \pi$, except when we consider two points $\zeta'=e^{i\theta'}$,
$\zeta''=e^{i\theta''}$, in which case we may obviously choose $\vert
\theta ' \vert, \ \vert \theta'' \vert \leq 2\pi$ with $0\leq \vert
\theta ' -\theta'' \vert \leq \pi$. Let now $X(\zeta)=\left(X_1
(\zeta), \dots, X_n (\zeta)\right)$ be an $\R^n$-valued mapping which
is of class $\mathcal{ C}^{1,\alpha}$ on the unit circle $\partial
\Delta$. We define its $\mathcal{ C}^{1, \alpha}$-norm precisely by
\def\theequation{6.8}\begin{equation}
\left\{
\aligned
\vert \vert X \vert \vert_{\mathcal{ C}^{1,\alpha}(\partial 
\Delta)}:= 
& \
\sup_{\vert \theta \vert \leq \pi}\, 
\left\vert X\left(e^{i\theta} \right) \right\vert+ 
\sup_{\vert \theta \vert \leq \pi}\, 
\left\vert
\frac{d X\left(e^{i\theta}\right)}{d \theta}
\right\vert+ \\
& \ \ \ \ \
+\sup_{0 < \left\vert \theta' - \theta''\right\vert \leq \pi} \ 
\frac{
\left\vert\frac{dX(e^{i\theta'})}{d\theta}-
\frac{dX(e^{i\theta''})}{d\theta}
\right\vert
}
{\left\vert \theta' - \theta''\right\vert^\alpha},
\endaligned\right.
\end{equation}
and we define its $\mathcal{ C}^1$-norm $\vert \vert X \vert
\vert_{\mathcal{ C}^1(\partial \Delta)}$ by keeping only the first two
terms. Let $h$ be as in Lemma~6.5.

\def\thelemma{{\bf 6.9}}\begin{lemma} Under the above assumptions, if
moreover $\left\vert X \left(e^{i\theta}\right) \right\vert \leq\rho$
for all $\theta$ with $\vert \theta \vert \leq \pi$, where $\rho \leq
\rho_1$, then we have the following three estimates{\rm :}
\def\theequation{6.10}\begin{equation}
\small
\left\{
\aligned
\vert \vert h(X) \vert \vert_{\mathcal{ 
C}^{1,\alpha}(\partial \Delta)}
\leq 
& \
\sup_{\vert x \vert \leq \rho} \, 
\vert h(x) \vert +
\left(
\sum_{k=1}^n\, 
\sup_{\vert x \vert \leq \rho}\, 
\left\vert h_{x_k}(x) \right\vert
\right)\cdot \vert \vert 
X \vert \vert_{\mathcal{ C}^1(\partial \Delta)}+\\
& \
+\left(
\sum_{k_1,k_2=1}^n\, \sup_{\vert x\vert\leq \rho}
\left\vert h_{x_{k_1}x_{k_2}} (x) \right\vert
\right)\cdot \pi^{1-\alpha} \cdot 
\left[
\vert \vert X \vert \vert_{\mathcal{ C}^1(\partial \Delta)}
\right]^2+ \\
& \
+
\left(
\sum_{k=1}^n\, 
\sup_{\vert x \vert\leq \rho} \, 
\vert h_{x_k} (x) \vert
\right)\cdot \vert \vert X 
\vert \vert_{\mathcal{ C}^{1,\alpha}
(\partial \Delta)}. \\
\sum_{k=1}^n \, 
\left\vert \left\vert h_{x_k}(X) \right\vert
\right\vert_{\mathcal{ C}^\alpha
(\partial \Delta)}\leq 
& \
\sum_{k=1}\, \sup_{\vert x \vert \leq \rho} \, 
\left\vert h_{x_k} (x) \right\vert + \\
& \ \
+
\left(
\sum_{k_1,k_2=1}^n\, 
\sup_{\vert x \vert \leq \rho} \, 
\left\vert h_{x_{k_1}x_{k_2}}(x) \right\vert
\right)\cdot \pi^{1-\alpha} \cdot
\vert \vert X \vert \vert_{\mathcal{ C}^1(\partial 
\Delta)},
\\
\sum_{k_1,k_2=1}^n \, 
\left\vert 
\left\vert
h_{x_{k_1}x_{k_2}}(X)
\right\vert
\right\vert_{\mathcal{ C}^\alpha(\partial \Delta)} \leq
& \
\sum_{k_1,k_2=1}^n\, 
\sup_{\vert x \vert \leq \rho} \,
\left\vert
h_{x_{k_1}x_{k_2}}(x)
\right\vert + \\
& \ \
+\vert\vert
h \vert \vert_{\mathcal{ C}^{2,\alpha}(
\{\vert x \vert \leq \rho\})}\cdot
\left(
\vert \vert X \vert \vert_{\mathcal{C}^1(\partial \Delta)}
\right)^\alpha.
\endaligned\right.
\end{equation}
\end{lemma}

\proof
We summarize the computations. Applying the definition~\thetag{6.8},
using the chain rule for the calculation of
$dh\left(X\left(e^{i\theta}\right)\right)/ d\theta$, and using the
trivial inequality $\left\vert a'b'-a''b''\right\vert \leq \left\vert
a' \right\vert \cdot \left\vert b'-b''\right\vert+ \left\vert b''
\right\vert\cdot \left\vert a' -a'' \right\vert$, we may majorize
\def\theequation{6.11}\begin{equation}
\small
\left\{
\aligned
\vert \vert h(X) \vert \vert_{\mathcal{ 
C}^{1,\alpha}(\partial \Delta)}
\leq 
& \
\sup_{\vert \theta \vert \leq \pi}\, 
\left\vert h(X(e^{i\theta})) \right\vert + 
\left(\sum_{k=1}^n\, 
\sup_{\vert \theta \vert\leq \pi}\, 
\left\vert h_{x_k}(X(e^{i\theta})) \right\vert
\right)\cdot \max_{1\leq k\leq n}\, 
\sup_{\vert \theta \vert \leq \pi}\, 
\left\vert
\frac{dX_k(e^{i\theta})}{d\theta}
\right\vert+ \\
& \
\sup_{0 < \vert \theta' -\theta'' 
\vert\leq \pi} \, \sum_{k=1}^n\, 
\frac{
\left\vert
h_{x_k}(X(e^{i\theta'}))-
h_{x_k}(X(e^{i\theta''}))
\right\vert
}{\vert \theta' -\theta ''\vert^\alpha} \cdot
\max_{1\leq k\leq n}\, 
\sup_{\vert \theta ' \vert \leq \pi}\, 
\left\vert
\frac{dX_k(e^{i\theta'})}{d\theta}
\right\vert + \\
& \
\left(
\sum_{k=1}^n\, \sup_{\vert \theta''\vert\leq \pi}\, 
\left\vert h_{x_k}(e^{i\theta''}) \right\vert
\right) \cdot \left(
\max_{1\leq k\leq n}\, 
\sup_{0< \vert \theta' -\theta''\vert\leq \pi}\,
\frac{
\left\vert
\frac{dX_k(e^{i\theta'})}{d\theta}-
\frac{dX_k(e^{i\theta''})}{d\theta}
\right\vert
}
{\vert \theta' -\theta '' \vert^\alpha}\right),
\endaligned\right.
\end{equation}
which yields the first inequality of~\thetag{ 6.10} after
using~\thetag{ 6.2} for the second line of~\thetag{ 6.11} and the
trivial majoration $\left\vert \theta' - \theta''
\right\vert^{1-\alpha}\leq \pi^{1-\alpha}$. The second and the third
inequalities of~\thetag{6.10} are established similarly, which
completes the proof.
\endproof

The following direct consequence will be strongly used in 
Section~7 below.

\def\thelemma{{\bf 6.12}}\begin{lemma} Under the above assumptions,
suppose that there exist constants $c_1>0$, $K_2>0$ with $c_1 K_2 \leq
\rho_1$ such that for each $c\in\R$ with $0\leq c \leq c_1$, there
exists $X_c\in\mathcal{ C}^{1,\alpha}(\partial \Delta,\R^n)$ with
$\vert \vert X_c \vert \vert_{\mathcal{ C}^{1,\alpha}(\partial
\Delta)}\leq c \cdot K_2$. Then there exists a constant $K_3>0$ such
that the following three estimates hold{\rm :}
\def\theequation{6.13}\begin{equation}
\left\{
\aligned
\vert \vert h(X_c)\vert 
\vert_{{\mathcal C}^{1,\alpha}(\partial 
\Delta)} \leq 
& \
c^{2+\alpha} \cdot K_3, \\
\sum_{k=1}^n\, 
\vert \vert h_{x_k}(X_c) 
\vert \vert_{\mathcal{ C}^\alpha(
\partial \Delta)} \leq 
& \
c^{1+\alpha} \cdot K_3, \\
\sum_{k_1,k_2=1}^n\, 
\left\vert \left\vert
h_{x_{k_1}x_{k_2}} (X_c)
\right\vert \right\vert_{\mathcal{ 
C}^\alpha(\partial \Delta)} \leq
& \
c^\alpha \cdot K_3.
\endaligned\right.
\end{equation}
\end{lemma}

\proof
Applying Lemmas~6.4 and~6.9, we see that it suffices to choose 
\def\theequation{6.14}\begin{equation}
K_3:=\max \left(
K_1K_2^{2+\alpha}(3+\pi^{1-\alpha}), \ 
K_1K_2^{1+\alpha}(1+\pi^{1-\alpha}), \ 
2K_1K_2^\alpha
\right),
\end{equation}
which completes the proof.
\endproof

Up to now, we have introduced three positive constants $K_1$, $K_2$,
$K_3$. In Sections~7, 8 and~9 below, we shall introduce further
positive constants $K_4$, $K_5$, $K_6$, $K_7$, $K_8$, $K_9$, $K_{10}$,
$K_{11}$, $K_{12}$, $K_{13}$, $K_{14}$, $K_{15}$, $K_{16}$, $K_{17}$,
$K_{18}$ and $K_{19}$, whose precise value will not be important.

\section*{\S7.~Families of analytic discs half-attached 
to maximally real submanifolds}

\subsection*{7.1.~Preliminary}
Let $E\subset \C^n$ be an arbitrary subset and let $A:
\overline{\Delta}\to \C^n$ be a continuous mapping, holomorphic in
$\Delta$. If $\partial^+\Delta:= \{\zeta\in\partial \Delta: {\rm Re}\,
\zeta \geq 0\}$ denotes the {\sl positive half-boundary} of $\Delta$,
we say that $A$ is {\sl half-attached} to $E$ if $A(\partial^+
\Delta)\subset E$. Such analytic discs which are glued in part to a
geometric object were studied by S.~Pinchuk in~\cite{p} to establish a
boundary uniqueness principle about continuous functions on a
maximally real submanifold of $\C^n$ which extend holomorphically to a
wedge. Further works on the CR edge of the wedge theorem using discs
partly attached to generic submanifolds were achieved by
R.~Ayrapetian~\cite{a} and by A.~Tumanov~\cite{tu2}.

In this section, we shall construct local families of analytic discs
$Z_{c,x,v}^1(\zeta): \overline{\Delta} \to \C^n$, where $c\in \R^+$ is
small, where $x\in\R^n$ is small and where $v\in\R^n$ is small, which
are half-attached to a $\mathcal{ C}^{2,\alpha}$-smooth maximally real
submanifold $M^1$ of $\C^n$, which satisfy $Z_{c,0,v}^1(1)\equiv
p_1\in M^1$, such that the boundary point $Z_{c,x,v}^1(1)$ covers a
neighborhood of $p_1$ in $M^1$ as $x$ varies ($c$ and $v$ being fixed)
and such that the tangent vector $\frac{\partial Z_{c,0,v}^1}{\partial
\theta}(1)$ at the fixed point $p_1$ covers a cone in
$T_{p_1}M^1$. These families will be used in Sections~8 and~9 below
for the final steps in the proof of the main Proposition~5.12. With
this choice, when $x$ varies, $v$ varies and $\zeta$ varies (but $c$
is fixed), the set of points $Z_{c,x,v}^1(\zeta)$, covers a thin wedge
of edge $M^1$ at $p_1$. Similar families of analytic discs were
constructed in~\cite{ber} to reprove S.~Pinchuk's boundary uniqueness
theorem, with $M^1$ of class $\mathcal{ C}^\infty$, using a method
(implicit function theorem in Banach spaces) which in the case where
$M$ is of class $\mathcal{ C}^{\kappa,\alpha}$ necessarily induces a
loss of smoothness, yielding families of analytic discs which are only
of class $\mathcal{ C}^{\kappa-1,\alpha}$. Since we want our families
to be of class at least $\mathcal{ C}^2$ and since $M$ is only of
class $\mathcal{ C}^{2,\alpha}$, we shall have to proceed differently.

To summarize symbolically the structure of the desired family:
\def\theequation{7.2}\begin{equation}
Z_{c,x,v}^1(\zeta): \left\{
\aligned
c= 
& \
\text{\rm small scaling factor}, \\
x=
& \
\text{\rm translation parameter}, \\
v=
& \
\text{\rm rotation parameter}, \\
\zeta=
& \
\text{\rm unit disc variable}. \\
\endaligned
\right.
\end{equation} 
We shall begin our constructions in the ``flat'' case where the
maximally real submanifold $M^1$ coincides with $\R^n$ and then
perform a pertubation argument, using the scaling parameter $c$ in an
essential way.

\subsection*{7.3.~A family of analytic
discs sweeping $\R^n \subset \C^n$ with prescribed first order jets}
We denote the coordinates over $\C^n$ by $z= x+iy= (x_1 + iy_1, \dots,
x_n+ iy_n)$. Let $c\in \R$ with $c \geq 0$ be a ``scaling factor'',
let $n \geq 2$, let $x= (x_1, \dots, x_n) \in \R^n$, let $v= (v_1,
\dots, v_n) \in \R^n$ and consider the algebraically parametrized
family of analytic discs defined by
\def\theequation{7.4}\begin{equation} 
B_{c,x,v}(s+it):=
\left(x_1+cv_1(s+it),\dots, 
x_n+cv_n(s+it)\right),
\end{equation}
where $s+it \in \C$ is the holomorphic variable. For $c \neq 0$, the
map $B_{c, x, v}$ embeds the complex line $\C$ into $\C^n$ and sends
$\R$ into $\R^n$ with {\it arbitrary first order jet at $0$}: center
point $B_{c, x, v}(0) =x$ and tangent direction $\partial B_{c, x, v}
(s) /\partial s \vert_{ s=0}= cv$.

To localize our family of analytic discs, we restrict the
map~\thetag{7.4} to the following specific set of values: $0\leq c\leq
c_0$ for some $c_0>0$; $\vert x \vert \leq c$~; $\vert v \vert \leq
2$~; and $\vert s+it \vert \leq 4$. To localize $\R^n$, we shall
denote $M^0:= \{x\in\R^n: \,\vert x \vert \leq \rho_0\}$, where
$\rho_0 >0$, and we notice that $B_{c,x,v}(\{\vert s+it \vert \leq
4\}) \subset M^0$ for all $c$, all $x$ and all $v$ provided that $c_0
\leq \rho_0/9$.

We then consider the mapping $(s+it) \longmapsto B_{c,x, v} (s+it)$ as
a local (nonsmooth) analytic disc defined on the rectangle
$\{s+it\in\C: \, \vert s \vert \leq 4, 0\leq t\leq 4\}$ whose bottom
boundary part $B_{c,x, v}([-4,4])$ is a small real segment contained
in $\R^n$.

\subsection*{7.5.~A useful conformal equivalence}
Next, we have to get rid of the corners of the rectangle $\{s+it\in
\C: \, \vert s \vert \leq 4, 0\leq t \leq 4\}$. We proceed as
follows. In the complex plane equipped with coordinates $s+it$, let
$\mathcal{D}(i\sqrt{3},2)$ be the open disc of center $i\sqrt{3}$ and
of radius $2$. Let $\mu: (-2,2)\to [0,1]$ be an {\it even}\,
$\mathcal{ C}^\infty$-smooth function satisfying $\mu(s)=0$ for $0\leq
s\leq 1$; $\mu(s)>0$ and $d\mu(s)/ds>0$ for $1< s < 2$; and $\mu(s)=
\sqrt{3}-\sqrt{4-s^2}$ for $\sqrt{3} \leq s < 2$. The simply connected
domain $C^+\subset \{t>0\}$ which is represented in {\sc Figure~1}
just below may be formally defined as
\def\theequation{7.6}\begin{equation}
\left\{
\aligned
C^+\cap \{t\geq \sqrt{3}-1\}
:= \mathcal{D}(i\sqrt{3},2)\cap \{t\geq
\sqrt{3}-1\}, \\
C^+\cap \{0<t< \sqrt{3}-1\}
:= \{s+it\in\C: \, t>
\mu(s)\}.
\endaligned\right.
\end{equation}

\bigskip
\begin{center}
\input conformal.pstex_t
\end{center}

\noindent
Let $\Psi: \Delta\to C^+$ be a conformal equivalence (Riemann's
theorem). Since the boundary $\partial C^+$ is $\mathcal{
C}^\infty$-smooth, the mapping $\Psi$ extends as a $\mathcal{
C}^\infty$-smooth diffeomorphism $\partial \Delta \to \partial C^+$.
Remind that $\partial^+ \Delta:=\{ \zeta \in \C: \vert \zeta \vert =1,
\ {\rm Re}\, \zeta \geq 0\}$ is the positive half-boundary of
$\Delta$. Then after a reparametrization of $\Delta$, we can (and we
shall) assume that $\Psi(\partial^+\Delta)=[-1,1]$, $\Psi (1) =0$ and
$\Psi(\pm i)=\pm 1$. It follows that $d\Psi \left(e^{i \theta}\right)
/d \theta$ is a positive real number for all
$e^{i\theta}\in\partial^+\Delta$. Although the precise shape of $C^+$
and the specific expression of $\Psi$ will not be crucial in the
sequel, it will be convenient to fix them once for all.

\subsection*{7.7.~Flat families of half-attached analytic discs}
Thanks to $\Psi$, we can define a family of small analytic
discs which are half-attached to the ``{\sl flat}'' maximally real
manifold $M^0\equiv \{x\in \R^n: \ \vert x \vert \leq \rho_0
\}$ as follows
\def\theequation{7.8}\begin{equation}
Z_{c,x,v}^0(\zeta):= 
B_{c,x,v} \left (\Psi (\zeta) \right)=
\left(
x+cv\Psi(\zeta)
\right).
\end{equation}
We then have $Z_{c, x, v}^0( \partial^+ \Delta) \subset M^0$ and $Z_{
c, x, v }^0 (1) = x$. Notice that every disc $Z_{c,x,v}
\left(\overline{\Delta}\right)$ 
is contained in a single complex line. Starting
with a maximally real submanifold of $\C^n$, as in Proposition~5.12, but
dealing with the ``flat'' maximally real submanifold $M^0 \equiv
\R^n$, we first construct a ``flat model'' of the desired family of
analytic disc.

\def\thelemma{{\bf 7.9}}\begin{lemma} 
Let $M^0=\{x\in \R^n: \ \vert x \vert \leq \rho_0\}$ be the ``flat''
local maximally real submanifold defined above, let $p_0 \equiv 0 \in
M^0$ denote the origin and let $v_0\in T_{p_0}M^0$ be a tangent vector
with $\vert v^0\vert =1$. Then there exists a constant $\Lambda_0>0$
and there exists a $\mathcal{ C }^\infty$-smooth family $A_{c,x,v}^0
(\zeta)$ of analytic discs defined for $c\in \R$ with $0 \leq c \leq
c_0$ for some $c_0>0$ satisfying $c_0\leq \rho_0/9$, for $x\in \R^n$
with $\vert x \vert \leq c$ and for $v\in \R^n$ with $\vert v \vert
\leq c$ which enjoy the following six properties{\rm :}

\smallskip
\begin{itemize}
\item[{\bf ($\mathbf{1_0}$)}] \
$A_{c,0,v}^0(1)=p_0=0$ for all $c$ and all $v$.
\item[{\bf ($\mathbf{2_0}$)}] \
$A_{c,x,v}^0: \overline{\Delta}\to \C^n$ is
an embedding and
$\left\vert A_{c,x,v}^0(\zeta) \right \vert \leq c\cdot \Lambda_0$ 
for all $c$, all $x$, all $v$ and all $\zeta$.
\item[{\bf ($\mathbf{3_0}$)}] \
$A_{c,x,v}^0(\partial^+\Delta)
\subset M^0$ for all $c$, all $x$ and all $v$.
\item[{\bf ($\mathbf{4_0}$)}] \
$\frac{\partial A_{c,0,0}^0}{\partial 
\theta}(1)$ is a positive multiple of 
$v_0$ for all $c\neq 0$.
\item[{\bf ($\mathbf{5_0}$)}] \
For all $c$, all $v$ and all $e^{i\theta}\in \partial^+\Delta$, the
mapping $x\longmapsto A_{c,x,v}^0\left(e^{i\theta}\right) \in M^0$ is
of rank $n$.
\item[{\bf ($\mathbf{6_0}$)}] \
For all $e^{i \theta} \in \partial^+ \Delta$, all $c \neq 0$ and all
$x$, the mapping $v \longmapsto \frac{ \partial A_{c, x, v}^0}{
\partial \theta }\left(e^{ i\theta}\right)$ is of rank $n$ at $v=0$.
Consequently, the positive half-lines $\R^+ \cdot \frac{ \partial
A_{c,0,v}}{ \partial \theta} (1)$ describe an open infinite cone
containing $v_0$ with
vertex $p_0$ in $T_{p_0} M^0$ when $v$ varies.
\end{itemize}
\end{lemma}

\proof
Proceeding as in the proof of Lemma~5.37, we can find a new affine
coordinate system centered at $p_0$ and stabilizing $\R^n$, which we
shall still denote by $(z_1,\dots,z_n)$, 
in which the vector $v_0$
has coordinates $(0,\dots,0,1)$. In this coordinate system, we then
construct the family $Z_{c,x,v}^0(\zeta)$ as in~\thetag{7.8} above and
we define the desired family simply as follows:
\def\theequation{7.10}\begin{equation}
A_{c,x,v}^0(\zeta):= Z_{c,x,v_0+v}^0(\zeta),
\end{equation}
where we restrict the variations of the parameter $v$ to $\vert v
\vert \leq c$. Notice that every disc $A_{c,x,v}^0
\left(\overline{\Delta}\right)$ is contained in a single complex line.
All the properties are then elementary consequences of the explicit
expression~\thetag{7.8} of $Z_{c,x,v}^0(\zeta)$.

Finally, we notice that it follows from properties {\bf
($\mathbf{5_0}$)} and {\bf ($\mathbf{6_0}$)} that the set of points
$A_{c,x,v}^0(\zeta)$, where $c>0$ is fixed, where $x$ varies, where $v$
varies and where $\zeta$ varies covers a local wedge of edge $M^0$ at
$p_0$. The proof of Lemma~7.9 is complete.
\endproof

\subsection*{7.11.~Curved families of half-attached analytic discs}
Our main goal in this section is to obtain a statement similar to
Lemma~7.9 after replacing the ``flat'' maximally real submanifold
$M^0\cong \R^n$ by a ``curved'' $\mathcal{ C}^{2, \alpha}$-smooth
maximally real submanifold $M^1$. We set up a formulation which will
be appropriate for the achievement of the end of the proof of
Proposition~5.12 in the next Sections~8 and~9. In particular, we shall
have to shrink the family of half-disc $Z_{c, x,v}^1(\zeta)$ which we
will construct as a perturbation of the family $Z_{c, x,v}^0(\zeta)$
in \S7.50 below, and we shall construct discs of size $\leq c^2 \cdot
\Lambda_1$ for some constant $\Lambda_1>0$, instead of requiring that
their size is $\leq c \cdot \Lambda_1$, which would be the property
analogous to {\bf ($\mathbf{2_0}$)}. Also, we shall loose the
$\mathcal{ C}^{ 2, \alpha-0}$-smoothness with respect to the scaling
parameter $c$.

\def\thelemma{7.12}\begin{lemma}
Let $M^1$ be $\mathcal{ C}^{2,\alpha}$-smooth maximally real
submanifold of $\C^n$, let $p_1\in M^1$ and let $v_1\in T_{p_1} M^1$
be a tangent vector with $\vert v_1 \vert =1$. Then there exists
a positive constant $\Lambda_1>0$ and there exists $c_1\in \R$ with
$c_1>0$ such that for every $c\in \R$ with $0 < c \leq c_1$, there
exists a family $A_{x,v:c}^1 (\zeta)$ of analytic discs defined for
$x\in \R^n$ with $\vert x \vert \leq c^2$ and for $v\in \R^n$ with
$\vert v \vert \leq c$ which is $\mathcal{ C }^{2,\alpha-0}$-smooth
with respect to $(x,v,\zeta)$ and which enjoys the following six
properties{\rm :}

\smallskip
\begin{itemize}
\item[{\bf ($\mathbf{1_1}$)}] \
$A_{0,v:c}^1(1)=p_1$ for all $v$.
\item[{\bf ($\mathbf{2_1}$)}] \
$A_{x,v:c}^1: \overline{\Delta}\to \C^n$ is
an embedding and
$\left\vert A_{x,v:c}^1 (\zeta) \right \vert \leq c^2\cdot \Lambda_1$ 
for all $x$, all $v$ and all $\zeta$.
\item[{\bf ($\mathbf{3_1}$)}] \
$A_{x,v:c}^1(\partial^+\Delta)
\subset M^1$ for all $x$ and all $v$.
\item[{\bf ($\mathbf{4_1}$)}] \
$\frac{\partial A_{0,0:c}^1}{\partial 
\theta}(1)$ is a positive multiple of 
$v_1$ for all $c\neq 0$.
\item[{\bf ($\mathbf{5_1}$)}] \ 
The mapping $x\longmapsto A_{x,0:c}^1(1)\in M^1$ is of
rank $n$.
\item[{\bf ($\mathbf{6_1}$)}] \
The mapping $v \longmapsto \frac{ \partial A_{0, v:c}^1}{ \partial
\theta }\left(e^{ i\theta}\right)$ is of rank $n$ at $v=0$.
Consequently, as $v$ varies, the positive half-lines $\R^+ \cdot
\frac{ \partial A_{0,v:c}^1 }{ \partial \theta} (1)$ describe an open
infinite cone containing $v_1$ with vertex $p_1$ in $T_{p_1} M^1$ and
the set of points $A_{x,v:c}^1(\zeta)$, as $\vert x\vert \leq c$,
$\vert v\vert \leq c$ and $\zeta\in \Delta$ vary, covers a wedge of
edge $M^1$ at $(p_1,Jv_1)$.
\end{itemize}
\end{lemma}

In Figure~16 drawn in Section~8 after Lemma~8.3 below, we have drawn
the property that the tangent direction $\frac{ \partial A_{0,v: c}^1}{
\partial \theta} (1)$ describes an open cone in $T_{p_1}M^1$ with
vertex $p_1$. The remainder of Section~3 is devoted to complete the
proof of Proposition~7.12.

\subsection*{7.13.~Perturbed family of analytic discs half-attached to 
a maximally real submanifold} Thus, let $M^1\subset \R^n$ be a locally
defined maximally real $\mathcal{ C}^{2,\alpha}$-smooth submanifold
passing through the origin. We can assume that $M^1$ is represented by
$n$ Cartesian equations
\def\theequation{7.14}\begin{equation}
y_1=h_1(x_1,\dots,x_n),\cdots\cdots, y_n=h_n(x_1,\dots,x_n),
\end{equation} 
where $z_k= x_k+i y_k\in\C$, for $k= 1,\dots, n$, where $\vert x \vert
\leq \rho_1$ for some $\rho_1 >0$, where $h= h(x)$ is of class
$\mathcal{ C}^{2, \alpha}$ in $\{\vert x \vert < 2\rho_1\}$, and
where, importantly, $h_j(0)=h_{j,x_k}(0)= h_{j, x_{k_1}
x_{k_2}}(0)=0$, for all $j, k, k_1,k_2 =1, \dots,n$. As
in~\thetag{6.5}, we set $K_1:=\vert\vert h \vert \vert_{ \mathcal{
C}^{2,\alpha}(\{ \vert x \vert \leq \rho_1\})}$. Also, we can assume
that $v_1=(0,\dots,0,1)$.

Our goal is to show that we can produce a $\mathcal{
C}^{2,\alpha-0}$-smooth (remind $\mathcal{ C}^{2,\alpha-0} \equiv
\bigcap_{ \beta < \alpha} \, \mathcal{ C}^{2, \beta}$) family of
analytic discs $Z_{c,x,v}^1(\zeta)$ which is half-attached to $M^1$
and which is sufficiently close, in $\mathcal{ C}^2$ norm, to the
original family $Z_{c,x,v}^0(\zeta)$. After having constructed the
family $Z_{c,x,v}^1(\zeta)$, we shall define the desired family
$A_{x,v:c}^1(\zeta)$.

Let $d\in\R$ with $0\leq d \leq 1$ and let the maximally real
submanifold $M^d$ (like ``$M$ {\it d}eformed'') be defined precisely
as the set of $z= x+iy \in \C^n$ with $\vert x \vert \leq \rho_1$
which satisfy the $n$ Cartesian equations
\def\theequation{7.15}\begin{equation}
y_1=d\cdot h_1(x_1,\dots,x_n),\cdots\cdots, \, 
y_n=d\cdot h_n(x_1,\dots,x_n).
\end{equation}
Notice that $M^0 \equiv\{ x\in \R^n: \ \vert x \vert \leq \rho_1\}$ is
essentially the same piece $M^0$ of $\R^n$ as in Lemma~7.9 (which
contains the $M^0$ of Lemma~7.9 if we choose $\rho_0\leq \rho_1$) and
notice that $\left. M^d\right \vert_{d=1} \equiv M^1$. Even better, we
shall construct for each $d$ with $0 \leq d\leq 1$ a one-parameter
family of analytic dics $Z_{c,x, v}^d (\zeta)$ which is of class at
least $\mathcal{ C}^{2, \alpha-0}$ with respect to all variables and
which is half-attached to $M^d$, by proceeding as follows.

First of all, the analytic disc $Z_{c,x, v}^d (\zeta)=: X_{c,x,
v}^d (\zeta)+ iY_{c,x, v}^d( \zeta)$ is half-attached to $M^d$ if
and only if
\def\theequation{7.16}\begin{equation}
Y_{c,x,v}^d(\zeta)=d\cdot h\left(
X_{c,x,v}^d(\zeta)
\right), \ \ \ \ \ 
{\rm for} \ \zeta\in\partial^+\Delta.
\end{equation}
Furthermore, $Y_{c,x, v}^d$ should be a harmonic conjugate of
$X_{c,x, v}^d$. However, the condition~\thetag{7.16} does not give any
relation between $X_{ c,x, v}^d$ and $Y_{c,x, v}^d$ on the
negative part $\partial^-\Delta$ of the unit circle. To fix this
point, we shall assign the following more
complete equation
\def\theequation{7.17}\begin{equation}
Y_{c,x,v}^d(\zeta)= d\cdot h\left(
X_{c,x,v}^d(\zeta)\right)+Y_{c,x,v}^0(\zeta), 
\ \ \ \ \
\text{\rm for {\it all}} \ \zeta\in\partial\Delta,
\end{equation}
which coincides with~\thetag{7.16} for $\zeta\in
\partial^+ \Delta$, since we have $Z_{c, x, v}^0( \partial^+
\Delta) \subset \R^n$ by construction ({\it cf.}~\thetag{7.8}).
 
As in~\cite{tu2}, \cite{tu3}, \cite{mp1}, \cite{mp3}, we denote by
$T_1$ the Hilbert transform (harmonic conjugate operator) on $\partial
\Delta$ vanishing at $1$, namely $(T_1 X) (1)=0$, whence
$T_1(T_1(X))=-X+X(1)$. By Privalov's theorem, for every integer
$\kappa \geq 0$ and every $\alpha\in\R$ with $0< \alpha < 1$, its norm
$\vert \vert \vert T_1 \vert \vert \vert_{ \kappa, \alpha}$ as an
operator $\mathcal{ C}^{ \kappa, \alpha}( \partial \Delta, \R^n) \to
\mathcal{ C}^{\kappa, \alpha} (\partial \Delta, \R^n)$ is finite and
explodes as $\alpha$ tends either to $0$ or to $1$. Also, we shall
require that $X_{c,x, v}^d (1) =x$, whence $Y_{ c,x, v}^d (1)= d
\cdot h(x)$.

With this choice, the mapping $\zeta \mapsto Y_{c,x, v}^d (\zeta)$
should necessarily coincide with the harmonic conjugate $\zeta \mapsto
\left[T_1 X_{c,x, v}^d\right] 
(\zeta)+ d\cdot h(x)$ (this property is already
satisfied for $d=0$) and we deduce that $X_{c, x,v}^d( \zeta)$
should satisfy the following Bishop type equation
\def\theequation{7.18}\begin{equation}
X_{c,x, v}^d(\zeta) =
-T_1\left[
d \cdot h\left(
X_{c,x,v}^d
\right)
\right](\zeta)+X_{c,x,v}^0(\zeta), 
\ \ \ \ \ 
\text{\rm for all} \ \zeta\in\partial\Delta.
\end{equation}
Conversely, if $X_{c,x, v}^d$ is a solution of this functional
equation, then setting $Y_{c,x, v}^d(\zeta):= T_1 X_{c,x, v}^d
(\zeta)+ d\cdot h(x)$, it is easy to see that the analytic disc $Z_{c,x,
v}^d (\zeta):= X_{c,x,v}^d (\zeta)+ iY_{c,x, v}^d(\zeta)$ is
half-attached to $M^d$ and more precisely, satisfies the
equation~\thetag{7.16}.

Applying now Theorem~1.2 of \cite{tu3}, we deduce that if the given
positive number $c_1$ is sufficiently small, and if $c$ satisfies
$0\leq c \leq c_1$, there exists a unique solution $X_{c,x,
v}^d(\zeta)$ to~\thetag{7.18} which is of class $\mathcal{ C}^{2,
\alpha}$ with respect to $\zeta$ and of class $\mathcal{ C}^{2,\alpha
-0}$ with
respect to all variables $(c,x, v, \zeta)$ with $0\leq c \leq c_1$,
$\vert x\vert \leq c$, $\vert v \vert \leq 2$ and $\zeta\in
\overline{\Delta}$. We shall now estimate the difference $\vert \vert
Z_{c,x,v}^d-Z_{c,x,v}^0\vert \vert_{ \mathcal{ C}^{1,\alpha} (\partial
\Delta)}$ and prove that it is bounded by a constant times
$c^{2+\alpha}$. In particular, if $c_1$ is sufficiently small, this
will imply that $Z_{c,x,v}^d$ is nonconstant.

\subsection*{7.19.~Size of the solution 
$X_{c,x,v}^d(\zeta)$ in $\mathcal{ C}^{1,\alpha}$ norm} Following the
beginning of the proof of Theorem~1.2 in~\cite{tu3}, we introduce the
mapping
\def\theequation{7.20}\begin{equation}
F: \ X(\zeta)\longmapsto X_{c,x,v}^0(\zeta)
-T_1\left[
d\cdot h(X)
\right](\zeta)
\end{equation}
from a neighborhood of $0$ in $\mathcal{ C}^{1,\alpha}( \partial
\Delta, \R^n)$ to $\mathcal{ C}^{1,\alpha}(\partial \Delta, \R^n)$,
and then, as in~\cite{b}, we introduce a Picard iteration processus by
defining $X\{0\}_{c,x,v}^d(\zeta):= X_{c,x,v}^0(\zeta)$ and for every
integer $\nu\geq 0$
\def\theequation{7.21}\begin{equation}
X\{\nu+1\}_{c,x,v}^d(\zeta):= 
F\left(X\{\nu\}_{c,x,v}^d(\zeta)\right).
\end{equation}
In a first moment, A.~Tumanov proves in~\cite{tu3} that the sequence
$\left( X\{\nu \}_{c,x, v}^d (\zeta) \right)_{ \nu \in \N}$ converges
towards the unique solution $X_{c,x, v}^d(\zeta)$ of~\thetag{7.18} in
$\mathcal{ C}^{1, \alpha}( \partial \Delta)$. Admitting this
convergence result, we need to extract the supplementary information
that $\left \vert \left \vert X_{c,x, v}^d \right \vert\right\vert_{
\mathcal{ C}^{1, \alpha}( \partial \Delta)} \leq c\cdot K_2$ for some
positive constant $K_2$, which will play the role of the constant
$K_2$ of Lemma~6.12.

To get this information, we observe that by construction
({\it cf.}~\thetag{7.8}) there exists a constant $K_4>0$ such that
\def\theequation{7.22}\begin{equation}
\left\vert \left\vert X_{c,x,v}^0\right\vert 
\right\vert_{\mathcal{
C}^{2,\alpha}(\partial \Delta)}\leq c \cdot K_4.
\end{equation}
Also, we set $K_5:= K_1 (3+ \pi^{1-\alpha}) \vert \vert \vert T_1
\vert \vert \vert_{ \mathcal{ C}^{1, \alpha}( \partial \Delta)}$.

\def\thelemma{{\bf 7.23}}\begin{lemma}
With these notations, if 
\def\theequation{7.24}\begin{equation}
c_1\leq\min\left( 
\frac{\rho_1}{2K_4}, \ 
\left(
\frac{1}
{2^{2+\alpha}\, K_4^{1+\alpha} \, K_5}
\right)^{\frac{1}{1+\alpha}}\right),
\end{equation}
then the solution of~\thetag{7.18} satisfies $\left\vert X_{c,x, v}^d
\left( e^{i \theta }\right) \right\vert \leq \rho_1$ for all
$e^{i\theta}\in\partial \Delta$ and there exists a constant $K_2>0$
such that
\def\theequation{7.25}\begin{equation}
\left
\vert \left \vert X_{c,x, v}^d 
\right \vert \right \vert_{ \mathcal{
C}^{1, \alpha} ( \partial \Delta)} \leq c \cdot K_2. 
\end{equation}
In fact, it
suffices to choose $K_2:= 2 K_4$.
\end{lemma}

\proof
Indeed, using Lemmas~6.4 and~6.9, if $X \in \mathcal{ C}^{1, \alpha}(
\partial \Delta, \R^n)$ satisfies $\left\vert 
X\left( e^{i \theta }\right) \right\vert \leq
\rho_1$ for all $e^{i\theta}\in\partial\Delta$ and $\vert \vert X \vert
\vert_{ \mathcal{ C}^{1, \alpha} (\partial \Delta )} \leq c \cdot
2K_4$ for all $c \leq c_1$, where $c_1$ is as in~\thetag{7.24}, we
may estimate (remind $0 \leq d \leq 1$)
\def\theequation{7.26}\begin{equation}
\left\{
\aligned
\vert \vert F(X) \vert \vert_{\mathcal{ C}^{1, 
\alpha}(\partial \Delta)} 
& \
\leq
\left\vert \left\vert X_{c,x,v}^0 \right\vert 
\right\vert_{\mathcal{ C}^{1,\alpha}
(\partial \Delta)}+
\vert \vert \vert T_1 \vert \vert \vert_{\mathcal{ C}^{1, \alpha}
(\partial \Delta)} \cdot
\vert \vert h(X) \vert \vert_{\mathcal{ C}^{1, \alpha}
(\partial \Delta)} \\
& \
\leq 
c\cdot K_4+
\vert \vert \vert T_1 \vert \vert \vert_{\mathcal{ C}^{1, \alpha}
(\partial \Delta)} \cdot K_1
( c \cdot 2K_4)^{2+\alpha}(
3+\pi^{1-\alpha}) \\
& \
= 
c\cdot \left(
K_4+c^{1+\alpha} 2^{2+\alpha} K_4^{2+\alpha} K_5
\right) \\
& \
\leq 
c\cdot (K_4+c_1^{1+\alpha} 
2^{2+\alpha}K_4^{2+\alpha} K_5) \\
& \
\leq c\cdot 2K_4.
\endaligned\right.
\end{equation}
Notice that from the last inequality, it also follows that $\left\vert
F\left( X\left( e^{i\theta} \right) \right)\right\vert \leq \rho_1$
for all $e^{i\theta}\in\partial \Delta$. Consequently, the processus
of successive approximations~\thetag{7.21} is well defined for each
$\nu\in\N$ and from the inequality~\thetag{7.26}, we deduce that the
limit $X_{c,x,v}^d$ satisfies the desired estimate $\left\vert
\left\vert X_{c,x, v}^d \right\vert \right\vert_{ \mathcal{ C}^{1,
\alpha}( \partial \Delta )}\leq c \cdot 2K_4$, which completes the
proof.
\endproof

\def\thecorollary{{\bf 7.27}}\begin{corollary}
Under the above assumptions, there exists a constant
$K_6>0$ such that
\def\theequation{7.28}\begin{equation}
\left\vert \left\vert X_{c,x,v}^d - X_{c,x,v}^0 
\right\vert \right\vert_{
\mathcal{ C}^{1,\alpha}(\partial \Delta)} \leq
c^{2+\alpha} \cdot K_6.
\end{equation}
\end{corollary}

\proof
We estimate
\def\theequation{7.29}\begin{equation}
\left\{
\aligned
\left\vert \left\vert X_{c,x,v}^d - X_{c,x,v}^0 
\right\vert \right\vert_{
\mathcal{ C}^{1,\alpha}(\partial \Delta)} 
& \
\leq
\vert \vert \vert T_1 \vert \vert \vert_{\mathcal{ C}^{1,\alpha}
(\partial \Delta)} \cdot 
\left\vert \left\vert h
\left(X_{c,x,v}^d\right) \right\vert 
\right\vert_{\mathcal{ C}^{1,\alpha}
(\partial \Delta)} \\
& \ 
\leq 
\vert \vert \vert T_1 \vert \vert \vert_{\mathcal{ C}^{1,\alpha}
(\partial \Delta)} \cdot 
K_1 (c\cdot 2K_4)^{2+\alpha}(3+\pi^{1-\alpha}) \\
& \
\leq
c^{2+\alpha} \cdot K_5 (2K_4)^{2+\alpha}.
\endaligned\right.
\end{equation}
so that it suffices to set 
$K_6:= K_5(2K_4)^{2+\alpha}$.
\endproof

\subsection*{7.30.~Smallness of the deformation in $\mathcal{ C}^2$
norm} As was already noticed (and admitted), the solution
$X_{c,x,v}^d(\zeta)$ is in fact $\mathcal{ C}^{2,\alpha}$-smooth with
respect to $\zeta$ and $\mathcal{ C}^{2,\alpha-0}$-smooth with respect
to all variables $(d,c,x,v,\zeta)$. We can therefore differentiate
twice Bishop's equation~\thetag{7.18}. First of all, if $X\in\mathcal{
C}^{2,\alpha-0}( \partial \Delta, \R^n)$, we remind the commutation
relation $\frac{\partial }{\partial \theta} (TX)= T\left(
\frac{\partial X}{\partial \theta} \right)$, whence
\def\theequation{7.31}\begin{equation}
\frac{\partial}{\partial \theta} \left(
T_1X
\right)= 
T\left(
\frac{\partial X}{\partial \theta}
\right),
\end{equation}
since $T_1X= TX-TX(1)$. We may then compute the first order 
derivative of~\thetag{7.18} with respect 
to $\theta$:
\def\theequation{7.32}\begin{equation}
\left\{
\aligned
\frac{\partial }{\partial \theta} 
X_{c,x,v}^d\left(e^{i\theta}\right) -
& \
\frac{\partial }{\partial \theta}
X_{c,x,v}^0\left(e^{i\theta}\right)=
-T\left[
d\cdot \sum_{l=1}^n\, 
\frac{\partial h}{\partial x_l} \left(X_{c,x,v}^d\right) \
\frac{\partial X_{l;c,x,v}^d}{\partial \theta}
\right]\left(e^{i\theta}\right).
\endaligned\right.
\end{equation}
and then its second order partial derivatives $\partial^2 / \partial
v_k \partial \theta$, for $k=1,\dots,n$, without writing the argument
$e^{i\theta}$:
\def\theequation{7.33}\begin{equation}
\left\{
\aligned
\frac{\partial^2 X_{c,x,v}^d}{\partial v_k\partial \theta}
-
\frac{\partial^2X_{c,x,v}^0}{\partial v_k
\partial \theta}
= 
& \
-T\left[
d\cdot \sum_{l_1,l_2=1}^n\,
\frac{\partial^2 h}{\partial x_{l_1} \partial x_{l_2}}\left(
X_{c,x,v}^d\right) \
\frac{\partial X_{l_1; c,x,v}^d}{\partial v_k} \,
\frac{\partial X_{l_2; c,x,v}^d}{\partial \theta} + \right. \\
& \
\left. \ \ \ \ \ \ \ \ \ \
+d\cdot \sum_{l=1}^n\, 
\frac{\partial h_j}{\partial x_l} \left(
X_{c,x,v}^d
\right) \
\frac{\partial^2 X_{l;c,x,v}^d}{\partial v_k 
\partial \theta}
\right].
\endaligned\right.
\end{equation}
Let now $K_2$ be as in~\thetag{7.25} and let 
$K_3$ be as in Lemma~6.12, applied to 
$X_{c,x,v}^d(\zeta)$.

\def\thelemma{7.34}\begin{lemma}
If in addition to the inequality~\thetag{7.24}, the constant
$c_1$ satisfies the inequality
\def\theequation{7.35}\begin{equation}
c_1\leq \left(
\frac{1}{2K_3 \vert \vert \vert 
T \vert \vert \vert_{
\mathcal{ C}^\alpha(\partial \Delta)}}
\right)^{\frac{1}{1+\alpha}},
\end{equation}
then there exists a positive constant $K_7>0$
such that for all $d$, all $c$, all $x$, all $v$, and for
$k=1,\dots,n$, the following two estimates hold{\rm :}
\def\theequation{7.36}\begin{equation}
\left\{
\aligned
\left\vert\left\vert
\frac{\partial^2 X_{c,x,v}^d}{
\partial v_k \partial \theta}-
\frac{\partial^2 X_{c,x,v}^0}{
\partial v_k \partial \theta}
\right\vert\right\vert_{\mathcal{ 
C}^\alpha(\partial \Delta)}\leq 
& \
c^{2+\alpha} \cdot K_7, \\
\left\vert\left\vert
\frac{\partial^2 X_{c,x,v}^d
}{\partial \theta^2}-
\frac{\partial^2 X_{c,x,v}^0
}{\partial \theta^2}
\right\vert\right\vert_{\mathcal{ 
C}^\alpha(\partial \Delta)}\leq 
& \
c^{2+\alpha} \cdot K_7.
\endaligned\right.
\end{equation}
\end{lemma}

\proof
We check only the first inequality, the proof of the second being
totally similar. According to Lemma~1.6 in~\cite{ tu3}, 
there exists a solution 
$\frac{\partial^2 X_{c,x,v}^d}{
\partial v_k \partial \theta}$
to the linearized Bishop equation~\thetag{ 7.33},
hence it suffices to make an estimate.

Introducing for the second line of~\thetag{7.33} a
new simplified notation $\mathcal{ R}:= -T\left[ d\cdot
\sum_{l_1,l_2=1}^n\, \frac{\partial^2 h}{ \partial x_{l_1} \partial
x_{l_2} }\left( X_{c,x,v}^d \right) \ \frac{\partial X_{l_1;
c,x,v}^d}{\partial v_k} \, \frac{\partial X_{l_2; c,x,v}^d}{ \partial
\theta}\right]$ and setting further obvious simplifying changes of
notation, we can rewrite~\thetag{7.33} more concisely as
\def\theequation{7.37}\begin{equation}
\mathcal{ X}^d-\mathcal{ X}^0= 
\mathcal{ R}- T\left[
d\cdot \mathcal{ H} \mathcal{ X}^d
\right].
\end{equation}
Here, thanks to the inequality $\left\vert \left\vert X_{c,x, v}^d
\right\vert \right\vert_{\mathcal{ C}^{1, \alpha} (\partial
\Delta)}\leq c\cdot K_2$ already established in Lemma~7.23 and thanks
to Lemma~6.12, we know that the vector $\mathcal{ R}\in \mathcal{
C}^{\alpha }( \partial \Delta, \R^n)$ and the matrix $\mathcal{ H}\in
\mathcal{ C}^{1, \alpha}( \partial \Delta, \mathcal{ M}_{n\times
n}(\R))$ are small and more precisely, they satisfy the following two
estimates
\def\theequation{7.38}\begin{equation}
\left\{
\aligned
\vert \vert \mathcal{ R} \vert \vert_{
\mathcal{ C}^\alpha(\partial \Delta)}\leq
& \
c^{2+\alpha} \cdot 
\vert \vert \vert T \vert \vert \vert_{\mathcal{ C}^\alpha
(\partial \Delta)} K_3 (K_2)^2
\\
\vert \vert
\mathcal{ H} \vert \vert_{\mathcal{ C}^\alpha(
\partial \Delta)}\leq 
& \
c^{1+\alpha} \cdot K_3.
\endaligned\right.
\end{equation}
We can rewrite~\thetag{7.37} under the form
\def\theequation{7.39}\begin{equation}
\mathcal{ X}^d -\mathcal{ X}^0 =
\mathcal{ S} -T\left[
d\cdot \mathcal{ H}(\mathcal{ X}^d-
\mathcal{ X}^0)
\right],
\end{equation}
with $\mathcal{ S}:= \mathcal{ R}-T \left[ d\cdot \mathcal{ H}
\mathcal{ X}^0 \right]$. Using the inequality $\left\vert \left\vert
\mathcal{ X}^0 \right\vert \right\vert_{ \mathcal{ C}^\alpha( \partial
\Delta)} \leq c\cdot K_4$ which is a direct consequence
of~\thetag{7.22} and taking the previous estimates~\thetag{7.38}
into account, we deduce the inequality
\def\theequation{7.40}\begin{equation}
\vert \vert \mathcal{ S}
\vert \vert_{\mathcal{ C}^\alpha(\partial \Delta)}\leq
c^{2+\alpha}\cdot
\vert \vert \vert T \vert \vert \vert_{\mathcal{ C}^\alpha
(\partial \Delta)}\left[
K_3(K_2)^2+
K_3K_4
\right].
\end{equation} 
Taking the $\mathcal{ C}^\alpha(\partial \Delta)$ norm of both sides
of~\thetag{7.40}, we deduce the estimate
\def\theequation{7.41}\begin{equation}
\left\{
\aligned
\left\vert \left\vert 
\mathcal{ X}^d - \mathcal{ X}^0 \right\vert 
\right\vert_{
\mathcal{ C}^\alpha(\partial \Delta)} \leq
& \ 
c^{2+\alpha} \cdot \frac{
\vert \vert \vert T \vert \vert \vert_{\mathcal{ C}^\alpha
(\partial \Delta)}\left[
K_3(K_2)^2+
K_3K_4
\right]
}{
1- c^{1+\alpha}
\cdot \vert \vert \vert T \vert \vert \vert_{\mathcal{ C}^\alpha
(\partial \Delta)} K_3} \\
\leq
& \
c^{2+\alpha} \cdot 2\vert 
\vert \vert T \vert \vert \vert_{\mathcal{ C}^\alpha
(\partial \Delta)}\left[
K_3(K_2)^2+
K_3K_4
\right]
\endaligned\right.
\end{equation}
where we use the assumption~\thetag{7.35} on $c_1$
to obtain the second second inequality. It suffices
to set $K_7:= 2\vert 
\vert \vert T \vert \vert \vert_{\mathcal{ C}^\alpha
(\partial \Delta)}\left[
K_3(K_2)^2+
K_3K_4
\right]$, which completes the proof.
\endproof

\subsection*{7.42.~Adjustment of the tangent vector}
Let $v_1\in T_{p_1}M^1$ with $\vert v_1 \vert =1$, as in Lemma~7.12.
Coming back to the first family $Z_{c,x, v}^0 (\zeta)$ defined
by~\thetag{7.8}, we observe that
\def\theequation{7.43}\begin{equation}
\left\{
\aligned
\frac{\partial Z_{j;c,0,v_1}^0}{\partial x_k}
(1)=
& \
\delta_k^j, 
\ \ \ \ \ 
j,k=1,\dots,n, \\
\frac{\partial^2 Z_{j;
c,0,v_1}^0}{\partial v_k 
\partial \theta}(1)=
& \
c\, \frac{\partial \Psi}{\partial \theta}
\left(e^{i\theta}\right) \, 
\delta_k^j, 
\ \ \ \ \ 
j,k=1,\dots,n.
\endaligned\right.
\end{equation}
From now on, we shall set $d=1$ and we shall only consider the family
$Z_{c,x,v}^1(\zeta)$. Thanks to the estimates~\thetag{7.28}
and~\thetag{7.36}, we deduce that if $c_1$ is sufficiently small, then
for all $c$ with $0< c \leq c_1$, the two Jacobian matrices
\def\theequation{7.44}\begin{equation}
\left(
\frac{\partial Z_{j;c,0,v_1}^1}{\partial x_k}(1)
\right)_{1\leq j,k\leq n} \ \ \ \ \
{\rm and} \ \ \ \ \ \
\left(
\frac{\partial^2 Z_{j;c,0,v_1}^1}{\partial v_k\partial \theta}(1)
\right)_{1\leq j,k\leq n}
\end{equation}
are invertible. It would follow that if we would set $A_{x, v:c}^1
(\zeta):= Z_{c, x, v_1 +v}^1(\zeta)$, similarly as in~\thetag{7.10},
then the disc $A_{x, v:c}^1 (\zeta)$ would satisfy the two rank
properties {\bf ($\mathbf{ 5_1}$)} and {\bf ($\mathbf{ 6_1}$)} of
Lemma~7.12. However, the tangency condition {\bf ($\mathbf{ 4_1}$)}
would certainly not be satisfied, because as $d$ varies from $0$ to
$1$, the disc $Z_{c,x,v}^d(\zeta)$ undergo a nontrivial deformation.

Consequently, for every $c$ with $0< c \leq c_1$, 
we have to adjust the ``cone
parameter'' $v$ in order to maintain the tangency 
condition.

\def\thelemma{7.45}\begin{lemma}
For every $c$ with $0< c \leq c_1$, there exists a
vector $v(c)\in \R^n$ such that
\def\theequation{7.46}\begin{equation}
\frac{\partial Z_{c,0,v_1
+v(c)}^1}{\partial \theta} 
(1)=
\frac{\partial Z_{c,0,
v_1}^0}{\partial \theta}(1)=c \cdot
\frac{\partial \Psi}{\partial \theta}(1) \cdot v_1.
\end{equation}
Furthermore, there exists a constant $K_8>0$ such that
$\vert v(c) \vert \leq c^{1+\alpha} \cdot K_8$.
\end{lemma}

\proof
Unfortunately, we cannot apply the implicit function theorem, because
the mapping $Z_{c,x,v}^1$ is identically zero when $c=0$, so we have
to proceed differently. First, we set
\def\theequation{7.47}\begin{equation}
C_1:= \frac{\partial \Psi}{\partial \theta}(1), \ \ \ \ \ 
{\rm and}
\ \ \ \ \
C_2:= \vert\vert
\Psi \vert\vert_{\mathcal{ C}^2\left(\overline{\Delta}\right)}.
\end{equation}
The constant $C_2$ will be used only in Section~8 below.
Choose $K_8\geq \frac{2K_6}{C_1}$. According to the explicit
expression~\thetag{7.8}, the set of points
\def\theequation{7.48}\begin{equation}
\left\{
\frac{\partial X_{c,0,v_1+v}^0}{\partial \theta}(1)\in \R^n: \ 
\vert v \vert \leq c^{1+\alpha} \cdot K_8
\right\}
\end{equation}
covers a cube in $\R^n$ centered at the point $\frac{ \partial X_{
c,0, v_1 }^0}{\partial \theta}(1)$ of radius $c^{2+ \alpha} \cdot C_1
K_8$.
Thanks to the estimate~\thetag{7.28}, we deduce
that the (deformed) set of points 
\def\theequation{7.49}\begin{equation}
\left\{
\frac{\partial X_{c,0,v_1+v}^1}{\partial \theta}(1)\in \R^n: \ 
\vert v \vert \leq c^{1+\alpha} \cdot K_8
\right\}
\end{equation}
covers a cube in $\R^n$ centered at the same
point $\frac{ \partial X_{
c,0, v_1 }^0}{ \partial \theta}(1)$, but of radius 
\def\theequation{7.50}\begin{equation}
c^{2+ \alpha} \cdot C_1
K_8-c^{2+\alpha} \cdot K_6\geq c^{2+\alpha} \cdot K_6.
\end{equation}
Consequently, there exists at least one $v(c) \in \R^n$ with $\vert
v(c) \vert \leq c^{ 1+\alpha} \cdot K_8$ such that~\thetag{7.46}
holds, which completes the proof. 
\endproof

\subsection*{7.51.~Construction of the family $A_{x, v:c}^1 (\zeta)$}
We can now complete the proof of the main Lemma~7.12 of the present
section. First of all, with $\Psi (\zeta)$ as in \S7.5 ans in {\sc
Figure~14}, we consider the composed conformal mapping
\def\theequation{7.52}\begin{equation}
\zeta \longmapsto c\Psi(\zeta) \longmapsto 
\frac{i-c\Psi(\zeta)}{i+c\Psi(\zeta)}=:
\Phi_c(\zeta).
\end{equation}
The image $\Phi_c ( \zeta)$ of the unit disc is a small domain
contained in $\Delta$ and concentrated near $1$. More precisely,
assuming that $c$ satifies $0< c \leq c_1$ with $c_1<<1$ as in the
previous paragraphs, and taking account of the definition of $\Psi
(\zeta)$, it can be checked easily that $\Phi_c(1) =1$, that $\Phi_c
(\partial^+ \Delta)$ is contained in $\{e^{ i\theta} \in
\partial^+\Delta: \ \vert \theta \vert < 10c\}$, and that
\def\theequation{7.53}\begin{equation}
\Phi_c\left(\overline{\Delta}\backslash 
\partial^+\Delta\right)\subset 
\{ \zeta \in \Delta : \ 
\vert \zeta -1 \vert < 8c\} \subset 
\{\rho e^{i\theta} \in \Delta : \
\vert \theta \vert < 10c, \
1-10c < \rho < 1\}.
\end{equation}
the second inclusion being trivial. Here is an illustration:

\bigskip
\begin{center}
\input sub-half-disc.pstex_t
\end{center}

We can now define the final desired family of analytic discs, 
writing the parameter $c$ after a semi-colon, since we lose the
$\mathcal{ C}^{2,\alpha-0}$-smoothness with respect to 
$c$ after the application of Lemma~7.45:
\def\theequation{7.54}\begin{equation}
A_{x,v:c}^1(\zeta):= Z_{c,x,v_1+v(c)+v}^1\left(
\Phi_c(\zeta)\right).
\end{equation}
We restrict the variation of the parameters $x$ to $\vert x \vert \leq
c^2$ and $v$ to $\vert v \vert\leq c$. Property {\bf ($\mathbf{
4_1}$)} holds immediately, thanks to the choice of $v(c)$. Properties
{\bf ($\mathbf{1_1}$)}, {\bf ($\mathbf{ 3_1}$)}, {\bf ($\mathbf{
5_1}$)} and {\bf ($\mathbf{ 6_1}$)} as well as the embedding property
in {\bf ($\mathbf{2_1}$)} are direct consequences of the similar
properties~\thetag{7.44} satisfied by $Z_{ c,x, v_1+ v(c)+v}^1 (\zeta)$,
using the chain rule and the nonvanishing of the partial derivative
$\frac{ \partial \Phi_c }{ \partial \theta}(1)$. The size estimate in
{\bf ($\mathbf{2_1}$)} follows from~\thetag{7.25}, from~\thetag{7.28},
from the restriction of the domains of variation of $x$ and of $v$ and
from~\thetag{7.53}. This completes the proof of Lemma~7.12.
\hfill 
\qed

\section*{\S8.~Geometric properties of families of half-attached 
analytic discs}

\subsection*{8.1.~Preliminary}
By Lemma~7.12, for every $c$ with $0 < c \leq c_1$, the family of
half-attached analytic discs $A_{x, v:c}^1 (\zeta)$ covers a local
wedge of edge $M^1$ at $p_1$. However, not only we want the family
$A_{x, v:c}^1$ to cover a local wedge of edge $M^1$ at $p_1$, but we
certainly want to remove the point $p_1$ by means of the continuity
principle, under the assumptions of the main Proposition~5.12, a final
task which will be achieved in Section~9 below. Consequently, in each
one of the three geometric situations {\bf ($\mathbf{I_1}$)}, {\bf
($\mathbf{I_2}$)} and {\bf (II)} which we have normalized in
Lemma~5.37 above, we shall firstly deduce from the tangency condition
{\bf ($\mathbf{4_1}$)} of Lemma~7.12 that the blunt half-boundary
$A_{0,0:c}^1( \partial^+ \Delta \backslash \{1\})$ is contained in the
open side $(H^1)^+$ (this is why we have normalized in Lemma~5.37 the
second order terms of the supporting hypersurface $H^1$ in order that
$(H^1)^+$ is strictly concave; the reason why we require that
$A_{0,0:c}^1(\partial^+\Delta\backslash \{1\})$ is contained in
$(H^1)^+$ will be clear in Section~9 below). Secondly, we shall show
that for all $x$ with $\vert x \vert \leq c^2$, the disc interior
$A_{x,0:c}(\Delta)$ is contained in the local half-wedge $\mathcal{
HW}_1^+$ in the cases {\bf ($\mathbf{I_1}$)}, {\bf ($\mathbf{I_2}$)}
and is contained in the wedge $\mathcal{ W}_2$ in case {\bf (II)}.

\subsection*{8.2.~Geometric disposition of the
discs with respect to $H^1$ and to $\mathcal{ HW}_1^+$ or to
$\mathcal{ W}_2$} We remember that the positive $c_1$ of Lemmas~7.12,
7.23 and 7.34 was shrunk explicitely, in terms of the constants
$K_1,K_2,K_3,\dots$. In this section, we shall again shrink $c_1$
several times, but without mentioning all the similar explicit
inequalities which will appear. The precise statement of the main
lemma of this section, which is a continuation of Lemma~7.12, is as
follows; whereas we can essentially gather the three cases in the
formal statement of the lemma, it is necessary to treat them
separately in the proof, because the normalizations of
Lemma~5.37 differ.

\def\thelemma{8.3}\begin{lemma}
Let $M$, let $M^1$, let $p_1$, let $H^1$, let $v_1$, let $(H^1)^+$,
let $\mathcal{ HW }_1^+$ or let $\mathcal{ HW }_2$ and let a
coordinate system $z= (z_1, \dots, z_n)$ vanishing at $p_1$ be as in
Case {\bf ($\mathbf{I_1}$)}, as in Case {\bf ($\mathbf{ I_2 }$)} or as
in Case {\bf (II)} of Lemma~5.37. Choose as a local one-dimensional
submanifold $T^1 \subset M^1$ transversal to $H^1$ in $M^1$ and
passing through $p_1$ the submanifold $T_1:=\{\left( x_1,0, \dots,0) +
ih(x_1, 0, \dots, 0) \right)\}$ in Case {\bf
($\mathbf{I_1}$)} and the submanifold $T_1:=\{\left( 0, \dots,0,x_n) +
ih(0, \dots, 0,x_n) \right)\}$ in Cases {\bf
($\mathbf{I_2}$)} and {\bf (II)}. For every $c$ with $0 < c \leq
c_1$, let $A_{x,v:c}^1(\zeta)$ be the family of analytic discs
satisfying properties {\bf ($\mathbf{1_1}$)}, {\bf ($\mathbf{2_1}$)},
{\bf ($\mathbf{3_1}$)}, {\bf ($\mathbf{4_1}$)}, {\bf ($\mathbf{5_1}$)}
and {\bf ($\mathbf{6_1}$)} of Lemma~7.12. Shrinking $c_1$ if
necessary, then for every $c$ with $0< c\leq c_1$, the following three
further properties hold
\begin{itemize}
\item[{\bf ($\mathbf{7_1}$)}] $A_{0, 0:c}^1( \partial^+ \Delta
\backslash\{ 1\}) \subset (H^1)^+$.
\item[{\bf ($\mathbf{8_1}$)}]
$A_{x,0:c}^1(\partial^+\Delta)$ is contained
in $(H^1)^+$ for all $x$ such that the point 
$A_{x,0:c}^1(1)$ belongs to 
$T^1\cap (H^1)^+$.
\item[{\bf ($\mathbf{9_1}$)}]
$A_{x,v:c}^1(\overline{\Delta}\backslash\partial^+\Delta)$ is
contained in the half-wedge
$\mathcal{ HW}_1^+$ or in the wedge $\mathcal{ W}_2$ for all
$x$ and all $v$.
\end{itemize}
\end{lemma}

\proof
For the three new properties {\bf ($\mathbf{ 7_1}$)}, {\bf ($\mathbf{
8_1}$)}, and {\bf ($\mathbf{ 9_1}$)}, we study thoroughly only Case
{\bf ($\mathbf{I_1}$)}, because the other two cases can be treated in
a totally similar way. {\sc Figure~16} just below illustrates
properties {\bf ($\mathbf{ 7_1}$)} and {\bf ($\mathbf{ 8_1}$)} and
also properties {\bf ($\mathbf{ 1_1}$)}, {\bf ($\mathbf{ 5_1}$)} and
{\bf ($\mathbf{ 6_1}$)} of Lemma~7.12.

\bigskip
\begin{center}
\input half-boundaries.pstex_t
\end{center}

Intuitively, the reason why property {\bf ($\mathbf{ 7_1}$)} holds
true is clear: the open set $(H^1)^+$ is strictly concave and the
small segment $A_{0,0:c}^1(\partial^+\Delta)$ is tangent to $H^1$ at
$p_1$; also, the reason why property {\bf ($\mathbf{ 8_1}$)} holds
true is equally clear: when $x$ varies, the small segments
$A_{x,0:c}^1(\partial^+\Delta)$ are essentially translated (inside
$M^1$) from $p_1$ by the vector $x\in\R^n$; and finally, the reason
why property {\bf ($\mathbf{ 8_1}$)} holds true has a simple geometric
interpretation: if the scaling parameter $c_1$ is small enough, the
small analytic disc $A_{x,v:c}^1(\overline{\Delta})$ is essentially a
slightly deformed small part of the straight complex line
$\C\cdot(v_1+Jv_1)$, where the half-wedge $\mathcal{ HW}_1^+$ or the
wedge $\mathcal{ W}_2$ is directed by the vector $Jv_1$ according to
Lemma~5.37. The next paragraphs are devoted to some elementary
estimates which will establish these properties rigorously.

Firstly, let us prove property {\bf ($\mathbf{ 7_1}$)} in Case {\bf
($\mathbf{I_1}$)}. According to Lemma~5.37, the vector $v_1$ is given
by $(0,1,\dots,1)$ and the side $(H^1)^+\subset M^1$ is defined by
$x_1>g(x')=-x_2^2-\cdots-x_n^2+ \widehat{ g}(x')$, where the
$\mathcal{ C}^{2,\alpha}$-smooth function $\widehat{ g}(x')$ vanishes
to second order at the origin, thanks to the normalization
conditions~\thetag{5.40}. By Lemma~6.4, the remainder $\widehat{
g}(x')$ then satisfies an inequality of the form $\vert \widehat{
g}(x') \vert \leq K_9 \cdot \vert x' \vert^{2+\alpha}\leq K_9 \cdot
\left( \sum_{j=2}^n\, x_j^2 \right)^{\frac{ \alpha+2}{2}}$, for some
constant $K_9>0$. Since the strictly concave open subset
$(\widetilde{ H}^1)^+$ of $M^1$ with $\mathcal{ C}^{2,\alpha}$-smooth
boundary defined by the inequality $x_1>-x_1^2-\cdots-x_n^2 + K_9
\cdot \left( \sum_{j=2}^n \, x_j^2 \right)^{\frac{2+\alpha}{2}}$ is
contained in $(H^1)^+$, it suffices to prove property {\bf ($\mathbf{
7_1}$)} with $(H^1)^+$ replaced by $(\widetilde{ H}^1)^+$.

By construction, the disc boundary $A_{0,0:c}(\partial\Delta)$ is
tangent at $p_1$ to $H^1$, hence also to $\widetilde{
H}^1$. Intuitively, it is clear that the blunt disc half-boundary
$A_{0,0:c}(\partial^+\Delta\backslash\{1\})$ should then be contained
in the strictly concave open subset $(\widetilde{ H}^1)^+$, 
{\it see}\, {\sc Figure~16} above. 

To proceed rigorously, we shall come back to the definition~\thetag{
7.53} which yields $A_{0, 0:c}^1( \zeta) \equiv Z_{c,0, v_1+ v(c)}^1
\left( \Phi_c(\zeta)\right)$, with the tangency condition~\thetag{
7.46} satisfied. First of all, denoting the $n$ components of $v(c)$
by $(v_1(c), \dots, v_n(c))$, we may compute the second order
derivatives of the similar discs attached to $M^0$:
\def\theequation{8.4}\begin{equation}
\left\{
\aligned
\frac{\partial^2 Z_{1;c,0,
v_1+v(c)}^0}{\partial \theta^2}(1)= 
& \
c\cdot \frac{\partial^2\Psi}{
\partial \theta^2}\left(e^{i\theta}\right)\cdot
v_1(c), \\
\frac{\partial^2 Z_{j;c,0,
v_1+v(c)}^0}{\partial \theta^2}(1)= 
& \
c\cdot \frac{\partial^2\Psi}{\partial 
\theta^2}\left(e^{i\theta}\right)\cdot
(1+v_j(c)), \ \ \ \ \ \ \ 
j=2,\dots,n.
\endaligned\right.
\end{equation}
Using the definition~\thetag{ 7.47}, the inequality
$\vert v(c) \vert \leq c^{1+\alpha} \cdot K_8$ and the 
second estimate~\thetag{7.36}, we deduce that
\def\theequation{8.5}\begin{equation}
\left\{
\aligned
\left\vert
\frac{\partial^2 Z_{1;c,0,
v_1+v(c)}^1}{\partial \theta^2}(1)\right\vert \leq
& \
c^{2+\alpha}\cdot K_7+c^{2+\alpha}\cdot C_2 K_8
=: c^{2+\alpha} \cdot 2K_{10} 
\\
\left\vert
\frac{\partial^2 Z_{j;c,0,
v_1+v(c)}^1}{\partial \theta^2}(1)\right\vert \leq
& \
c\cdot 2C_2, 
\ \ \ \ \ \ \
j=2,\dots,n.
\endaligned\right.
\end{equation}
Applying then Taylor's integral formula $F (\theta)=F(0) +\theta \cdot
F'(0)+ \int_0^\theta \, (\theta- \theta') \cdot \partial_\theta
\partial_\theta F(\theta') \cdot d \theta'$ to $F(\theta):= X_{1;c,
0,v_1 +v(c)}^1\left(e^{i\theta}\right)$ and afterwards to $F(\theta):=
X_{j;c, 0,v_1 +v(c)}^1\left(e^{i\theta}\right)$ for $j=2,\dots,n$,
taking account of the tangency conditions
\def\theequation{8.6}\begin{equation}
\frac{\partial X_{1;c,0,v_1+v(c)}^1}{\partial \theta}(1)=0, 
\ \ \ \ \ \ \ \ \ \ \ \ 
\frac{\partial X_{j;c,0,v_1+v(c)}^1}{\partial \theta}(1)=
c\cdot C_1, \ \ \ \ \ \ \ 
j=2,\dots,n,
\end{equation}
(a simple rephrasing of~\thetag{7.46}) 
and using the inequalities~\thetag{8.5}, 
we deduce that
\def\theequation{8.7}\begin{equation}
\left\{
\aligned
\left\vert
X_{1;c,0,v_1+v(c)}^1\left(e^{i\theta}\right)
\right\vert \leq
& \
\theta^2\cdot c^{2+\alpha} \cdot K_{10}, 
\\
\left\vert
X_{j;c,0,v_1+v(c)}^1\left(e^{i\theta}\right) -\theta \cdot c \cdot C_1
\right\vert \leq
& \
\theta^2 \cdot
c\cdot C_2, 
\ \ \ \ \ 
j=2,\dots,n.
\endaligned\right.
\end{equation}
Recall that
\def\theequation{8.8}\begin{equation}
x_1>\widetilde{ g}(x'):= -x_2^2-\cdots-x_n^2+K_9\left(
\sum_{j=2}^n\, x_j^2
\right)^{\frac{2+\alpha}{2}}
\end{equation}
denotes the equation of $(\widetilde{ H}^1)^+$. We now claim that if
$c_1$ is sufficiently small, then for every $\theta$ with $0<\vert
\theta \vert < 10c$, we have
\def\theequation{8.9}\begin{equation}
X_{1;c,0,v_1+v(c)}^1\left(e^{i\theta}\right) >
\widetilde{ g}\left(
X_{2;c,0,v_1+v(c)}^1\left(e^{i\theta}\right),\dots\dots,
X_{n; c,0,v_1+v(c)}^1\left(e^{i\theta}\right)
\right).
\end{equation}
Since $\Phi_c( \partial^+\Delta)$ is 
contained in $\{e^{ i\theta}\in \partial^+
\Delta: \ \vert \theta \vert <
10c\}$, this will imply the desired inclusion
for proving {\bf ($\mathbf{ 7_1}$)}:
\def\theequation{8.10}\begin{equation}
\left\{
\aligned
A_{x,v:c}^1(\partial^+\Delta\backslash \{1\})=
& \
Z_{c,0,v_1+v(c)}^1\left(
\Phi_c(\partial^+\Delta\backslash\{1\})
\right)\subset \\
\subset 
& \
Z_{c,0,v_1+v(c)}^1\left(
\{e^{i\theta}\in \partial^+\Delta: \ 0 <\vert \theta \vert \leq 10c\}
\right)
\subset (\widetilde{ H}^1)^+.
\endaligned\right.
\end{equation}
To prove the claim, we notice a minoration 
of the left hand side of~\thetag{8.9}, using~\thetag{8.7}
\def\theequation{8.11}\begin{equation}
X_{1;c,0,v_1+v(c)}^1\left(e^{i\theta}\right) 
\geq -\theta^2\cdot c^{2+\alpha}
\cdot K_{10}.
\end{equation}
On the other hand, using two
inequalities which are direct consequences of
the second line of~\thetag{8.7}, provided that 
$10c_1\cdot C_2\leq \frac{C_1}{2}$:
\def\theequation{8.12}\begin{equation}
\left\{
\aligned
\left\vert
X_{j;c,0,v_1+v(c)}^1\left(e^{i\theta}\right)
\right \vert \leq 
& \
\vert \theta\vert \cdot c \cdot \left(
C_1+ \vert \theta \vert \cdot C_2\right)
\leq \vert \theta \vert 
\cdot c\cdot \frac{ 3C_1}{2}, \\
\left[
X_{j;c,0,v_1+v(c)}^1
\right]^2\geq
& \
\theta^2\cdot c^2\cdot 
(C_1-\vert \theta\vert \cdot C_2)^2 \geq
\theta^2\cdot c^2 \cdot
\frac{ C_1^2}{4}, 
\endaligned\right.
\end{equation}
for $j=2,\dots,n$, we deduce the following
majoration of the 
right hand side of~\thetag{8.9}
\def\theequation{8.13}\begin{equation}
\left\{
\aligned
\widetilde{ g}
& \
\left(
X_{2;c,0,v_1+v(c)}^1\left(e^{i\theta}\right),\dots\dots,
X_{n;c,0,v_1+v(c)}^1\left(e^{i\theta}\right)
\right)= \\
& \
=-\sum_{j=2}^n\, \left[
X_{j;c,0,v_1+v(c)}^1
\right]^2+K_9\left(
\sum_{j=2}^n\, 
\left[
X_{j;c,0,v_1+v(c)}^1\left(e^{i\theta}\right)\right]^2
\right)^{\frac{2+\alpha}{2}} \\
& \
\leq
-\theta^2\cdot c^2 \cdot \frac{C_1^2}{4}
(n-1)+
\vert \theta \vert^{2+\alpha} \cdot c^{2+\alpha} \cdot
\left(\frac{(n-1)9C_1^2}{4}\right)^{\frac{2+\alpha}{2}} K_9 \\
& \
\leq 
-\theta^2\cdot c^2\left(
\frac{ C_1^2}{4}(n-1)-
c^\alpha\cdot \left(\frac{(n-1)9C_1^2}{4}\right)^{
\frac{2+\alpha}{2}} K_9
\right).
\endaligned\right.
\end{equation}
Thanks to the minoration~\thetag{8.11} and to the
majoration~\thetag{8.13}, in order that the inequality~\thetag{8.9}
holds for all $\theta$ with $0< \vert \theta \vert \leq 10c$, it
suffices that the right hand side of~\thetag{8.11} be greater than the
last line of~\thetag{8.13}. By writing this strict inequality
and clearing the factor $\theta^2\cdot c^2$, we see
that it suffices that
\def\theequation{8.14}\begin{equation}
-K_{10}\cdot c^\alpha > -\left( \frac{ C_1^2}{4}
(n-1)-c^\alpha\cdot
\left(
\frac{(n-1)9C_1^2}{4}
\right)^{\frac{2+\alpha}{2}}K_9
\right),
\end{equation}
or equivalently
\def\theequation{8.15}\begin{equation}
c_1 <
\left(
\frac{
\frac{ C_1^2}{4}(n-1)
}{
K_{10}+\left(\frac{ (n-1)9C_1^2}{4}\right)^{
\frac{2+\alpha}{2}} K_9
}
\right)^{\frac{1}{\alpha}}.
\end{equation}
This completes the proof of property 
{\bf ($\mathbf{7_1}$)}.

Secondly, let us prove property {\bf ($\mathbf{ 8_1}$)} in Case {\bf
($\mathbf{I_1}$)}, proceeding similarly. As above, we come back to
the definition $A_{x, 0 :c}^1 (\zeta):= Z_{c, x,v_1 +v(c)}^1 \left(
\Phi_c (\zeta) \right)$ and we remind that $A_{x, 0: c}^1 (1)= Z_{c,
x,v_1+ v(c)}^1(1)= x+ih(x)$, which follows by putting $d=1$ and
$\zeta=1$ in~\thetag{ 7.18}. Thanks to the inclusion $\Phi_c(
\partial^+ \Delta)\subset \{e^{i\theta}\in \partial^+\Delta : \ \vert
\theta \vert < 10c\}$, it suffices to prove that the segment $Z_{c, x,
v_1+v(c)}\left( \left\{ e^{i\theta}: \ \vert \theta \vert < 10c
\right\} \right)$ is contained in the open side $( \widetilde{ H}^1)^+
\subset (H^1 )^+$ defined by the inequation~\thetag{8.8}, if the point
$x+ i h(x)$ belongs to the transversal half-submanifold $T^1 \cap
(H^1)^+$, namely if $x=(x_1, 0, \dots,0)$ with $x_1>0$. In the sequel,
we shall denote the disc $Z_{c,x, v_1+ v(c)}^1(\zeta)$ by $Z_{c, x_1,
x',v_1+ v(c) }^1( \zeta)$, emphasizing the decomposition $x= (x_1,
x')\in \R \times \R^{n-1}$, and we shall also use the convenient
notation
\def\theequation{8.16}\begin{equation}
Z_{c,x_1,x',v_1+v(c)}^{'1}\left(\rho e^{i\theta}\right):=
\left(
Z_{2;c,x_1,x',v_1+v(c)}^1\left(\rho e^{i\theta}\right),\dots\dots, 
Z_{n;c,x_1,x',v_1+v(c)}^1\left(\rho e^{i\theta}\right)
\right).
\end{equation}
So, we have to show that for all $c$ with $0 < c \leq c_1$, for all
$x_1$ with $0 < x_1 \leq c^2$ and for all $\theta$ with $\vert \theta
\vert < 10c$, then the following strict inequality holds true
\def\theequation{8.17}\begin{equation}
X_{1;c,x_1,0,v_1+v(c)}^1\left(e^{i\theta}\right) > 
\widetilde{ g}\left(
X_{c;x_1,0,v_1+v(c)}^{1'}\left(e^{i\theta}\right)
\right),
\end{equation}
First of all, coming back to the family of discs attached to $M^0$, we
see by differentiating~\thetag{ 7.8} twice with respect to $x_1$ that
$\frac{\partial^2 Z_{ c,x_1,0,v_1+v(c)}^0}{\partial x_1^2}(\zeta)
\equiv 0$. Next, by differentiating twice Bishop's equation~\thetag{
7.18} with respect to $x_1$ and by reasoning as in Lemma~7.34, we
deduce the estimate
\def\theequation{8.18}\begin{equation}
\left\vert
\left\vert
\frac{\partial^2 Z_{c,x_1,0,v_1+v(c)}^1}{
\partial x_1^2}
\right\vert
\right\vert_{\mathcal{ C}^\alpha(\partial \Delta)}
\leq c^{2+\alpha} \cdot K_7,
\end{equation}
say, with the same constant $K_7>0$ as in Lemma~7.34, after enlarging
it if necessary. Applying then Taylor's integral formula $F(x_1)=
F(0) +x_1 \cdot \partial_{x_1}F(0)+ \int_0^{ x_1}(x_1- \widetilde{
x}_1) \cdot \partial_{x_1} \partial_{ x_1} F( \widetilde{ x}_1) \cdot
d\widetilde{ x}_1$ to the function $F(x_1 ):= X_{1;c, x_1,0, v_1+
v(c)}^1 \left(e^{ i\theta}\right)$, we deduce the minoration
\def\theequation{8.19}\begin{equation}
X_{1;c,x_1,0,v_1+v(c)}^1\left(e^{i\theta}\right)\geq 
X_{1;c,0,0,v_1+v(c)}^1\left(e^{i\theta}\right)+
x_1 \cdot
\frac{\partial X_{1;c,0,0,v_1+v(c)}^1}{
\partial x_1}\left(e^{i\theta}\right)
- x_1^2 \cdot c^{2+\alpha} \cdot \frac{K_7}{2}.
\end{equation}
On the other hand, by differentiating Bishop's equation~\thetag{7.18}
with respect to $x_1$ at $x=0$, the derivative $\partial_{x_1} x$
yields the vector $(1,0,\dots,0)$ and we obtain
\def\theequation{8.20}\begin{equation}
\small
\left\{
\aligned
\frac{\partial X_{c,0,0,v_1+v(c)}}{\partial x_1}
\left(e^{i\theta} \right)=
& \
-T_1\left[
\sum_{l=1}^n\, 
\frac{\partial h}{\partial x_l}
\left(
X_{c,0,0,v_1+v(c)}^1(\cdot)
\right) \, 
\frac{\partial X_{l;c,0,0,v_1+v(c)}^1}{\partial x_1}(\cdot)
\right](e^{i\theta})+ \\
& \ \
+
(1,0,\dots,0).
\endaligned\right.
\end{equation}
Using then the second inequality~\thetag{ 6.13} and
the estimate~\thetag{ 7.25}, we deduce from~\thetag{ 8.20}
\def\theequation{8.21}\begin{equation}
\left\{
\aligned
\left\vert\left\vert
\frac{\partial X_{1;c,0,0,v_1+v(c)}^1}{\partial x_1}(\cdot)-1
\right\vert\right\vert_{
\mathcal{ C}^\alpha(\partial \Delta)}\leq
& \
c^{2+\alpha} \cdot \vert\vert\vert
T_1 \vert \vert \vert_{\mathcal{ C}^\alpha(\partial \Delta)} 
K_2 K_3, \\
\left\vert\left\vert
\frac{\partial X_{j;c,0,0,v_1+v(c)}^1}{\partial x_1}(\cdot)
\right\vert\right\vert_{
\mathcal{ C}^\alpha(\partial \Delta)}\leq
& \
c^{2+\alpha} \cdot \vert\vert\vert
T_1 \vert \vert \vert_{\mathcal{ C}^\alpha(\partial \Delta)} 
K_2 K_3, \ \ \ \ \ \ \ 
j=2,\dots,n.
\endaligned\right.
\end{equation}
Thanks to the first line of~\thetag{8.21}, we can refine the
minoration~\thetag{8.19} by replacing the first order partial
derivative $\frac{ \partial X_{1;c, 0,0, v_1+v(c)}^1 }{ \partial x_1
}\left (e^{i \theta} \right)$ in the right hand side of~\thetag{8.19}
by the constant $1$, modulo an error term and also, we can use the
trivial minoration $-x_1^2\geq -x_1$, which yields a new, more
interesting minoration of the form
\def\theequation{8.22}\begin{equation}
X_{1;c,x_1,0,v_1+v(c)}^1\left(e^{i\theta}\right)\geq 
X_{1;c,0,0,v_1+v(c)}^1\left(e^{i\theta}\right)+
x_1 - x_1 \cdot c^{2+\alpha} \cdot K_{11}, 
\end{equation}
for some constant $K_{ 11}>0$. On the other hand, using the
inequalities $\vert \partial_{ x_j} \widetilde{ g} (x') \vert \leq
\vert x' \vert+ K_9 \cdot \vert x' \vert^{1+\alpha} \cdot \left(1+
\frac{ \alpha }{2}\right) (n-1)^{ \frac{ \alpha }{2}}$ for $j=
2,\dots,n$, using the estimate~\thetag{ 7.25} and using~\thetag{ 6.2},
we deduce an inequality of the form
\def\theequation{8.23}\begin{equation}
\widetilde{ g}\left(
X_{c,x_1,0,v_1+v(c)}^{'1}\left(e^{i\theta}\right)
\right)\leq 
\widetilde{ g}\left(
X_{c,0,0,v_1+v(c)}^{'1}\left(e^{i\theta}\right)
\right)+ x_1 \cdot c \cdot K_{12}, 
\end{equation}
for some constant $K_{12}>0$. Finally, putting together the two
inequalities~\thetag{ 8.22} and~\thetag{ 8.23}, and using the
following inequality, which is an immediate consequence
of the strict inequality~\thetag{ 8.9}:
\def\theequation{8.24}\begin{equation}
X_{1;c,0,0,v_1+v(c)}\left(e^{i\theta}\right)
\geq \widetilde{ g}\left(
X_{c,0,0,v_1+v(c)}^{'1}\left(
e^{i\theta}
\right)
\right),
\end{equation}
valuable for all $\theta$ with $\vert \theta \vert < 10c$, we deduce
the desired inequality~\thetag{ 8.17} as follows:
\def\theequation{8.25}\begin{equation}
\left\{
\aligned
X_{1;c,x_1,0,v_1+v(c)}^1\left(e^{i\theta} \right) 
\geq 
& \
X_{1;c,0,0,v_1+v(c)}^1\left(e^{i\theta} \right)+
x_1-x_1\cdot c^{2+\alpha} \cdot K_{11} \\
\geq 
& \ \widetilde{ g}\left(
X_{c,0,0,v_1+v(c)}^{'1}\left(
e^{i\theta}
\right)
\right)+ x_1-x_1\cdot c^{2+\alpha} \cdot K_{11} \\
\geq 
& \
\widetilde{ g}\left(
X_{c,x_1,0,v_1+v(c)}^{'1}\left(
e^{i\theta}
\right)
\right)+ x_1-x_1\cdot c \cdot K_{11}-x_1\cdot c \cdot K_{12} \\
> 
& \
\widetilde{ g}\left(
X_{c,x_1,0,v_1+v(c)}^{'1}\left(
e^{i\theta}
\right)
\right),
\endaligned\right.
\end{equation}
for all $x_1$ with $0< x_1 \leq c^2$, all $\theta$ 
with $\vert \theta \vert < 10c$ and all $c$ with 
$0 < c \leq c_1$, 
provided
\def\theequation{8.26}\begin{equation}
c_1 \leq \frac{ 1/2}{K_{11}+K_{12}}.
\end{equation}
This completes the proof of property {\bf ($\mathbf{ 8_1}$)}.

Thirdly, let us prove property {\bf ($\mathbf{ 9_1}$)} in Case {\bf
($\mathbf{I_1}$)}. The half-wedge $\mathcal{ HW}_1^+$ is defined by
the $n$ inequalities of the last two lines of~\thetag{5.38}, where
$a_2+\cdots+a_n=1$. For notational reasons, it will be convenient to
set $a_1:=1$ and to write the first inequality defining $\mathcal{
HW}_1^+$ simply as $\sum_{j=1}^n\, a_j y_j> \psi(x,y')$.

Because $\Phi_c \left( \overline{ \Delta} \backslash \partial^+ \Delta
\right)$ is contained in the open sector $\{\rho e^{i \theta} \in
\overline{ \Delta}: \ \vert \theta \vert < 10c, \ 1-10c < \rho < 1$,
taking account of the definition~\thetag{7.53} of $A_{x,v
:c}^1(\zeta)$, in order to check property {\bf ($\mathbf{ 9_1}$)}, it
clearly suffices to show that $Z_{c,x, v_1+ v(c)+ v}^1 \left( \left\{
\rho e^{i\theta}\in \Delta : \ 1-10c < \rho < 1, \ \vert \theta \vert
< 10c \right\} \right)$ is contained in $\mathcal{ HW}_1^+$, which
amounts to establish that for all $x$ with $\vert x \vert \leq c^2$,
all $v$ with $\vert v \vert \leq c$, all $\rho e^{i \theta}$ with
$1-10c < \rho < 1$ and with $\vert \theta \vert < 10c$, the following
two collections of strict inequalities hold true
\def\theequation{8.27}\begin{equation}
\left\{
\aligned
\sum_{k=1}^n\, a_k 
Y_{j;c,x,v_1+v(c)+v}^1\left(\rho e^{i\theta}\right) 
> 
& \
\psi\left(
X_{c,x,v_1+v(c)+v}^1\left(\rho e^{i\theta}\right), \
Y_{c,x,v_1+v(c)+v}^{'1}\left(\rho e^{i\theta}\right)
\right), \\
Y_{j; c,x,v_1+v(c)+v}^1\left(
\rho e^{i\theta}
\right) > 
& \
\varphi_j\left(
X_{c,x,v_1+v(c)+v}^1\left(
\rho e^{i\theta}
\right), \ 
Y_{1;c,x,v_1+v(c)+v}^1\left(
\rho e^{i\theta}
\right)
\right),
\endaligned\right.
\end{equation}
for $j=2,\dots,n$, 
provided $c_1$ is sufficiently small, where we use the
notation~\thetag{ 8.16}.

We first treat the collection of $(n-1)$ strict inequalities
in the second line of~\thetag{ 8.27}.
First of all, by differentiating~\thetag{ 7.8}
twice with respect to $\theta$, we obtain
\def\theequation{8.28}\begin{equation}
\frac{\partial^2 Z_{c,x,v_1+v(c)+v}^0}{\partial \theta^2}\left(
e^{i\theta}
\right)= c\cdot \frac{
\partial^2\Psi}{\partial \theta^2}
\left(
e^{i\theta}
\right) \cdot
\left[
v_1+v(c)+v
\right].
\end{equation}
Using the second estimate~\thetag{7.36}, we deduce that
there exists a constant $K_{13}>0$ such that
\def\theequation{8.29}\begin{equation}
\left\vert
\frac{\partial^2 Z_{c,x,v_1+v(c)+v}^1}{
\partial \theta^2}\left(e^{i\theta}\right)
\right\vert\leq c\cdot K_{13}.
\end{equation}
Using the inequality~\thetag{ 6.2}, using~\thetag{ 8.29}, and then
taking account of the inequalities $\vert \theta \vert < 10c$, $\vert
x \vert \leq c^2$ and $\vert v \vert < c$, we deduce the following
inequality
\def\theequation{8.30}\begin{equation}
\left\{
\aligned
\left\vert
\frac{\partial Z_{c,x,v_1+v(c)+v}^1}{\partial \theta}
\left( e^{i\theta}
\right)-
\frac{\partial Z_{c,0,v_1+v(c)}^1}{\partial \theta} (1)
\right\vert
\leq 
& \
c \cdot \left(
\vert \theta \vert + \vert x \vert + \vert v \vert
\right) \\
\leq 
& \
c^2\cdot K_{14},
\endaligned\right.
\end{equation}
for some constant $K_{14} >0$. On the other hand, by
differentiating~\thetag{ 7.8} with respect to $\theta$ at $\theta=0$
and applying the inequality~\thetag{ 7.28}, we obtain
\def\theequation{8.31}\begin{equation}
\left\vert
\frac{\partial Z_{c,0,v_1+v(c)}^1}{\partial \theta} (1)
-c\cdot C_1\cdot (0,1,\dots,1)
\right\vert\leq c^{2+\alpha}\cdot K_6, 
\end{equation}
where $C_1=\frac{ \partial \Psi }{ \partial \theta}(1)$, as defined
in~\thetag{7.47}. We remind that for every $\mathcal{ C}^1$-smooth
function $Z$ on $\overline{ \Delta}$ which is holomorphic in $\Delta$,
we have $i \frac{ \partial }{ \partial \theta} Z ( e^{i\theta})=
-\frac{ \partial }{\partial \rho} Z( e^{i\theta})$. Consequently, we
deduce from~\thetag{ 8.30} the following first (among three)
interesting inequality
\def\theequation{8.32}\begin{equation}
\left\vert
-\frac{\partial Z_{c,x,v_1+v(c)+v}^1}{\partial \rho}
\left(
e^{i\theta}
\right)- c\cdot C_1 \cdot (0,i,\dots,i)
\right\vert\leq 
c^2 \cdot K_{15}, 
\end{equation}
for some constant $K_{15}>0$.

According to the definition~\thetag{ 7.8}, 
we may compute 
\def\theequation{8.33}\begin{equation}
\frac{ \partial^2 Z_{c,x,v_1+v(c)+v}^0}{\partial \rho^2}
\left(
\rho e^{i\theta}\right)= c\cdot 
\frac{\partial^2 \Psi}{\partial \rho^2}
\left(
\rho e^{i\theta} 
\right)\cdot
(v_1+v(c)+v)
\end{equation}
By reasoning as in the proof of Lemma~7.34, we may obtain an
inequality similar to~\thetag{ 7.36}, with the second order partial
derivative $\partial^2/ \partial \theta^2$ replaced by the second
order partial derivative $\partial^2 / \partial \rho^2$. Putting this
together with~\thetag{ 8.33}, we deduce that there exists a constant
$K_{16} >0$ such that
\def\theequation{8.34}\begin{equation}
\left\vert
\frac{\partial^2 Z_{c,x,v_1+
v(c)+v}^1}{\partial \rho^2}
\left(
\rho e^{i\theta}
\right)
\right\vert\leq 
c \cdot 2K_{16},
\end{equation}
for some constant $K_{16}>0$. Applying then Taylor's integral formula
$F (\rho) = F(1) + (\rho -1) \cdot \partial_\rho F(1)+ \int_1^\rho (
\rho- \widetilde{ \rho})\cdot \partial_\rho \partial_\rho F(
\widetilde{ \rho}) \cdot d\widetilde{ \rho}$ to the functions $F(
\rho):= Y_{k;c, x,v_1 +v(c) +v}^1\left( \rho e^{ i\theta} \right)$ for
$k=1,\dots,n$, we deduce the second interesting collection of
inequalities
\def\theequation{8.35}\begin{equation}
\left\{
\aligned
{}
&
\left\vert
Y_{k;c,x,v_1+v(c)+v}^1
\left(
\rho e^{i\theta}
\right)
-Y_{k;c,x,v_1+v(c)+v}^1\left(
e^{i\theta} 
\right) - \right. \\
& \ \ \ \ \ \ \ \ \ \ \ \ \ \ \ \ \ \ \ \ \ \ \ \ 
\left.
-(\rho-1)\cdot
\frac{\partial Y_{k;c,x,v_1+v(c)+v}^1}{\partial \rho }
\left(
e^{i\theta}
\right)
\right\vert 
\leq (1-\rho)^2 \cdot c \cdot K_{16},
\endaligned\right.
\end{equation}
for $k=1,\dots,n$.

On the other hand, thanks to the normalizations of the functions
$\varphi_j (x, y_1)$ given in~\thetag{ 5.40}, namely $\varphi_j (0)=
\partial_{ x_k} \varphi_j (0)= \partial_{ y_1} \varphi_j (0)=0$, $j=2,
\dots, n$, $k=1, \dots, n$, we see that, possibly after increasing the
constant $K_1 >0$ of Lemma~6.4, we have inequalities of the form
\def\theequation{8.36}\begin{equation}
\left\{
\aligned
{}
&
\sum_{k=1}^n\, 
\left\vert
\varphi_{j,x_k}(x,y_1)
\right\vert+
\left\vert
\varphi_{j,y_1}(x,y_1)
\right\vert\leq (\vert x \vert 
+ \vert y_1 \vert) \cdot
K_1, \\
&
\left\vert \varphi_j(x,y_1) - 
\varphi_j\left(\widetilde{ x}, 
\widetilde{ y}_1\right)
\right\vert \leq 
\left(
\left\vert 
x-\widetilde{ x}
\right\vert+
\left\vert
y_1- \widetilde{ y}_1 
\right\vert
\right)\cdot
\left(
\sum_{k=1}^n\, 
\sup_{\vert x \vert, \, 
\vert y_1 \vert \leq c\cdot K_2}
\left\vert
\varphi_{j,x_k}(x,y_1)
\right\vert+ \right. \\
& \ \ \ \ \ \ \ \ \ \ \ \ \ \ \
\left. 
+\sup_{\vert x \vert, \, 
\vert y_1 \vert \leq c\cdot K_2}
\left\vert
\varphi_{j,y_1}(x,y_1)
\right\vert
\right),
\endaligned\right.
\end{equation}
for $j=2, \dots, n$, provided $\vert x\vert, \, \vert \widetilde{ x}
\vert, \, \vert y_1 \vert, \, \vert \widetilde{ y_1} \vert \leq c\cdot
K_2$. On the other hand, computing $\frac{ \partial Z_{c,x, v_1+
v(c)+v}^0}{\partial \rho} \left( \rho e^{ i\theta} \right)$
in~\thetag{ 7.8}, using~\thetag{ 7.25}, \thetag{ 7.28} and an
inequality of the form~\thetag{ 6.2}, we deduce that there exists a
constant $K_{17}>0$ such that
\def\theequation{8.37}\begin{equation}
\left\vert
Z_{c,x,v_1+v(c)+v}^1\left(
\rho e^{i\theta} 
\right)-
Z_{c,x,v_1+v(c)+v}^1\left(
e^{i\theta} 
\right) 
\right\vert\leq
(1-\rho) \cdot c \cdot K_{17}.
\end{equation}
Finally, using the inequality $\vert Z_{ c,x, v_1+ v(c)+v}^1\left(
\rho e^{i \theta} \right) \vert \leq c \cdot K_2$ obtained in~\thetag{
7.25}, using the collection of inequalities~\thetag{ 8.36} and using
the inequality~\thetag{ 8.37}, we may deduce the third (and last)
interesting inequality for $j=2, \dots, n$:
\def\theequation{8.38}\begin{equation}
\left\{
\aligned
{}
&
\left\vert
\varphi_j\left(
X_{c,x, v_1+v(c)+v}^1\left(
\rho e^{i\theta}
\right), \
Y_{1;c,x, v_1+v(c)+v}^1
\left(\rho e^{i\theta} \right)
\right)- \right. \\
&
\ \ \ \ \ \ \ \ \ \ \ \ \ 
\left.
-\varphi_j\left(
X_{c,x, v_1+v(c)+v}^1\left(
e^{i\theta}
\right), \
Y_{1;c,x, v_1+v(c)+v}^1
\left(e^{i\theta} \right)
\right)
\right\vert\leq \\
&
\left(
\left\vert
X_{c,x,v_1+v(c)+v}^1\left(
\rho e^{i\theta}
\right)- X_{c,x,v_1+v(c)+v}^1\left(
e^{i\theta}\right)
\right\vert+ \right. \\
& \ \ \ \ \ \ \ \ \ \ \ \ \ \ \ \ \ \
\left.
+\left\vert
Y_{1;c,x, v_1+v(c)+v}^1\left(
\rho e^{i\theta}
\right)- Y_{1;c,x, v_1+v(c)+v}^1\left(
e^{i\theta}\right)
\right\vert
\right)\cdot \\
& \ \ \ \ \ \ \ \ \ \ \ \ \ \ \ \
\cdot 
\left(
\sum_{k=1}^n\, 
\sup_{\vert x \vert, \, 
\vert y_1 \vert \leq c\cdot K_2}
\left\vert
\varphi_{j,x_k}(x,y_1)
\right\vert
+\sup_{\vert x \vert, \, 
\vert y_1 \vert \leq c\cdot K_2}
\left\vert
\varphi_{j,y_1}(x,y_1)
\right\vert
\right) \leq \\
& 
\leq (1-\rho) \cdot c^2 \cdot K_{18}, 
\endaligned\right.
\end{equation}
for some constant $K_{18}>0$.

We can now complete the proof of the collection of inequalities in the
second line of~\thetag{ 8.27}. As before, 
let $c$ with $0 < c \leq c_1$, let $\rho$ with $10c < \rho
< 1$, let $\theta$ with $\vert \theta \vert < 10 c$, 
let $x$ with $\vert x \vert \leq c^2$, let $v$ with
$\vert v \vert \leq c$ and let $j=2, \dots, n$. Starting
with~\thetag{ 8.35}, using~\thetag{ 8.32}, using the fact that
$Z_{c,x,v_1+ v(c)+v}^1 \left( \partial^+\Delta \right) \subset M^1
\subset M$ and using~\thetag{ 8.38}, we have
\def\theequation{8.39}\begin{equation}
\left\{
\aligned
{}
&
Y_{j;c,x,v_1+v(c)+v}^1\left(
\rho e^{i\theta} 
\right) \geq \\
& \ \ \ \
\geq Y_{j;c,x,v_1+v(c)+v}^1
\left(
e^{i\theta} 
\right)+ \\
& \ \ \ \ \ \ \ 
+(\rho-1) \cdot 
\frac{\partial Y_{j;c,x,
v_1+v(c)+v}^1}{\partial
\rho}\left(
e^{i\theta}
\right)- (1-\rho)^2\cdot c \cdot K_{16} \geq \\
& \ \ \ \ \
\geq 
Y_{j;c,x,v_1+v(c)+v}^1
\left(
e^{i\theta} 
\right)+
(1-\rho)\cdot
c\cdot C_1
-(1-\rho) \cdot c^2 \cdot K_{15} - 
(1-\rho)^2\cdot c \cdot K_{16} \\
& \ \ \ \ \
= 
\varphi_j\left(
X_{1;c,x,v_1+v(c)+v}^1\left(
e^{i\theta} 
\right), \
Y_{1;c,x,v_1+v(c)+v}^1\left(
e^{i\theta}
\right)
\right)+
(1-\rho)\cdot
c\cdot C_1 - \\
& \ \ \ \ \ \ \ \ \ \ \ \ \ \ \ \ 
\ \ \ \ \ \ \ \ \ \ \ \ \ \ \ \
\ \ \ \ \ \ \ \ \ \ \ \ \ \ \ \
\ \ \ \ \ \ \ \ \ \ \ \ \ \ \ \
-(1-\rho)\cdot c^2 
\cdot K_{15} - (1-\rho)^2\cdot c \cdot K_{16} \\
& \ \ \ \ \ 
\geq 
\varphi_j\left(
X_{1;c,x,v_1+v(c)+v}^1\left(
\rho e^{i\theta} 
\right), \
Y_{1;c,x,v_1+v(c)+v}^1\left(
\rho e^{i\theta}
\right)
\right)+
(1-\rho)\cdot
c\cdot C_1 - \\
& \ \ \ \ \ \ \ \ \ \ \ \ \ \ \ \ 
\ \ \ \ \ \ \ \ \ \ \ \ \ \ \ \
\ \ \ \ \
-(1-\rho) \cdot c^2 \cdot K_{15} - (1-\rho)^2\cdot c \cdot K_{16}-
(1-\rho)\cdot c^2\cdot K_{18} \\
& \ \ \ \ \ 
\geq 
\varphi_j\left(
X_{1;c,x,v_1+v(c)+v}^1\left(
\rho e^{i\theta} 
\right), \
Y_{1;c,x,v_1+v(c)+v}^1\left(
\rho e^{i\theta}
\right)
\right)+
(1-\rho)\cdot
c\cdot \left[C_1 - \right.\\
& \ \ \ \ \ \ \ \ \ \ \ \ \ \ \ \ 
\ \ \ \ \ \ \ \ \ \ \ \ \ \ \ \
\ \ \ \ \ \ \ \ \ \ \ \ \ \ \ \
\ \ \ \ \ \ \ \ \ \
\left.
-c \cdot K_{15} - 10c \cdot K_{16}-
c\cdot K_{18}\right]\\
& \ \ \ \ \
\geq 
\varphi_j\left(
X_{1;c,x,v_1+v(c)+v}^1\left(
\rho e^{i\theta} 
\right), \
Y_{1;c,x,v_1+v(c)+v}^1\left(
\rho e^{i\theta}
\right)
\right)+
(1-\rho)\cdot
c\cdot \frac{ C_1}{2} \\
& \ \ \ \ \ 
> 
\varphi_j\left(
X_{1;c,x,v_1+v(c)+v}^1\left(
\rho e^{i\theta} 
\right), \
Y_{1;c,x,v_1+v(c)+v}^1\left(
\rho e^{i\theta}
\right)
\right),
\endaligned\right.
\end{equation}
provided that
\def\theequation{8.40}\begin{equation}
c_1 \leq 
\frac{ C_1/2}{K_{15}+10K_{16}+K_{18}}.
\end{equation}
This yields the collection of
inequalities in the second line of~\thetag{ 8.27}.

For the first inequality~\thetag{ 8.27}, we proceed similarly. Recall
that $v_1= (0,1, \dots,1)$, that $a_1=1$ and that $a_2+ \cdots+a_n=1$.
Since $Z_{c, x, v_1+ v(c)+v }^1( \partial^+\Delta)\subset M^1 \subset
N^1$, we have for all $\theta$ with $\vert \theta \vert \leq
\frac{\pi}{2}$ the following relation
\def\theequation{8.41}\begin{equation}
\sum_{k=1}^n\, 
a_k\, 
Y_{k;c,x,v_1+v(c)+v}^1\left(
e^{i\theta}
\right)= \psi\left(
X_{c,x,v_1+v(c)+v}^1\left(
e^{i\theta}
\right), \ 
Y_{c,x,v_1+v(c)+v}^{'1} \left(
e^{i\theta}\right)
\right).
\end{equation}
Using that $\psi$ vanishes to order one at the origin by the
normalization conditions~\thetag{ 5.40} and proceeding as in the
previous paragraph concerning the functions $\varphi_j$, we obtain an
inequality similar to~\thetag{ 8.38}:
\def\theequation{8.42}\begin{equation}
\left\{
\aligned
{}
&
\left\vert
\psi\left(X_{c,x,v_1+v(c)+v}^1\left(
\rho e^{i\theta} 
\right), \
Y_{c,x,v_1+v(c)+v}^{'1}\left(
\rho e^{i\theta}
\right)\right)- \right . \\
& \ \ \ \ \ \ \ \ \ \ \ \ \ \ \ \ \ \
\left.
-\psi\left(X_{c,x,v_1+v(c)+v}^1\left(
e^{i\theta} 
\right), \
Y_{c,x,v_1+v(c)+v}^{'1}\left(
e^{i\theta}
\right)\right)
\right\vert\leq 
(1-\rho)\cdot c^2 \cdot K_{19}, 
\endaligned\right.
\end{equation}
for some constant $K_{19}>0$.

As before, let $c$ with $0 < c \leq c_1$, let $\rho$ with $10c < \rho
< 1$, let $\theta$ with $\vert \theta \vert < 10 c$, let $x$ with
$\vert x \vert \leq c^2$ and let $v$ with $\vert v \vert \leq c$. Using
then~\thetag{ 8.35}, \thetag{ 8.32}, \thetag{ 8.41} and~\thetag{
8.42}, we deduce the desired strict inequality
\def\theequation{8.43}\begin{equation}
\small
\left\{
\aligned
{}
&
\sum_{k=1}^n\, 
a_k \, 
Y_{k; c,x,v_1+v(c)+v}^1\left(
\rho e^{i\theta}
\right)\geq 
\sum_{k=1}^n\, 
a_k \, 
Y_{k; c,x,v_1+v(c)+v}^1\left(
e^{i\theta}
\right) +\\
& \ \ \ \ \ \ \ \ 
+(1-\rho) \left[
\sum_{k=1}^n\, 
a_k\left(
-\frac{\partial Y_{k;c,x,v_1+v(c)+v}^1}{\partial
\rho}\left(
e^{i\theta}
\right)
\right)
\right]-
(1-\rho)^2 \cdot c \cdot \left( \sum_{k=1}^n\, a_k\right)K_{16} \\
& \ \ \ \ \ \ \ \ \ 
\geq 
\sum_{k=1}^n\, 
a_k \, 
Y_{k; c,x,v_1+v(c)+v}^1\left(
e^{i\theta}
\right) +
(1-\rho) \left[
\sum_{j=2}^n\,
a_j\cdot 
c\cdot C_1 - \sum_{k=1}^n\, a_k \cdot c^2 \cdot
K_{15} 
\right] - \\
& \ \ \ \ \ \ \ \ \ \ \ \ \ \ \ \ 
\ \ \ \ \ \ \ \ \ \ \ \ \ \ \ \
-(1-\rho)^2\cdot c \cdot 2K_{16} \\
& \ \ \ \ \ \ \ \ \ 
\geq 
\sum_{k=1}^n\, 
a_k \, 
Y_{k; c,x,v_1+v(c)+v}^1\left(
e^{i\theta}
\right) +
(1-\rho)\cdot c \cdot C_1- \\
& \ \ \ \ \ \ \ \ \ \ \ \ \ \ \ \ \
\ \ \ \ \ \ \ \ \ \ \ \ \ \ \ \ 
\ \ \ \ \ \ \ \ \ \ \ \ \ \ \ \
-(1-\rho)\cdot c^2
\cdot 2K_{15}- (1-\rho)^2 \cdot c \cdot 2K_{16} \\
& \ \ \ \ \ \ \ \ \ 
=
\psi\left(
X_{c,x,v_1+v(c)+v}^1\left(
e^{i\theta}
\right), \ 
Y_{c,x,v_1+v(c)+v}^{'1} \left(
e^{i\theta}\right)
\right) + \\
& \ \ \ \ \ \ \ \ \ \ \ \ \ \ \ \ \
\ \ \ \ \ \ \ \ \ \ \ \ \ \ \ \ 
\ \ \ \ \ \ \ \ \ \ \ \ \ \ \ \
+ (1-\rho) \cdot c \cdot \left[
C_1-c \cdot 2K_{15}- 10c \cdot 2K_{16}
\right] \\
& \ \ \ \ \ \ \ \ \ 
\geq 
\psi\left(
X_{c,x,v_1+v(c)+v}^1\left(
\rho e^{i\theta}
\right), \ 
Y_{c,x,v_1+v(c)+v}^{'1} \left(
\rho e^{i\theta}\right)
\right) + \\
& \ \ \ \ \ \ \ \ \ \ \ \ \ \ \ \ \
\ \ \ \ \ \ \ \ \ \ \ \ \ \ \ \ 
\ \ \ \ \ \ \ \ \ \ \ \ \ \ \ \
+ (1-\rho) \cdot c \cdot \left[
C_1-c \cdot 2K_{15}- 10c \cdot 2K_{16}-c\cdot K_{19}
\right] \\
& \ \ \ \ \ \ \ \ \ 
\geq 
\psi\left(
X_{c,x,v_1+v(c)+v}^1\left(
\rho e^{i\theta}
\right), \ 
Y_{c,x,v_1+v(c)+v}^{'1} \left(
\rho e^{i\theta}\right)
\right) +
(1-\rho) \cdot c \cdot \frac{ C_1}{2} \\
& \ \ \ \ \ \ \ \ \ 
>
\psi\left(
X_{c,x,v_1+v(c)+v}^1\left(
\rho e^{i\theta}
\right), \ 
Y_{c,x,v_1+v(c)+v}^{'1} \left(
\rho e^{i\theta}\right)
\right), 
\endaligned\right.
\end{equation}
provided
\def\theequation{8.44}\begin{equation}
c_1\leq \frac{ C_1/2}{2K_{15}+20K_{16}+K_{19}}
\end{equation}
This yields the first inequality of~\thetag{ 8.27} and
completes the proof of {\bf ($\mathbf{ 9_1}$)} in 
Case {\bf ($\mathbf{ I_1}$)}.

The proof of Lemma~8.3 is complete, because by following the
normalizations of Lemma~5.37 and by formulating analogous inequalities,
cases {\bf ($\mathbf{ I_2}$)} and {\bf (II)} are achieved in a totally
similar way.
\endproof

\section*{\S9.~End of proof of Theorem~1.2': application
of the continuity principle}

\subsection*{9.1.~Preliminary}
In this section, we shall now complete the proof of Proposition~5.12
(at last!), hence the proof of Theorem~3.19 {\bf (i)} (which is a
direct consequence of Proposition~5.12, as was explained in Section~5)
and hence also the proof of Theorem~1.2', modulo supplementary
arguments postponed to \S9.27 below. By means of a deformation
$A_{x,v,u:c}^1(\zeta)$ (we add a real parameter
$u$) of the family of analytic discs
$A_{x,v:c}^1(\zeta)$ satisfying properties {\bf ($\mathbf{ 1_1}$)} to
{\bf ($\mathbf{ 9_1}$)} of Lemmas~7.12 and~8.3, and by means of the
continuity principe, we shall show that, in Cases {\bf ($\mathbf{
I_1}$)} and {\bf ($\mathbf{ I_2}$)}, there exists a local wedge
$\mathcal{ W}_{p_1}$ of edge $M$ at $p_1$ to which all holomorphic
functions in $\mathcal{ O}\left( \Omega \cup \mathcal{ HW}_1^+\right)$
extend holomorphically and we shall show that in Case {\bf (II)},
there exists a neighborhood $\omega_{p_1}$ of $p_1$ in $\C^n$ to which
all holomorphic functions in $\mathcal{ O}\left( \Omega \cup \mathcal{
W}_2\right)$ extend holomorphically. To organize this last main
step of the proof of Proposition~5.12, we shall consider jointly
Cases~{\bf ($\mathbf{ I_1}$)}, {\bf ($\mathbf{ I_2}$)} and then
afterwards Case~{\bf (II)} separately in \S9.22 below.

\subsection*{9.2.~Isotopies of analytic
discs and continuity principle} To begin with, we shall formulate a
convenient version of the continuity principle. If $E\subset \C^n$ is
an arbitrary subset, we denote by
\def\theequation{9.3}\begin{equation}
\mathcal{V}_{\C^n}(E,\rho)=\bigcup_{p\in E}\, 
\{z\in \C^n: \vert z - p \vert < \rho\}
\end{equation}
the union of polydiscs of radius $\rho>0$ centered at points of $E$.
We then have the following lemma, extracted from~\cite{m2}, which
applies to families of analytic discs $A_\tau(\zeta)$ which are
embeddings of $\overline{\Delta}$ into $\C^n$:

\def\thelemma{9.4}\begin{lemma}
\text{\rm (\cite{m2}, Proposition~3.3)}
Let $\mathcal{ D}$ be a nonempty domain in $\C^n$ and let $A_\tau:
\overline{\Delta} \to \C^n$ be a one-parameter family of analytic
discs, where $\tau\in\R$ satisfies $0\leq \tau \leq 1$. Assume that
there exist constants $c_\tau$ and $C_\tau$ with $0 < c_\tau < C_\tau$
such that 
\def\theequation{9.5}\begin{equation} c_\tau \vert
\zeta_1 - \zeta_2 \vert < \vert A_\tau(\zeta_1) 
- A_\tau (\zeta_2) \vert <
C_\tau \vert \zeta_1 - \zeta_2 \vert,
\end{equation}
for all distinct points $\zeta_1,\zeta_2 \in \overline{ \Delta}$ and
all $0\leq \tau \leq 1$. Assume that $A_1(\overline{ \Delta} )
\subset \mathcal{ D}$, set $\rho_\tau := \inf\{ \vert t - A_\tau
(\zeta)\vert: \, t\in \partial \mathcal{ D},\, \zeta \in \partial
\Delta\}$, namely $\rho_\tau$ is the polydisc distance between $A_\tau
(\partial \Delta)$ and $\partial \mathcal{ D}$, assume
$\rho_\tau>0$ for all $\tau$ with $0\leq \tau \leq 1$,
and set $\sigma_\tau
:= \rho_\tau c_\tau / 2C_\tau$. Then for every holomorphic function
$f \in \mathcal{O} (\mathcal{ D})$, and for every $\tau\in [0, 1]$,
there exists a holomorphic function $F_\tau \in \mathcal{ O}\left(
\mathcal{ V}_{\C^n} (A_\tau\left( \overline{ \Delta}\right),
\sigma_\tau)\right)$ such that $F_\tau =f$ in $\mathcal{ V}_{\C^n}
(A_\tau ( \partial \Delta), \sigma_\tau)\subset \mathcal{ D}$.
\end{lemma}

Two analytic discs $A', A'': \overline{\Delta} \to \C^n$ which are of
class $\mathcal{ C}^1$ over $\overline{\Delta}$ and holomorphic in
$\Delta$ and which are both {\it embeddings}\, of $\overline{\Delta}$
into $\C^n$ are said to be {\sl analytically isotopic} if there exists
a $\mathcal{ C}^1$-smooth family of analytic discs $A_\tau: \overline{
\Delta}\to \C^n$ which are of class $\mathcal{ C}^1$ over
$\overline{\Delta}$ and holomorphic in $\Delta$ such that $A_0=A'$,
such that $A_1=A''$ and such that $A_\tau$ is an {\it embedding}\, of
$\overline{ \Delta}$ into $\C^n$ for all $\tau$ with $0\leq
\tau \leq 1$.

\subsection*{9.6.~Translations of $M^1$ in $M$}
According to Lemma~5.37, in Case~{\bf ($\mathbf{ I_1}$)}, the
one-codimensional submanifold $M^1\subset M$ is given by the equations
$y'=\varphi'(x,y_1)$ and $x_1=g(x')$. If $u\in \R$ is a small real
parameter, we may define a ``translation'' $M_u^1$ of $M^1$ in $M$ by
the $n$ equations
\def\theequation{9.7}\begin{equation}
M_u^1: \ \ \ \ \ 
y'=\varphi'(x,y_1), \ \ \ \ \ 
x_1=g(x')+u.
\end{equation}
Clearly, we have $M_0^1 \equiv M^1$, we have $M_u^1 \subset (M^1 )^+$
if $u>0$ and we have $M_u^1 \subset (M^1)^-$ if $u <0$. We may perturb
the family of analytic discs $Z_{c, x,v}^d ( \zeta)$ attached to $M^1$
and satisfying Bishop's equation~\thetag{ 7.18} by requiring that it
is attached to $M_u^1$. We then obtain a new family of analytic discs
$Z_{c,x, v,u}^d (\zeta)$ which is half-attached to $M_u^1$ and which
is of class $\mathcal{ C}^{2,\alpha-0}$ with respect to all variables
$(c,x, v, u, \zeta)$, thanks to the stability under perturbation of
the solutions to Bishop's. For $u=0$, this solution coincides with the
family $Z_{c, x, v}^d ( \zeta)$ constructed in \S7.13. Using a similar
definition as in~\thetag{7.53}, namely setting $A_{x,v, u:c}^1
(\zeta):= Z_{c,x, v_1+v(c) +v,u}^1 \left( \Phi_c (\zeta) \right)$, we
obtain a new family of analytic discs which coincides, for $u=0$, with
the family of analytic discs $A_{x, v: c}^1 (\zeta)$ of Lemmas~7.12
and~8.3. Similarly, in Case~{\bf ($\mathbf{ I_2}$)}, taking account
of the normalizations stated in Lemma~5.37, we can also construct an
analogous family of analytic discs $A_{x,v,u:c}^1(\zeta)$. From now on,
we shall fix the scaling parameter $c$ with $0 < c \leq c_1$, so that
the nine properties {\bf ($\mathbf{ 1_1}$)} to {\bf ($\mathbf{ 9_1}$)}
of Lemmas~7.12 and~8.3 are satisfied.

\subsection*{9.8.~Definition of a local wedge of edge $M$ at $p_1$
in Cases~($\mathbf{ I_1}$) and~($\mathbf{I_2}$)} First of all, in
Cases~{\bf ($\mathbf{I_1}$)} and~{\bf ($\mathbf{ I_2}$)}, we shall
restrict the variation of the parameter $v$ to a certain
$(n-2)$-dimensional linear subspace $V_2$ of $T_{p_1} \R^n \cong \R^n$
as follows. By hypothesis, the vector $v_1$ does not belong to the
characteristic direction $T_{p_1}M^1 \cap T_{p_1}^cM$, so the real
vector space $\left( \R \cdot v_1 \right) \oplus \left( T_{p_1}M^1
\cap T_{p_1 }^cM \right) \subset T_{ p_1} M^1$ is $2$-dimensional. We
choose an arbitrary $(n-2)$-dimensional real vector subspace
$V_2\subset T_{p_1}M^1$ which is a supplementary in $T_{p_1}M^1$ to
$\left( \R\cdot v_1\right) \oplus \left( T_{p_1} M^1 \cap T_{ p_1}^cM
\right)$ and we shall let the parameter $v$ vary
only in $V_2$.

From the rank properties {\bf
($\mathbf{ 5_1}$)} and {\bf ($\mathbf{ 6_1}$)} of Lemma~7.12 and from
the definition of $V_2$, we deduce that there exists $\varepsilon>0$
small enough with $\varepsilon < < c^2$ such that the mapping
\def\theequation{9.9}\begin{equation}
(x,v,u,\rho)\longmapsto 
A_{x,v,u:c}^1(\rho)
\end{equation}
is a {\it one-to-one immersion}\, from the open set $\{(x,v,u,\rho)\in
\R^n\times V_2 \times \R\times \R: \ \vert x \vert < \varepsilon,
\ \vert v \vert < \varepsilon, \ \vert u \vert < \varepsilon, \
1-\varepsilon < \rho < 1\}$ into $\C^n$. This property will be
important for uniqueness of the holomorphic extension in our
application of the continuity principle to be conducted in Lemma~9.20
below. Furthermore, shrinking $\varepsilon >0$ if necessary, we can
insure that the open subset
\def\theequation{9.10}\begin{equation}
\left\{
\aligned
{}
&
\mathcal{ W}_{p_1} :=\left\{
A_{x,v,u:c}^1(\rho)\in \C^n: 
(x,v,u,\rho)\in
\R^n\times V_2 \times \R\times \R, \right. \\
& \ \ \ \ \ \ \ \ \ \ \ \ \ \ \ \ \ \ \ \ \
\ \ \ \ \ \ \ \ \ \ \ \ \ \ \ \ \ \ \ \ \ \ \ \
\left.
\vert x \vert < \varepsilon,
\ \vert v \vert < \varepsilon, \ 
\vert u \vert < \varepsilon, \
1-\varepsilon < \rho < 1
\right\}
\endaligned\right.
\end{equation}
is a local wedge of edge $M$ at $(p_1, Jv_1)$, with $\mathcal{ W}_{
p_1} \cap M=\emptyset$.

Let the closed subset $C$ with $p_1 \in C$ and $C \backslash \{ p_1\}
\subset (H^1)^-$, let the neighborhood $\Omega$ of $M \backslash C$ in
$\C^n$, let the half-wedge $\mathcal{ HW }_{p_1 }^+$ be as in
Proposition~5.12, and let the sub-half-wedge $\mathcal{ HW}_1^+\subset
\mathcal{ HW}_{ p_1 }^+$ be as in \S5.14 and Lemma~5.37. 
In Cases {\bf ($\mathbf{ I_1
}$)} and {\bf ($\mathbf{ I_2}$)}, we shall consider the envelope of
holomorphy of the open subset $\Omega \cup \mathcal{ HW }_1^+ $. We
shall prove in the next paragraphs that, after possibly shrinking it a
little bit, its envelope of holomorphy contains the wedge $\mathcal{ W}_{ p_1 }$.

\subsection*{9.11.~Boundaries of analytic discs}
Since we want to apply the continuity principle Lemma~9.4, we have to
show that most discs $A_{x, v, u:c}^1 ( \zeta)$ have their boundaries
in $\Omega \cup \mathcal{ HW}_1^1+$. 
To this aim, it will be useful to decompose the boundary
$\partial \Delta$ in three closed parts $\partial \Delta = \partial^1
\Delta \cup \partial^2 \Delta \cup \partial^3 \Delta$, where
\def\theequation{9.12}\begin{equation}
\left\{
\aligned
\partial^1\Delta := 
& \
\left\{
e^{i\theta} \in \partial\Delta : \ 
\vert \theta \vert \leq \pi/2 -\varepsilon
\right\} \subset \partial^+\Delta, \\
\partial^2\Delta := 
& \
\left\{
e^{i\theta} \in \partial\Delta : \ 
\pi/2 +\varepsilon \leq
\vert \theta \vert \leq \pi
\right\} \subset \partial^-\Delta, \\
\partial^3\Delta := 
& \
\left\{
e^{i\theta} \in \partial\Delta : \ 
\pi/2-\varepsilon \leq 
\vert \theta \vert \leq \pi/2 +\varepsilon
\right\} \subset \partial\Delta, \\
\endaligned\right.
\end{equation}
where $\varepsilon$ with $0< \varepsilon < < c^2$ is as in \S9.8 just
above. This decomposition is illustrated in the left hand side of {\sc
Figure~17} below. Next, we observe that the two points $A_{0,
0,0:c}^1(i)$ and $A_{ 0,0, 0:c}^1(-i)$ belong to $(H^1 )^-\subset M
\backslash C\subset \Omega$, hence there exists a fixed open
neighborhood of these two points which is contained in $\Omega$. We
shall denote by $\omega^3$ such a (disconnected) neighborhood, for 
instance the union of two small open polydiscs centered at 
these two points. To
proceed further, we need a crucial geometric information about the
boundaries of the analytic discs $A_{x,v,u:c}^1 (\zeta)$ with $u\neq
0$.

\def\thelemma{9.13}\begin{lemma}
Under the assumptions of Cases~{\bf ($\mathbf{ I_1}$)} and~{\bf
($\mathbf{ I_2}$)} of Proposition~5.12, after
shrinking $\varepsilon >0$ if necessary, then 
\def\theequation{9.14}\begin{equation}
A_{x,v,u:c}^1\left(\partial \Delta
\right)\subset \Omega \cup 
\mathcal{ HW}_1^+,
\end{equation}
for all $x$ with $\vert x \vert < \varepsilon$, for all $v$ with
$\vert v \vert < \varepsilon$ and for all {\rm nonzero} $u\neq 0$ with
$\vert u \vert < \varepsilon$. 
\end{lemma}

\proof 
Firstly, since $A_{0, 0, 0:c}^1 (\pm i) \in \omega^3$, it
follows just by continuity of the family $A_{x, v, u: c}^1(\zeta)$
that, after possibly shrinking $\varepsilon >0$, the closed arc $A_{x,
v,u: c}^1 \left( \partial^3 \Delta \right)$ is contained in
$\omega^3$, for all $x$ with $\vert x \vert < \varepsilon$, for all
$v$ with $\vert v \vert < \varepsilon$ and for all $u$ with $\vert u
\vert < \varepsilon$. Secondly, since $A_{0, 0, 0:c}^1( \partial^2
\Delta) \subset A_{0, 0, 0: c}^1 \left( \partial^-\Delta \backslash
\{i,-i\}\right) \subset \mathcal{ HW }_1^+$, then by property {\bf
($\mathbf{ 9_1 }$)} of Lemma~8.3, it follows just by continuity of the
family $A_{x, v, u:c}^1(\zeta)$ that, after possibly shrinking
$\varepsilon >0$, the closed arc $A_{x, v, u:c}^1 \left( \partial^2
\Delta \right)$ is contained in $\mathcal{ HW}_1^+$, for all $x$ with
$\vert x \vert < \varepsilon$, for all $v$ with $\vert v \vert <
\varepsilon$ and for all $u$ with $\vert u \vert < \varepsilon$.
Thirdly, it follows from the inclusion $A_{x, v, u: c}^1 \left(
\partial^1 \Delta \right) \subset A_{x, v,u:c}^1 \left( \partial^+
\Delta \right) \subset M_u^1$ and from the inclusion $M_u^1 \subset
\Omega$ for all $u\neq 0$ that, after possibly shrinking
$\varepsilon>0$, the closed arc $A_{x, v, u:c}^1 \left( \partial^1
\Delta \right)$ is contained in $\Omega$, for all $x$ with $\vert x
\vert < \varepsilon$, for all $v$ with $\vert v \vert < \varepsilon$
and for all $u$ with $\vert u \vert < \varepsilon$ and $u \neq 0$. 
This completes the proof of Lemma~9.13.
\endproof

\bigskip
\begin{center}
\input three-boundary.pstex_t
\end{center}

In addition to this lemma, we notice that it follows immediately from
properties {\bf ($\mathbf{8_1}$)} and {\bf ($\mathbf{9_1}$)} of
Lemma~8.3 that $A_{x,0,0:c}^1(\overline{\Delta})$ is contained in
$\Omega \cup \mathcal{ HW}_1^+$ for all $x$ with $\vert x \vert <
\varepsilon$ such that $A_{x,0,0:c}^1(1)\in T^1\cap (H^1)^+$. This
property and Lemma~9.13 are illustrated in the right hand side of {\sc
Figure~17} just above.

\subsection*{9.15.~Analytic isotopies}
Next, in Case {\bf ($\mathbf{ I_1}$)}, we fix some $x_0= (x_{1; 0},0,
\dots, 0) \in \R^n$ with $0< x_{1;0} < \varepsilon$. Then $A_{x_0, 0,
0:c }^1 (1)= x_0+ i h(x_0)$ belongs to $T^1 \cap (
H^1)^+$. Analogously, in Cases {\bf ($\mathbf{ I_2}$)}, we fix some
$x_0=(0,\dots,0, x_{ n;0})\in \R^n$ with $0< x_{n;}0 <
\varepsilon$. Then in this second case, the point $A_{x_0, 0,
0:c}^1(1)=x_0+ ih(x_0)$ also belongs to $T^1 \cap (H^1 )^+$. We fix
the disc $A_{x_0, 0,0:c}^1( \zeta )$, which satisfies $A_{x_0,
0,0:c}^1( \overline{ \Delta })\subset \Omega \cup \mathcal{ HW 
}_1^+$.

\def\thelemma{9.16}\begin{lemma}
In Cases~{\bf ($\mathbf{ I_1})$} and~{\bf ($\mathbf{ I_2}$)}, every
disc $A_{x, v, u:c}^1 (\zeta)$ with $\vert x \vert < \varepsilon$,
$\vert v \vert < \varepsilon$, $\vert u \vert < \varepsilon$ and
$u\neq 0$ is analytically isotopic to the disc $A_{x_0, 0,
0:c}^1(\zeta)$, with the boundaries of the analytic discs of the
isotopy being all contained in $\Omega \cup \mathcal{ HW}_1^+$.
\end{lemma}

\proof
Indeed, since the set $\{u=0\}$ is a hyperplane, there clearly exists
a $\mathcal{ C}^{2, \alpha- 0}$-smooth curve $\tau \mapsto (x (\tau),v
(\tau),u (\tau))$ in the parameter space which joins a given arbitrary
point $(x^*, v^*,u^*)$ with $u^*\neq 0$ to the point $(x_0, 0,0)$
without meeting the hyperplane $\{u= 0\}$, except at its endpoint
$(x_0, 0,0)$. According to the previous Lemma~9.13, each boundary
$A_{x (\tau), v (\tau), u (\tau) :c}^1( \partial \Delta)$ is then
automatically contained in $\Omega \cup \mathcal{ HW }_1^+$, which
completes the proof.
\endproof

\subsection*{9.17.~holomorphic extension 
to a local wedge of edge $M$ at $p_1$} In Cases {\bf ($\mathbf{
I_1}$)} and~{\bf ($\mathbf{ I_2}$)}, we define the
following $\mathcal{ C}^{2,\alpha-0}$-smooth connected
hypersurface of $\mathcal{ W}_{p_1}$:
\def\theequation{9.18}\begin{equation}
\left\{
\aligned
\mathcal{ M}_{p_1}:= 
\left\{
A_{x,v,0:c}^1(\rho): \ 
(x,v,\rho)\in \R^n\times V_2 \times \R, \right. \\
\left.
\vert x \vert < \varepsilon, \ 
\vert v \vert < \varepsilon, \
1-\varepsilon < \rho < \varepsilon
\right\}, 
\endaligned\right.
\end{equation}
together with the following closed subset of 
$\mathcal{ M}_{p_1}$:
\def\theequation{9.19}\begin{equation}
\left\{
\aligned
\mathcal{ C}_{p_1} := 
& \
\left\{
A_{x,v,0:c}^1(\rho): \ 
(x,v,\rho)\in \R^n\times V_2 \times \R, \right. \\
& \ \ \ \ \ \ \ \ \ \ \ \ \ \ \ \ \ \ \ \ \ \ \ \ \ \
\left.
A_{x,v,0:c}^1(\partial^+ \Delta) \cap 
C= \emptyset, \ 
\vert x \vert < \varepsilon, \ 
\vert v \vert < \varepsilon, \
1-\varepsilon < \rho < \varepsilon
\right\}.
\endaligned\right.
\end{equation}
Since $A_{x, 0, 0: c}^1( \partial^+ \Delta)$ is contained in $(
H^1)^+$ for all $x$ such that $A_{ x, 0, 0}^1(1) \in T^1 \cap (
H^1)^+$, the closed subset $\mathcal{ C}_{ p_1}$ of $\mathcal{ M}_{
p_1}$ is a {\it proper}\, closed subset of $\mathcal{ M}_{ p_1}$.
The following figure provides a geometric illustration.

\bigskip
\begin{center}
\input wedge-half-wedge.pstex_t
\end{center}

We can now state the main lemma of this section, which will complete
the proof of Proposition~5.12 in Cases~{\bf ($\mathbf{ I_1}$)}
and~{\bf ($\mathbf{ I_2}$)}.

\def\thelemma{9.20}\begin{lemma}
After possibly shrinking $\Omega$ in a small neighborhood of $p_1$ and
after possibly shrinking $\varepsilon >0$, the set $\mathcal{ W}_{
p_1} \cap \left[ \Omega \cup \mathcal{ HW}_1^+ \right]$ is connected
and for every holomorphic function $f \in \mathcal{ O}\left( \Omega
\cup \mathcal{ HW }_1^+ \right)$, there exists a holomorphic function
$F \in \mathcal{ O} \left( \Omega \cup \mathcal{ HW }_1^+ \cup
\mathcal{ W}_{p_1} \right)$ such that $F \vert_{ \Omega \cup \mathcal{
HW}_1^+}=f$.
\end{lemma}

\proof
Remind that $\varepsilon < < c^2$ and remind that the wedge $\mathcal{
W}_{p_1}$ with $\mathcal{ W}_{p_1}\cap M= \emptyset$ in the two cases
is of size ${\rm O}( \varepsilon)$. Since $C$ is contained in
$(H^1)^-\cup \{ p_1\}\subset M^1$, we observe that $M \backslash C$ is
locally connected at $p_1$. Since the half-wedge $\mathcal{ HW }_1^+$
defined in Lemma~5.37 by simple inequalities is of size ${\rm O} (
\delta_1)$, if moreover $\varepsilon < < \delta_1$, after shrinking
$\Omega$ if necessary in a smal neighborhood of $p_1$ whose size
is ${\rm
O}(\varepsilon)$, it follows that we can assume that $\mathcal{
W}_{p_1} \cap \left[ \Omega \cup \mathcal{ HW}_1^+ \right]$ is
connected. However, in {\sc Figure~18} above, because we draw $M$ as
if it were one-dimensional, the intersection $\mathcal{ W}_{p_1} \cap
\left[ \Omega \cup \mathcal{ HW }_1^+ \right]$ appears to be
disconnected, which is a slight incorrection.

Let $f$ be an arbitrary holomorphic function in $\mathcal{ O} \left(
\Omega \cup \mathcal{ HW}_1^+ \right)$. Thanks to the isotopy
Lemma~9.16 and thanks to the continuity principle Lemma~9.4, we deduce
that $f$ extends holomorphically to a neighborhood in $\C^n$ of every
disc $A_{x, v, u: c}( \overline{ \Delta})$ whose boundary $A_{ x, v,
u:c} (\partial \Delta)$ is contained in $\Omega \cup \mathcal{
HW}_1^+$. Using the fact that the mapping~\thetag{ 9.9} is one-to
one, we deduce that we can extend $f$ uniquely by means of Cauchy's
formula at points of the form $A_{x,v,u:c}^1(\rho)$ with 
such values $\vert x
\vert < \varepsilon$, $\vert v \vert < \varepsilon$, $\vert u \vert <
\varepsilon$ and $1-\varepsilon < \rho < \varepsilon$ for 
which $A_{x,v,u:c}^1(\partial \Delta) \subset
\Omega \cup \mathcal{ HW}_1^+$, simply as follows
\def\theequation{9.21}\begin{equation}
f\left(
A_{x,v,u:c}^1(\rho)
\right):=\int_{\partial \Delta} \, 
\frac{ f\left(
A_{x,v,u:c}^1(\widetilde{ \zeta})
\right)}{\widetilde{ \zeta} -\rho} \, 
d \widetilde{ \zeta}.
\end{equation}
With this definition, we extend $f$ holomorphically and uniquely to
the domain $\mathcal{ W}_{ p_1} \backslash \mathcal{ C}_{p_1}$, where
$\mathcal{ C}_{p_1}$ is the proper closed subset, defined by~\thetag{
9.21}, of the $\mathcal{ C}^{2,\alpha-0}$-smooth hypersurface
$\mathcal{ M}_{ p_1 }$ defined by~\thetag{ 9.20}. Let $F \in \mathcal{
O} \left(\mathcal{ W}_{p_1} \backslash \mathcal{ C}_{p_1} \right)$
denote this holomorphic extension. Since $\mathcal{ W}_{p_1} \cap
\left[ \Omega \cup \mathcal{ HW}_1^+ \right]$ is connected, it follows
that $\left[ \mathcal{ W}_{p_1} \backslash \mathcal{ C}_{p_1}\right]
\cap \left[ \Omega \cup \mathcal{ HW }_1^+ \right]$ is also
connected. By Lemma~9.4, the function $f$ and its holomorphic
extension $F$ coincide in a neighborhood of every boundary
$A_{x,v,u:c}^1 (\partial \Delta)$ which is contained in the domain
$\Omega \cup \mathcal{ HW}_1^+$. From the analytic continuation
principle, we deduce that there exists a well-defined function, still
denoted by $F$, which is holomorphic in $\left[ \mathcal{ W}_{p_1}
\backslash \mathcal{ C}_{p_1}\right] \cap \left[ \Omega \cup \mathcal{
HW}_1^+ \right]$ and which extends $f$, namely $F\vert_{\Omega \cup
\mathcal{ HW}_1^+}=f$.

To conclude the proof of Proposition~5.12 in Cases {\bf ($\mathbf{
I_1}$)} and~{\bf ($\mathbf{ I_2}$)}, it suffices to extend $F$
holomorphically through the closed subset $\mathcal{ C}_{p_1}$ of the
connected hypersurface $\mathcal{ M}_{p_1}$ in the domain $\mathcal{
W}_{p_1} \subset \C^n$. Since $n\geq 2$, we notice that we are
exactly in the situation of Theorem~1.4 in the special, much simpler
case where the generic submanifold $M$ is replaced by a domain of
$\C^n$. It may therefore appear to be quite satisfactory to have
reduced Theorem~1.2' to the CR dimension $\geq 2$ version Theorem~1.4,
in an open subset of $\C^n$ (notice however that Theorem~1.4 as well
as Theorems~1.2 and~1.2' are stated in positive codimension, since the
case where $M$ is replaced by a domain of $\C^n$ is relatively trivial
in comparison). The removability of such a proper closed subset
contained in a connected hypersurface of a domain in $\C^n$ is known,
follows from~\cite{j4} and is explicitely stated and proved as
Lemma~2.10, p.~842 in~\cite{ mp1}. However, we shall still provide
another different geometric proof of this simple removability result
in Lemma~10.10 below, using fully the techniques developed in the
previous sections.

The proofs of Lemma~9.20 together with Cases~{\bf ($\mathbf{I_1}$)}
and~{\bf ($\mathbf{ I_2}$)} of Proposition~5.12 are complete now.
\endproof

\subsection*{9.22.~End of proof of Proposition~5.12 in 
Case {\bf (II)}} According to Lemma~5.37, in Case {\bf (II)}, the
one-codimensional totally real submanifold $M^1 \subset M$ is given by
the equations $y'=\varphi'(x,y_1)$ and $x_n=g(x'')$. If $u\in \R$ is
a small real parameter, we may define a ``translation'' $M_u^1$ of
$M^1$ in $M$ by the equations
\def\theequation{9.23}\begin{equation}
y'=\varphi'(x,y_1), \ \ \ \ \ \ \ \ 
x_n=g\left(x''\right)+u.
\end{equation}
Similarly as in \S9.6, we may construct a family of analytic discs
$A_{x,v,u:c}^1(\zeta)$ half-attached to $M_u^1$. We then we fix a
small scaling parameter $c$ with $0 < c \leq c_1$ so that properties
{\bf ($\mathbf{ 1_1}$)} to {\bf ($\mathbf{ 9_1}$)} of Lemmas~7.12
and~8.3 hold true. Similarly as in \S9.8, we shall restrict 
the variation
of the parameter $v$ to an arbitrary $(n-1)$-dimensional subspace
$V_1$ of $T_{p_1}M^1\cong \R^n$ which is supplementary to the real
line $\R\cdot v_1$ in $T_{p_1}M^1$. If $\varepsilon>0$ is small
enough with $\varepsilon < < c^2$, it follows that the mapping
\def\theequation{9.24}\begin{equation}
(x,v,u,\rho) \longmapsto 
A_{x,v,u:c}^1(\rho) 
\end{equation}
is a {\it one-to-one immersion}\, from the open set $\{(x,v,\rho)\in
\R^n\times V_1 \times \R\times \R: \ \vert x \vert < \varepsilon, \
\vert v \vert < \varepsilon, \ 1-\varepsilon < \rho < 1\}$ into
$\C^n$. Thanks to the choice of the linear subspace $V_1$, shrinking
$\varepsilon >0$ if necessary, it follows that for every $u$ with
$\vert u \vert < \varepsilon$, the open subset
\def\theequation{9.25}\begin{equation}
\left\{
\aligned
{}
&
\mathcal{ W}_u^1 :=\left\{
A_{x,v,u:c}^1(\rho)\in \C^n: 
(x,v,\rho)\in
\R^n \times \R\times \R, \right. \\
& \ \ \ \ \ \ \ \ \ \ \ \ \ \ \ \ \ \ \ \ \
\ \ \ \ \ \ \ \ \ \ \ \ \ \ \ \ \ \ \ \ \ \ \ \
\left.
\vert x \vert < \varepsilon,
\ \vert v \vert < \varepsilon, \ 
1-\varepsilon < \rho < 1
\right\}
\endaligned\right.
\end{equation}
is a local wedge of edge $M_u^1$. Clearly, this wedge 
$\mathcal{ W}_u^1$ depends $\mathcal{ C}^{2,\alpha-0}$-smoothly 
with respect to $u$. 

Using the fact that in Case~{\bf (II)} we have
\def\theequation{9.26}\begin{equation}
\frac{\partial A_{0,0,0:c}^1}{\partial \theta}(1)=
v_1=(1,0,\dots,0)\in T_{p_1}M^1\cap T_{p_1}^cM,
\end{equation}
one can prove that Lemma~9.13 holds true with $\mathcal{ HW}_1^+$
replaced by $\mathcal{ W}_2$ in~\thetag{ 9.14} and also that
Lemma~9.16 holds true, again with $\mathcal{ HW}_1^+$ replaced by
$\mathcal{ W}_2$. Similarly as in the proof of Lemma~9.20, applying
then the continuity principle and using the fact that the
mapping~\thetag{ 9.24} is one-to-one, after possibly shrinking
$\Omega$ in a neighborhood of $p_1$, and shrinking $\varepsilon>0$, we
deduce that for each $u\neq 0$, there exists a holomorphic function
$F\in \mathcal{ O}\left(\Omega \cup \mathcal{ W}_2 \cup
\mathcal{ W}_u^1 \right)$ with 
$F_{\vert \Omega \cup \mathcal{ W}_2}=f$.

To conclude the proof of Proposition~5.12 in Case {\bf (II)}, it
suffices to observe that for every fixed small $u$ with $
-\varepsilon < < u < 0$, the wedge $\mathcal{ W}_u^1$ contains in fact
a neighborhood $\omega_{p_1}$ of $p_1$ in $\C^n$.

The proofs of Proposition~5.12 and of Theorem~3.19 {\bf (i)} 
are complete now. \qed

\subsection*{9.27.~End of proof of Theorem~1.2'}
In order to derive Theorem~1.2' from Theorem~3.19 {\bf (i)}, we now
remind the necessity of supplementary arguments about the stability of
our constructions under deformation. Coming back to the strategy
developped in \S3.16, we had a wedge $\mathcal{ W}_1$ attached to $M
\backslash C_{\rm nr}$. Using a partition of unity, we may introduce a
one-parameter $\mathcal{ C}^{ 2, \alpha}$-smooth family of generic
submanifolds $M^d$, $d\in \R$, $d\geq 0$, with $M^0 \equiv M$, with
$M^d$ containing $C_{\rm nr}$ and with $M^d \backslash C_{ \rm nr}$
contained in $\mathcal{ W}_1$. In the proof of Theorem~3.19 {\bf (i)},
thanks to this deformation, the wedge $\mathcal{ W}_1$ was replaced by
a neighborhood $\Omega$ of $M\backslash C_{\rm nr}$ in $\C^n$.

In Sections~4 and~5, we constructed an important semi-local half-wedge
$(\mathcal{ HW}_\gamma^+)^d$ attached to a one-sided neighborhood of
$(M^1)^d$ in $M^d$ along a characteristic segment $\gamma^d$ of
$M^d$. Now, we make the crucial claim that, after possibly adapting
the deformation $M^d$, we may achieve that the geometric extent of
this semi-local half-wedge be uniform as $d>0$ tends to zero, namely
$(\mathcal{ HW}_\gamma^+)^d$ tends to a semi-local half-wedge
$(\mathcal{ HW}_\gamma^+)^0$ attached to a one-sided neighborhood of
$M^1$ in $M$ along $\gamma$, as $d$ tends to zero. Indeed, in
Section~4 we have constructed a family of analytic discs $(\mathcal{
Z}_{t,\chi,\nu:s}(\zeta))^d$ ({\it cf.}~\thetag{ 4.61}) which covers
the half-wedge $(\mathcal{ HW}_\gamma^+)^d$. Thanks to the stability
of Bishop's equation under $\mathcal{ C}^{2,\alpha}$-smooth
perturbations, the deformed family $(\mathcal{
Z}_{t,\chi,\nu:s}(\zeta))^d=: \mathcal{ Z}_{t,\chi,\nu:s}^d(\zeta)$ is
also of class $\mathcal{ C}^{2,\alpha-0}$ with respect to the
parameter $d$. We remind that for every $d>0$, the family $\mathcal{
Z}_{t,\chi,\nu:s}^d(\zeta)$ was in fact constructed by means of a
family $\widehat{ Z}_{ r_0,t,\tau,\chi,\nu:s}^d(\zeta)$ obtained by
solving Bishop's equation~\thetag{ 4.40}, where we now add the
parameter $d$ in the function $\Phi'$. In order to construct the
semi-local attached half-wedge, we have used the rank property stated
in Lemma~4.34. This rank property relied on the possibility of
deforming the disc $\widehat{ Z}_{r_0,t:s}(\zeta)$ near the point
$\widehat{ Z}_{r_0,t:s}^d(-1)$ in the open neighborhood
$\Phi_s(\Omega) \equiv \Phi_s\left(\mathcal{ W}_1 \right)$ of $
\Phi_s\left(M^d \right)$. As $d>0$ tends to zero, if $M^d$ tends to
$M$, the size of the neighborhood $\Phi_s\left(\mathcal{ W}_1 \right)$
shrinks to zero, hence it could seem that the we have no control on
the semi-local attached half-wege $(\mathcal{ HW}_\gamma^+)^d$ as
$d>0$ tends to zero. Fortunately, since the points $\widehat{
Z}_{r_0,0:s}^d(-1)$ in a neighborhood of which we introduce the
deformations~\thetag{ 4.30} are at a uniformly positive distance
$\delta >0$ from $\gamma$, we may choose the deformation $M^d$ of $M$
to tend to $M$ as $d$ tends to zero only in a small neighborhood of
$\gamma$, whose size is small in comparison to this distance
$\delta$. By smoothness with respect to $d$ of the family $\mathcal{
Z}_{t,\chi,\nu:s}^d (\zeta)$, we then deduce that the semi-local
half-wedge $(\mathcal{ HW}_\gamma^+)^d$ tends to a nontrivial
semi-local half-wedge $(\mathcal{ HW}_\gamma^+)^0$ as $d$ tends to
zero, which proves the claim.

Next, again thanks to the stability of Bishop's equation under
perturbation, all the constructions of Sections~5, 6, 7, 8 and~9 above
may be achieved to depend $\mathcal{ C}^{2, \alpha-0}$-smoothly, hence
uniformly, with respect to $d$. Importantly, we observe that if the
deformation $M^d$ is chosen so that $M^d$ tends to $M$ only in a small
neighborhood of $p_1$ of size $ < < \varepsilon$, then the shrinking
of $\varepsilon$ which occurs in Lemma~9.13 may be achieved to be
uniform as $d$ tends to zero, because the part
$A_{x,v,u:c}(\partial^3\Delta)$ stays in a uniform compact subset of
$\Omega$, as $d$ tends to zero. At the end of the proof of
Proposition~5.12, we then obtain univalent holomorphic extension to a
local wedge $\mathcal{ W}_{p_1}^d$ of edge $M^d$ or to a neighborhood
$\omega_{ p_1 }^d$ of $M^d$ in $\C^n$, and they tend smoothly to a
wedge $\mathcal{ W }_{ p_1}^0$ of edge $M$ at $p_1$ or to a
neighborhood $\omega_{ p_1 }$ of $p_1$ in $\C^n$.

The proof of Theorem~1.2' is complete.
\qed

\section*{\S10.~Three proofs of Theorem~1.4}

\subsection*{10.1.~Preliminary}
Theorem~1.4 may be established by means of the processus of
minimalization of generic submanifods developed by B.~J\"oricke
in~\cite{j2}, as is done effectively in~\cite{j4} in the hypersurface
case and then in~\cite{p2} in arbitrary codimension. In this section,
we shall suggest three more different proofs of Theorem~1.4. As
explained in Section~3, it suffices to treat $\mathcal{
W}$-removability, and essentially to prove Proposition~3.22. The first
proof appears already in~\cite{ m2} and also in~\cite{p1}. The second
proof consists in repeating some of the constructions of the previous
sections, using the fact that $M^1$ is of positive CR dimension in
order to simplify substantially the reasonings. The third proof
consists of a slicing argument showing that {\it Theorem~1.4 is in
fact a logical consequence of Theorem~1.2'}. In fact, because these
three proofs are already written elsewhere or very close to the
constructions developed in the previous sections, we shall only
provide summaries here.

\subsection*{10.2.~Normal deformation of analytic discs attached 
to $M^1$} Firstly, in the situation of Proposition~3.22, because $M^1$
is of positive CR dimension, we can construct a small analytic disc $A
(\zeta)$ attached to $M^1$ which satisfies $A(1) = p_1 \in M^1$ and $A
(\partial \Delta \backslash \{1\}) \subset (H^1)^+$. As in~\cite{m2},
\cite{ mp1}, using normal deformations of $A$ near $A(-1)$, we may
include $A$ in a $\mathcal{ C}^{ 2, \alpha-0}$-smooth parametrized
family $A_v (\zeta)$ of analytic discs attached to $M^1$, where $v \in
\R^{d+1}$ is small, so that the rank at $v=0$ of the mapping $v
\longmapsto -\frac{ \partial A_v}{ \partial \rho}(1) \in T_{ p_1} \C^n
\ {\rm mod} \ T_{p_1} M^1 \cong \R^{ d+1}$ is maximal equal to $(d
+1)$, the codimension of $M^1$ in $\C^n$. For this, we use a
deformation lemma which is essentially due to A.~Tumanov~\cite{tu3},
which appears as Lemma~2.7 in~\cite{ mp1} and which was already used
above (with a supplementary parameter $s$) in Lemma~4.34. Then we add
a ``translation'' parameter $x\in \R^{ 2m+d-1}$, getting a family of
analytic discs $A_{x,v} (\zeta)$ with $A_{x,v} (\partial
\Delta)\subset M^1$ so that the rank at $x=0$ of the mapping $x
\longmapsto A_{x,0} (1) \in M^1$ is maximal equal to $(2m+d-1)$, the
dimension of $M^1$. Finally, we introduce some ``translations''
$M_u^1$ of $M^1$ in $M$, where $u \in \R$, and we obtain a family $A_{
x, v,u} (\zeta)$ of analytic discs attached to $M_u^1$. For $u\neq 0$,
since the discs $A_{x,v,u}(\zeta)$ are attached to $M_u^1$, their
boundaries are contained in $M\backslash C$. Applying the
approximation Lemma~4.8, we deduce that for every $u\neq 0$, all
holomorphic functions in the open set $\Omega$ of Proposition~3.22
(which contains $M\backslash C$) extend holomorphically to the local
wedge $\mathcal{ W }_u := \{A_{x, v,u} (\rho): \ \vert x \vert <
\varepsilon, \ \vert v \vert < \varepsilon, \ 1-\varepsilon < \rho <
1\}$ of edge $M_u^1$. To control the univalence of the holomorphic
extension, it suffices to shrink $\Omega$ a little bit in a
neighborhood of $p_1$. To conclude the proof, one observes as
in~\cite{ j4}, \cite{ cs} that the union $\bigcup_{ u \neq 0} \,
\mathcal{ W}_u^1$ contains a local wedge of edge $M$ at $p_1$. This
processus is sometimes called ``sweeping out by wedges''. We notice
that this proof is geometrically much simpler than the proof of
Theorem~1.2' achieved in the previous sections.

\subsection*{10.3.~Half-attached analytic discs}
Secondly, we may generalize our constructions achieved in the previous
sections from the CR dimension $m=1$ case to the CR dimension $m\geq
2$ case, as follows. Let $M$, $M^1$, $p_1$, $H^1$ and $C$ be as in
Proposition~3.22. The fact that $M^1$ is of positive CR dimension
will provide substantial geometric simplifications for essentially two
reasons. Indeed, as the hypersruface $H^1\subset M^1$ is generic and
at least of real dimension $n$, there exists a (in fact infinitely many)
maximally real submanifold $K^1$ passing through $p_1$ which is
contained in $H^1$. Then
\begin{itemize}
\item[{\bf (i)}]
Applying the considerations of Section~7, we may construct families of
small analytic discs $A_{x,v} (\zeta)$ which are half-attached to
$K^1$ and which cover a local wedge of edge $K^1$ at $p_1$, when the
translation parameter $x$ and the rotation parameter $v$ vary. We
notice that this would be impossible in the case $m=1$, because in
this case $H^1$ is of real dimension $(n-1)$, hence does {\it not}\,
contain any maximally real submanifold.
\item[{\bf (ii)}]
We can even prescribe the direction $\frac{ \partial A_{ 0,0} }{
\partial \theta} (1)$ as an arbitrary given nonzero vector $v_1\in T_{
p_1}K^1$. Again, this would be impossible in the CR dimension $m=1$
case.
\end{itemize}
Let $\mathcal{ HW }_1^+$ be a local half-wedge of edge $(M^1)^+$ at
$p_1$, whose construction is suggested in \S4.64. Since $K^1$ is
generic, we may choose a nonzero vector $v_1\in T_{ p_1}K^1$ with the
property that $\mathcal{ HW }_1^+$ is directed by $J
v_1$. Generalizing Lemma~8.3, we see that $A_{x,v} \left( \overline{
\Delta } \backslash \partial^+ \Delta \right)$ is contained in
$\mathcal{ HW}_1^+$. We notice that in the CR dimension $m=1$ case,
$H^1$ is {\it not}\, generic, and we remember that the choice of a
special point $p_1$ to be removed locally and the choice of a
supporting hypersurface $H^1 \subset M^1$ was much more subtle,
because we had to insure that there exists a vector $v_1 \in T_{
p_1}H^1$ such that $\mathcal{ HW }_1^+$ is directed by $J v_1$.

Next, we can translate $K^1$ in $M^1$ by means of a small parameter
$t\in \R^{d-1}$ and then $M^1$ in $M$ by means of a small parameter
$u\in \R$. By stability of Bishop's equation, we get a family of
analytic discs $A_{x,v,t,u}(\zeta)$ half-attached to the translation
$K_{t,u}^1$ of $K^1$. Nine properties analogous to properties {\bf
($\mathbf{ 1_1}$)} to {\bf ($\mathbf{ 9_1}$)} of Lemmas~7.12 and~8.3
are then satisfied and we conclude the proofs of Proposition~3.22 and
of Theorem~1.4 in essentially the same way as in Section~9
above. We shall not write down all the details.

We notice that this second strategy of proof is much more complicated
than the first (known) proof summarized in \S10.2 just above.

\subsection*{10.4.~Slicing argument: reduction of Theorem~1.4 to 
Theorem~1.2'} Thirdly, we may provide a new proof of the central
Theorem 1.4 valid in CR dimension $m\geq 2$, in order to illustrate
how our results in CR dimension 1 can be applied to the general case
via slicing techniques. We shall see that Theorem~1.4 is a logical 
consequence of Theorem~1.2'.

Let $M$, $M^1$ and $C$ be as in Theorem~1.4. To begin with, we shall
treat the three notions of removability (CR-, $L^{\sf p}$- and
$\mathcal{ W}$-) commonly. However, we remind that CR-removability is
immediately reduced to $\mathcal{ W}$-removability thanks to
Lemma~3.5, hence it suffices to consider only $L^{\sf p}$- and
$\mathcal{ W}$-removability.

Arguing by contradiction, we see as in the
previous parts of this paper that we lose essentially nothing if we
consider $C$ to be the minimal nonremovable subset. Also, we may
assume holomorphic extension to a wedge attached to $M\backslash C$. 
As explained in \S3.16 it is enough to remove one single point of $C$.

As in \S3.16 ({\it see} especially the statement of Proposition~3.22),
we can show that there is a $\mathcal{ C}^{2,\alpha}$-smooth
hypersurface $H$ of $M$ which is generic in $\C^n$, transversal to
$M^1$, and which has the following property: $H$ contains some point
$p_1\in C$, and we can choose a small neighborhood of $H$ in $M$ in
which $H$ has two connected open sides $H^+$ and $H^-$ such that $C$
is contained in $H^-\cup\{p_1\}$ locally in a neighborhood of
$p_1$. Notice that the hypersurface $H^1$ constructed in Lemma~3.21 is
simply the intersection of $M^1$ with such a hypersurface $H$.

Next, in a neighborhood of $p_1$, we construct a local foliation of
$M$ by generic submanifolds of CR dimension 1 (hence of dimension
$n+1$) as follows. We notice that $\dim_\R H=2m+d-1$. We first choose
a local $\mathcal{ C}^{2,\alpha}$-smooth foliation of $H$ by generic
submanifolds of CR dimension 1 (hence of real dimension $n+1$) which
we denote by $\widetilde{ M}_s$, where the transversal parameters $s=
(s_1,\ldots,s_{m-2})\in \I_{m-2}(\rho_1)$ belong to a cube in
$\R^{m-2}$ of radius some $\rho_1>0$. Afterwards, we extend this
foliation to a local $\mathcal{ C}^{2,\alpha}$-smooth foliation
$\widetilde{ M}_{s,t}$ of a neighborhood of $p_1$ in $M$, where $s\in
\I_{m-2} (\rho_1)$, where $t\in \I(\rho_1)$ and where $\widetilde{
M}_{ s,0}\equiv \widetilde{ M}_s$; if $m=2$, we notice that the
parameter $s$ disappears and that we have only one real parameter $t$.
Also, we can assume that $\widetilde{ M}_{s,t}$ is contained in 
$H^+$ if and only if $t>0$.

By genericity, the submanifolds $\widetilde{ M}_{s,t}$ may be chosen
in addition to be transversal to the one-codimensional submanifold
$M^1\subset M$ containing the singularity $C$ with of course $p_1\in
\widetilde{ M}_{0,0}$. Then for all parameters $s$ and $t$ the
intersections $\widetilde{ M}_{s,t}^1:=M^1\cap \widetilde{ M}_{s,t}$
are maximally real submanifolds of $\C^n$. Considering $\widetilde{
M}_{s,t}^1$ as a maximally real one-codimensional submanifold of
$\widetilde{ M}_{ s,t}$, a characteristic foliation is induced on each
$\widetilde{ M}_{s,t}^1$. After contraction around $p_1$ we can
assume that these characteristic foliations all have trivial topology:
their leaves are the level sets of an $\R^{m-1}$-valued submersion. As
the intersection $\widetilde{ M}_{0,0}^1 \cap C$ is the singleton
$\{p_1\}$ (by construction of $H$), there exists $\varepsilon>0$ such
that for all $s$ with $\vert s \vert < \varepsilon$ and all $t$ with
$\vert t \vert < \varepsilon$, the closed subsets
$C_{s,t}:=\widetilde{ M}_{s, t}^1\cap C$ are compact in $\widetilde{
M}_{s, t}^1$. Notice that $C_{s,t}$ is even empty if 
$t>0$, because $C$ is contained in $H^-\cup \{p_1\}$.
The following simple fact shows that $C_{s,t}$ satisfies
the nontransversality condition $\mathcal{ F}^c_{
\widetilde{ M}_{s,t}^1}\{C_{s,t}\}$ of
Theorem~1.2', for all $s$ with $\vert s \vert < \varepsilon$ and all
$t$ with $-\varepsilon < t \leq 0$.

\def\thelemma{10.5}\begin{lemma}
Let $\mathcal{ F }_\mathcal{ M}$ be a $\mathcal{ C}^{ 1,
\alpha}$-smooth foliation by curves on some $m$-dimensional $\mathcal{
C}^{ 2, \alpha}$-smooth real manifold $\mathcal{ M}$ defined by a
surjective $\mathcal{ C}^{ 1, \alpha}$-smooth submersion $F: \mathcal{
M} \rightarrow \I_{ m-1}( \rho_1)$. Then every compact set $\mathcal{
C} \subset \mathcal{ M}$ satisfies the nontransversality condition
$\mathcal{ F}_\mathcal{ M}\{ \mathcal{ C} \}$.
\end{lemma}

\proof
Let $\mathcal{ C}'$ be an arbitrary compact subset of $\mathcal{ C}$.
As $\mathcal{ C}'$ is compact, there exists the smallest $\rho_2<
\rho_1$ with $\mathcal{ C}' \subset F^{-1}( \overline{ \I_{ m-1}(
\rho_2)})$. The semi-local projection $\pi_{ \mathcal{ F}_{\mathcal{
M}}}$ along the leaves of $\mathcal{ F }$ may of course be identified
with $F$. Thus, $\pi_{ \mathcal{ F}_{ \mathcal{ M } }} ( \mathcal{
C}')$ is contained in $\overline{ \I_{m-1 }( \rho_2)}$ and meets the
boundary $\partial \I_{ m-1} (\rho_2 )$ of the cube $\I_{m-1}
(\rho_2)$. Also, by compactness, the set
$\mathcal{ C}'$ cannot contain a fiber of $F$ in the whole. 
This completes the proof.
\endproof

We can now show that Theorem~1.4 is a logical consequence of
Theorem~1.2'. In fact, we cannot insure that the generic submanifolds
$\widetilde{ M}_{s,t}$ of CR dimension $1$ defined above are all
globally minimal, hence it seems that Theorem~1.2' itself does not
apply. However, we notice that the wedge attached to $M\backslash C$
restricts to a wedge attached to $\widetilde{ M}_{s,t}\backslash
C_{s,t}$, for all $s$ and $t$. Hence, we can observe that everything
that was needed in the proof of Theorem~1.2' was the existence of a
wedge attached to $M\backslash C$ to which holomorphic extension is
already assumed. One can even formulate a slightly more general
version of Theorem~1.2', where the global minimality assumption is
replaced by the assumption of holomorphic extension to a wedge
attached to $M\backslash C$. Of course, with this more general
assumption, $M$ may consists of several CR orbits, but thanks to
Lemma~3.5 about stability of CR orbits, one may check that the proof
of the main Theorem~3.19 {\bf (i)} and of the main Proposition~5.12
remain unchanged, in the case where holomorphic extension is assumed
in a wedgelike domain over $M\backslash C$ and not only in wedgelike
domains attached to the CR orbits of $M\backslash C$. 

Thus this slight generalization of Theorem~1.2' together with the
observation made in Lemma~10.5 just above yield that for all $s$ with
$\vert s \vert < \varepsilon$ and for all $t$ with $-\varepsilon < t
\leq 0$, the closed subset $C_{ s,t}$ is $\mathcal{ W}$-removable in
the generic submanifold $\widetilde{ M}_{s, t}$. We deduce that for
every $(s, t)$ with $\vert s \vert < \varepsilon$ and $-\varepsilon <
t \leq 0$, we get holomorphic extension from the given restricted
wedge attached to $\widetilde{ M}_{ s, t} \backslash C_{ s, t}$ into
an open wedge $\widetilde{ \mathcal{ W}}_{ s,t}$ attached $\widetilde{
M}_{ s, t}$. Notice that this does not immediately achieves the proof,
since the direction of the $\widetilde{ \mathcal{ W }}_{s, t}$ need
not depend continuously on $(s, t)$. In fact, the proof of the slight
generalization of Theorem~1.2' contains arguments (for example the
localization near a very special point) which do not depend nicely on
external parameters. Hence the attached wedges $\widetilde{ \mathcal{
W}}_{s,t}$ may well be completely unrelated.

To overcome this difficulty we proceed in the following way, already
argued in~\cite{ m2}, \cite{ mp1} (Lemma~2.7) in slightly different
contexts. We first construct a regular family $A_{ x, v}( \zeta)$ of
analytic discs attached to $\widetilde{ M}_{ 0, 0 }\cup \Omega$ whose
size is small in comparison to the basis of the wedge $\widetilde{
\mathcal{ W} }_{0, 0}$ and which sweep out a local wedge $\mathcal{ W}
(A_{ x, v})$ of edge $\widetilde{ M}_{ 0, 0}$ at $p_1$. Here, the
parameter $x\in \R^{ n+ 1}$ corresponds to translations in
$\widetilde{ M}_{ 0, 0}$ and $v \in \R^{ n-2}$ to normal deformations
in a neighborhood of the point $A_{ 0, 0} (-1)\in \Omega$. Deforming
this family thanks to the flexibility of Bishop's equation, we
construct a family $A_{x, v,s,t} (\zeta)$ attached to $\widetilde{
M}_{ s,t} \cup \Omega$, still sweeping out a local wedge $\mathcal{ W
}(A_{ x, v, s, t})$ of edge $\widetilde{ M}_{ s, t}$. This family is
of class $\mathcal{ C}^{2,\alpha-0}$ with respect to all
parameters. Using the $\mathcal{ W}$-removability of $C_{s,t}$, we can
introduce for every $(s, t)$ a one-parameter deformation $\widetilde{ M
}_{ s,t}^d$ of $\widetilde{ M}_{s,t}$ which is contained in the
attached wedge $\widetilde{ W}_{ s,t}$ whenever $d>0$ and which
coincides with $\widetilde{ M}_{s,t}$ when $d=0$. Thanks to the
flexibility of Bishop's equation with parameters, we get a deformed
family $A_{x,v, s, t: d}(\zeta)$ of analytic discs. Since the wedges
$\widetilde{ W}_{ s,t}$ are {\it a priori}\, unrelated, we loose the
smoothness with respect to all variables, including $d$. Fortunately,
by an application of the continuity principle, for every $d>0$, we
deduce holomorphic extension to the wedge generated by the family
$A_{x, v,s, t:d} (\zeta)$. If we let $d$ tend to zero, fixing $(s,t)$,
we obtain univalent holomorphic extension to the wedge generated by
the family $A_{x,v,s,t}(\zeta)$. Finally, as $(s,t)$ varies, the
wedges $\mathcal{ W}( A_{x,v,s,t})$ varies smoothly and covers a local
wedge $\mathcal{ W}$ of edge $M$ at $p_1$. By the continuity
principle, we may verify that we obtain univalent holomorphic
extension to $\mathcal{ W}$.

Secondly, we explain how $L^{\sf p}$-removability of the
point $p_1\in C$ in $M$ follows logically from the $L^{\sf
p}$-removability of every $C_{s,t}$ in $\widetilde{ M}_{s,t}$. The
main argument relies on the following simple but useful fact: if $M$
is a generic CR submanifold of $\C^n$ and if $N\subset M$ is a lower
dimensional submanifold which is itself a generic CR submanifold of
$\C^n$ of positive CR dimension, then differentiable CR functions on
$M$ obviously restrict to CR functions on $N$. More generally, a
foliated version of this observation with lower regularity assumptions
is as follows.

\def\thelemma{10.6}\begin{lemma}
Let $M\subset\C^n$ be a generic submanifold of class $\mathcal{
C}^{2,\alpha}$ and of CR dimension $m\geq 2$.

\begin{itemize}
\item[{\bf (a)}]
If $p_1\in M$ and
if $M$ carries a local $\mathcal{ C }^{ 2, \alpha}$-smooth foliation
by a family $N_u$ of generic submanifolds of CR dimension $1$ where
$u \in \R^{m-1}$ is a small parameter and where
$p_1 \in N_0$, then for every CR function $f
\in L_{ loc}^{\sf p} (M)$, ${\sf p}\geq 1$, and for almost every $u\in
\R^{ m-1}$, its restriction $f|_{ N_u}$ is an $L_{ loc}^{\sf p} ( N_u)$
function which is CR on $N_u$.
\item[{\bf (b)}]
Conversely, if $p_1 \in M$ and if $M$ carries $m$ local $\mathcal{
C}^{2, \alpha}$-smooth foliations by families $N_{u_j}^j$,
$j=1, \dots,m$, $u_j\in \R^{m-1}$, of generic submanifolds of CR
dimension $1$ satisfying $p_1\in N_0^j$ for $j=1,\dots,m$ and
\def\theequation{10.7}\begin{equation}
T_{p_1}N_0^1+\cdots + T_{p_1}N_0^m =T_{p_1}M, 
\end{equation}
then a function $f\in L_{ loc}^{ \sf p }(M)$ is CR in a neighborhood of
$p_1$ if and only if for all $j =1, \dots,m$ and for almost every
$u_j \in \R^{ m-1}$, its restriction $f\vert_{ N_{ u_j}^j}$ is CR on
$N_{ u_j }^j$.
\end{itemize}
\end{lemma}

\proof 
Of course, property {\bf (a)} only makes sense for an everywhere
defined representative of $f$ and the nullset of excluded parameters
$t$ depends on the choice of the representative of $f$.

To establish {\bf (a)}, we choose a small box-neighborhood $U\cong
\I_{ m-1}( \rho_1)\times \widetilde{ N}$, where $\I_{ m-1}(\rho_1)$ is
a cube of some positive radius $\rho_1 >0$ in $\R^{ m-1}$, such that
every plaque $\{ \upsilon \}\times \widetilde{N}$ is an open subset of
some leaf $N_u$. By the $L^{\sf p}$ version of the
approximation theorem ({\it cf.}~\cite{j5}, \cite{p1}, \cite{mp1}),
the restriction $f|_U$ is the limit in the $L^{ \sf p}$ norm of the
restrictions of holomorphic polynomials $(P_\nu )_{ \nu \in
\N}$. Thanks to Fubini's theorem, we deduce that for almost every
$\upsilon \in \I_k ( \rho_1)$, the restriction $P_\nu |_{ \{\upsilon
\} \times \widetilde{ N }}$ converges in $L^{\sf p}$ norm to $f|_{ \{
\upsilon \} \times \widetilde{ N }}$. Hence for such parameters
$\upsilon$, the restriction $f \vert_{\{\upsilon \} \times \widetilde{
N}}$ is CR, which completes the proof of {\bf (a)}.

To establish {\bf (b)}, we observe first that the ``only if'' part is
a direct consequence of {\bf (a)}. To prove the ``if'' part, we may
introduce for every $j=1,\dots,m$ and for every $u_j\in \R^{m-1}$ a
$(0,1)$ vector field $\overline{L}_{u_j}^j$ tangent to $N_{u_j}^j$ and
$\mathcal{ C}^{1,\alpha}$-smoothly parameterized by $u_j$. The
geometric assumption~\thetag{ 10.7} entails that the $m$ vector fields
$\overline{L}_{u_1}^1,\dots,\overline{L}_{u_m}^m$ generate the CR
bundle $T^{0,1}M$ in a neighborhood of $p_1$. By assumption, the
$L_{loc}^{\sf p}$ function $f$ is annihilated in the distributional
sense by these $m$ vector field, hence it is CR. This completes the
proof of {\bf (b)}.
\endproof

We can now prove that the $L^{ \sf p}$-removability of $C$ in $M$
follows from an application of Theorem~1.2'. Let a function $f\in
L^{\sf p}_{ loc}(M)$ which is CR on $M \backslash C$. Coming back to
the construction of the submanifolds $\widetilde{ M}_{s,t}$ achieved
in the paragraphs before Lemma~10.5, it is clear that for almost every
$(s,t) \in \R^{ m-1}$, the restriction $f\vert_{ \widetilde{
M}_{s,t}}$ is $L_{ loc}^{ \sf p}$-integrable. More generally,
proceeding as in the paragraph before Lemma~10.5, we may construct $m$
such families $\widetilde{ M}_{j; s_j,t_j}$ for $j=1, \dots, m$ with
$p_1\in \widetilde{ M}_{ j; 0,0}$ and
\def\theequation{10.8}\begin{equation}
T_{p_1}\widetilde{ M}_{1;s_1,t_1}+\cdots+
T_{p_1}\widetilde{ M}_{m;s_m,t_m}=
T_{p_1}M,
\end{equation}
without changing the conclusion that the corresponding closed subsets
$C_{s_j, t_j}^j$ are $L^{\sf p}$-removable in $\widetilde{ M }_{j;
s_j,t_j}$ for $j =1, \dots, m$. Applying Lemma~10.7 just above, we
finally deduce that $f$ is CR on $M$, as desired.

This completes the description of the reduction of
Theorem~1.4 to Theorem~1.2' via a slicing argument.

\subsection*{10.9.~Version of Theorem~1.4 in an open subset 
of $\C^n$} To conclude this section, we remind that in the end of the
proof of Proposition~5.12 in Cases {\bf ($\mathbf{ I_1}$)} and {\bf
($\mathbf{ I_2}$)}, we came down to the removability of a proper
closed subset $\mathcal{ C}_{p_1}$ of a one-codimensional submanifold
$\mathcal{ M}_{p_1}$ of an open subset of $\C^n$ ($n\geq 2$), namely
the wedge $\mathcal{ W}_{p_1}$ (remind {\sc Figure~18} above), which
amounts exactly to prove Theorem~1.4 in the case where the generic
submanifold $M$ is replaced by an open subset of $\C^n$. We may
formulate this result as the following lemma. To our knowledge, its
first known proof is given in~\cite{ j4}. Here, we provide a slightly
different proof, using half-attached analytic discs.

\def\thelemma{10.10}\begin{lemma}
Let $\mathcal{ D}\subset \C^n$ be a domain, let $\mathcal{ M}^1
\subset \mathcal{ D}$ be a connected $\mathcal{ C}^{2,\beta}$-smooth
hypersurface with $0 < \beta < 1$ and let $\mathcal{ C}$ be a proper
closed subset of $\mathcal{ M}^1$ which does not contain any CR orbit
of $\mathcal{ M}^1$. Then for every holomorphic function $f \in
\mathcal{ O} \left( \mathcal{ D} \backslash \mathcal{ C} \right)$,
there exists a holomorphic function $F \in \mathcal{ O}(\mathcal{ D})$
such that $F \vert_{\mathcal{ D} \backslash \mathcal{ C }} = f$.
\end{lemma}

\proof
We summarize the proof, which anyway is very similar to the proof of
Proposition~3.22 delineated in \S10.3 above. Reasoning by
contradiction and constructing a supporting hypersurface, we come down
to the local removability of a single point $p_1$ in a geometric
situation analogous to the one described in Proposition~3.22, with $M$
replaced by the domain $\mathcal{ D}$, with $M^1$ replaced by
$\mathcal{ M}^1$ and with a generic $\mathcal{ C}^{2,\beta}$-smooth
submanifold $\mathcal{ H}^1\subset \mathcal{ M}^1$ such that, locally
in a neighborhood of $p_1$, we have $\mathcal{ C}\subset (\mathcal{
H}^1)^-\cup \{p_1\}$, where we use the same notation $\mathcal{ C}$
for the smallest non-removable subset of the original $\mathcal{ C}$.

Notice that $\mathcal{ H}^1$ is of codimension $2$. Let $\mathcal{
K}^1\subset \mathcal{ H}^1$ be a maximally real submanifold passing
through $p_1$. We may translate $\mathcal{ K}^1$ in $\mathcal{ M}^1$
by means of a small parameter $t\in \R^{n-1}$ and then $\mathcal{
M}^1$ in $\mathcal{ D}$ by means of a small parameter $u\in \R$.
Following Section~7, we then construct a small family of analytic
discs $\mathcal{ A}_{x,v,t,u}^1(\zeta)$ half-attached to the
``translations'' $\mathcal{ K}_{t,u}^1$. Nine properties analogous to
properties {\bf ($\mathbf{ 1_1}$)} to {\bf ($\mathbf{ 9_1}$)} of
Lemmas~7.12 and~8.3 are then satisfied and we conclude the proof 
in essentially the same way as in Section~9 above.
\endproof

\section*{\S11~$\mathcal{ W}$-removability implies
$L^{\sf p}$-removability}

\subsection*{11.1.~Preliminary}
This section is devoted to prove Lemma~3.15 about $L^{\sf
p}$-removability of the proper closed subset $C\subset M^1$, granted
it is $\mathcal{ W}$-removable. More generally, we shall establish the
$L^{\sf p}$-removability of certain proper closed subsets $\Phi$ of
$M$ that are nullsets with respect to the Lebesgue measure of $M$.

As a preliminary, we remind that if $M'$ is a globally minimal
$\mathcal{ C}^{2,\alpha}$-smooth generic submanifold of $\C^n$ of CR
dimension $m \geq 1$ and of codimension $d= n-m\geq 1$, there exists a
wedge $\mathcal{ W}'$ attached to $M'$ constructed by means of
analytic discs successively glued to $M'$ and to conelike submanifolds
attached to $M'$ consisting of parametrized families of pieces of
analytic discs. By means of the approximation theorem of~\cite{ bt},
one deduces classically that continuous CR functions on $M'$ extend
holomorphically to $\mathcal{ W}'$, and continuously to $M'\cup
\mathcal{ W}'$.

For the holomorphic extension of the $L_{loc}^{\sf p}$ CR functions to
a wedge attached to $M'$, some supplementary routine, though not
obvious, work has to be achieved. Firstly, using a convolution with
Gauss' kernel as in~\cite{bt}, one shows that on a $\mathcal{
C}^2$-smooth generic submanifold $M'$ of $\C^n$, every $L_{loc}^{\sf
p}$ CR function on $M'$ is locally the limit, in the $L^{\sf p}$ norm,
of a sequence of polynomials ({\it see} Lemma~3.3 in~\cite{j5}). In
the case where $M'$ is a hypersurface, studied in~\cite{j5}, the wedge
$\mathcal{ W}'$ is in fact a one-sided neighborhood attached to $M'$,
which we will denote by $\mathcal{ S}'$. The theory of Hardy spaces on
the unit disc transfers to parameterized families of small analytic
discs glued to $M'$ which cover local one-sided neighborhoods of a
hypersurface, provided the boundaries of these discs foliate an open
subset of $M'$. Using in an essential way L.~Carleson's imbedding
theorem, B.~J\"oricke established in~\cite{j5} that every $L^{\sf
p}_{loc}$ CR function on a globally minimal $\mathcal{ C}^2$-smooth
hypersurface $M'$ extends holomorphically in the Hardy space $H^{\sf
p}(\mathcal{ S}')$ of holomorphic functions defined in the one-sided
neighborhood $\mathcal{ S}'$, with
$L^{\sf p}$ boundary values on the hypersurface $M'$. 
In his thesis~\cite{p1}, the second author
of the present paper has built the theory in higher codimension,
introducing the Hardy space $H^{\sf p}(\mathcal{ W}')$ of functions
holomorphic in the wedge $\mathcal{ W}'$ attached to $M'$, with
$L^{\sf p}$ boundary values on the edge $M'$.

At present, these background statements about holomorphic
extendability of $L_{loc}^{\sf p}$ CR functions on globally minimal
generic submanifolds may be reproved in a more elegant way than by
going through the rather complicated technology dispersed in the
articles~\cite{tu1}, \cite{tu2}, \cite{m1}, \cite{j2}, thanks to a
simplification of the wedge extendability theorem obtained recently by
the second author of this paper, which treats in an unified way local
and global minimality. We refer the reader to the work in
preparation~\cite{p3} for a substantial cleaning of the theory.

\subsection*{11.2.~$L^{\sf p}$-removability of nullsets}
Let us say that a subset $\Phi$ of a $\mathcal{ C}^{2,\alpha}$-smooth
generic submanifold is {\sl stably $\mathcal{ W}$-removable} if it is
$\mathcal{ W}$-removable on every compactly supported sufficiently
small $\mathcal{ C}^{2,\alpha}$-smooth deformation $M^d$ of $M$
leaving $\Phi$ fixed. In the situations of Theorems~1.2' and~1.4, the
assumptions of Lemma~11.3 just below are satisfied with $\Phi=C$,
taking account of the fact that we have already established the
$\mathcal{ W}$-removability of $C$ and that for logical reasons only,
the closed set $C$ in the statements of Theorems~1.2' and~1.4 is
obviously stably removable.

\def\thelemma{11.3}\begin{lemma}
Let $M$ be a $\mathcal{ C}^{2, \alpha}$-smooth generic submanifold of
$\C^n$ of CR dimension $m \geq 1$ and of codimension $d= n-m\geq 1$,
hence of dimension $(2m+ d)$, let $\Phi \subset M$ be a nonempty
proper closed subset whose $(2m +d)$-dimensional Hausdorff measure is
equal to zero. Assume that $M \backslash \Phi$ is globally minimal and
let $\mathcal{ W}$ be a wedge attached to $M\backslash \Phi$ such that
every function in $L_{loc}^{\sf p}(M) \cap CR(M\backslash \Phi)$
extends holomorphically as a function in the Hardy space $H^{\sf
p}(\mathcal{ W})$. If $\Phi$ is stably 
$\mathcal{ W}$-removable, then $\Phi$
is $L^{\sf p}$-removable.
\end{lemma}

Before giving the proof,
let us summarize intuitively the reason why this strong $L^{\sf
p}$-removability result Lemma~11.3 holds. Indeed, let $f\in L_{loc
}^{\sf p}(M) \cap CR (M \backslash \Phi)$. As soon as wedge extension
over points of $\Phi$ is known, thanks to the fact that we can deform
$M$ over $\Phi$ in the wedgelike domain, thus erasing the singularity
$\Phi$, we get a $L_{loc}^{ \sf p}$ CR function $f^d$ on the deformed
manifold $M^d$, without singularities anymore, and in additition, we
can let the deformation $M^d$ tend to $M$ with a uniform with $L^{\sf
p}$ control of the extension $f^d$, which therefore tends to a CR
extension of $f$ through $\Phi$.

\proof
First of all, we remind that for every ${\sf p}$ with $1\leq {\sf p}
\leq \infty$, the space $L_{loc}^{\sf p}(M)$ is contained in
$L_{loc}^1(M)$. We claim that it follows that $\Phi$ is $L^{\sf
p}$-removable for every ${\sf p}$ with $1\leq {\sf p} \leq \infty$ if
and only if $\Phi$ is $L^1$-removable. Indeed, suppose that $\Phi$ is
$L^1$-removable, namely for every function $f\in L_{loc}^1(M)\cap {\rm
CR}(M\backslash \Phi)$, and every $\mathcal{ C}^1$-smooth
$(n,m-1)$-form with compact support, we have $\int_M \, f \cdot
\overline{\partial} \psi =0$. In particular, since $L_{loc}^{\sf p}$
is contained in $L_{loc}^1$ by H\"older's inequality, this property
holds for every function $g\in L_{loc}^{\sf p}(M) \cap CR(M\backslash
\Phi)$, hence $\Phi$ is $L^{\sf p}$-removable, as
claimed. Consequently, it suffices to show that $\mathcal{
W}$-removability implies $L^1$-removability.

Let $f\in L_{loc}^1(M\backslash \Phi)\cap L^1(M)$ be an arbitrary
function. The goal is to show that $f$ is in fact CR on $\Phi$. Of
course, it suffices to show that $f$ is CR locally at every point of
$\Phi$. So, we fix an arbitrary point $q\in \Phi$. If $\psi$ is an
arbitrary $(n,m-1)$-form of class $\mathcal{ C}^1$ supported in a
sufficiently small neighborhood of $q$, we have to prove that $\int_M
\, f \cdot \overline{\partial} \psi =0$.

We may also fix a small open polydisc $\mathcal{ V}_q$ centered at
$q$. We shall first argue that we can assume that the $L_{loc}^1$
function $f$ is holomorphic in a neighborhood of $\left(M\backslash
\Phi\right)\cap \mathcal{ V}_q$ in $\C^n$. Indeed, since $M \backslash
\Phi$ is globally minimal, there exists a wedge $\mathcal{ W}$
attached to $M\backslash \Phi$ such that every $L_{loc}^1$ CR function
on $M \backslash \Phi$, and in particular $f$, extends holomorphically
as a function which belongs to the Hardy space $H^1(\mathcal{ W})$. By
slightly deforming $\left(M\backslash \Phi\right) \cap \mathcal{ V}_q$
into $\mathcal{ W}$ along Bishop discs glued to $M\backslash \Phi$,
keeping $\Phi$ fixed, using the theory of Hardy spaces in wedges
developed in~\cite{p1}, we may obtain the following deformation result
with $L^1$ control, a statement which is a particular case of
Proposition~1.16 in~\cite{mp1}.

\def\theproposition{11.4}\begin{proposition}
For every $\varepsilon > 0$, every $\beta <\alpha$, there exists a
small $\mathcal{ C}^{ 2, \beta}$-smooth deformation $M^d$ of $M$ with
support contained in $\overline{ \mathcal{ V }_q}$ and there exists a
function $f^d \in L_{ loc}^1 \left(M^d \right)\cap CR \left( M^d
\backslash \Phi \right)$, such that
\begin{itemize}
\item[{\bf (1)}]
$M^d \cap \mathcal{V}_q \supset \Phi \cap \mathcal{ V}_q\ni q$.
\item[{\bf (2)}]
$\left(M^d\backslash \Phi \right)\cap \mathcal{ V}_q 
\subset \mathcal{ W}\cap \mathcal{ V}_q$.
\item[{\bf (3)}]
$f^d$ is holomorphic in the neighborhood
$\mathcal{ W}\cap \mathcal{ V}_q$ of $(M^d\backslash \Phi)\cap 
\mathcal{ V}_q$ in $\C^n$.
\item[{\bf (4)}]
$M\cap \mathcal{ V}_q$ and $M^d\cap \mathcal{ V}_q$ are graphed over
the same $(2m+d)$ linear real subspace and $\left\vert\left\vert M^d\cap
\mathcal{ V}_q - M \cap \mathcal{ V}_q
\right\vert\right\vert_{\mathcal{ C}^{2,\beta}} \leq \varepsilon$.
\item[{\bf (5)}]
The volume forms of $M\cap \mathcal{ V}_q$ and of $M^d \cap \mathcal{
V}_q$ may be identified and $\left\vert f-f^d
\right\vert_{L^1(M\cap \mathcal{ V}_q)} \leq \varepsilon$. 
\end{itemize}
\end{proposition}

Since it will suffice to have a control of the
deformation $M^d$ only in $\mathcal{ C}^2$ norm, we shall
replace $\mathcal{ C}^{2,\beta}$ and $\mathcal{ C}^{1,\beta}$
by $\mathcal{ C}^2$ and $\mathcal{ C}^1$ in the sequel.

Let us be more explicit about conditions {\bf (4)} and {\bf
(5)}. Without loss of generality, we can assume that in coordinates
$(z,w)= (x+iy, u+iv)\in \C^m\times \C^d$ centered at $q$, we have
$T_qM= \{v=0\}$, hence the generic submanifolds $M$ and $M^d$ are
represented locally by vectorial equations $v=\varphi(x,y,u)$ and $v=
\varphi^d (x,y,u)$, where $\varphi$ and $\varphi^d$ are defined in the
real cube $\I_{2m+d}(2\rho_1)$, for some small $\rho_1>0$ and that
$\mathcal{ V}_q$ is the polydisc $\Delta_n(\rho_1)$ of radius
$\rho_1$. Then condition {\bf (4)} simply means that $\vert \vert
\varphi^d - \varphi \vert \vert_{\mathcal{ C}^2 (\I_{2m+d}
(\rho_1) )}\leq \varepsilon$ and condition {\bf (5)} is clear if we
choose $dxdydu$ as the volume form on $M$ and on $M^d$.

Suppose that for every $\varepsilon>0$ and for every 
deformation $M^d$, we can show that
the function $L_{loc}^1$ function $f^d$ on $M^d$ is in fact CR over
$M^d\cap \Delta_n(\rho_1)$. Then we claim that $f$ is CR in a
neighborhood of $q$.

Indeed, to begin with, let us denote by
$\overline{L}_1,\dots,\overline{L}_m$ a basis of $(0,1)$ vector fields
tangent to $M$, having coefficients depending on the first order
derivatives of $\varphi$. More precisely, in slightly abusive matrix
notation, we can choose the basis $\overline{L}:=
\frac{\partial}{\partial \bar z}+2(i-\varphi_u)^{-1} \, \varphi_{\bar
z} \frac{\partial}{\partial \bar w}$. Let us denote this basis
vectorially by $\overline{L}=\frac{\partial}{\partial \bar z}+ A \,
\frac{\partial}{\partial \bar w}$. To compute the formal adjoint of
$\overline{L}$ with respect to the local Lebesgue measure $dx dy du$
on $M$, we choose two $\mathcal{ C}^1$-smooth functions $\psi$, $\chi$
of $(x,y,u)$ with compact support in $\I_{2m+d}(\rho_1)$. Then the
integration by part $\int \overline{L} (\psi) \cdot \chi \cdot dxdydu=
\int\, \psi\cdot {}^T \overline{L}
(\chi)\cdot dx dy du$ yields the explicit
expression ${}^T\overline{L} (\chi):= -\overline{L}(\chi)-A_{\bar w}
\cdot \chi$ of the formal adjoint of $\overline{L}$.

It follows immediately that if we denote by ${}^T(\overline{L}^d)$ the
formal adjoint of the basis of CR vector fields tangent to $M^d$, then
we have an estimate of the form $\vert \vert {}^T(\overline{L}^d) -
{}^T(\overline{L})\vert \vert_{ \mathcal{ C}^1} \leq {\rm C} \cdot
\varepsilon$, for some constant $C>0$. Recall that $f^d$ is assumed to
be CR in $M^d\cap \Delta_n(\rho_1)$. Equivalently, we have $\int \,
f^d \cdot \, {}^T(\overline{L}^d)(\psi)\cdot dxdydu=0$ for every
$\mathcal{ C}^1$-smooth function $\psi$ with compact support in the
cube $\I_{2m+d}(\rho_1)$. Then we deduce that (some explanation
follows)
\def\theequation{11.5}\begin{equation}
\left\{
\aligned
{}
&
\left\vert
\int\, f\cdot {}^T\overline{L}(\psi) \cdot dxdydu \right\vert=
\left\vert
\int\, 
\left[
f\cdot {}^T\overline{L}(\psi)-
f^d\cdot {}^T(\overline{L}^d)(\psi)
\right] \cdot dxdydu \right\vert \\
& \
\leq 
\left\vert
\int\, 
\left[
f\cdot {}^T\overline{L}(\psi)-f\cdot {}^T(\overline{L}^d)(\psi)+
f\cdot {}^T(\overline{L}^d)(\psi)-
f^d\cdot {}^T(\overline{L}^d)(\psi)
\right] \cdot dxdydu \right\vert
\\
& \
\leq 
C_1(\psi) \cdot \varepsilon \cdot 
\int_{\I_{2m+d}(\rho_1)} \, 
\vert f \vert \cdot dxdydu + C_2(\psi) \cdot
\int_{\I_{2m+d}(\rho_1)}\, 
\vert f - f^d \vert \cdot dxdydu
\\
& \
\leq C(\psi,f,\rho_1)\cdot \varepsilon,
\endaligned\right. 
\end{equation}
taking account of property {\bf (5)} of Proposition~11.4 for the
passage from the third to the fourth line, where $C(\psi,f,\rho_1)$ is
a positive constant depending only on $\psi$, $f$ and $\rho_1$. As
$\varepsilon$ was arbitrarily small, it follows that $\int\, f\cdot
{}^T\overline{L}(\psi) \cdot dxdydu =0 $ for every $\psi$, namely $f$
is CR on $M\cap \Delta_n(\rho_1)$, as was claimed.

It remains to show that $f^d$ is CR on $M^d\cap \Delta_n(\rho_1)$.
First of all, we need some observations. For every compactly
supported small deformation $M^d$ stabilizing $\Phi$, the wedge
$\mathcal{ W}$ attached to $M\backslash \Phi$ is still a wedge
attached to $M^d\backslash \Phi$. In addition, this wedge contains a
neighborhood of $\left(M^d \backslash \Phi\right)\cap
\Delta_n(\rho_1)$ in $\C^n$ by property {\bf (3)} of Proposition~11.4.
As $\Phi$ was supposed to be stably removable, it follows that there
exists a wedge $\mathcal{ W}_1$ attached to $M^d$ (including points of
$\Phi$) to which holomorphic functions in $\mathcal{ W}$ extend
holomorphically.

Consequently, replacing $M^d\cap \Delta_n(\rho_1)$ by $M$, we are led
to prove the following lemma, which, on the geometric side, is totally
similar to Lemma~11.3, except that the wedge $\mathcal{ W}$ attached
to $M\backslash \Phi$ appearing in the formulation of Lemma~11.3 is
now replaced by a neighborhood $\Omega$ of $M\backslash \Phi$ in
$\C^n$.

\def\thelemma{11.6}\begin{lemma} 
Let $M$ be a $\mathcal{ C}^{2, \alpha}$-smooth generic submanifold of
$\C^n$ of CR dimension $m \geq 1$ and of codimension $d= n-m\geq 1$, 
let $\Phi \subset M$ be a nonempty
proper closed subset whose $(2m +d)$-dimensional Hausdorff measure is
equal to zero. Let $\Omega$ be a neighborhood of $M\backslash \Phi$ in
$\C^n$ and let $\mathcal{ W}_1$ be a wedge attached to $M$, including
points of $\Phi$.
Let $f\in L_{loc}^1(M)$ and assume that its restriction to
$M\backslash \Phi$ extends as a holomorphic function $f'\in \mathcal{
O}(\Omega \cup \mathcal{ W}_1)$. Then $f$ is CR all over $M$.
\end{lemma}

\proof
It suffices to prove that $f$ is CR at every point of $\Phi$. Let
$q\in \Phi$ be arbitrary and let $\mathcal{ W}_q$ be a local wedge of
edge $M$ at $q$ which is contained in $\mathcal{ W}_1$. Without loss
of generality, we can assume that in coordinates
$(z,w)=(x+iy,u+iv)\in\C^m\times \C^d$ vanishing at $q$ with
$T_qM=\{v=0\}$, the generic submanifold $M$ is represented locally in
the polydisc $\Delta_n(\rho_1)$ by $v=\varphi(x,y,u)$ for some
$\mathcal{ C}^{2,\alpha}$-smooth $\R^d$-valued mapping $\varphi$
defined on the real cube on $\I_{2m+d}(\rho_1)$. First of all, we
construct a family of analytic discs half attached to $M$ whose
interior is contained in the local wedge $\mathcal{ W}_q\subset
\mathcal{ W}_1$.

\def\thelemma{11.7}\begin{lemma}
There exists a family of analytic
discs $A_s(\zeta)$, with $s\in \R^{2m+d-1}$, $\vert s \vert \leq
2\delta$ for some $\delta>0$, and $\zeta \in\overline{\Delta}$, which
is of class $\mathcal{ C}^{2,\alpha-0}$ with respect to all variables,
such that
\begin{itemize}
\item[{\bf (1)}]
$A_0(1)=q$.
\item[{\bf (2)}]
$A_s(\overline{\Delta})\subset \Delta_n(\rho_1)$. 
\item[{\bf (3)}]
$A_s(\Delta) \subset \mathcal{ W}_q\cap \Delta_n(\rho_1)$.
\item[{\bf (4)}]
$A_s(\partial^+\Delta)\subset M$.
\item[{\bf (5)}]
$A_s(i)\in M\backslash \Phi$ and $A_s(-i)\in M\backslash\Phi$
for all $s$. 
\item[{\bf (6)}]
The mapping
$[-2\delta,2\delta]^{2m+d-1}\times [-\pi/2,\pi/2]
\ni (s,\theta)\longmapsto 
A_s(e^{i\theta})\in M$ is an embedding onto 
a neighborhood of $q$ in $M$.
\item[{\bf (7)}]
There exists $\rho_2>0$ such that the image of
$[- \delta, \delta]^{ 2m+d-1} \times [-\pi /4,\pi /4]$ through 
this mapping contains $M \cap \Delta_n (\rho_2)$.
\end{itemize}
\end{lemma} 

\proof
Let $M^1$ be a $\mathcal{ C}^{2,\alpha}$-smooth maximally real
submanifold of $M$ passing through $q$ such that $M^1\cap \Phi$ is of
zero measure with respect to the Lebesgue measure of $M^1$. Let $t\in
\R^d$ and include $M^1$ in a parametrized family of maximally real
submanifolds $M_t^1$ which foliates a neighborhood of $q$ in $M$.
Starting with a family of analytic discs $A_{c,x,v}^1(\zeta)$ which
are half-attached to $M^1$ as constructed in Lemma~7.12 above, 
we first choose the rotation
parameter $v_0$ and a sufficiently small scaling factor $c_0$ in order
that $A_{c_0,0,v_0}^1(\pm i)$ does not belong to $\Phi$. In fact, this
can be done for almost every $(c_0,v_0)$, because the
mapping $(c,v)\mapsto A_{c,0,v}^1(\pm i)$ is of rank $n$ at
every point $(c,v)$ with $c\neq 0$ and $v\neq 0$. In addition, we
adjust the rotation parameter $v_0$ in order that the vector $Jv_0$
points inside a proper subcone of the cone which defines the wedge
$\mathcal{W}_q$. If the scaling parameter $c$ is sufficiently small,
this implies that $A_{c_0,0,v_0}^1(\Delta)$ is contained in $\mathcal{
W}_q\cap \Delta_n(\rho_1)$, as in Lemma~8.3 above. The translation
parameter $x$ runs in $\R^n$ and we may select a $(n-1)$-dimensional
parameter subspace $x'$ which is transversal in $M^1$ to the half
boundary $A_{c_0, 0,v_0}^1(\partial^+ \Delta)$. With such a choice,
there exists $\delta>0$ such that the mapping $[-2 \delta,
2\delta]^{n-1}\times [-\pi/2, \pi/2 ]\ni (x', \theta) \longmapsto
A_{c_0, x', v_0}^1(e^{i\theta})$ is a diffeomorphism onto a
neighborhood of $q$ in $M^1$. Finally, using the stability of Bishop's
equation under perturbations, we can deform this family of discs by
requiring that it is half attached to $M_t^1$, thus obtaining a family
$A_s(\zeta) := A_{c_0, x',v_0, t}^1( \zeta)$ with $s:= (x',t) \in\R^{
2m+d-1}$. Shrinking $\delta$ if necessary, we can check as
in the proof of Lemma~8.3 {\bf ($\mathbf{9_1}$)} that
condition {\bf (5)} holds. This completes the proof.
\endproof

Let now $f \in L_{ loc}^1 (M)$ and let $f' \in \mathcal{ O} ( \Omega
\cup \mathcal{ W }_1)$. Thanks to the foliation propery {\bf (6)} of
Lemma~11.7, it follows from Fubini's theorem that for almost every
translation parameter $s$, the mapping $e^{ i\theta} \mapsto f \left(
A_s (e^{ i\theta })\right)$ defines a $L^1$ function on $\partial^+
\Delta$. In addition, the restriction of the function $f' \in
\mathcal{ O}( \Omega\cup \mathcal{ W}_1)$ to the disc $A_s( \Delta)
\subset \mathcal{ W}_q \subset \mathcal{ W}_1$ yields a holomorphic
function $f' \left( A_s( \zeta) \right)$ in $\Delta$.

\def\thelemma{11.8}\begin{lemma}
For almost every $s$ with $\vert s \vert \leq 2 \delta$, the function
$f' \left( A_s (\zeta) \right)$ belongs to the Hardy space $H^1 (
\Delta )$.
\end{lemma}

\proof
Indeed, for almost every $s$, the restriction $f\left (A_s(e^{
i\theta})\right)$ belongs to $L^1(\partial^+ \Delta)$. We can also
assume that for almost every $s$, the intersection $\Phi \cap
A_s(\partial^+ \Delta)$ is of zero one-dimensional measure. By the
assumption of Lemma~11.6, the restriction of $f\circ A_s$ and of
$f'\circ A_s$ to $\partial^+\Delta \backslash \Phi$ coincide. Recall
that $\partial^-\Delta=\{\zeta\in \partial \Delta: \, {\rm Re}\, \zeta
\leq 0\}$. Since $A_s( \pm i)$ does not belong to $\Phi$ and since
$A_s \left(e^{i\theta}\right)$ belongs to $\mathcal{ W}_q$ for all
$\theta$ with $\pi/2 < \vert \theta \vert \leq \pi$, it follows that
$f \circ A_s \vert_{\partial^+ \Delta}$ and $f' \circ A_s \vert_{
\partial^-\Delta}$ (which is holomorphic in a neighborhood of
$\partial^-\Delta$ in $\C$) match together in a function which is
$L^1$ on $\partial \Delta$. Let us denote this function by $f_s$.
Furthermore, $f_s$ extends holomorphically to $\Delta$ as $f'\circ A_s
\vert_\Delta$. Consequently, $f'\circ A_s\vert_\Delta$ belongs to the
Hardy space $H^1( \Delta)$, which proves the lemma.
\endproof

Since we have now established that the boundary value of
$f'$ on $M\backslash \Phi$ along the family
of discs $A_s(\zeta)$ coincides with $f$, we can
now denote both functions by the same letter $f$. 

For $\varepsilon\geq 0$ small, let now $\chi_\varepsilon \left(s,e^{i
\theta} \right)$ be a $\mathcal{ C}^2$-smooth function on $[ -2\delta,
2\delta]\times \partial \Delta$ which equals $\varepsilon$ for $\vert
s \vert \leq \delta$ and for $\theta \in [-\pi/4,\pi/4]$ and which
equals $0$ if either $\pi/2 \leq \vert \theta \vert \leq \pi$ or
$\vert s \vert \geq 2\delta/3$. We may require in addition that $\vert
\vert \chi_\varepsilon \vert \vert_{\mathcal{ C}^2}\leq
\varepsilon$. We define a deformation $M^\varepsilon$ of $M$ compactly
supported in a neighborhood of $q$ by pushing $M$ inside $\mathcal{
W}_q$ along the family of discs $A_s(\zeta)$ as follows:
\def\theequation{11.9}\begin{equation}
M^\varepsilon:= \{A_s\left([1-\chi_\varepsilon (s,e^{i\theta})]\, 
e^{i\theta}\right): \, 
\vert \theta \vert \leq \pi/2, \ 
\vert s \vert \leq 2\delta\}.
\end{equation}
Notice that $M^\varepsilon$ coincides with $M$ outside a small
neighborhood of $q$. Then we have $\vert \vert M^\varepsilon - M \vert
\vert_{\mathcal{ C}^2}\leq C\cdot \varepsilon$, for some constant
$C>0$ which depends only on the $\mathcal{ C}^2$ norms of $A_s(\zeta)$
and of $\chi_\varepsilon(s,e^{i\theta})$. If the radius $\rho_2$ is as
in Property {\bf (7)} of Lemma~11.7 above, the deformation
$M^\varepsilon \cap \Delta_n(\rho_2)$ is entirely contained in
$\mathcal{ W}_q$ and since $f$ is holomorphic in $\mathcal{ W}_q$,
its restriction to $M^\varepsilon \cap \Delta_n (\rho_2)$ is obviously
CR. An illustration is provided in the right hand side of the
following figure.

\bigskip
\begin{center}
\input M-deformation-discs.pstex_t
\end{center}

As in~\cite{j5}, \cite{p1}, \cite{mp1}, we notice that for every $s$
and every $\varepsilon$, the one-dimensional Lebesgue measure on the
arc 
\def\theequation{11.10}\begin{equation}
\Gamma_{\varepsilon,s}:=\{[1-\chi_\varepsilon
(s,e^{i\theta})]e^{i\theta}\in \Delta: \vert
\theta \vert \leq \pi\}
\end{equation}
is a Carleson measure. Thanks to the geometric uniformity of these
arcs $\Gamma_{\varepsilon,s}$, it follows from an inspection of the
proof of L.~Carleson's imbedding theorem that there exists a (uniform)
constant $C$ such that for all $s$ with $\vert s \vert \leq 2\delta$
and all $\varepsilon$, one has the estimate
\def\theequation{11.11}\begin{equation}
\int_{\Gamma_{\varepsilon,s}}\, 
\left\vert f\left(
A_s\left(
\left[
1-\chi_\varepsilon(s,e^{i\theta})
\right]\, e^{i\theta}
\right)
\right)\right\vert \cdot d\theta
\,
\leq \, C 
\int_{\partial\Delta}\, 
\vert f \vert \cdot d\theta.
\end{equation} 

We are now ready to complete the proof of Lemma~11.6. Let
$\pi_{x,y,u}$ denote the projection parallel to the $v$-space onto the
$(x, y,u)$-space. The mapping $(s, \theta) \mapsto \pi_{x, y,u }\left(
A_s ( \theta) \right)$ may be used to define new coordinates in a
neighborhood of the origin in $\C^m\times \R^d$, an open subset above
which $M$ and $M^\varepsilon$ are graphed. We shall now work with
these coordinates. With respect to the coordinates $(s,\theta)$, on
$M$ and on $M^\varepsilon$, we have formal adjoints ${}^T \overline{
L}$ and ${}^T( \overline{L }^\varepsilon)$ of the basis of CR vector
fields with an estimation of the form $\left\vert \left\vert {}^T(
\overline{ L }^\varepsilon)- {}^T \overline{ L} \right \vert \right
\vert_{ \mathcal{ C}^1}\leq C \cdot \varepsilon$, for some constant
$C>0$. Let now $\psi= \psi (s,\theta)$ be $\mathcal{ C}^1$-smooth
function with compact support in the set $\{\vert s \vert < \delta, \,
\vert \theta \vert \leq \pi/4\}$. By construction, the subpart of
$M^\varepsilon$ defined by $\widetilde{ M }^\varepsilon:= \{A_s
\left([1-\chi_\varepsilon (s,e^{ i \theta })]\, e^{ i\theta} \right):
\, \vert \theta \vert \leq \pi/4, \ \vert s \vert \leq \delta\}$ is
contained in the wedge $\mathcal{ W}_q$, hence the restriction of the
holomorphic function $f\in \mathcal{ W}_q$ to $\widetilde{
M}^\varepsilon$ is obviously CR on $\widetilde{ M}^\varepsilon$.

For simplicity of notation, we shall denote $f \left(A_s(e^{ i\theta})
\right)$ by $f_s(\theta)$ and $f\left( A_s \left( \left[ 1-
\chi_\varepsilon(s, e^{i\theta }) \right]\, e^{ i\theta} \right)
\right)$ by $f_s^\varepsilon ( \theta)$. Since by construction
for every $\varepsilon >0$, the $L^1$ function $(s,\theta)\mapsto
f_s^\varepsilon (\theta)$ is annihilated in the distributional sense
by the CR vector fields $\overline{ L}^\varepsilon$ on $\widetilde{
M}^\varepsilon$, we may compute (not writing the arguments $(s,
\theta)$ of $\psi$)
\def\theequation{11.12}\begin{equation}
{\small
\left\{
\aligned
{}
&
\left\vert
\int_{\vert s \vert \leq \delta}\,\int_{\vert
\theta \vert \leq \pi/4}\, 
f_s(\theta)\cdot 
{}^T\overline{L}(\psi) \cdot dsd\theta \right\vert=\\
& \
=\left\vert
\int_{\vert s \vert \leq \delta}\,\int_{\vert
\theta \vert \leq \pi/4}\, 
\left[
f_s(\theta)
\cdot {}^T\overline{L}(\psi)-
f_s^\varepsilon(\theta)
\cdot {}^T(\overline{L}^\varepsilon)(\psi)
\right] \cdot dsd\theta \right\vert \\
& \
\leq 
\left\vert
\int_{\vert s \vert \leq \delta} \left(
\int_{\vert \theta \vert \leq \pi/4}\,
\left[
f_s(\theta)\cdot {}^T\overline{L}(\psi)-f_s(\theta)\cdot 
{}^T(\overline{L}^\varepsilon)(\psi)+ \right.\right.\right.\\
& \ \ \ \ \ \ \ \ \ \ \ \ \ \ \ \ \ \ \ \ \ \ \ \ \ \ \ \ \ \ \
\ \ \ \ \ \ \
\left.\left.\left.
+f_s(\theta)\cdot {}^T(
\overline{L}^\varepsilon)(\psi)-
f_s^\varepsilon(\theta)
\cdot {}^T(\overline{L}^\varepsilon)(\psi)
\right] \cdot d\theta \right)\cdot ds \right\vert
\\
& \
\leq 
C_1(\psi) \cdot \varepsilon \cdot 
\int_{\vert s \vert \leq \delta}\,\int_{\vert
\theta \vert \leq \pi/4}\, 
\vert f_s(\theta) \vert \cdot dsd\theta + \\
& \ \ \ \ \ \ \ \ \ \ \ \ \ \ \ \ \ \ \ \ \ \ \ \ \ \ \ \ 
+C_2(\psi) \cdot
\int_{\vert s \vert \leq \delta}\,\int_{\vert
\theta \vert \leq \pi/4}\,
\vert f_s(\theta) - f_s^\varepsilon(\theta) \vert \cdot dsd\theta
\\
& \
\leq C_1(\psi,f,\delta)\cdot \varepsilon+
C_2(\psi,\delta)\cdot \max_{\vert s \vert \leq \delta}\, \int_{\vert
\theta \vert \leq \pi/4} \, 
\vert f_s(\theta) - f_s^\varepsilon(\theta) \vert \cdot dsd\theta.
\endaligned\right. 
}
\end{equation}
However, thanks to the estimate~\thetag{11.11} and thanks to
Lebesgue's dominated convergence theorem, the last integral tends to
zero as $\varepsilon$ tends to zero. It follows that the integral in
the first line of~\thetag{11.12} can be made arbitrarily small, hence
it vanishes. This proves that $f$ is CR in a neighborhood
of $q$ and completes the proof of Lemma~11.6.
\endproof

The proof of Lemma~11.3 is complete.
\endproof

\section*{\S12.~Proofs of Theorem~1.1 and of Theorem~1.3}

\subsection*{12.1.~Tree of separatrices linking hyperbolic points}
Let $M\subset \C^2$ be a globally minimal $\mathcal{
C}^{2,\alpha}$-smooth hypersurface, let $S\subset M$ be a $\mathcal{
C}^{2,\alpha}$-smooth open surface (without boundary) and let
$K\subset S$ be a proper compact subset of $S$. Assume that $S$ is
totally real outside a discrete subset of complex tangencies which are
hyperbolic in the sense of E.~Bishop. Since we aim to remove the
compact subset $K$ of $S$, we can shrink the open surface
$S$ around $K$ in order that $S$ contains {\it 
only finitely many}\, such hyperbolic
complex tangencies, which we shall
denote by $\{h_1,\dots,h_\lambda\}$, where
$\lambda$ is some integer, possibly zero. Furthermore, we can
assume that $\partial S$ is of class $\mathcal{ C}^{2,\alpha}$. As a
corollary of the qualitative theory of planar vector fields, due to
H.~Poincar\'e and I.~Bendixson, we know that
\begin{itemize}
\item[{\bf (i)}]
The hyperbolic points $h_1, \dots, h_\lambda$ are singularities of the
characteristic foliation $\mathcal{ F}_S^c$.
\item[{\bf (ii)}]
Incoming to every hyperbolic point $h_1,\dots,h_\lambda$, there
are exactly four $\mathcal{ C}^{2,\alpha}$-smooth open {\sl
separatrices} (to be defined precisely below).
\item[{\bf (iii)}]
After perturbing slightly the boundary $\partial S$ if necessary,
these separatrices are all transversal to $\partial S$ and the union of
all separatrices together with all hyperbolic points makes a {\sl
finite tree without cycles in $S$} (to be defined below).
\end{itemize} 

Precisely, by an (open) {\sl separatrix}, we mean a $\mathcal{
C}^{2,\alpha}$-smooth curve $\tau: (0,1) \to S$ with $\frac{ d
\tau}{ds}(s) \in T_{\tau(s)}S \cap T_{\tau(s)}^c M\backslash \{0\}$
for every $s\in (0,1)$, namely its tangent vectors are all nonzero and
characteristic, such that one limit point, say $\lim_{s\to 0} \,
\tau(s)$ is a hyperbolic point, and the other $\lim_{s\to 1} \,
\tau(s)$ either belong to the boundary $\partial S$ or is a second
hyperbolic point.

From the local study of saddle phase diagrams ({\it cf.}~\cite{
ha}), we get in addition:

\begin{itemize}
\item[{\bf (iv)}] 
There exists $\varepsilon >0$ and for every $l = 1, \dots, \lambda$,
there exist two curves $\gamma_l^1, \gamma_l^2 : (- \varepsilon,
\varepsilon ) \to S$ which are of class $\mathcal{ C }^{ 1, \alpha}$,
{\it not more}, with $\gamma_l^i (0) = h_l$ and $\frac{ d \gamma_l^1
}{ dt} (s) \in T_{ \gamma_l (s)} S \cap T_{ \gamma_l(s) }^c M
\backslash \{0\}$ for every $s \in (- \varepsilon, \varepsilon)$ and
for $i= 1,2$, such that the four open segments $\gamma_l^1 (
-\varepsilon,0)$, $\gamma_l^1 ( 0, \varepsilon )$, $\gamma_l^2 (
-\varepsilon, 0)$ and $\gamma_l^2 ( 0, \varepsilon )$ cover the four
pieces of open separatrices incoming at $h_l$.
\end{itemize}

Let $\tau_1,\dots,\tau_\mu: (0, 1) \to S$ denote all the separatrices of
$S$, where $\mu$ is some integer, possibly equal to zero. By the
{\sl finite hyperbolic tree $T_S$ of $S$}, we mean:
\def\theequation{12.2}\begin{equation}
T_S := \{h_1,\dots,h_\lambda\} 
\bigcup_{1\leq k\leq \mu} \ \tau_k(0,1).
\end{equation} 
We say that $T_S$ has {\sl no cycle} if it does not contain any subset
homeomorphic to the unit circle. For instance, in the case where
$S\equiv D$ is diffeomorphic to a real disc (as in the assumptions of
Theorem~1.1), its hyperbolic tree $T_D$ necessarily has no
cycle. However, in the case where $S$ is an annulus (for instance),
there is a trivial example of a characteristic foliation with two
hyperbolic points and a circle in the hyperbolic tree.

\subsection*{12.3.~Hyperbolic decomposition in the disc case}
Let the real disc $D$ and the compact subset $K \subset D$ be as in
Theorem~1.1. As in \S12.1 just above, we shrink $D$ slightly and
smooth out its boundary, so that its hyperbolic tree
$T_D$ is finite and has
no cycle. We may decompose $D$ as the disjoint union
\def\theequation{12.4}\begin{equation}
D= T_D \cup D_o,
\end{equation}
where the complement of the hyperbolic tree $D_o: = D \backslash T_D$
is an open subset of $D$ entirely contained in the totally real part
of $D$. Then $D_o$ has finitely many connected
components $D_1,\dots,D_\nu$, the {\sl hyperbolic sectors of $D$}.
Then, for $j=1,\dots,\nu$, we define the proper closed subsets $C_j:=
D_j \cap K$ of $D_j$ as illustrated in the left hand side of the
following figure.

\bigskip
\begin{center}
\input hyperbolic.pstex_t
\end{center}

Again from H.~Poincar\'e and I.~Bendixson's 
theory, we know that for every
component $D_j$ (in which the characteristic foliation is
nonsingular), the proper closed subset $C_j$ satisfies the nontransversality
condition $\mathcal{ F}_{ D_j}^c\{ C_j\}$ formulated in Theorem~1.2.
In {\sc Figure~20} just above, we have drawn the characteristic curves
only for the two sectors $D_4$ and $D_6$. One may observe that
$\mathcal{ F}_{D_4}\{ C_4\}$ and $\mathcal{ F}_{D_6}\{ C_6\}$ hold
true. Also, $K\cap T_D$ is a proper closed subset of the hyperbolic
tree of $D$.

\subsection*{12.5.~Global minimality of some complements}
Before proceeding to the deduction of Theorems~1.1 and~1.3 from
Theorem~1.2, we must verify that the complement $M\backslash K$ is
also globally minimal. Here, we state a generalization of Lemma~3.5 to
the case where some hyperbolic complex tangencies are allowed. Its
proof is not immediate.

\def\thelemma{12.6}\begin{lemma}
Let $M$ be a $\mathcal{ C}^{2,\alpha}$-smooth hypersurface in $\C^2$
and let $S\subset M$ be $\mathcal{ C}^{2,\alpha}$-smooth surface which
is totally real outside a discrete subset of hyperbolic complex
tangencies. Assume that the hyperbolic tree $T_S$ of $S$ has no cycle.
Then for every compact subset $K\subset S$ and for
an arbitrary point $p\in
M\backslash K$, its CR orbit in $M\backslash K$ coincides with its CR
orbit in $M$, minus $K$, namely
\def\theequation{12.7}\begin{equation}
\mathcal{ O}_{CR}(M\backslash K, p)= 
\mathcal{ O}_{CR}(M, p) \backslash K.
\end{equation}
\end{lemma}

\proof
Of course, we may assume that $S$ coincides with the shrinking of a
slightly larger surface and has finitely many hyperbolic points
$\{h_1,\dots,h_\lambda\}$, as described in \S12.1, with the same
notation. Let $K_{T_S}:= K \cap T_S$ be the track of $K$ on the
hyperbolic tree $T_S$. Since the intersection of $K_{T_S}$ with any
open separatrix may in general coincide
with any arbitrary closed subset of an
interval, in order to fix ideas, it will be convenient to deal with an
enlargement $\overline{K}$ of $K_{T_S}$, simply defined by
filling the possible holes of $K_{T_S}$ in $T_S$: more
precisely, $\overline{K}$
should contain all hyperbolic points together with all separatrices
joining them and for every separatrix $\tau_k(0,1)$ with right limit
point $\lim_{s\to 1} \, \tau_k(s)$ belonging to the boundary of $S$,
we require that $\overline{K}$ contains the segment $\tau_k[0, r_1]$,
where $r_1<1$ is close enough to $1$ in order that $\overline{ K}$
effectively contains $K_{T_S}$. 

Obviously, from the inclusions
\def\theequation{12.8}\begin{equation}
K_{T_S} \subset \overline{ K} \subset K,
\end{equation}
we deduce that for every point $p\in M\backslash K$, we have
the reverse inclusions
\def\theequation{12.9}\begin{equation}
\mathcal{ O}_{CR}(M\backslash K,p) \subset
\mathcal{ O}_{CR}(M\backslash \overline{ K},p) \subset
\mathcal{ O}_{CR}(M\backslash K_{T_S},p).
\end{equation}
The main step in the proof of Lemma~12.6 will be to establish the
following two assertions, implying the third, desired
assertion, already stated as~\thetag{12.7}.

\begin{itemize}
\item[{\bf (A1)}]
For every point $q\in M\backslash \overline{ K}$, we have
$\mathcal{ O}_{CR}(M\backslash \overline{ K}, q)= 
\mathcal{ O}_{CR}(M, q) \backslash \overline{K}$.
\item[{\bf (A2)}]
For every point $r\in M\backslash K_{T_S}$, we have $\mathcal{
O}_{CR}(M\backslash K_{T_S}, r)= \mathcal{ O}_{CR}(M, r) \backslash
K_{T_S}$.
\end{itemize}

Indeed, taking these two assertions for granted, let us conclude the
proof of Lemma~12.6. Let $p\in M\backslash K$ and decompose $K$ as a
disjoint union $K= K_{T_S}\cup C'$, where $C':= K\backslash K_{T_S}$
is a relatively closed subset of the hypersurface $M':= M\backslash
K_{T_S}$. Notice that $C'$ is contained in the totally real 
part of $S$.
Again thanks to foliation theory, we see that
the assumption that $K_{T_S}$ does not contain any cycle
entails that $C'$ does not contain maximal characteristic lines of the
totally real surface $S\backslash K_{T_S}$. Consequently, all the
assumptions of Lemma~3.5 are satisfied, hence by applying it to $p$, we
deduce that $\mathcal{ O}_{CR}(M' \backslash C', p) = \mathcal{
O}_{CR}(M', p) \backslash C'$. By developing in length this identity
between sets, we get
\def\theequation{12.10}\begin{equation}
\aligned
\mathcal{ O}_{CR}(M\backslash K, p)
& \ 
= 
\mathcal{ O}_{CR}\left((M\backslash K_{T_S})
\backslash C', p\right) \\
& \
=
\left[\mathcal{ O}_{CR}(M\backslash K_{T_S},p)
\right]\backslash C' \\
& \ 
= 
\left[\mathcal{ O}_{CR}(M, p)\backslash K_{T_S}
\right]\backslash C' \\
& \ 
=
\mathcal{ O}_{CR}(M,p) \backslash K,
\endaligned
\end{equation}
where, for the passage from the second to the third line, we
use {\bf (A2)}. This is~\thetag{ 12.7}, as desired. 

Thus, the (main) remaining task is to establish the 
assertion {\bf (A2)}, with {\bf (A1)} being a preliminary step.

First of all, we show how to deduce {\bf (A2)} from {\bf (A1)}. Pick
an arbitrary point $r \in M\backslash K_{ T_S}$. Of course, we have
the trivial inclusion $\mathcal{ O}_{ CR}(M \backslash K_{ T_S},r)
\subset \mathcal{ O}_{CR} (M, r) \backslash K_{ T_S}$ and we want an
equality. As $\overline{ K}$ contains $K_{ T_S}$, we have either $r
\in M \backslash \overline{ K}$ (first case) or $r \in \overline{ K}
\backslash K_{T_S}$ (second case), {\it see}\, {\sc Figure~21} just
below, where $r$ is located in $\overline{ K}\backslash K_{T_S}$.

\bigskip
\begin{center}
\input kts.pstex_t
\end{center}

Since the second case makes it impossible to apply {\bf (A1)}, we need
to find another point $r'$ in the CR orbit of $r$ in $M\backslash
K_{T_S}$ such that $r'$ belongs to $M\backslash \overline{ K}$. This
is elementary. We make a dichotomy: either $r=h_l$ is a hyperbolic
point or it belongs to an open separatrix $\tau_k(0,1)$. If $r=h_l\in
M\backslash K_{T_S}$ is a hyperbolic point, we may use one of the two
$\mathcal{ C}^{2,\alpha}$-smooth complex tangent curves $\gamma_l^1$
or $\gamma_l^2$ passing through $h_l$ and running in $M\backslash
K_{T_S}$ to join $r$ with another point which belongs to an open
separatrix and which obviously lies in the same CR orbit $\mathcal{
O}_{CR }( M \backslash K_{T_S}, r)$. Hence, we may assume that $r\in
\overline{ K} \backslash K_{T_S}$ belongs to an open separatrix. Since
$S$ is now maximally real near $r$, we may choose a $T^cM$-tangent
vector field $Y$ defined in a neighborhood of $r$ which is transversal
to $S$ at $r$. Then for all $\delta>0$ small enough, the point $r':
=\exp (\delta Y)(r)$ is outside $S$, hence does not belong to
$\overline{ K}$ and clearly lies in the same CR orbit $\mathcal{ O
}_{CR }( M \backslash K_{ T_S}, r)$.

In summary, when $r \in \overline{ K} \backslash K_{ T_S}$, we have
exhibited a point $r'\in \mathcal{ O }_{CR}(M \backslash K_{T_S}, r)$
with $r '\in M \backslash \overline{ K}$ so that it suffices now to
show that for every point $r\in M \backslash \overline{ K}$, we have
$\mathcal{ O}_{ CR}(M \backslash K_{T_S}, r)= \mathcal{ O}_{ CR} (M,r)
\backslash K_{ T_S}$.

Using a trivial inclusion and
applying {\bf (A1)}, we deduce that
\def\theequation{12.11}\begin{equation}
\mathcal{ O}_{CR}(M\backslash K_{T_S}, r) \supset
\mathcal{ O}_{CR}(M\backslash \overline{ K}, r) = 
\mathcal{ O}_{CR}(M, r) \backslash \overline{ K}.
\end{equation}
Unfortunately, there may well exist points of $\mathcal{ O}_{CR}(M,
r)$ belonging to $\overline{ K} \backslash K_{T_S}$, so that it
remains to show that {\it every}\, point $r'\in \mathcal{ O}_{CR}(M,
r)\cap [ \overline{ K} \backslash K_{T_S}]$ also belong to $\mathcal{
O}_{ CR}(M \backslash K_{T_S}, r)$. Again, this last step is
elementary and totally analogous to the above argument: we first claim
that we can join such a point $r'$ to a point $r'''\in M\backslash
\overline{ K}$ by means of a piecewise smooth CR curve running in
$M\backslash K_{T_S}$. Indeed, if $r'=h_l$ is a hyperbolic point, we
may first use one of the two $\mathcal{ C}^{2, \alpha}$-smooth complex
tangent curves $\gamma_l^1$ or $\gamma_l^2$ passing through $h_l$ and
running in $M \backslash K_{ T_S}$ to join $r'$ with another nearby
point $r''$ which belongs to an open separatrix. If $r'$ already
belongs to an open separatrix, we simply set $r'':= r'$. Since $S$ is
now maximally real near $r''$, we may choose a $T^cM$-tangent vector
field $Y$ defined in a neighborhood of $r''$ which is transversal to
$S$ at $r''$. Then for all $\delta>0$ small enough, the point
$r''':=\exp( \delta Y )(r'')$ satisfies $r''' \not \in S$, whence
$r''' \not \in \overline{ K}$. Of course, by choosing $r''$ and $r'''$
sufficiently close to $r'$, it follows that the piecewise smooth CR
curve joining them does not meet $K_{T_S}$. We deduce that $r'\in
\mathcal{O}_{CR} (M\backslash K_{ T_S}, r''')$.

Since by assumption $r' \in \mathcal{ O}_{CR}(M, r)$, it follows that
$r'''\in \mathcal{ O}_{CR}(M,r)$ and then $r'''\in \mathcal{ O}_{
CR}(M,r)\backslash \overline{ K}$. By means of the
supclusion~\thetag{ 12.11}, we deduce that $r'''\in \mathcal{
O}_{CR}(M\backslash K_{T_S},r)$.

Finally, from the two relations $r'''\in \mathcal{ O}_{CR}(M\backslash
K_{T_S},r)$ and $r'\in \mathcal{O}_{CR} (M\backslash K_{ T_S}, r''')$,
we conclude immediately that $r'\in \mathcal{ O}_{CR}(M\backslash
K_{T_S}, r)$. We have thus shown that the supclusion in~\thetag{
12.11} is an equality.

This completes the deduction of {\bf (A2)}
from {\bf (A1)}.

It remains now to establish {\bf (A1)}. We remind that in 
Section~3, we derived Lemma~3.5 from Lemma~3.7.
By means of a totally similar argument, which we shall
not repeat, one deduces {\bf (A1)} from the following
assertion. Remind that as $M$ is a hypersurface in 
$\C^2$, its CR orbits are of dimension either $2$ or
$3$.

\def\thelemma{12.12}\begin{lemma}
Let $M$, $S$, $T_S$ and $\overline{ K} \subset T_S$ be as above.
There exists a connected submanifold $\Omega$ embedded in $M$
containing the hyperbolic tree $T_S$ such that
\begin{itemize}
\item[{\bf (1)}]
$\Omega$ is a $T^c M$-integral manifold, namely $T_p^cM \subset
T_p \Omega$ for all $p \in \Omega$.
\item[{\bf (2)}]
$\Omega$ is contained in a single CR orbit of $M$.
\item[{\bf (3)}]
$\Omega\backslash \overline{ K}$ is also 
contained in a single CR orbit of $M \backslash
\overline{ K}$.
\end{itemize}
More precisely, $\Omega$ is an open neighborhood of $T_S$ if it is of
real dimension $3$ and a complex curve surrounding $T_S$ if it is of
dimension $2$.
\end{lemma}

\proof
We shall construct $\Omega$ by means of a flowing procedure, 
starting from a local piece of it. We start
locally in a neighborhood of a fixed point
$p_0 \in \overline{ K} \backslash \{ h_1, \dots, h_\lambda \}$,
whose precise
choice does not matter. Since $S$ is totally real in a neighborhood of
$p_0$, there exists a locally defined $T^cM$-tangent vector field $Y$
which is transversal to $S$ at $p_0$. Consequently, for $\delta >0$
small enough, the small segment $I_0:= 
\{ \exp(sY) (p_0): \, - \delta < s < \delta\}$ is transversal to 
$S$ at $p_0$ and moreover,
the
two half-segments
\def\theequation{12.13}\begin{equation}
I_0^\pm:= \left\{
\exp(s Y)(p_0): \, 
0 < \pm s < \delta 
\right\}
\end{equation}
lie in $M\backslash S$. Since $p_0$ belongs to some $T^cM$-tangent
open separatrix $\tau_k(0,1)$, there exists a $\mathcal{ C}^{1,
\alpha}$-smooth vector field $X$ defined in a neighborhood of $p_0$ in
$M$ which is tangent to $S$ and whose integral curve passing through
$p_0$ is a piece of $\tau_k(0,1)$. Since $Y$ is transversal to $S$ at
$p_0$, it follows that the set $\omega_0 := \{\exp(s_2 X)(\exp ( s_1
Y) (p_0)): \, -\delta < s_1, \, s_2 < \delta\}$ is a well-defined
$\mathcal{ C}^{1, \alpha}$-smooth codimension one small submanifold
passing through $p_0$ which is transversal to $S$ at $p_0$. Clearly,
we even have $T_{p_0}\omega_0 = T_{ p_0}^cM$. Thanks to the fact that
the flow of $X$ stabilizes $S$, we see that the integral curves $s_2
\mapsto \exp(s_2 X)(\exp (s_1 Y) (p_0))$ are contained in $M
\backslash S$ for every starting point $\exp (s_1 Y) (p_0)$ in the
segment $I_0$ which does not lie in $S$, namely for all $s_1 \neq
0$. We deduce that the two open halves of $\omega_0$ defined by
\def\theequation{12.14}\begin{equation}
\omega_0^\pm:=\{
\exp(s_2
X)(\exp (s_1 Y)(p_0)): \, 
0 < \pm s_1 < \delta, \ 
-\delta < s_2 < \delta
\}
\end{equation}
are contained in a single CR orbit of $M \backslash \overline{ K}$.

To begin with, assume that the CR orbit in $M \backslash \overline{
K}$ the point $q_0^+:= \exp\left(\frac{ \delta}{2} Y \right) (p_0)$
which belong to $\omega_0^+$, as drawn in {\sc Figure~22} below,
is of real dimension $2$. Afterwards, we shall treat the case where
its CR orbit in $M\backslash \overline{ K}$ is of dimension $3$.

Since $M$ is a hypersurface in $\C^2$ and since we have just proved
that the CR orbit $\mathcal{ O}_{ CR}(M \backslash \overline{ K},
q_0^+)$ already contains the $2$-dimensional half piece $\omega_0^+$,
we deduce that $\omega_0^+$ is a piece of complex curve whose boundary
$\partial \omega_0^+$ is (by construction) contained in the separatrix
$\tau_k (0,1)$. Since $\tau_k (0,1)$ is an embedded segment, we may
suppose from the beginning that the vector field $X$ is defined in a
neighborhood of $\tau_k (0,1)$ in $M$. Using then the flow of $X$, we
may easily prolong the small piece $\omega_0^+$ to get a semi-local
$\mathcal{ C }^{1, \alpha}$-smooth submanifold $\omega_k^+$ stretched
along $\tau_k(0,1)$, which constitutes its boundary. Again, this
piece $\omega_k^+$ is (by construction) contained in the CR orbit of
$q_0^+$ in $M\backslash \overline{K}$. By the fundamental stability
property of CR orbits under flows, we deduce that $\omega_k^+$ is in
fact a piece of complex curve with boundary $\tau_k(0,1)$.

Remind that by definition of separatrices, the point $\tau_k(0)$ is
always a hyperbolic point. There is a dichotomy: either $\tau_k(1)$ is
also a hyperbolic point or it lies in $\partial S$. If $\tau_k(1)$ is
a hyperbolic point, then by the definition of $\overline{ K}$, the
complete boundary $\partial \omega_k^+ =\tau_k(0,1)$ is contained in
$\overline{ K}$, hence it may {\it not}\, be crossed by means of a CR
curve running in $M \backslash \overline{ K}$. Since the piece
$\omega_k^+$ will be flowed all around $T_S$, our filling $\overline{
K}$ of $K_{T_S}$ was motivated by the desire of simplifying the
geometric situation without having to discuss whether $K_{T_S}$
contains or does not contain the whole segment $\tau_k (0,1)$, for
each $k=1, \dots, \mu$.

Before studying the case where $\tau_k (1) \in \partial S$, 
let us analyze the local situation in a neighborhood of the
hyperbolic point $\tau_k(0)=:h_l$, for some $l$ with $1\leq l\leq
\lambda$. 

As a preliminary, in order to understand clearly the situation, 
let us assume that the two characteristic curves $\gamma_l^1$ and
$\gamma_l^2$ passing through $h_l$ are of class $\mathcal{
C}^{ 2,\alpha}$, an assumption which would be satisfied if
we had assumed that $M$ and $S$ are of class $\mathcal{ C}^{3,
\alpha}$. By a straightening of $\gamma_l^1$ and of $\gamma_l^2$, 
which induces a loss of one derivative,
we easily show that there exist two linearly independent $\mathcal{
C}^{ 1,\alpha}$-smooth $T^cM$-tangent vector fields $X_1$ and $X_2$
whose integral curves issued from $h_l$ coincide with $\gamma_l^1$ and
$\gamma_l^2$. After possibly renumbering and reversing $X_1$ and $X_2$
and also reparametrizing $\gamma_l^1$ and
$\gamma_l^2$, we may assume that
$\gamma_l^1(s) =\exp(sX_1) (h_l)$ and that 
$\gamma_l^2(s) =\exp(sX_2) (h_l)$, for all small
$s>0$. Furthermore, we may assume that the direction of $X_2$ in a
neighborhood of $h_l$ is the same as the direction from 
$\partial \omega_k^+$ to $\omega_k^+$.

\bigskip
\begin{center}
\input tree-orbit.pstex_t
\end{center}

As in {\sc Figure~22} just above, let $\tau_j(0,1)$ be the separatrix
issued from $h_l$ in the positive direction of $X_2$. 
We may assume that $\tau_j (0) = h_l$.
Thanks to the
flow of the vector field $X_2$, we may now propagate the piece of
complex curve $\omega_k^+$ by stretching it along $\tau_j(0,1)$ in a
neighborhood of $h_l$. Using then the flow of a semi-locally defined
complex tangent vector field defined in a neighborhood of
$\tau_j(0,1)$, we may extend this local piece as a complex curve
$\omega_j^+$ with boundary $\tau_j(0,1)$. Finally, $\omega_k^+$ and
$\omega_j^+$ glue together as a complex curve with boundary
$\tau_k(0,1) \cup \{h_l\} \cup \tau_j(0,1)$ and corner $h_l$.

However, by {\bf (iv)} above, $\gamma_l^1$ and $\gamma_l^2$ are only
of class $\mathcal{ C}^{1, \alpha}$. Examples for which this
regularity is optimal are easily found. Straightening them is again
possible, but the vector fields $X_1$ and $X_2$ would be of class
$\mathcal{ C}^\alpha$, and we would lose the uniqueness of their
integral curves as well as the regularity of their flow. Consequently,
to prove that $\omega_k^+$ propagates along the second separatrix, with
its boundary contained in it, we must proceed differently~:
the proof is longer and we need one more diagram.

In {\sc Figure~23} just below, we draw the sadlle-looking surface $S$
in the $3$-dimensional space $M$~; the horizontal plane passing
through $h_l$ is thought to be the complex tangent plane $T_{
h_l}^cM$.

\bigskip
\begin{center}
\input saddle-orbit.pstex_t
\end{center}

Let us introduce two $T^cM$-tangent vector fields $X_1$ and $X_2$
defined in a neighborhood of $h_l$ with $X_1 (h_l)$ directed along
$\tau_k$ in the sense of increasing $s$ and $X_2 (h_l)$ directed along
$\tau_j$ in the sense of increasing $s$. Let $Z$ denote the vector
field $X_1 +X_2$, as shown in the top of the left hand side of {\sc
Figure~23} above. Using the flow of $Z$ we can begin by extending the
banana-looking piece $\omega_k^+$ of complex curve by introducing the
submanifold $\omega$ consisting of points
\def\theequation{12.15}\begin{equation}
\exp(s_2 Z ) (\tau_k(s_1)),
\end{equation}
where $0 < s_1 < \delta$ and $0 < s_2 < \delta$, for some small
$\delta > 0$. One checks that all these points stay in $M \backslash
S$, hence are contained in the same CR orbit as
$\omega_k^+$ in $M \backslash \overline{ K}$. By the stability
property of CR orbits, it follows of course that $\omega$ is a piece
of complex curve contained in $M$.

For $0 < s < \delta$, let $\mu(s) := \exp (sZ) (h_l)$ denote the CR
curve lying ``between'' $\tau_k$ and $\tau_j$ and which
constitutes a part of the boundary of $\omega$. Let $p$ be an
arbitrary point of this curve, close
to $h_l$. 

\def\thelemma{12.16}\begin{lemma}
The integral curve $s\mapsto
\exp (-s X_1)(p)$ of $-X_1$ issued from $p$ necessarily
intersects $S$ at a point $q$ close to $h_l$ and close to $\tau_j$
{\rm (}{\it cf.} {\sc Figure~23}{\rm )}.
\end{lemma}

\proof
First of all, we need some preliminary.

Thanks to the existence of a ``$1/8$ piece'' of complex curve $\omega$
with $h_l\in \overline{ \omega}$ which is contained in the
hypersurface $M$, we see that $M$ is necessarily Levi-degenerate at
$h_l$.

Next, we introduce local holomorphic coordinates $(z, w) = (x+ iy,
u+iv)\in \C^2$ vanishing at $h_l$ in which the hypersurface $M$ is
given as the graph $v= \varphi (x, y, u)$, where $\varphi$ is a
$\mathcal{ C}^{ 2, \alpha}$-smooth function. Since $M$ is Levi
degenerate at $h_l$, we may assume that $\varphi (x, y, u) \vert \leq
C \cdot \left(\vert x \vert + \vert y \vert + \vert u \vert \right)^{
2+ \alpha}$. We may also assume that the surface $S$, as a subset of
$M$, is represented by one supplementary equation of the form $u = h
(x,y)$, where the $\mathcal{ C}^{ 2,\alpha}$-smooth
function $h$ satisfies
\def\theequation{12.17}\begin{equation}
\left\{
\aligned
h(x,y)
& \
= z \bar z+ \gamma \, (z^2+ \bar z^2)+ 
{\rm O}\left( \vert z \vert^{ 2+ \alpha} \right) \\
& \
=
(2\gamma + 1) \, x^2 - 
(2\gamma -1) \, y^2+ {\rm O} \left( \vert z \vert^{ 2+ \alpha}
\right),
\endaligned\right.
\end{equation}
and where $\gamma >\frac{ 1}{2}$ is E.~Bishop's invariant. Then the
tangents at $h_l$ to the two half-separatrices $\tau_k$ and $\tau_j$
are given respectively by the linear (in)equations $x>0$, $y = -
\frac{ 2\gamma + 1}{ 2\gamma -1} \, x$, $u=0$ and $x< 0$, $y = -
\frac{ 2\gamma + 1}{ 2\gamma -1} \, x$, $u=0$. In {\sc Figure~23},
where we do not draw the axes, the $u$-axis is vertical, the $y$ axis
points behind $h_l$ and the $x$-axis is horizontal, from left to
right.

Expressing the two $T^cM$-tangent vector fields $X_1$ and $X_2$ in the
(natural) real coordinates $(x,y,u)$ over $M$, we may write them as
\def\theequation{12.18}\begin{equation}
\left\{
\aligned
X_1
& \
 = 
\frac{ \partial}{ \partial x} -
\left( 
\frac{ 2\gamma + 1}{ 2\gamma -1}\right) \, 
\frac{ \partial }{ \partial y} + 
A_1 (x,y,u)\, 
\frac{ \partial }{ \partial u}, \\ 
X_2 
& \
= 
-
\frac{ \partial}{ \partial x} -
\left( 
\frac{ 2\gamma + 1}{ 2\gamma -1}\right) \, 
\frac{ \partial }{ \partial y} + 
A_2 (x,y,u)\, 
\frac{ \partial }{ \partial u}. 
\endaligned\right.
\end{equation}
Since $\varphi$ vanishes to second order at $h_l$, the two
$\mathcal{ C}^{ 1,\alpha}$-smooth
coefficients $A_1$ and $A_2$ satisfy an estimate
of the form
\def\theequation{12.19}\begin{equation}
\left\vert 
A_1, A_2 (x,y,u) 
\right\vert < C \cdot \left(
\vert x \vert + \vert y \vert + 
\vert u \vert
\right)^{ 1+ \alpha}.
\end{equation}

Now, we come back to the integral curve of Lemma~12.16. It is
contained in the real $2$-surface passing through $h_l$ defined by
\def\theequation{12.20}\begin{equation}
\Sigma:= \{\exp(-s_2 X_1)
(\exp(s_1 Z)(h_l)): \, -\delta < s_1, s_2 <
\delta\},
\end{equation}
for some $\delta >0$. Because the vector fields $X_1$, $X_2$ and $Z=
X_1 + X_2$ have $\mathcal{ C}^{ 1,\alpha}$-smooth coefficients, the
surface $\Sigma$ is only $\mathcal{ C}^{ 1,\alpha}$-smooth in
general. In $M$ equipped with the three real coordinates $(x,y,u)$, we
may parametrize $\Sigma$ by a mapping of the form
\def\theequation{12.21}\begin{equation}
(s_1, s_2) \longmapsto 
\left(
s_2 - 
2s_1 \left(
\frac{ 2\gamma +1}{ 2\gamma -1}
\right), \ s_2
\left(
\frac{ 2\gamma+ 1}{ 2\gamma -1}
\right), \ 
u( s_1,s_2)
\right), 
\end{equation}
where $u$ is of class $\mathcal{ C}^{ 1,\alpha}$. It is clear that
$u(0)= \partial_{ s_1} u(0) = \partial_{ s_2} u( 0) =0$, so that
there is a constant $C$ such that
\def\theequation{12.22}\begin{equation}
\left \vert u (s_1,\, s_2) 
\right \vert < C \cdot
\left(
\vert s_1 \vert + \vert s_2 \vert
\right)^{ 1+ \alpha},
\end{equation}
since $u$ is of class $\mathcal{ C}^{ 1,\alpha}$. Furthermore, by
inspecting the flows appearing in~\thetag{ 12.20}, taking account of
the estimates~\thetag{ 12.19}, we claim that $u$ satisfies the
better estimate
\def\theequation{12.23}\begin{equation}
\vert u (s_1, s_2) \vert < C \cdot
\left(
\vert s_1 \vert + \vert s_2 \vert 
\right)^{ 2+ \alpha},
\end{equation}
for some constant $C>0$.
In other words, $\Sigma$ osculates the complex tangent plane $T_{
h_l}^cM$ to second order at $h_l$~: $\Sigma$ is more flat than $S$ at
$h_l$. One may check the estimate~\thetag{ 12.23} is sufficient to
establish Lemma~12.16, because the second jet of the saddle function
$h(x,y)$ does not vanish at the origin.

To prove the claim, we formulate the main argument as an independent
assertion. Mild modifications of this argument apply to our case, but
we shall not provide all the details.

Let $L_1 := \frac{ \partial }{ \partial x} + A_1 (x,y, u) \, \frac{
\partial }{ \partial u}$ and $L_2 := \frac{ \partial }{ \partial y} +
A_2 (x,y, u) \, \frac{ \partial}{ \partial u}$ be two vector fields
having $\mathcal{ C}^{ 1,\alpha }$-smooth coefficients
satisfying~\thetag{ 12.19}. Denote 
by $s_1 \longmapsto (s_1, \lambda (s_1), \mu(s_1))$
the integral curve of $L_1$ passing through the origin. It is
$\mathcal{ C}^{ 2,\alpha}$-smooth and
we have 
\def\theequation{12.24}\begin{equation}
\left \vert
\lambda (s_1)
\right \vert < C \cdot \vert s_1 \vert^{ 2+\alpha}
\ \ \ \ \ \ \
{\rm and} 
\ \ \ \ \ \ \
\left \vert
\mu (s_1)
\right \vert < C \cdot \vert s_1 \vert^{ 2+\alpha},
\end{equation}
for some constant $C>0$. Consider the composition of flows $\exp (s_2
L_2) (\exp (s_1 L_1)(0))$. We have to solve the system of ordinary
differential equations
\def\theequation{12.25}\begin{equation}
\frac{ dx}{ ds_2}= 0, 
\ \ \ \ \ \ \
\frac{ dy}{ ds_2} = 1, 
\ \ \ \ \ \ \
\frac{ du}{ ds_2}= A_2 (x,y,u)
\end{equation}
with initial conditions
\def\theequation{12.26}\begin{equation}
x(0)= s_1, 
\ \ \ \ \ \ \
y(0) = \lambda( s_1), 
\ \ \ \ \ \ \ 
u(0) = \mu(s_1).
\end{equation}
This yields $x(s_1,s_2) =s_1$, $y(s_1,s_2)= \lambda (s_1)+ s_2$ and
the integral equation
\def\theequation{12.27}\begin{equation}
u(s_1,s_2) = 
\mu(s_1)+ 
\int_0^{ s_2} \, 
A_2 \left(
s_1, s_2'+ \lambda (s_1), u(s_1, s_2')
\right) \, ds_2'.
\end{equation}
Since $u$ is at least $\mathcal{ C}^{ 1,\alpha}$-smooth and vanishes
to order $1$ at $(s_1,s_2)= (0,0)$, we know already that it
satisfies~\thetag{ 12.22}. Using~\thetag{ 12.19}, it is now
elementary to provide an upper estimate of the right hand side
of~\thetag{ 12.27} which yields the desired estimate~\thetag{ 12.23}.

The proof of Lemma~12.16 is complete.
\endproof

So, for various points $p= \mu(s)$ close to $h_l$ the intersection
points $q\in S$ exist. If all points $q$ belong to $\tau_j$, we are
done: the piece $\omega$ extends a $1/4$ piece of complex curve with
boundary $\tau_k \cup \tau_j$ near $h_l$ and corner $h_l$.

Assume therefore that one such point $q$ does not belong to $\tau_j$,
as drawn in the left hand side of {\sc Figure~23} above. Suppose that
$q$ lies above $\tau_j$, the case where $q$ lies under $\tau_j$ being
similar and in fact simpler. The characteristic curve $\gamma'\subset
S$ passing through $q$ stays above $\tau_j$ and is nonsingular.
Propagating the complex curve $\omega$ in $M\backslash \overline{ K}$
by means of the flow of $-X_1$, we deduce that there exists at $q$ a
local piece $\omega_q^+$ of complex curve with boundary contained in
$\gamma'$ which is contained in the same CR orbit as $\omega$. Using
then the flow of a CR vector having $\gamma'$ as an integral curve, we
can propagate $\omega_q^+$ along $\gamma'$, which yields a long thin
banana-looking complex curve with boundary in $\gamma'$. However, this
piece may remain too thin. Fortunately, thanks to the flow of $X_1 -
X_2$, we can extend it as a piece $\omega'$ of complex curve with
boundary $\gamma'$ which goes over $h_l$, with respect to a
complex projection onto $T_{h_l}^cM$, as illustrated in {\sc
Figure~23} above. We claim that this yields a contradiction.

Indeed, as $\omega$ and $\omega'$ are complex curves, they are locally
defined as graphs of holomorphic functions $g$ and $g'$ defined in
domains $D$ and $D'$ in the complex line $T_{h_l}^cM$. By
construction, there exists a point in $r \in D\cap D'$ at which the
values of $g$ and $g'$ are distinct. However, since by construction
$g$ and $g'$ coincide in a neighborhood of the CR curve joining $p$ to
$q$, they must coincide at $r$ because of the principle of analytic
continuation: this is a contradiction. In conclusion, the CR orbit
passes through the hyperbolic point $h_l$, in a neighborhood of which
it consists of a cornered complex curve with boundary $\tau_k\cup
\tau_j$.

We can now continue the proof. Since the hyperbolic tree $T_S$ does
not contain any cycle, by proceeding this way we claim that the small
piece of complex curve $\omega_0^+$ propagates all around $T_S$ and
matches up as a smooth complex curve $\Omega$ containing the
hyperbolic tree. Indeed, in the case where $\tau_k( 1)$ is not a
hyperbolic point, recall that we arranged at the beginning that
$\overline{ K}\cap \tau_k(0,1)= \tau_k(0,r_1]$, where $r_1<1$. It is
then crucial that when a limit point $\tau_k(1)$ belongs to $\partial
S$, we escape from $\overline{ K}$ and using a local CR vector field
$Y$ transversal to $S$, we may cross the separatrix $\tau_k(0,1)$ at a
point $\tau_k(r_2)$ where $r_2$ satisfies $r_1 < r_2 < 1$. Hence, we
pass to the other side of $S$ in $M$ and then, by means of a futher
flowing, we turn around to the other side of $\tau_k(0,1)$. Also, the
two pieces in either side of $\tau_k(0,1)$ match up at least
$\mathcal{ C}^{1,\alpha}$-smoothly. Then thanks to the stability
property of orbits under flows, we deduce that these two pieces match
up as a piece of complex curve containing $\tau_k(0,1)$ in its
interior.

We thus construct the complex curve $\Omega$ surrounding $T_S$, which
is obviously contained in a single CR orbit of $M$. Also, by
construction, $\Omega\backslash \overline{ K}$ is contained in a
single CR orbit of $M\backslash \overline{ K}$. Thus, we have
established Lemma~3.12 under the assumption that the CR orbit of
$q_0^+$ is two-dimensional.

Assume finally that the CR orbit of $q_0^+$ is $3$-dimensional. By a
similar propagation procedure, we easily construct a neighborhood
$\Omega$ in $M$ of the hyperbolic tree satisfying conditions {\bf
(1)}, {\bf (2)} and {\bf (3)} of Lemma~12.12. This complete its proof.
\endproof

The proof of Lemma~12.6 is complete.
\endproof

\subsection*{12.28.~Proofs of Theorems~1.1 and 1.3}
We can now prove Theorem~1.1. In fact, we shall directly prove the
more general version stated as Theorem~1.3 which implies Theorem~1.1
as a corollary, thanks to the geometric observations of \S12.3.

First of all, we notice that as $M\backslash K$ is globally minimal,
there exists a wedge attached to $M\backslash K$ to which continuous
CR functions on $M\backslash K$ extend holomorphically. Hence, the
CR-removability of $K$ is a consequence of its $\mathcal{
W}$-removability. Also, Lemma~11.3 shows that the 
$L^{\sf p}$-removability of $K$ is a consequence of
its $\mathcal{ W}$-removability. Consequently, it suffices
to establish that $K$ is $\mathcal{ W}$-removable in 
Theorem~1.3.

Let $T_S$ be the hyperbolic tree of (a suitable shrinking of) $S$,
which contains no cycle by assumption. Let $\omega_1$ be a one-sided
neighborhood of $M \backslash K$ in $\C^2$. Because the
nontransversality condition $\mathcal{ F }_{ S \backslash T_S}^c \{ K
\cap (S\backslash T_S) \}$ holds true by assumption, we may apply
Theorem~1.2 to the totally real surface $S \backslash T_S$ in the
globally minimal (thanks to Lemma~12.6) hypersurface $M \backslash K_{
T_S}$ to remove the proper closed subset $K\cap (S \backslash
T_S)$. We deduce that there exists a one-sided neighborhood $\omega_2$
of $M\backslash K_{ T_S}$ in $\C^2$ such that (after shrinking
$\omega_1$ if necessary), holomorphic functions in $\omega_1$ extend
holomorphically to $\omega_2$. Then we slightly deform $M$ inside
$\omega_2$ over points of $K\cap (S \backslash T_S)$. We obtain a
$\mathcal{ C}^{2, \alpha}$-smooth hypersurface $M^d$ with
$M^d\backslash K_{ T_S} \subset \omega_2$. Also, by stability of
global minimality under small perturbations, we can assume that $M^d$
is also globally minimal. By construction, we obtain holomorphic
functions in the neighborhood $\omega_2$ of $M^d \backslash K_{T_S}$
in $\C^2$.

Since $M$ and $M^d$ are of codimension $1$, the union of a one-sided
neighborhood $\omega^d$ of $M^d$ in $\C^2$ together with $\omega_2$
constitutes a one-sided neighborhood of $M$ in $\C^2$. To conclude the
proof of Theorem~1.1, it suffices therefore to show that the closed
set $K_{T_S}$ is $\mathcal{ W}$-removable.
The reader is referred to {\sc Figure~20} above for
an illustration.
 
Reasoning by contradiction (as for the proof of Theorem~1.2'), let
$K_{\rm nr}\subset K_{T_S}$ denote the smallest nonremovable subset of
$K_{T_S}$. If $K_{\rm nr}$ is empty, we are done, gratuitously.
Assume therefore that $K_{\rm nr}$ is nonempty. Let $T'$ be a
connected component of the minimal subtree of $T$ containing $K_{\rm
nr}$. By a {\sl subtree} of a tree $T$ defined as in~\thetag{12.2}
above, we mean of course a finite union of some of the separatrices
$\tau_k(0,1)$ together with all hyperbolic points which are endpoints
of separatrices. Since $T'$ does not contain any subset homeomorphic
to the unit circle, there exists at least one extremal branch of $T'$,
say $\tau_1(0,1)$ after renumbering, with $\tau_1(1)\in \partial
S$. To reach a contradiction, we shall show that at least one point of
the nonempty set $K_{\rm nr} \cap T'$ is in fact $\mathcal{
W}$-removable.

If the subtree $T'$ consists of the single branch $\tau_1 (0,1)$
together with the single elliptic point $\tau_1 (0)$, thanks to
properties {\bf (iii)} and {\bf (iv)} of \S12.1, we can enlarge a
little bit this branch by prolongating the curve $\tau_1 (0,1)$ to an
open $\mathcal{ C}^{ 2, \alpha}$-smooth Jordan arc $\tau_1
(-\varepsilon, 1+ \varepsilon)$, for some $\varepsilon>0$, with the
appendix $\tau_1(1,1+\varepsilon)$ outside $S$ (in the slightly larger
surface containing $S$). Then we remind that by a special case of
Theorem~4 (ii) of~\cite{mp1}, every proper closed subset of
$\tau_1(-\varepsilon, 1+\varepsilon)$ is $\mathcal{ W}$-removable. It
follows that the proper closed subset $K_{\rm nr}$ of the Jordan arc
$\tau_1(-\varepsilon, 1+\varepsilon)$ is removable, which yields the
desired contradiction in the case where $T'$ consists of a single
branch together with a single elliptic point.

If $T'$ consists of at least two branches, again with $\tau_1(1)\in
\partial S$, then applying Theorem~4 (ii) of~\cite{mp1}, we may at
least deduce that $K_{\rm nr} \cap \tau_1(0,1)$ is $\mathcal{
W}$-removable, since $K_{\rm nr}\cap \tau_1(0,1)$ is contained in
$\tau_1(0, r_1]$, for some $r_1<1$. But possibly, this set $K_{\rm nr}
\cap \tau_1(0,1)$ could be empty.

However, we claim that it is
nonempty. Indeed, otherwise, if $K_{ \rm nr} \cap \tau_1[0,1)$
consists of the single point $\tau_1(0)$, which is a hyperbolic point,
then $K_{ \rm nr}$ is in fact contained in the smaller subtree $T''$
defined by $T'' := T' \backslash \tau_1 (0,1)$ (here, we use that
$\tau_1 (0)$ is a hyperbolic point, hence there exists another branch
$\tau_k (0,1)$ with $\tau_k (1)= \tau_1 (0)$ or $\tau_k (0)=
\tau_1(0)$). This contradicts the assumption that $T'$ is the minimal
subtree containing $K_{\rm nr}$. Then $K_{\rm nr} \cap \tau_1 (0,1)$
is nonempty and removable, which contradict the assumption that
$K_{\rm nr}$ is the smallest nonremovable subset of $K$.

The proofs of Theorems~1.1 and~1.3 are complete.
\endproof

\section*{\S13.~Applications to the edge of the wedge theorem}

In this section we formulate three versions of the edge of the wedge
theorem for holomorphic and meromorphic functions, two of which are
based upon an application of our removable singularities theorems.
Let us begin with some definition.

\subsection*{13.1.~Preliminary}
Let $E$ be a generic submanifold of $\C^n$, which may be maximally
real. By a {\sl double wedge attached to $E$}, we mean a pair
$(\mathcal{ W}_1, \mathcal{ W }_2)$ of disjoint wedges attached to $E$
which admit a nowhere vanishing continuous vector field $v: E
\rightarrow T \C^n\vert_E / T E$ such that $ Jv ( p)$ points into
$\mathcal{ W }_1$ and $-Jv(p)$ into $\mathcal{ W }_2$, for every $p \in
E$.

In the case where $E=\R^n$, the classical edge of the wedge theorem
states that there exists a neighborhood $\mathcal{ D}$ of $E$ in
$\C^n$ such that every function which is continuous on 
$\mathcal{ W}_1 \cup E \cup 
\mathcal{ W }_2$ and holomorphic in $\mathcal{
W}_1 \cup \mathcal{ W }_2$ extends holomorphically to $\mathcal{ D}$.
Also, generalizations are known in the case where $f|_{ \mathcal{
W}_1}$ and $f|_{\mathcal{ W }_2}$ have coinciding distributional
boundary values on $E$.

The assumption about the matching up of boundary values along $E$ from
$\mathcal{ W}_1$ and from $\mathcal{ W}_1$ is really needed, even if
the two boundary values coincide on a thick subset of $E$. To support
this observation, consider the following elementary example: the
complex hyperplane $H := \{z_n=0\} \subset \C^n$ and the maximally
real plane $E:= \R^n \subset \C^n$ intersect transversally in the
$(n-1)$-dimensional totally real plane $C:= \{y=0, \ x_n =0\}$; the
pair of wedges $\mathcal{ W}_1 := \{y_1 >0, \dots, y_n >0\}$ and
$\mathcal{ W }_2:= \{y_1<0, \dots, y_n <0\}$ clearly form a double wedge
attached to $E$; the function $\exp(-1 /z_n)$ restricted to the two
wedges is holomorphic there, has coinciding boundary values on the
thick set $E \backslash C$, but does not extend holomorphically to a
neighborhood of $E$ in $\C^n$. Evidently, the envelope of holomorphy
of the union of $\mathcal{ W}_1 \cup \mathcal{ W}_2$ together with a
thin neighborhood of $E \backslash C$ in $\C^n$ does
{\it not}\, contain any neighborhood of $E$ in $\C^n$.

Thus, in order to apply our removability theorems (which are
essentially statements about envelopes of holomorphy), the first
question is how to impose coincidence of boundary values on the
edge. We shall first see that the fact that $C\subset E$ is exactly of
codimension one in the above example is the limiting case for the
obstruction to holomorphic extension. Since we want to treat
also meromorphic extension, let us remind some definitions.

\subsection*{13.2.~Meromorphic functions and envelopes}
Let $U$ be a domain in $\C^n$. A meromorphic function $f \in \mathcal{
M} ( U)$ is a collection of equivalence classes of quotients of
locally defined holomorphic functions. It defines a $P_1 (\C)$-valued
function, which is single-valued only on some Zariski dense open
subset $D_f \subset U$. More geometrically, we may represent $f$ by
the closure $\Gamma_f$ of its graph $f|_{ D_f}$ over $D_f$, which
always constitutes an irreducible $n$-dimensional complex analytic
subset of $U \times P_1 (\C)$ with surjective, almost everywhere
biholomorphic projection onto $U$ (equivalent definition). It is well
known that the {\sl indeterminacy set} of $f$, namely the set of $z
\in U$ over which the whole fiber $\{ z\} \times P_1 (\C)$ is
contained in $\Gamma_f$, is an analytic subset of $U$ of codimension
at least 2. It is the only set where $f$ is multivalued.

We shall constantly apply a theorem due to P.~Thullen
(generalized by S.~Ivashkovitch in~\cite{ i} in the
context of K\"ahler manifolds) according to which 
the envelope of holomorphy of a domain in $\C^n$
coincides with its envelope of meromorphy. As holomorphic
functions are meromorphic, we shall state
Lemma~13.4, Corollary~13.8 and Corollary~13.11 below
directly for meromorphic functions.

\subsection*{13.3.~Edge of the wedge theorem over a maximally 
real edge} Let $E \subset \C^n$ $(n \geq 2)$ be a real analytic
maximally real submanifold, let $( \mathcal{ W }_1, \mathcal{
W}_2)$ be a double wedge attached to $E$ and let $f_1$, $f_2$ be two
meromorphic functions in $\mathcal{ W}_1$, $\mathcal{ W}_2$.

We need an assumption which tames the behaviour of their indeterminacy
sets, as one tends towards the edge $E$ from either $\mathcal{ W}_1$
or $\mathcal{ W}_2$. It will be sufficient to impose a matching up of
their boundary values on the complement of a closed subset $C$ whose
$(n-1)$-dimensional Hausdorff measure vanishes. Let ${\sf H}_d$ denote
the $d$-dimensional Hausdorff measure.

\def\thelemma{13.4}\begin{lemma} If there is a closed subset $C\subset
E$ with ${\sf H }_{ n-1} (C) =0$ such that both $\overline{
\Gamma_{f_1 }}\cap [(E \backslash C) \times P_1(\C)]$ and $\overline{
\Gamma_{ f_2 }} \cap [( E \backslash C)\times P_1(\C)]$ coincide with
the graph of a {\rm (}single{\rm )} continuous mapping from $E
\backslash C$ to $ P_1 (\C)$, then there exists a neighborhood
$\mathcal{ D}$ of $E$ in $\C^n$ which depends only on $(\mathcal{ W}_1, 
\mathcal{ W}_2)$ and a meromorphic function 
\def\theequation{13.5}\begin{equation}
f\in \mathcal{ M}\left(\mathcal{ D} \cup
\mathcal{ W}_1 \cup \mathcal{ W}_2\right), 
\end{equation}
extending the $f_j$, namely such that $f|_{ \mathcal{ W }_j}
=f_j$, for $j= 1,2$.
\end{lemma}

\proof
First of all, the assumption of continuous coincidence of boundary
values enables us to apply the classical edge of the wedge theorem at
each point of $E \backslash C$. This yields a neighborhood $\mathcal{
D }_0$ of $E \backslash C$ in $\C^n$ and a meromorphic extension $f_0
\in \mathcal{ M}( \mathcal{ D}_0 \cup \mathcal{ W}_1 \cup \mathcal{
W}_2)$. We claim that the envelope of meromorphy of
$\mathcal{ D}_0 \cup \mathcal{ W}_1 \cup \mathcal{ W}_2$ contains a
neighborhood $\mathcal{ D}_1$ of $E$ in $\C^n$.

Indeed, this follows from a very elementary application of the
continuity principle. Let $p\in C$ be arbitrary. After a local
straightening, we may insure that $p$ is the origin, that $E=\R^n$, that
$\mathcal{ W}_1$ contains $\{y_1 >0, \dots, y_n >0 \}$ and that
$\mathcal{ W}_2$ contains $\{y_1 < 0, \dots, y_n < 0\}$.

Let us introduce the trivial family of analytic discs
\def\theequation{13.6}\begin{equation}
A_{c,x,v}( \zeta) :=
\left(x_1+c(1+v_1)\zeta,\dots,
x_n+(1+v_n)\zeta\right), 
\end{equation}
where $c >0$ is a sufficiently small fixed scaling factor, where $x
\in \R^n$ is a small translation parameter and where $v\in \R^n$ is a
small pivoting parameter. Clearly, $A_{c, x,v}( \partial^+ \Delta)$ is
contained in $\mathcal{ W }_1$ and $A_{c, x,v} (\partial^- \Delta)$ is
contained in $\mathcal{ W }_2$. However $A_{c, x,v}( \pm 1)$ may
encounter $C$.

First of all, using the submersiveness of the two mappings $v \mapsto
A_{c,0,v}(\pm 1)\in E$, we may find $v_0$ arbitrarily 
close to the origin in $\R^n$ such that $A_{c,0,v_0}(\pm
1)$ does not belong to $C$. It follows that for all small
translation vectors $q\in \C^n$, the disc boundary
$A_{c,0,v_0}(\partial \Delta)+q$ is contained in the domain $\mathcal{
D}_0\cup \mathcal{ W}_1 \cup \mathcal{ W}_2$.

Furthermore, because $C$ is of Hausdorff $(n-1)$-dimensional measure
zero, for almost all $x\in \R^n$, the segment $A_{c,x,v_0}([-1,1])$
does not meet $C$. It follows that for such $x$, the disc
$A_{c,x,v_0}(\overline{ \Delta})$ is contained in the domain
$\mathcal{ D}_0\cup \mathcal{ W}_1 \cup \mathcal{ W}_2$. We deduce
that every disc $A_{c,0,v_0}(\partial \Delta)+q$ is analytically
isotopic to a point in $\mathcal{ D}_0\cup \mathcal{ W}_1 \cup
\mathcal{ W}_2$. An application of the continuity principle yields
meromorphic extension to a neighborhood of 
$p=A_{c,0,0}(0)\in C$.

In sum, we have constructed a neighborhood $\mathcal{ D}_1$ of $E$ in
$\C^n$ and a meromorphic extension $f\in \mathcal{ M}( \mathcal{
D}_1)$. But $\mathcal{ D}_1$ is not independent of $(f_1,f_2)$, since
it depends on $C$. Fortunately, once we know meromorphic extension to
a neighborhood $\mathcal{ D}_1$ of $E$ in $\C^n$, we may reemploy the
analytic disc technique of the classical edge of the wedge theorem to
describe a neighborhood $\mathcal{ D}$ of $E$ in $\C^n$ which depends
only on $(\mathcal{ W}_1, \mathcal{ W}_2)$ ({\it see}\, the end of the
proof of Corollary~13.8 below for more arguments). This completes the
proof of Lemma~13.4.
\endproof

\subsection*{13.7.~Edge of the wedge theorem over an edge of positive
CR dimension} Let $M$ be a $\mathcal{ C}^{2,\alpha}$-smooth generic
submanifold of $\C^n$ of positive CR dimension and let $C$ be a proper
closed subset of $M$ such that $M$ and $M\backslash C$ are globally
minimal. In~\cite{ mp3}, Theorem~1.1, it was shown as a main theorem
that every such closed subset $C$ of $M$ is $CR$-, $\mathcal{ W}$- and
$L^{\sf p}$-removable. We may formulate the following application,
where, for simplicity, we assume local minimality at every point.

\def\thecorollary{13.8}\begin{corollary}
Let $E \subset \C^n$ $( n \geq 2)$ be a generic manifold of class
$\mathcal{ C }^{ 2, \alpha}$ of positive CR dimension which is locally
minimal at every point, let $(\mathcal{ W}_1, \mathcal{ W}_2)$ a
double wedge attached to $E$ and let two meromorphic functions $f_j\in
\mathcal{ M} ( \mathcal{ W}_j)$ for $j=1,2$. If there is a closed
subset $C\subset E$ with ${\sf H }_{ n-1} (C) =0$ such that both
$\overline{ \Gamma_{f_1 }} \cap [(E \backslash C) \times P_1( \C)]$
and $\overline{ \Gamma_{ f_2 }} \cap [( E \backslash C)\times P_1(
\C)]$ coincide with the graph of a {\rm (}single{\rm )} continuous
mapping from $E \backslash C$ to $ P_1 (\C)$, then there exists a
neighborhood $\mathcal{ D }$ of $E$ in $\C^n$ which depends only on
$(\mathcal{ W}_1, \mathcal{ W}_2)$ and a meromorphic function
\def\theequation{13.9}\begin{equation}
f\in \mathcal{ M} \left(\mathcal{ W}_1 \cup
\mathcal{ D} \cup \mathcal{ W}_2\right), 
\end{equation}
extending the $f_j$, namely such that $f|_{ \mathcal{ W }_j}
=f_j$, for $j = 1,2$.
\end{corollary}

\proof
Applying the classical edge of the wedge theorem, we get a meromorphic
extension $f_0 \in \mathcal{ M}( \mathcal{ W}_1 \cup \mathcal{ D}_0
\cup \mathcal{ W}_2)$, where $\mathcal{ D}_0$ is some open
neighborhood of $E \backslash C$ in $\C^n$. Next, we include $E$ in a
CR manifold $M$ with $M \subset \mathcal{ W}_1 \cup E \cup \mathcal{
W}_2$ and $\dim_\R M = 1 + \dim_\R E$, as shown in the following
figure.

\bigskip
\begin{center}
\input double-wedge.pstex_t
\end{center}

Of course, the domain $\mathcal{ W}_1 \cup \mathcal{ D}_0 \cup
\mathcal{ W}_2$ constitues a (rather thick) wedge attached to
$M\backslash C$. Since $E$ is locally minimal at every point and
since one CR tangential direction of $M$ is transversal to $E$, both
$M$ and $M \backslash \Sigma$ are both globally minimal. Applying
Theorem~1.1 in~\cite{ mp2}, we deduce that there exists a wedge
$\mathcal{ W}$ attached to $M$ which is contained in the envelope of
meromorphy of $\mathcal{ W}_1 \cup \mathcal{ D}_0 \cup \mathcal{
W}_2$. We then claim that there exists a neighborhood $\mathcal{ D}$
of $p$ in $\C^n$, which depends only on $(\mathcal{ W}_1, \mathcal{
W}_2)$ such that $\mathcal{ D}$ is contained in the envelope of
meromorphy of $\mathcal{ W}_1 \cup \mathcal{ W} \cup \mathcal{ W}_2$.

Indeed, by deforming slightly $M$ inside $\mathcal{ W}$ near $E$, we
get a $\mathcal{ C}^{2, \alpha}$-smooth generic submanifold $M^d
\subset \mathcal{ W}_1 \cup \mathcal{ W} \cup \mathcal{ W}_2$. Instead
of functions meromorphic in the disconnected open set $\mathcal{ W}_1
\cup \mathcal{ W}_2$, we now consider meromorphic functions in the
{\it connected}\, open set $\mathcal{ W}_1 \cup \mathcal{ W} \cup
\mathcal{ W}_2$, which is a neighborhood of $M^d$ in $\C^n$. Then by
following the proof of the edge of the wedge theorem given in~\cite{a}
and applying the continuity principle, one deduces meromorphic
extension to a neighborhood $\mathcal{ D }^d$ in $\C^n$ of the
deformed submanifold $M^d$. Since the size of $\mathcal{ W }_1$ and
the size of $\mathcal{ W }_2$ are uniform with respect to $d$, the
size of the domain $\mathcal{ D }^d$ is also uniform with respect to
$d$, as follows from the stability of the edge of the wedge theorem
established in~\cite{ a}, since it relies on E.~Bishop's equation.
Hence for $M^d$ sufficiently close to $M$, the domain $\mathcal{ D
}^d$ contains a neighborhood of $p$ in $\C^n$. This completes the
proof of Corollary~13.8.
\endproof

\subsection*{13.10.~Hartogs-Bochner phenomenon and edge of the wedge 
theorem} Next we turn to the question whether a Hartogs-Bochner
phenomenon holds in presence of a double wedge. More
precisely we ask when it is sufficient to require coincidence of
boundary values only outside some compact $K\subset E$. Let us first
look at a prototypical case where the answer is particularly neat,
thanks to Theorem~1.1. Obviously, the proof is totally similar to the
proof of Corollary~13.8 and will not be repeated.

\def\thecorollary{13.11}\begin{corollary}
Let $E \subset \C^2$ be an embedded real analytic totally real disc,
let $(\mathcal{ W }_1, \mathcal{ W }_2)$ a double wedge attached to
$E$ and let two meromorphic functions $f_j\in \mathcal{ M} ( \mathcal{
W}_j)$ for $j=1,2$. If there is a compact subset $K\subset E$ such
that both $\overline{ \Gamma_{f_1 }} \cap [(E \backslash K) \times
P_1( \C)]$ and $\overline{ \Gamma_{ f_2 }} \cap [( E \backslash
K)\times P_1( \C)]$ coincide with the graph of a {\rm (}single{\rm )}
continuous mapping from $E \backslash K$ to $ P_1 (\C)$, then there
exists a neighborhood $\mathcal{ D}$ of $E$ in $\C^2$ which depends
only on $(\mathcal{ W}_1,\mathcal{ W}_2)$ and a meromorphic function
\def\theequation{13.12}\begin{equation}
f\in \mathcal{ M}\left(\mathcal{ W}_1 \cup
\mathcal{ D} \cup \mathcal{ W}_2\right), 
\end{equation}
extending the $f_j$, namely such that $f|_{ \mathcal{ W }_j}
=f_j$, for $j = 1,2$.
\end{corollary}

In order to find the most general application of Theorems~1.2 and
1.2', we first remark that in {\sc Figure~24}, we have a considerable
freedom in the choice of the generic submanifolds $M$ of CR dimension
$1$ with $E\subset M \subset \mathcal{ W}_1 \cup E \cup \mathcal{ W}_2$,
depending on the aperture of $\mathcal{ W }_1$ and $\mathcal{ W
}_2$. Since $\dim_\R M = 1+ \dim_\R E$, the tangent space
to such an $M$ at
a point $p\in E$ is uniquely determined by some nonzero vector $v_p\in
T_p \C^n / T_p E$. In order that $M$ is locally contained in
$\mathcal{ W}_1 \cup E \cup \mathcal{ W}_2$, it is necessary and
sufficient that either $Jv$ points in $\mathcal{ W}_1$ and $-Jv$
points in $\mathcal{ W}_2$, or vice versa, depending on the
orientations of $\mathcal{ W}_1$ and $\mathcal{ W}_2$ with respect to
$E$. Without loss of generality, after a possible shrinking, we can
therefore assume that the cones of $\mathcal{ W}_1$ and of $\mathcal{
W}_2$ are exactly opposite to each other at every point of $E$;
indeed, it would be impossible to construct an $M$ locally contained
in $\mathcal{ W}_1 \cup E \cup \mathcal{ W}_2$ which satisfies $T_p M
= T_p E \oplus \R v_p$ at a point $p\in E$, in the case where the
vector $Jv_p$ points in the cone at $p$ of $\mathcal{ W}_1$ but $-J
v_p$ lies outside the cone at $p$ of $\mathcal{ W}_2$, or
vice versa.

Assuming $\mathcal{ W}_2$ to be opposite to $\mathcal{ W}_1$, let us
define an induced field of open cones $p\mapsto {\sf C }_p^{ \mathcal{
W }_1, \mathcal{ W}_2 }$ as follows: a nonzero vector $v_p \in T_p E
\backslash \{0\}$ belongs to ${\sf C }_p^{ \mathcal{ W }_1, \mathcal{
W}_2 }$ if either $Jv$ or $-Jv$ points into $\mathcal{ W }_1$. A
nowhere vanishing vector field $p \mapsto v (p)$ is said to be {\sl
directed by ${\sf C }_p^{ \mathcal{ W }_1, \mathcal{ W}_2 }$} if $v(p)
\in {\sf C}_p^{ \mathcal{ W }_1, \mathcal{ W}_2 }$ for every $p \in
E$. Clearly, for every vector field $p\mapsto v(p)$ directed by $
p\mapsto {\sf C }_p^{ \mathcal{ W }_1, \mathcal{ W}_2 }$, we may
construct a $\mathcal{ C}^{2,\alpha}$-smooth semi-local generic
submanifold $M$ containing $E$, contained in $\mathcal{ W}_1 \cup E
\cup \mathcal{ W}_2$ which satisfies $T_p M = T_p E \oplus J v(p)$ at
every point $p\in E$.

In the statement of Theorem~1.2', we defined the condition $\mathcal{
F}_{M^1}^c\{C\}$ with respect to some CR manifold $M$ containing the
totally real manifold $M^1$. But $M$ entered in the definition only
via the characteristic foliation induced on $M^1$. Hence its
reasonable to define a more general nontransversality property, by
replacing the characteristic foliation by the foliation induced by any
vector field directed by the field of cones $p\mapsto {\sf
C}_p^{\mathcal{ W}_1, \mathcal{ W}_2}$, as follows:

\smallskip
\begin{itemize}
\item[
$\mathcal{ F}_{ \mathcal{ W}_1, \mathcal{ W}_2 } \{C\}:$] For every
closed subset $C' \subset C$ there is a smooth vector field $ p
\mapsto v(p)$ directed by the field
of cones $p\mapsto {\sf C }_p^{ \mathcal{ W}_1, \mathcal{ W
}_2}$ such that there exists a simple $\mathcal{ C}^{2,
\alpha}$-smooth curve $\gamma': [-1,1] \to E$ whose range
$\gamma'([-1,1])$ is contained in a single integral curve of $p
\mapsto v (p)$ with $\gamma'(-1) \not \in C'$, $\gamma' (0) \in C'$
and $\gamma' (1) \not \in C'$, there exists a local $(n -
1)$-dimensional transversal $R \subset E$ to $\gamma'$ passing through
$\gamma' (0)$ and there exists a thin open neighborhood $V$ of
$\gamma' ( [-1, 1])$ in $E$ such that if $\pi: V \to R$ denotes the
semi-local projection parallel to the flow lines of $v$, then
$\gamma'(0)$ lies on the boundary, relatively to the topology of $R$,
of $\pi (C' \cap V)$.
\end{itemize}
\smallskip

The proper application of Theorem~1.2' is the following.
Its proof follows by a direct examination of the
proof of Theorem~1.2'.

\def\thecorollary{13.13}\begin{corollary}
Let $E \subset \C^n$ $(n \geq 2)$ be a real analytic maximally real
submanifold, let $(\mathcal{ W }_1, \mathcal{ W }_2)$ a double wedge
attached to $E$. Let $C$ be a proper closed subset of $E$ satisfying
the nontransversality property $\mathcal{ F}_{ \mathcal{ W }_1,
\mathcal{ W}_2} \{C\}$ above. Let two meromorphic functions $f_j\in
\mathcal{ M} ( \mathcal{ W}_j)$ for $j=1,2$ such that both $\overline{
\Gamma_{f_1 }} \cap [(E \backslash C) \times P_1( \C)]$ and
$\overline{ \Gamma_{ f_2 }} \cap [( E \backslash C)\times P_1( \C)]$
coincide with the graph of a {\rm (}single{\rm )} continuous mapping
from $E \backslash C$ to $ P_1 (\C)$, Then there exists a neighborhood
$\mathcal{ D}$ of $E$ in $\C^2$ which depends
only on $(\mathcal{ W}_1,\mathcal{ W}_2)$ and a meromorphic function
\def\theequation{13.14}\begin{equation}
f\in \mathcal{ M}\left(\mathcal{ W}_1 \cup
\mathcal{ D} \cup \mathcal{ W}_2\right), 
\end{equation}
extending the $f_j$, namely such that $f|_{ \mathcal{ W }_j} =f_j$,
for $j = 1,2$. 
\end{corollary}

\subsection*{13.15.~Further applications}
We now conclude this section by suggesting two applications of
Theorem~1.2' in higher dimensions, in the case where $M^1$ is not
everywhere totally real. However, we must mention that we did not try
to generalize the results of E.~Bishop to understand the local
geometry of complex tangencies of generic submanifolds of CR dimension
$1$ in $\C^n$, for $n\geq 3$. Consequently, our formulations
should be considered as mild generalizations of
Theorem~1.1 and~1.3.

Thus, let $M$ be a $\mathcal{ C }^{2, \alpha}$-smooth generic
submanifold of $\C^n$ {\rm (}$n \geq 2${\rm )} of CR dimension $1$,
let $M^1$ be a codimension one submanifold of $M$ which is maximally
real except at every point of some proper closed subset $E \subset
M^1$. Let $C$ be a proper closed subset of $M^1$. For simplicity, we
assume that $M$ is locally minimal at every point, an assumption which
insures that for every closed subset $\widetilde{ C}$ of $M$, 
both $M$ and $M \backslash \widetilde{ C}$ are globally
minimal.

Firstly, applying Theorem~1.2' to remove
$C \cap (M^1 \backslash E)$ and then Theorem~1.1 of~\cite{ mp2}
to remove $C \cap E$, we deduce the following. 

\def\thecorollary{13.16}\begin{corollary}
Assume that $E$ is of vanishing $(n - 1)$-dimensional Hausdorff
content, and that the nontransversality condition $\mathcal{ F}_{M^1
\backslash E}^c \{C \cap (M^1 \backslash E)\}$ holds. Then $C$ is
CR-, $\mathcal{ W}$- and $L^{ \sf p}$-removable.
\end{corollary}

Secondly, we may generalize the notion of hyperbolic tree and assume
that $E$ consists of finitely many compact submanifolds of codimension
$2$ in $M^1$ joined by a collection of finitely many codimension one
submanifolds of $M^1$ with boundaries in $E$ which are foliated by
characteristic curves. Under some easily found assumptions, one could
formulate a second corollary analogous to Theorem~1.3.

\section*{\S14.~An example of a nonremovable three-dimensional torus}

This final section is devoted to exhibit a crucial example of a closed
subset $C$ violating the main nontransversality condition $\mathcal{
F}_{M^1}^c\{ C\}$ of Theorem~1.2' such that $C$ is truly nonremovable.
In addition, we may require that $M$ and $M^1$ have the simplest
possible topology.

\def\thelemma{14.1}\begin{lemma}
There exists a triple $(M, M^1, C)$, where
\begin{itemize}
\item[{\bf (i)}]
$M$ is a $\mathcal{ C}^\infty$-smooth generic submanifold in
$\C^3$ of CR dimension $1$, diffeomorphic to a real $4$-ball{\rm ;}
\item[{\bf (ii)}]
$M^1$ is a $\mathcal{ C}^\infty$-smooth one-codimensional
submanifold of $M$ which is maximally real in $\C^n$ and diffeomorphic
to a real $3$-ball{\rm ;}
\item[{\bf (iii)}]
$C$ is a compact subset of $M^1$
diffeomorphic to a real three-dimensional torus which is everywhere
transversal to the characteristic foliation $\mathcal{ F}_{M^1}^c$, 
hence the nontransversality condition $\mathcal{ F}_{M^1}^c\{C\}$ of
Theorem~1.2' clearly does not hold{\rm ;}
\item[{\bf (iv)}]
$M$ of finite type $4$ in the sense of
T.~Bloom and I.~Graham at every point, hence globally minimal,
\end{itemize} 
such that $C$ is neither CR- nor 
$\mathcal{ W}$- nor $L^{\sf p}$-removable with respect to $M$.
\end{lemma}

By type 4 at a point $p\in M$, we mean of course that
the Lie brackets of the complex tangent bundle $T^cM$ up to length $4$
generate $T_pM$.

\subsection*{14.2.~The geometric recipe}
We first construct the $3$-torus $C$, then construct the maximally
real $M^1$ and finally define $M$ as a certain thickening of
$M^1$. The argument for insuring global minimality of $M$ involves
computations with Lie brackets and is postponed to the
end.

Firstly, in $\R^3 =\R^3\oplus i \{0\}\subset \C^3$ equipped with the
coordinates $(x_1,x_2,x_3)$, where $x_j={\rm Re}\, z_j$ for $j=1,2,3$,
pick the ``standard'' $2$-dimensional torus $T^2$ of Cartesian
equation
\def\theequation{14.3}\begin{equation}
\left(\sqrt{x_1^2+x_2^2}-2\right)^2+x_3^2=1.
\end{equation}
This torus is stable under the rotations directed by the $x_3$-axis;
its intersection with the $(x_1,x_3)$-plane consists of two circles of
radius $1$ centered at the points $x_1=2$ and $x_1=-2$; it bounds a
three-dimensional open ``full'' torus $T^3$; both $T^2$ and $T^3$ are
contained in the ball $B^3$ of radius $5$ centered at the origin.

It is better to drop the square root: one checks that the equations
of $T^2$ and $T^3$ are equally given by $T^2:=\{\rho =0\}$ and
$T^3:= \{\rho<0\}$, by means of the {\it polynomial}\, defining 
function
\def\theequation{14.4}\begin{equation}
\rho(x_1,x_2,x_3):=
(x_1^2+x_2^2+x_3^2+3)^2-16\, (x_1^2+x_2^2),
\end{equation}
which has nonvanishing differential at every point of $T^2$.
Consequently, the extrinsinc complexification of $T^2$, namely the
complex hypersurface defined by
\def\theequation{14.5}\begin{equation}
\Sigma:= \{(z_1,z_2,z_3)\in \C^3: \
\rho(z_1,z_2,z_3)=0\}
\end{equation}
cuts $\R^3$ along $T^2$ with the transversality property $T_x \R^3
\cap T_x\Sigma= T_x T^2$ for every point $x\in T^2$.

Secondly, according to G.~Reeb ({\it see}\, \cite{ cln}, pp.~25--27;
{\it see}\, also the figures there), by considering the space
$\R^3\equiv S^3 \backslash \{ \infty\}$ as a punctured
three-dimensional sphere $S^3$, one may glue a second
three-dimensional full torus $\widetilde{ T}^3$ to $T^3$ along $T^2$
with $\infty\in \widetilde{ T}^3$ and then construct a foliation of
$S^3$ by $2$-dimensional surfaces all of whose leaves, except one, are
diffeomorphic to $\R^2$, are contained in either $T^3$ or in
$\widetilde{ T}^3$ and are accumulating on $T^2$, and finally, whose
single compact leaf is the above $2$-torus $T^2$. This yields the
so-called {\it Reeb foliation} of $S^3$, which is $\mathcal{
C}^\infty$-smooth and orientable. Consequently, there exists a
$\mathcal{ C}^\infty$ smooth vector field $L= a_1(x)\, \partial_{x_1}+
a_2(x) \, \partial_{x_2} + a_3(x) \, \partial_{x_3}$ of norm $1$,
namely $a_1(x)^2+ a_2(x)^2+ a_3(x)^2=1$ for every $x\in \R^3$, which
is everywhere orthogonal (with respect to the standard Euclidean
structure) to the leaves of the Reeb foliation. Geometrically, the
integral curves of $L$ accumulate asymptotically on the two nodal
(central) circles of $T^3$ and of $\widetilde{ T}^3$.

The open ball $B^3\subset \R^3$ of radius $5$ centered at the origin
will be our maximally real submanifold $M^1$. The two-dimensionaly
torus $T^2$ will be our nonremovable closed set $C$. The integral
curves of the vector field $L$ will be our characteristic lines.
Since $L$ is orthogonal to $T^2$, these characteristic lines will of
course be everywhere transverse to $C$, so that $\mathcal{
F}_{M^1}^c\{C\}$ clearly does not hold.

Thirdly, it remains to construct the generic submanifold $M$ of
CR dimension $1$ containing $M^1$ and to check that $C$ will be
nonremovable.

First of all, we notice that $L$ provides the characteristic
directions of $M^1$ if and only if $T_x M= T_p\R^3\oplus \R \, J\,
L(x)$ for every point $x\in M^1\equiv B^3$. Consequently, all
submanifolds $M\subset \C^3$ obtained by slightly thickening $M^1$ in
the direction of $J\, L(x)$ will be convenient; in other words, only
the first jet of $M$ along $M^1$ is prescribed by our choice of the
characterisctic vector field $L$. Notice that all such thin strips
$M$ along $M^1$ will be diffeomorphic to a real $4$-ball.

The fact that $C$ is nonremovable for all such generic submanifolds
$M$ is now clear: the hypersurface $\Sigma=\{z\in \C^3 : \, \rho(z)
=0\}$ satisfying $T_x\Sigma= T_xT^2\oplus \R \, J\, T_xT^2$ for all
$x\in T^2$ and $L$ being transversal to $T^2$, we easily deduce the
transversality property $T_x \Sigma + T_x M= T_x\C^3$ for all $x\in
T^2$, a geometric property which insures that the holomorphic function
$1/\rho(z)$, which is CR on $M\backslash C$, does not extend
holomorphically to any wedge of edge $M$ at any point of $C$.
Intuitively, $T_x\Sigma/T_xM$ absorbs all the normal space 
$T_x\C^3/T_x M$ at every point $x\in T^2$, leaving no room
for any open cone.

Finally, to fulfill all the hypotheses of Theorem~1.2' (except of
course $\mathcal{ F}_{M^1}^c\{C\}$), we have to insure that $M$ is
globally minimal. We claim that by bending strongly the second and the
fourth order jet of $M$ along $M^1$ (without modifying the first order
jet which must be prescribed by $J\, L$), one may insure that $M$ is
of type 4 in the sense of T.~Bloom and I.~Graham at every point of
$M^1$; since being of finite type is an open property, it follows that
$M$ is finite type at every point provided that, as a strip, $M$ is
sufficiently thin along $M^1$. As is known, finite-typeness at every
point implies local minimality at every point which in turn implies
global minimality. This completes the recipe.

We would like to mention that by following a similar recipe, one may
construct an elementary example of a non-removable compact subset of a
generic submanifold of codimension one diffeomorphic to a $4$-ball
lying in a globally minimal hypersurface in $\C^3$ which is (also)
diffeomorphic to a $5$-ball ({\it cf.} \cite{ js}).

\subsection*{14.6.~Finite-typisation}
Thus, it remains to construct a generic submanifold $M\subset \C^3$ of
CR dimension $1$ satisfying $T_x M = T_xM^1 \oplus \R \, J \, L(x)$
for every $x\in M^1$, which is of {\it type $4$ at every point $x\in
M^1$}.

First of all, let us denote by $L=a_1(x)\, \partial_{x_1}+ a_2(x) \,
\partial_{x_2}+ a_3(x) \, \partial_{x_3}$ the unit vector field which
was constructed as a field orthogonal to the Reeb foliation: it is
defined over $\R^3$ and has $\mathcal{ C}^\infty$-smooth coefficients
satisfying $a_1(x)^2+ a_2(x)^2+ a_3(x)^2= 1$ for all $x\in\R^2$. The
two-dimensional quotient vector bundle $T\R^3 / (\R L)$ with
contractible base being necessarily trivial, it follows that we can
complete $L$ by two other $\mathcal{ C }^\infty$-smooth unit
vector fields $K^1$ and $K^2$ defined over $\R^3$ such that the triple
$(L(x), K^1(x), K^2(x))$ forms a direct orthonormal frame at every
point $x\in \R^3$. Let us denote the coefficients of $K^1$ and of
$K^2$ by
\def\theequation{14.7}\begin{equation}
\aligned
K^1 
& \
= \rho_1 \, \partial_{x_1}+ 
\rho_2\, \partial_{x_2} + 
\rho_3 \, \partial_{x_3},\\
K^2 
& \
=
r_1 \, \partial_{x_1}+ 
r_2 \, \partial_{x_2}+ 
r_3 \, \partial_{x_3},
\endaligned
\end{equation}
where $\rho_j$ and $r_j$ for $j=1,2,3$ are $\mathcal{
C}^\infty$-smooth functions of $x\in \R^3$ satisfying $\rho_1^2+
\rho_2^2+ \rho_3^2=1$ and $r_1^2+r_2^2+r_3^2=1$. In our case, $K^1$
and $K^2$ may even be constructed directly by means of a
trivialization of the bundle tangent to the Reeb foliation.

Let $P>0$ be a constant, which will be chosen later to be large. Since
by construction we have the two orthogonality relations $a_1 \rho_1 +
a_2 \rho_2 + a_3 \rho_3= 0$ and $a_1 r_1 + a_2 r_2 + a_3 r_3= 0$, it
follows that every generic submanifold $M_P\subset 
\C^3$ defined by the two
Cartesian equations
\def\theequation{14.8}\begin{equation}
\aligned
0 
& \
= \rho 
= 
y_1 \, \rho_1(x) + 
y_2 \, \rho_2(x) +
y_3 \, \rho_3(x) + 
P 
\left[
y_1^2+ y_2^2+y_3^2
\right], \\
0 
& \
= r 
= 
y_1\, r_1(x) + 
y_2 \, r_2(x) + 
y_3 \, r_3(x)+
P^3
\left[
y_1^4 + y_2^4 + y_3^4
\right]
\endaligned
\end{equation}
enjoys the property that the vector field $J L(x)= a_1(x) \,
\partial_{y_1}+ a_2(x) \, \partial_{x_2}+ a_3(x) \, \partial_{x_3}$ is
tangent to $M_P$ at every $x\in \R^3$. As desired, we deduce that
$T_x^c M = \R L(x) \oplus J \R L(x)$ for every $x\in \R^3$, a property
which insures that $\R L(x)$ is the characteristic direction of $M^1$
in $M_P$, independently of $P$.

To complete the final minimalization argument for the construction of
a nonremovable compact set $C:= T^2 \subset M^1\subset M$ which
appears in the Introduction, it suffices now to apply the following
lemma with $R= 5$. Though calculatory, its proof is totally
elementary.

\def\thelemma{14.9}\begin{lemma}
For every $R>0$, there exist $P>0$ sufficiently large
such that $M_P$ is of type $4$ at every point
$x\in \R^3$ with $x_1^2+ x_2^2+ x_3^2\leq R^2$. 
\end{lemma}

\proof
As above, let $M_P = \{z\in \C^3: \, \rho =r =0\}$. By writing the
tangency condition, one checks immediately that the one-dimensional
complex vector bundle $T^{1, 0}M_P$ is generated over $\C$ by the
vector field $\LL := A_1 \, \partial_{z_1}+ A_2 \, \partial_{ z_2}+
A_3\, \partial_{ z_3}$, with the explicit expressions
\def\theequation{14.10}\begin{equation}
\aligned
A_1:= 
& \ 
4\rho_{z_3}r_{z_2}- 4 \rho_{z_2} r_{z_3}, \\
A_2:= 
& \ 
4\rho_{z_1}r_{z_3}- 4 \rho_{z_3} r_{z_1}, \\
A_3:= 
& \ 
4\rho_{z_2}r_{z_1}- 4 \rho_{z_1} r_{z_2}. \\
\endaligned
\end{equation}
Using the expressions~\thetag{ 14.8} for $\rho$ and $r$, we see that
these three components restrict on $\{y=0\}$ as the Pl\"ucker
coordinates of the bivector $(K^1,K^2)$, namely
\def\theequation{14.11}\begin{equation}
\aligned
A_1 \vert_{y=0}=
& \
\rho_2 r_3- \rho_3 r_2=:
\Delta_{2,3}, \\
A_2\vert_{y=0}=
& \
\rho_3 r_1- \rho_1 r_3
=:\Delta_{3,1}, \\
A_3 \vert_{y=0}=
& \
\rho_1 r_2- \rho_2 r_1
=:\Delta_{1,2}.
\endaligned
\end{equation}
As $K^1$ and $K^2$ are of norm $1$ and orthogonal at every point, it
follows by direct computation that $\Delta_{2,3}^2 + \Delta_{3,1}^2+
\Delta_{ 1,2}^2 = 1$ and that the vector of coordinates $(\Delta_{
2,3}, \Delta_{3,1}, \Delta_{1,2})$ is orthogonal to both $K^1$ and
$K^2$. Moreover, as the orthonormal trihedron $(L(x), K^1(x), K^2(x))$
is direct at every point, we deduce that necessarily
\def\theequation{14.12}\begin{equation}
\Delta_{2,3} \equiv a_1, \ \ \ \ \
\Delta_{3,1} \equiv a_2, \ \ \ \ \
\Delta_{1,2} \equiv a_3.
\end{equation}

Next, we compute in length $A_1$, $A_2$ and $A_3$ using~\thetag{
14.8}. As their complete explicit development will not be crucial for
the sequel and as we shall perform with them differentiations and
linear combinations yielding relatively complicated expressions, let
us adopt the following notation: by $\mathcal{ R}^0$, we denote
various expressions which are polynomials in the jets of the functions
$\rho_1,\rho_2,\rho_3$ and $r_1,r_2,r_3$. Similarly, by $\mathcal{
R}^{I}$, by $\mathcal{ R}^{II}$, by $\mathcal{ R}^{III}$ and by
$\mathcal{ R}^{IV}$, we denote polynomials in the transverse variables
$(y_1,y_2,y_3)$ which are homogeneous of degree $1$, $2$, $3$ and $4$
and have as coefficients various expressions $\mathcal{ R}^0$.

Importantly, we make the convention that such expressions $\mathcal{
R}^0$, $\mathcal{ R}^I$, $\mathcal{ R}^{II}$, $\mathcal{ R}^{III}$ and
$\mathcal{ R}^{IV}$ should be totally independent of the constant
$P$. Consequently, if $P$ appears somehow, we shall write it as a
factor, as for instance in $P \, \mathcal{ R}^I$ or in $P^3 \,
\mathcal{ R}^{III}$.

With this convention at hand, we may develope~\thetag{ 14.10} using the
expressions~\thetag{ 14.8} by writing out only the terms which will be
useful in the sequel and by treating the rest as controlled
remainders. Let us detail the computation of $A_1$:
\def\theequation{14.13}\begin{equation}
\aligned
A_1 
& \
=
4 \, 
\left[
-\frac{ i}{2} \rho_3 - iP y_3 + \mathcal{ R}^I
\right] 
\left[
-\frac{ i}{2} r_2 - 2i P^3 y_2^3 + 
\mathcal{ R}^I
\right]- \\
& \ \ \ \ \ \ \ 
- 
4\, 
\left[
-\frac{ i}{2} \rho_2 - iP y_2 + \mathcal{ R}^I
\right] 
\left[
- \frac{ i}{2} r_3 - 2i P^3 y_3^3 + 
\mathcal{ R}^I
\right] \\
& 
\
=
-\rho_3 r_2 -4 P^3 \rho_3y_2^3+ \mathcal{ R}^I - 
2P r_2y_3+ P^4 \mathcal{ R}^{IV}+ P \mathcal{ R}^I + 
\mathcal{ R}^I + P^3\mathcal{ R}^{IV}+ 
\mathcal{ R}^{II} \\
& \ \ \ \ \
+
\rho_2 r_3 + 4 P^3 \rho_2 y_3^3+
\mathcal{ R}^I +
2P r_3y_2+ P^4 \mathcal{ R}^{IV}+ P \mathcal{ R}^I + 
\mathcal{ R}^I + P^3\mathcal{ R}^{IV}+ 
\mathcal{ R}^{II} \\
& \
=
\rho_2 r_3 - \rho_3 r_2 + 
2P r_3 y_2 - 2P r_2 y_3 + 
4P^3 \rho_2 y_3^3 - 4 P^3 \rho_3 y_2^3 + \\
& \ \ \ \ \
+
\mathcal{ R}^I + \mathcal{ R}^{II}+
P \mathcal{ R}^{II}+
P^3 \mathcal{ R}^{IV}+
P^4 \mathcal{ R}^{IV}.
\endaligned
\end{equation}
In the development, before simplification, we firstly write out in 
lines 3 and 4 all
the $9\times 2$ terms of the two product: for instance, the third term
of the first product, namely 
$4(-\frac{i}{ 2} \rho_3 )(\mathcal{ R}^I)$, yields
a term $\mathcal{ R}^I$ whereas the fifth term $4(-iPy_3)(-2iP^3
y_2^3)$ yields a term $P^4 \mathcal{ R}^{IV}$; secondly, we simplify
the obtained sum: by our convention, $\mathcal{ R}^{I}+ \mathcal{
R}^I=\mathcal{ R}^I$, whereas $\mathcal{ R}^I+ P\mathcal{ R}^I$ cannot
be simplified, since the large constant $P$ will be chosen later.
With these technical explanations at hand, we shall not provide any
intermediate detail for the further computations, whose rules are
totally analogous. For $A_1$, $A_2$ and $A_3$, we obtain
\def\theequation{14.14}\begin{equation}
\left\{
\aligned
A_1 
& \
=
\rho_2 r_3 - \rho_3 r_2 + 
2P r_3 y_2 - 2P r_2 y_3 + 
4P^3 \rho_2 y_3^3 - 4 P^3 
\rho_3 y_2^3 + \\
& \ \ \ \ \
+
\mathcal{ R}^I + \mathcal{ R}^{II}+
P \mathcal{ R}^{II}+
P^3 \mathcal{ R}^{IV}+
P^4 \mathcal{ R}^{IV}, \\
A_2
& \
=
\rho_3 r_1 - \rho_1 r_3 + 
2P r_1 y_3 - 2P r_3 y_1 + 
4P^3 \rho_3 y_1^3 - 4 P^3 
\rho_1 y_3^3 + \\
& \ \ \ \ \
+
\mathcal{ R}^I + \mathcal{ R}^{II}+
P \mathcal{ R}^{II}+
P^3 \mathcal{ R}^{IV}+
P^4 \mathcal{ R}^{IV}, \\
A_3
& \
=
\rho_1 r_2 - \rho_2 r_1 + 
2P r_2 y_1 - 2P r_1 y_2 + 
4P^3 \rho_1 y_2^3 - 4 P^3 
\rho_2 y_1^3 + \\
& \ \ \ \ \
+
\mathcal{ R}^I + \mathcal{ R}^{II}+
P \mathcal{ R}^{II}+
P^3 \mathcal{ R}^{IV}+
P^4 \mathcal{ R}^{IV}.
\endaligned\right.
\end{equation}

Now that we have written the complex vector field $\LL$ and its
coefficients $A_1$, $A_2$ and $A_3$, in order to establish Lemma~14.9,
it suffices to choose $P>0$ sufficiently large in order that the four
complex vector fields
\def\theequation{14.15}\begin{equation}
\overline{ \LL}\, \vert_{y=0} , 
\ \ \ \ \
\LL \vert_{y=0}, 
\ \ \ \ \
\left[
\overline{ \LL}, \LL
\right]\, \vert_{y=0}, 
\ \ \ \ \
\left[
\overline{\LL}, 
\left[
\overline{\LL}, 
\left[
\overline{ \LL}, \LL
\right]\right]\right]
\vert_{y=0}
\end{equation}
are linearly independent at every point $x \in \R^3$ with $x_1^2 +
x_2^2 + x_3^2 \leq R^2$. At the end of the proof, we shall explain why
we cannot insure type $3$ at every point, namely why the consideration
of $\left[ \overline{\LL}, \left[ \overline{ \LL}, \LL
\right]\right]\vert_{y=0}$ instead of the length four last Lie bracket
in~\thetag{ 14.15} would fail.

As promised, we shall now summarize all the subsequent
computations. As we aim to restrict the last Lie bracket to $\{y=0\}$
which is of length four and whose coefficients involve derivatives of
order at most three of the coefficients $A_1$, $A_2$ and $A_3$, we can
already neglect the last two remainders $P^3\mathcal{ R}^{IV}$ and
$P^4 \mathcal{ R}^{IV}$ in~\thetag{ 14.14}. In other words, we can
consider $A^1$, $A^2$ and $A^3 \ {\rm mod} (IV)$. Similarly, in the
computation of the Lie bracket
\def\theequation{14.16}\begin{equation}
\left[\overline{\LL},
\LL\right]=:
C_1 \, \partial_{z_1}+
C_2 \, \partial_{z_2}+
C_3 \, \partial_{z_3}- 
\overline{C_1} \, \partial_{\bar z_1}- 
\overline{C_2}\, \partial_{\bar z_2} - 
\overline{C_3} \, \partial_{\bar z_3},
\end{equation}
before restriction to $\{y=0\}$, we can restrict our task to
developing the coefficients
\def\theequation{14.17}\begin{equation}
\aligned
C_1 
& \
:=
\overline{A_1} A_{1,\bar z_1}+
\overline{A_2} A_{1,\bar z_2}+ 
\overline{A_3} A_{1,\bar z_3}, \\
C_2
& \
:=
\overline{A_1} A_{2,\bar z_1}+
\overline{A_2} A_{2,\bar z_2}+ 
\overline{A_3} A_{2,\bar z_3}, \\
C_3
& \
:=
\overline{A_1} A_{3,\bar z_1}+
\overline{A_2} A_{3,\bar z_2}+ 
\overline{A_3} A_{3,\bar z_3} \\
\endaligned
\end{equation}
only modulo order $(III)$, which yields by means of
the expressions~\thetag{ 14.14}
\def\theequation{14.18}\begin{equation}
\aligned
C_1 \ {\rm mod} \, 
(III) 
& \
\equiv 
-iP \rho_1 +
6i P^3 a_3 \rho_2 y_3^2 -
6i P^3 a_2 \rho_3 y_2^2 +
\mathcal{ R}^0 +
\mathcal{ R}^I+ \\
& \ \ \ \ \ \ 
+
P\mathcal{ R}^I+
P^2 \mathcal{ R}^I
+\mathcal{ R}^{II}+
P\mathcal{ R}^{II} +
P^2 \mathcal{ R}^{II}, \\ 
C_2 \ {\rm mod} \, 
(III) 
& \
\equiv 
-iP \rho_2 +
6i P^3 a_1 \rho_3 y_1^2 -
6i P^3 a_3 \rho_1 y_3^2 +
\mathcal{ R}^0 +
\mathcal{ R}^I+ \\
& \ \ \ \ \ \ 
+
P\mathcal{ R}^I+
P^2 \mathcal{ R}^I
+\mathcal{ R}^{II}+
P\mathcal{ R}^{II} +
P^2 \mathcal{ R}^{II}, \\
C_3 \ {\rm mod} \, 
(III) 
& \
\equiv 
-iP \rho_3 +
6i P^3 a_2 \rho_1 y_2^2 -
6i P^3 a_1 \rho_2 y_1^2 +
\mathcal{ R}^0 +
\mathcal{ R}^I+ \\
& \ \ \ \ \ \ 
+
P\mathcal{ R}^I+
P^2 \mathcal{ R}^I
+\mathcal{ R}^{II}+
P\mathcal{ R}^{II} +
P^2 \mathcal{ R}^{II}. 
\endaligned
\end{equation} 
We must mention the use of natural rule hold for computing the partial
derivatives $A_{j, \bar z_k}$: we have for instance $\partial_{\bar
z_k} \left( \mathcal{ R}^{II} \right)= \mathcal{ R}^I+ \mathcal{
R}^{II}$. Also, we have used the hypothesis that $(L(x), K^1(x),
K^2(x))$ provides a direct orthonormal frame at every $x\in \R^3$,
which yields in particular the three relations
\def\theequation{14.19}\begin{equation}
a_2 r_3 - a_3 r_2=-\rho_1, \ \ \ \ \ 
a_3 r_1-a_1 r_3= -\rho_2, \ \ \ \ \ 
a_1r_2-a_2r_1= -\rho_3.
\end{equation}
After mild computation, the coefficients 
$F_1$, $F_2$ and $F_3$ of the length four Lie bracket
\def\theequation{14.20}\begin{equation}
\left[
\overline{\LL}, 
\left[
\overline{\LL}, 
\left[
\overline{\LL}, \LL
\right]
\right]
\right]=
F_1 \, \partial_{z_1}+
F_2 \, \partial_{z_2}+
F_3 \, \partial_{z_3}+
G_1 \, \partial_{\bar z_1}+
G_2 \, \partial_{\bar z_2}+
G_3 \, \partial_{\bar z_3}
\end{equation}
are given, after restriction to $\{y=0\}$, by 
\def\theequation{14.21}\begin{equation}
\aligned
F_1\vert_{y=0} 
& \
=
3i P^3 a_2^3\rho_3 -3i P^3 a_3^3 \rho_2 +
\mathcal{ R}^0+ 
P\mathcal{R}^0 +
P^2\mathcal{ R}^0, \\
F_2\vert_{y=0} 
& \
=
3i P^3 a_3^3\rho_1 -3i P^3 a_1^3 \rho_3 +
\mathcal{ R}^0+ 
P\mathcal{R}^0 +
P^2\mathcal{ R}^0, \\
F_3\vert_{y=0} 
& \
=
3i P^3 a_1^3\rho_2 -3i P^3 a_2^3 \rho_1 +
\mathcal{ R}^0+ 
P\mathcal{R}^0 +
P^2\mathcal{ R}^0, \\
\endaligned
\end{equation}

We can now complete the proof of Lemma~14.9. In the basis
$(\partial_{z_1}, \partial_{z_2}, \partial_{z_3}, \partial_{\bar z_1},
\partial_{\bar z_2}, \partial_{\bar z_3})$, the $4\times 6$ matrix
associated with the four vector fields~\thetag{ 14.15} (without
mentioning $\vert_{y=0}$)
\def\theequation{14.22}\begin{equation}
\left(
\begin{array}{cccccc}
0 & 0 & 0 & a_1 & a_2 & a_3 \\
a_1 & a_2 & a_3 & 0 & 0 & 0 \\
C_1 & C_2 & C_3 & -\overline{C_1} 
& - \overline{C_2} & -\overline{C_3} \\
F_1 & F_2 & F_3 & G_1 & G_2 & G_3
\end{array}
\right)
\end{equation}
has rank four at a point $x\in\R^3$ if and only if the $3\times 3$
determinant in the left low corner is nonvanishing, namely if and only
if the developped expression
\def\theequation{14.23}\begin{equation}
\aligned
{}
& \
\left\vert
\begin{array}{ccc}
a_1 & a_2 & a_3 \\
-iP\rho_1+\mathcal{ R}^0 &
-iP\rho_2+\mathcal{ R}^0 &
-iP\rho_3+\mathcal{ R}^0 \\
3iP^3a_2^3\rho_3-3iP^3a_3^3\rho_2+ &
3iP^3a_3^3\rho_1-3iP^3a_1^3\rho_3+ &
3iP^3a_1^3\rho_2-3iP^3a_2^3\rho_1+ \\ 
+\mathcal{ R}^0+P\mathcal{ R}^0+P^2\mathcal{ R}^0 &
+\mathcal{ R}^0+P\mathcal{ R}^0+P^2\mathcal{ R}^0 &
+\mathcal{ R}^0+P\mathcal{ R}^0+P^2\mathcal{ R}^0 
\end{array}
\right\vert \\
& \
=
3P^4 \left(
r_3 [a_1^3 \rho_2 -a_2^3 \rho_1]+
r_2[a_3^3\rho_1 -a_1^3 \rho_3]+
r_1 [a_2^3 \rho_3 -a_3^3 \rho_2]
\right)+ \\
& \ \ \ \ \ \
+
\mathcal{ R}^0+
P\mathcal{ R}^0+
P^2 \mathcal{ R}^0+
P^3\mathcal{ R}^0 +
P^4\mathcal{ R}^0 \\
& \
=
3P^4 \left(a_1^4+a_2^4+a_3^4\right)+
\mathcal{ R}^0+
P\mathcal{ R}^0+
P^2 \mathcal{ R}^0+
P^3\mathcal{ R}^0 +
P^4\mathcal{ R}^0
\endaligned
\end{equation}
is nonvanishing. 

At this point, the conclusion of the lemma is now an immediate
consequence of the following trivial assertion: {\it Let $a_1$, $a_2$
and $a_3$ be $\mathcal{ C }^\infty$-smooth functions on $\R^3$
satisfying $a_1 (x)^2+ a_2 (x)^2+a_3(x)^2=1$ for all $x \in \R^3$ and
let $\mathcal{ R }_0^0$, $\mathcal{ R}_1^0$, $\mathcal{ R }_2^0$,
$\mathcal{ R }_3^0$ and $\mathcal{ R}_4^0$ be $\mathcal{ C
}^\infty$-smooth functions on $\R^3$. For every $R>0$, there exists a
constant $P>0$ large enough so that the function
\def\theequation{14.24}\begin{equation}
3P^4 \left(
a_1^4+a_2^4+a_3^4\right)+
\mathcal{ R}_0^0+P\mathcal{ R}_1^0+
P^2\mathcal{ R}_2^0 +
P^3 \mathcal{ R}_3^0+
P^4\mathcal{ R}_4^0 
\end{equation}
is positive at every $x\in \R^3$ with $x_1^2 + x_2^2+ x_3^2 \leq R^2$}.

If we had put $y_1^3+ y_2^3+ y_3^3$ instead of $y_1^4+ y_2^4+ y_3^4$
in the second equation~\thetag{14.8}, we would have considered the
length three Lie bracket $\left[ \overline{\LL}, \left[
\overline{\LL}, \LL \right] \right] \vert_{y=0}$ instead of the length
four Lie bracket in~\thetag{14.15}, and hence instead of the quartic
$a_1^4+ a_2^4+ a_3^4$ in~\thetag{ 14.24}, we would have obtained the
cubic $a_1^3+a_2^3+a_3^3$, a function which (unfortunately) vanishes,
for instance if $a_1(x) =\frac{ 1}{\sqrt{ 2}}$, $a_2(x)= -\frac{ 1}{
\sqrt{ 2}}$ and $a_3(x) =0$. We notice that in our example, this value
of $(a_1,a_2,a_3)$ is indeed attained at the point $x \in T^2$ of
coordinates $(\frac{ 3}{\sqrt{2}}, -\frac{ 3}{ \sqrt{ 2}}, 0)$, whence
the necessity of passing to type 4.
The proof of Lemma~14.9 is complete.
\endproof

\vfill
\end{document}